%% file: local_acceleration_uiNEPv.tex
\documentclass[letterpaper]{article}
\usepackage[margin=1in]{geometry}
\usepackage{pgfplots,tikz,float,placeins}
\pgfplotsset{compat=1.18}
\usetikzlibrary{plotmarks}
\usetikzlibrary{positioning}
\usepgfplotslibrary{external} 
\usepackage{amsmath,amsfonts,amssymb,amsthm,xcolor,xspace} 
\usepackage{algorithm,algorithmic,url,hyperref,bm} 
\usepackage{enumerate}
\usepackage{cleveref}
\usepackage{multirow,array,makecell}
\usepackage{booktabs}

\newtheorem{assumption}{Assumption}
\newtheorem{theorem}{Theorem}
\newtheorem{lemma}{Lemma}

\newtheorem{definition}{Definition}
\newtheorem{remark}{Remark}

\newtheorem{example}{Example}

\newcommand{\twoplotscale}{0.496\textwidth}
\newcommand{\threeplotscale}{0.328\textwidth}

\newcommand{\UU}{\mathbb{U}}
\newcommand{\RR}{\mathbb{R}}
\newcommand{\CC}{\mathbb{C}}
\newcommand{\DD}{\mathbf{D}}
\newcommand{\LL}{\mathcal{L}}

\newcommand{\ALIGN}[1]{[\![{#1}]\!]}
\newcommand{\ALIGNM}[1]{[\![{#1}]\!]_{_M}}

\newcommand{\Log}{\mathrm{Log}}
\newcommand{\Exp}{\mathrm{Exp}}
\newcommand{\vect}{\mathrm{vec}}
\newcommand{\TR}{\operatorname{tr}}
\DeclareMathOperator*{\argmin}{arg\,min}
\def\range{\mathrm{range}}

\def\Re{\mathrm{Re}}

\def\Gr{\mathcal{G}}
\def\AAm{AA($m$)\xspace}
\def\X{\mathcal{X}}
\def\Y{\mathcal{Y}}
\def\T{\mathcal{T}}
\def\R{\mathcal{R}}
\def\lo{\mathbf{o}}
\graphicspath{{Figs/}}

\def\dX{E}
\def\dL{D}

\definecolor{myalign}{rgb}{0.49400,0.18400,0.55600}%
\definecolor{mygrassaa}{rgb}{0.00000,1.00000,1.00000}%
\definecolor{myscf}{rgb}{1.00000,0.00000,1.00000}%
\definecolor{mygnt}{rgb}{0.63500,0.07800,0.18400}%

\colorlet{gray}{.}

\title{Anderson Acceleration of SCF Iteration and Newton's Method for 
Solving Nonlinear Eigenvector Problems
}

\author{
	Ding Lu\thanks{
			Department of Mathematics, University of Kentucky,
			Lexington, KY, USA.
			Email: \texttt{Ding.Lu@uky.edu}.
		}
	\and 
	Tom Werner\thanks{
		Institute for Numerical Analysis,
		Technische Universit\"at Braunschweig,
		38106 Braunschweig, Germany.
		Email: \texttt{tom.werner@tu-bs.de}.
	}
    \and
    Zhaojun Bai\thanks{
            Department of Computer Science and Department of Mathematics, 
            University of California, Davis, 
            Davis, CA, USA.
            Email: \texttt{zbai@ucdavis.edu}.
        }
}

\date{\empty}

\begin{document}

\maketitle

\begin{abstract}
The self-consistent field (SCF) iteration is a gateway method for solving
nonlinear eigenvector problems, also known as nonlinear eigenvalue problems
with eigenvector dependency (NEPv). However, the convergence of SCF can be slow
in practice. In this paper, we study two local acceleration techniques
for SCF applied to unitarily invariant NEPv: Anderson acceleration and
Newton acceleration. A central difficulty in developing these techniques
is the non-uniqueness of the eigenvector basis, which complicates direct
application of standard Anderson acceleration and introduces singularity
into the correction equation of the standard Newton formulation.
For Anderson acceleration, we present three formulations based on
unitarily invariant parameterization, basis alignment, and Grassmannian
geometry, and establish regularity properties of the associated SCF maps.
For Newton acceleration, we introduce basis-aligned updates and derive
a correction equation whose operator is nonsingular under suitable
assumptions, and prove local quadratic convergence.
Numerical experiments on NEPv arising from computational physics and
data analysis confirm the theoretical findings and demonstrate the
effectiveness of the proposed acceleration techniques.
\end{abstract}

\noindent
\textbf{Keywords:}
  SCF iteration, nonlinear eigenvalue problem, eigenvector nonlinearity, Anderson acceleration, Newton's method 


\tableofcontents
\newpage

\section{Introduction}
We consider the following nonlinear eigenvalue
problem with eigenvector dependency (NEPv): find a matrix
$X\in\CC^{n,k}$ with orthonormal columns (i.e., $X^HX=I_k$)
and a square matrix $\Lambda\in\CC^{k,k}$ satisfying 
\begin{equation}\label{eq:nepv}
    H(X)X=X\Lambda,
\end{equation}
where $H:\CC^{n,k}\to\CC^{n,n}$ is a given Hermitian matrix-valued function.
Here, the eigenvalues of $\Lambda\equiv X^HH(X)X$ correspond to a subset of $k$
eigenvalues of $H(X)$. We assume that they are the $k$ smallest throughout the paper.
The largest or those at other prescribed positions can be treated analogously.
Typically, $k\ll n$.
In addition, we assume that $H$ is right {\em unitarily invariant}, i.e.,
\begin{equation}\label{eq:uinepv}
    H(X)=H(XQ),\quad
    \text{for all }Q\in\CC^{k,k}\text{ with }Q^HQ=I_k.
\end{equation}
This implies that any solution $(X,\Lambda)$ would give rise
to an equivalent class of solution $(XQ,Q^H\Lambda Q)$.
Therefore, a solution $X$ is intrinsically the subspace $\range(X)$.

NEPv~\eqref{eq:nepv} arises in a wide range of applications.
They have long played important roles in quantum physics
and chemistry through the Hartree--Fock and Kohn--Sham
equations~\cite{BaoD2004,Bretonnet2017,CanLB2000,EchA2007,JarKM2014,Martin2004}.
More recently, they have attracted growing interest in machine learning,
with applications including sparse PCA, trace ratio optimization, and Rayleigh
quotient optimization;
see, e.g.,~\cite{BaiLV2018,Benner:2025,Lu:2025b,WanWL2025,WanZL2023,ZhaL2014,ZhaL2015,ZhaZL2017}.
The widespread use of NEPv stems from their fundamental role
as first-order optimality conditions for optimization
with orthogonality constraints, and the invariance
condition~\eqref{eq:uinepv} reflects that the solution sought is a subspace rather than a particular basis in those
applications; see, e.g.,~\cite{Li:2024}.

A standard approach to solving~\eqref{eq:nepv} is the
self-consistent field (SCF) iteration.
Starting from an initial matrix $X_0$ with orthonormal columns,
SCF repeatedly solves the linear eigenvalue problem%
\begin{equation}\label{eq:scf}
    H(X_i)X_{i+1}=X_{i+1}\Lambda_{i+1},
    \qquad i=0,1,2,\dots,
\end{equation}
where the columns of $X_{i+1}$ are orthonormal eigenvectors
corresponding to the $k$ smallest eigenvalues of $H(X_i)$,
and $\Lambda_{i+1}$ is Hermitian but need not be diagonal.
We note that the iterate matrix $X_{i+1}$ is generally not unique, 
with different choices related by an orthonormal rotation $X_{i+1}Q$,
provided that the $k$ smallest eigenvalues are separated from the rest of the spectrum.
By the unitary invariance~\eqref{eq:uinepv}, these choices
give the same $H(X_{i+1})$ and thus do not affect the subsequent iterations.
Hence, the SCF~\eqref{eq:scf} is an iteration on subspaces 
$\{\range(X_i)\}_{i=0}^\infty$, 
with each $X_i$ serving as a representative orthonormal basis.

The SCF iteration was originally developed for electronic
structure calculations~\cite{Hartree:1928,Roothaan:1951}.
Under suitable assumptions, local convergence can be
established, typically at a linear rate \cite{BaiLL2022,CaiZBL2018,Liu:2014,Upadhyaya:2021a,Yang:2009}.
In practice, slow convergence or failure to converge
may occur and has motivated various techniques to stabilize and accelerate SCF,
such as damping, level shifting, and mixing of iterates \cite{ChuDLS2021,Kerker1981,Pulay1980,Pulay1982,SauH1973,ThoOYJSH2004}.
However, many of these techniques have been developed
for specific applications in computational physics,
such as density functional theory, where they exploit
particular functional structures of $H(X)$,
and they do not readily apply to
general NEPv~\eqref{eq:nepv} in a computationally efficient way.

The purpose of this paper is to investigate acceleration
techniques of the SCF iteration~\eqref{eq:scf} 
for unitarily invariant NEPv~\eqref{eq:nepv} in its most general form.
We focus on two types of acceleration:
Anderson acceleration (AA), which combines recent iterates, 
and Newton acceleration, which uses
derivative information of $H(X)$ to achieve superlinear acceleration. 

AA is a widely used technique
for accelerating fixed-point iterations~\cite{Anderson1965}.
The basic idea is to form a refined approximation through a linear combination
(or mixing) of recent iterates, using suitably chosen coefficients;
see, e.g.,~\cite{Anderson2019,Evans:2020, FanS2009,TotK2015,WalN2011} for details.
Applying AA to the SCF~\eqref{eq:scf} requires
extra care because the iterates $X_i$ are not uniquely
determined, and directly mixing them depends on the arbitrary unitary rotations $X_iQ$.
In electronic structure calculations, AA is closely related
to Pulay mixing, also known as direct inversion in the
iterative subspace (DIIS)~\cite{Pulay1980,Pulay1982}.
For these specialized NEPv, the SCF can be
expressed in terms of quantities such as the charge density
$\mathrm{diag}(XX^H)$, density matrix $XX^H$, or Fock
matrix $H(X)$, which are independent of the chosen basis;
see, e.g.,~\cite{GarS2012,HuY2010,KudSC2002,Pulay1982}.
For example, for the Kohn-Sham NEPv that arises in density functional theory, 
where $H(X)$ is fully determined by the charge density $\rho(X)=\mathrm{diag}(XX^H)\in\RR^n$,
the SCF can be recast as a fixed-point iteration $\rho_{i+1}=g(\rho_i)$ 
in charge densities, and AA can then be applied directly to this iteration;
see, e.g.,~\cite{Saad:2025}.
However, for general NEPv~\eqref{eq:nepv} such a compact density-based
representation of $H(X)$ may not be readily available,
while mixing the full matrices $H(X_i)\in\CC^{n,n}$ can be costly for
large $n$.
This calls for an in-depth investigation of AA in these settings, 
particularly into how to work economically with tall-skinny basis matrices $X_i$.

On the other hand, Newton's method provides another means of locally
accelerating SCF by using derivative information to achieve
superlinear convergence under suitable assumptions.
For NEPv~\eqref{eq:nepv}, applying Newton's acceleration
also requires care 
because the solution basis is not uniquely determined.
Most existing Newton-type methods are derived for real-variable NEPv 
with $X\in\RR^{n,k}$ and $H(X)\in\RR^{n,n}$.
A direct approach applies Newton's method to the equations
$H(X)X-X\Lambda=0$ and $X^TX=I_k$;
see, e.g.,~\cite{GaoYM2007,Werner2025}.
However, the unitary invariance~\eqref{eq:uinepv} implies
that a solution $(X,\Lambda)$ is not isolated, since
$(XQ,Q^T\Lambda Q)$ remains a solution for every unitary $Q$,
which leads to a singular Jacobian (derivative) operator 
that creates difficulties for standard Newton analysis and numerical
solution~\cite{Werner2025}.
To avoid this singularity, the approach in~\cite{JarU2022} replaces
the orthogonality constraint $X^TX=I_k$ by a normalization condition 
$C^TX=I_k$ with a fixed matrix $C\in\RR^{n,k}$.
This formulation assumes that the NEPv satisfies $H(XQ)=H(X)$ for every nonsingular
$Q\in\RR^{k,k}$, which is stronger than
the unitary invariance~\eqref{eq:uinepv}.
Along a different line, Riemannian Newton methods on the
Grassmann manifold have been developed for NEPv arising
from energy minimization, which removes the basis non-uniqueness
by working with subspaces; see, e.g.,~\cite{Altmann:2024,Zhao:2015}.
These derivations, however, rely on an underlying energy
functional and differential-geometric tools.
To our knowledge, a Newton acceleration for general
NEPv~\eqref{eq:nepv}, based on an elementary derivation 
and a correction operator that is nonsingular under suitable assumptions,
has yet to be developed.

\paragraph{Contributions.}
The main contributions of this work are threefold:
\begin{itemize}
\item
For Anderson acceleration, we present three formulations
that resolve the basis non-uniqueness of SCF iteration: (i) a parameterization-based
framework that exploits the structure of $H(X)$,
(ii) an alignment-based AA, and (iii) an adaptation of Riemannian AA.
The first extends the density-based approach to general
NEPv, while the latter two work with $n\times k$ basis
matrices without requiring a parameterization of $H(X)$. 
The regularity properties of the associated SCF maps are established. 

\item
For Newton acceleration, we introduce the notion of a
basis-aligned update to characterize local corrections
to an approximate solution of the NEPv.
This naturally removes basis non-uniqueness and leads to 
a Newton equation for correction through elementary linear algebra operations. 
It also provides a unified treatment to real and complex NEPv.
Under suitable assumptions, we establish the nonsingularity
of the correction operator and local quadratic convergence.

\item
Through extensive numerical experiments on NEPv arising
from applications in computational physics and machine learning,
we confirm the theoretical findings and demonstrate the effectiveness and computational benefits of the proposed acceleration techniques. 
\end{itemize}

\paragraph{Organization.}
The rest of the paper is organized as follows.
\Cref{sec:aa} reviews Anderson acceleration and presents
three variants for SCF, 
including parameterization-based, alignment-based, and Riemannian
formulations of AA.
\Cref{sec:newton} derives the basis-aligned Newton correction,
establishes its local quadratic convergence, and discusses its implementation.
\Cref{sec:experiments} presents numerical experiments
on problems from computational physics and data science.
Concluding remarks and future directions are given in~\Cref{sec:conclusion}.

\paragraph{Notation.}
Throughout this paper, we use standard notation from matrix analysis. 
We denote by $\RR^{n,k}$ and $\CC^{n,k}$ the spaces of
real and complex $n\times k$ matrices, respectively,
with $\RR^n=\RR^{n,1}$ and $\CC^n=\CC^{n,1}$.
The set of complex orthonormal matrices is denoted by 
$\UU^{n,k}:=\{X\in\CC^{n,k}:X^HX=I_k\}$.
The Grassmann manifold of $k$-dimensional subspaces of $\CC^n$
(or $\RR^n$ in the real case) is denoted by $\Gr(n,k)$.
For a matrix $A$, $A^T$ and $A^H$ denote its transpose
and conjugate transpose, respectively, 
$\range(A)$ denotes its column space,
and $\vect(A)$ denotes the vector obtained by stacking the columns of $A$.
For square $A$, $\mathrm{tr}(A)$ denotes its trace.
The identity matrix of size $n$ is denoted by $I_n$,
with the subscript omitted when its size is clear.
We use $\mathrm{diag}(A)$ for the vector of diagonal
entries of $A$ and $\mathrm{Diag}(v)$ for the diagonal
matrix with diagonal $v$.
For a Hermitian $A\in \CC^{n,n}$,
its eigenvalues are ordered as $\lambda_1(A)\leq \cdots\leq \lambda_n(A)$.
For a block matrix $\mathcal A=[A_1,\dots,A_m]$,
and a vector $\beta=(\beta_1,\dots,\beta_m)^T\in\RR^m$,
we define the block product
$\mathcal A\odot\beta :=\sum_{j=1}^m\beta_jA_j$.
Unless otherwise stated, $\|\cdot\|$ denotes the 
spectral norm for matrices and $\|\cdot\|_F$ denotes the Frobenius norm.
For complex matrices, we use the (real) Frobenius inner product
$\langle A,B\rangle:=\Re(\TR(A^HB))$.
Differentiability of functions with complex variables
is understood with respect to the real and imaginary parts
of their entries,
and we write $\DD F(X)[E]$ for the Fr\'echet derivative
of $F$ at $X$ along the direction $E$.
We use the standard big-$O$ and little-$o$ notation:
$f=\mathcal O(g)$ means $|f|\leq C|g|$ for some constant $C$
as $g\to0$,  and $f=\lo(g)$ means $|f|/|g|\to0$.
Other notations will be explained as needed.


\section{Anderson Acceleration for NEPv} \label{sec:aa}
We consider Anderson Acceleration (AA) of the SCF iteration~\eqref{eq:scf}.
Recall that applying AA to SCF requires care because of the unitary
invariance~\eqref{eq:uinepv} of the NEPv.
Indeed, each iterate $X_{i+1}$ by the SCF~\eqref{eq:scf}
is not uniquely defined: 
For any unitary matrix $Q\in\UU^{k,k}$, 
the matrix $X_{i+1}\leftarrow X_{i+1} Q$ 
is also a valid solution.
Thus, at the level of basis matrices, SCF is a set-valued
fixed-point iteration 
\begin{equation}\label{eq:scffp}
X_{i+1}\in g(X_i),
\end{equation}
where 
\begin{equation}\label{eq:scffpmap1}
g(X)\!:=\! \left\lbrace Y: \text{$Y\in\UU^{n,k}$ is an 
	orthonormal eigenbasis for 
	$k$ smallest eigenvalues of $H(X)$}\right\rbrace.
\end{equation}
Consequently, the fixed-point residual $r(X)=g(X)-X$,
as used in AA, is not well defined, 
since it depends on the particular choice of $Y\in g(X)$ and need
not vanish even at a true solution $X^*$ of the NEPv~\eqref{eq:nepv}.
In this section, we resolve this ambiguity.
After a brief review of AA for the 
fixed-point iteration, we will present three approaches, 
including 
a parameterization-based AA framework that exploits the functional 
structure of the NEPv for computational efficiency, 
an alignment-based AA that works directly with basis matrices,
and a Riemannian AA that uses geometric tools on the Grassmann manifold.

\subsection{Anderson Acceleration for Fixed-Point Iterations}
In this section, we briefly review AA for standard
fixed-point iterations; further details on the method
and its convergence properties are available in, e.g.,~\cite{Anderson1965,Anderson2019, Evans:2020,FanS2009,Garner:2026,TotK2015,WalN2011}.
Consider the fixed-point equation 
\begin{equation}\label{eq:fpeq}
   x = g(x), 
\end{equation}
where $g\colon \RR^n\to \RR^n$ is continuous.
Given an initial $x_0$, the corresponding fixed-point iteration is 
\begin{equation}\label{eq:fpgx}
x_{i+1}=g(x_i), \quad i=0,1\dots.
\end{equation}
If the iteration converges, its limit $x^*:=\lim_{i\to \infty} x_i$ 
satisfies~\eqref{eq:fpeq}.
Convergence, however, is not guaranteed and depends on
both the mapping $g$ and the initial approximation $x_0$.
For initial $x_0$ sufficiently close to a fixed
point $x^*$, convergence is ensured if $g$ is continuously
differentiable near $x^*$ and the spectral radius
of the Jacobian matrix $g'(x^*)$ (largest absolute value of its eigenvalues) 
is strictly less than one; see, e.g.,~\cite{Kelley:1995}.
In such cases, the convergence is generally linear and can be slow
when this spectral radius is close to one.

AA seeks to improve the local convergence
of fixed-point iteration by mixing recent
iterates~\cite{Anderson1965}.
Specifically, \AAm scheme uses the $m+1$ most recent
iterates $x_i,\dots,x_{i-m}$ to construct
\begin{subequations} \label{eq:aa}
\begin{equation}\label{eq:aastep}
 x_{i+1}:=\sum_{\ell =i-m}^{i} \alpha_\ell \cdot g(x_\ell),
\end{equation}
where the coefficients $\{\alpha_\ell\}$ are chosen to 
minimize the norm of a linear combination of the corresponding residuals:
\begin{equation}\label{eq:aacoeff}
\min_{\sum_{\ell=i-m}^{i}\alpha_\ell =1} \left\|\sum_{\ell=i-m}^{i} \alpha_\ell \cdot r(x_{\ell})\right\|
\qquad \text{with} \quad r(x):=g(x) - x.
\end{equation}    
\end{subequations}
The motivation for~\eqref{eq:aacoeff} follows from a local linearization
of the residual. 
If $\{x_\ell\}$ are sufficiently close to $x_i$, and the coefficients
$\{\alpha_\ell\}$ with $\sum_{\ell=i-m}^{i}\alpha_\ell =1$ 
are bounded absolutely by some constant (independent of $\{x_\ell\}$), 
then the first-order Taylor's expansion about $x_i$ gives 
\begin{equation}\label{eq:resd}   
r\left( \sum_{\ell=i-m}^{i} \alpha_\ell x_{\ell}\right)
=
\sum_{\ell = i-m}^{i} \alpha_\ell \cdot r(x_{\ell})
+ 
\lo\left(\max_{\ell=i-m:i-1} \|x_{\ell} - x_i\|\right).
\end{equation}
Thus, the least-squares problem~\eqref{eq:aacoeff}
approximately minimizes the residual norm of the mixed
approximation $x_{i+1}':=\sum_{\ell=i-m}^{i}\alpha_\ell x_\ell$.
The update~\eqref{eq:aastep} then replaces each $x_\ell$
by its fixed-point image $g(x_\ell)$, incorporating
a further fixed-point correction into the mixed approximation.
In practice, one may also damp the update~\eqref{eq:aastep} by setting
$x_{i+1}\leftarrow (1-\gamma)x'_{i+1}+\gamma x_{i+1}$
with a damping parameter $\gamma\in[0,1]$; 
see, e.g.,~\cite{Saad:2025}.
For simplicity of exposition, we take $\gamma=1$ (no damping) 
throughout this paper.

In implementation, the constraint $\sum_{\ell=i-m}^{i}\alpha_\ell =1$
on the mixing coefficients in the least squares~\eqref{eq:aacoeff}
can be eliminated by expressing \AAm in an equivalent
difference form; see, e.g.,~\cite{FanS2009}.
Specifically, the update~\eqref{eq:aastep} becomes
\begin{subequations}\label{eq:aa:inc}
\begin{equation}\label{eq:aastep:inc}
    x_{i+1}=g(x_i)
    -\sum_{\ell=i-m}^{i-1}\beta_\ell
    \bigl(g(x_{\ell+1})-g(x_\ell)\bigr),
\end{equation}
where the coefficients $\{\beta_\ell\}$ solve the unconstrained
least-squares problem
\begin{equation}\label{eq:aacoeff:inc}
    \min_{\beta}
    \left\|
    r(x_i)-\sum_{\ell=i-m}^{i-1}\beta_\ell
    \bigl(r(x_{\ell+1})-r(x_\ell)\bigr)
    \right\|.
\end{equation}
\end{subequations}
The equivalence between~\eqref{eq:aa:inc} and~\eqref{eq:aa}
follows from the one-to-one correspondence
$\alpha_{i-m}=\beta_{i-m}$,
$\alpha_\ell=\beta_\ell-\beta_{\ell-1}$
for $i-m+1\leq\ell\leq i-1$, and
$\alpha_i=1-\beta_{i-1}$.
The standard least-squares problem~\eqref{eq:aacoeff:inc} 
can then be solved  by direct QR factorization
or iterative methods such as LSQR;
see, e.g.,~\cite{Bjorck:2024,GolV2013}.


%
\begin{algorithm}[!htbp]
\caption{AA($m$) in difference form}\label{alg:aag}
\begin{algorithmic}[1]
\REQUIRE Initial $x_0\in\RR^{n}$, $g:\RR^n\rightarrow\RR^n$, $m\in\mathbb{N}$, error tolerance $\tau>0$ 
\ENSURE Approximate fixed-point $x^*\in\RR^{n}$ of $g$.
\FOR{$i=0,1,2,\dots$}
\STATE Compute $y_i=g(x_i)$ and $r_i=y_i-x_i$
\STATE {\bf if} $\|r_i\|<\tau$, {\bf then return} $x^*=y_i$ \label{alg:aag:line:res}
\STATE Set $m_i=\min\{m,i\}$
\IF{$m_i=0$}
    \STATE Set $x_{i+1}=y_i$
\ELSE
    \STATE Form $\Delta X_i=[\Delta x_{i-m_i},\dots,\Delta x_{i-1}]$ and
    $\Delta R_i=[\Delta r_{i-m_i},\dots,\Delta r_{i-1}]$, where
    $\Delta x_\ell=x_{\ell+1}-x_\ell$ and $\Delta r_\ell=r_{\ell+1}-r_\ell$
    \STATE \label{alg:aag:i:ls}
    Compute $\beta^{(i)}\in\RR^{m_i}$ by solving
    $\beta^{(i)}=\argmin_{\beta\in\RR^{m_i}}
    \|r_i-\Delta R_i\beta\|$
    \STATE Set $x_{i+1}= \widetilde x_{i+1} +\widetilde r_{i+1}$,
	where 
	$\widetilde x_{i+1} := x_i - \Delta X_i\beta^{(i)}$
	and 
	$\widetilde r_{i+1} := r_i - \Delta R_i\beta^{(i)}$
	\label{alg:aag:i:xp1}
\ENDIF
\ENDFOR
\end{algorithmic}
\end{algorithm}
Algorithm~\ref{alg:aag} summarizes \AAm in difference form,
where we apply AA from the first iteration using the
available history length $m_i=\min\{m,i\}$.
Throughout this paper, we use a prescribed, fixed window
size $m$.
For greater flexibility, one may also consider adaptive
truncation~\cite{ChuDLS2021} or switching between Anderson
mixing and fixed-point steps~\cite{BarL2026,Lupo2019}.
We stop the algorithm when the fixed-point residual norm is below 
an error tolerance $\tau$.
\subsection{Parameterization-based Anderson Acceleration} \label{sec:aapara}

In this section, we describe a general framework for applying
AA to SCF through a parameterization
of the NEPv~\eqref{eq:nepv} that is independent of the
chosen orthonormal basis.

\subsubsection{Unitary-Invariant Parameterization}
Our development is based on the notion of
unitarily-invariant parameterization of the NEPv~\eqref{eq:nepv} given below.  
\begin{definition} [ui--parameterization]
\label{def:uiparam}
	We call $\rho\colon \CC^{n,k} \to \CC^{p}$ a 
	unitarily-invariant parameterization of the NEPv~\eqref{eq:nepv},
	if 
	\begin{equation}
	\text{$\rho(XQ) \equiv \rho(X)$, \quad for all $X\in\UU^{n,k}$ and
	$Q\in\UU^{k,k}$},
	\end{equation}
	and there exists a Hermitian matrix-valued function 
	$\widehat H \colon \CC^p \to \CC^{n,n}$ such that 
	\begin{equation}\label{eq:parah}
	H(X) = \widehat H(\rho(X)), \quad\text{for all $X\in \UU^{n,k}$.}
	\end{equation}
\end{definition}

There is always a ui-parameterization of the NEPv~\eqref{eq:nepv}, although 
it is generally not unique.
For example, since the Hermitian matrix $H(X)$ is completely defined by 
its lower triangular entries (including the diagonal),
a trivial ui-parameterization is 
$\rho(X):=\mbox{tril}(H(X))\in\CC^p$ with $p=\frac{n(n+1)}{2}$,
and $\widehat H(\cdot)$ reconstructs the full matrix $H(X)$ from
its lower triangular part.
In practice, since $H(\cdot)$ has a structure such as sparsity or other
functional structure, one may construct a much more compact
ui-parameterization with $p\ll n^2$ or even $p\ll n$ as shown in 
the following examples. 

\begin{example}\label{rmk:parameterization}
\rm
Here are three practically relevant examples of ui-parameterization of NEPv.

\begin{enumerate}[(i)]
\item 
{\em Sparsity:}
Suppose that $H(X)=A+B(X)$, where $A$ is a constant matrix 
and $B(X)$ is a sparse Hermitian matrix with a fixed 
sparsity pattern $\Pi$ with $p$ non-zero entries
in the upper (or lower) triangle.
Then, we can define $\rho(X)\in\CC^p$ as a collection of the non-zero entries
in $B(X)$, 
and $\widehat H$ updates $A$ according to the given 
sparsity pattern $\Pi$ using the entries of $\rho(X)$. 

\item
{\em Functional structure through charge density:}
A Kohn-Sham NEPv from density functional theory~\cite{Bretonnet2017,CaiZBL2018} involves
\[
H(X) = 
A +  \mbox{Diag}(L^\dag \rho(X)) +  \mbox{Diag}(\mu_{xc}(\rho(X))^T\bf 1),
\]
where $\rho\colon \mathbb R^{n\times k}\to \mathbb R^n$ 
with $\rho(X) = \mbox{diag}(XX^T)$ is the charge density, $\mu_{xc}\colon \mathbb R^{n}\to \mathbb R^{n\times n}$ 
is a given function accounting for exchange energy and $L^\dag$ is the pseudo-inverse of a discrete Laplacian $L$.
Since the dependence on $X$ is entirely through $\rho(X)$, we have the
ui-parameterization~\eqref{eq:parah}: 
\[
\widehat H(\rho) =  A +  \mbox{Diag}(L^\dag \rho) + \mbox{Diag}(\mu_{xc}(\rho)^T\bf 1),
\]
where $\rho(X)\in\RR^n$ has $p=n$ components. 
Recall that this charge-density parameterization underlies
density-based Pulay mixing (DIIS) for SCF; see, e.g.,~\cite{Saad:2025}.

\item
{\em Functional structure through Rayleigh quotients:}
Many practical NEPv involve $H(X)$ that depends on $X$ only through a small number of
Rayleigh quotients $\mbox{tr}(X^H M_j X)$, 
where $M_j\in\mathbb C^{n\times n}$ are Hermitian matrices; see e.g.~\cite{Bai:2024,Benner:2025,Lu:2025b}.
In these cases, a ui-parameterization is given by 
$\rho(X):=[\rho_1(X),\ldots,\rho_p(X)]^T$, where $\rho_j(X):=\mbox{tr}(X^H M_jX)$, 
with $H(X) = \widehat{H}(\rho(X))$ for some $\widehat H\colon\mathbb R^p\to\mathbb C^{n\times n}$.
As a representative example, consider the sum-of-trace-ratio
NEPv~\cite{Li:2024}, where 
\[
H(X) = \frac{1}{\rho_2(X)}\cdot M_1 
     - \frac{\rho_1(X)}{\rho_2(X)^2}\cdot M_2 + M_3,
\]
with the ui-parameterization given by 
$\rho(X)=[\rho_1(X),\rho_2(X)]^T\in\RR^p$ with $p=2$,  and 
\[
\widehat H(\rho) = \frac{1}{\rho_2}\cdot M_1
                 - \frac{\rho_1}{\rho_2^2}\cdot M_2 + M_3.
\]
\end{enumerate}

\end{example}

\subsubsection{SCF in Parameterized Form}
Based on a ui-parameterization $\rho$ of the NEPv, 
the SCF~\eqref{eq:scf} can be recast as an iteration 
in terms of the parameters $\rho_i :=\rho(X_i)$: 
Given an initial $\rho_0=\rho(X_0)$, 
compute iteratively
\begin{equation}\label{eq:rhoscf}
	\rho_{i+1} = \rho(X_{i+1})
	\quad \text{with}\quad
	\widehat H(\rho_i) X_{i+1} = X_{i+1}\Lambda_{i+1} 
    \quad \text{for} \quad i=0,1,\dots,
\end{equation}
where $X_{i+1}$ contains the $k$ eigenvectors  of $\widehat
H(\rho_i)\equiv H(X_i)$ corresponding to the $k$ smallest eigenvalues.
In the form of fixed-point iteration, we can write~\eqref{eq:rhoscf} as 
\begin{equation}\label{eq:fprho}
\rho_{i+1} = g(\rho_i), 
\end{equation}   
with the fixed-point mapping $g:\CC^p\to\CC^p$ defined as
\begin{equation} \label{eq:scfpf}
    g(y) := \lbrace \rho(Y):Y\in\UU^{n,k} \text{ is an eigenbasis 
    of the $k$ smallest eigenvalues of $\widehat H(y)$}\rbrace.
\end{equation}
Note that $g(y)$ may have multiple values when the $k$-dimensional
eigenspace associated with the $k$ smallest eigenvalues is not uniquely determined.
Such ambiguity, however, does not occur under the eigengap condition
$\lambda_k<\lambda_{k+1}$ for the matrix $\widehat H(y)$.
The following lemma formalizes the well-definedness and local
differentiability of the parameter-space fixed-point mapping $g(y)$ under
this condition.
\begin{lemma}\label{lem:param_scf_smooth}
Let $\rho$ be a ui-parameterization of
the NEPv~\eqref{eq:nepv},  with associated
matrix-valued function $\widehat H$ satisfying~\eqref{eq:parah}.  
Suppose that, at $y\in\CC^p$, the $k$ smallest eigenvalues of
$\widehat H(y)$ are separated
from the rest of the spectrum, i.e.,
$\lambda_k(\widehat H(y))<\lambda_{k+1}(\widehat H(y))$,
and $Y\in\UU^{n,k}$ is a corresponding eigenbasis.
Then
\begin{enumerate}[(i)]
\item\label{i:lem:param_scf_smooth}
the fixed-point mapping $g$ in~\eqref{eq:scfpf} is single-valued at
$y$ (i.e., $g(y)$ is uniquely defined);
\item\label{ii:lem:param_scf_smooth}
the mapping $g$ is locally continuous (and differentiable) at $y$,
if $\rho$ and $\widehat H$ are continuous (and differentiable) at $Y$
and $y$, respectively, 
both in terms of the real and imaginary parts of their arguments.
\end{enumerate}
\end{lemma}

\begin{proof}
For item~\eqref{i:lem:param_scf_smooth}:
The eigengap condition implies that the invariant subspace associated with the
$k$ smallest eigenvalues of $\widehat H(y)$ is unique
(see, e.g.,~\cite{GolV2013}).  Hence, any two orthonormal bases of this
subspace differ only by a right unitary factor.
By the unitary invariance of $\rho$, the value $\rho(X)$ is independent of
the choice of such a basis, and therefore, $g(y)$ is single-valued.

For item~\eqref{ii:lem:param_scf_smooth}:
If $\widehat H$ is continuous (and differentiable), 
then the entries of the Hermitian matrix
$\widehat H(y)$ are continuous (and differentiable)
functions of the real and imaginary
parts of $y$.  
By classical eigenvalue perturbation theory for Hermitian
matrices and the gap condition,
there is, locally around $y$, a choice of orthonormal
eigenbasis $Y(y)$ for the $k$ smallest eigenvalues
that is continuous (and differentiable); see,
e.g.,~\cite[Chap 2.5]{Kato:1995}.
Since $\rho$ is
also continuous (and differentiable) by assumption, 
the composition $g(y)=\rho(Y(y))$ is locally
continuous (and differentiable).
\end{proof}

Conditions for the local convergence of SCF have been 
established in~\cite{BaiLL2022,CaiZBL2018,Liu:2014}.
Under these conditions, the fixed point map $g$ in~\eqref{eq:scfpf} is locally contractive. 
The regularity condition in~\Cref{lem:param_scf_smooth},
together with existing local convergence analyses for AA (e.g.,~\cite{Evans:2020,TotK2015,WalN2011}),
provides a basis for establishing the local convergence of Anderson accelerated SCF~\eqref{eq:fprho}.

\subsubsection{Anderson Acceleration in Parameter Space}
By~\Cref{lem:param_scf_smooth}, the parameterized
SCF iteration~\eqref{eq:fprho} produces a uniquely
defined sequence $\{\rho_i\}$, provided that
$\lambda_k(\widehat H(\rho_i))<
\lambda_{k+1}(\widehat H(\rho_i))$ at each step.
This removes the ambiguity associated with the choice
of eigenbasis in~\eqref{eq:scf} and allows AA to be
applied directly in parameter space.
The resulting adaptation of \AAm in~\Cref{alg:aag}
is summarized in~\Cref{alg:aarho}.

\begin{algorithm}[!htbp]
\caption{AA($m$) for parameterized SCF in difference form}
\begin{algorithmic}[1]
\REQUIRE Initial $X_0\in\UU^{n,k}$, ui-parameterization $\rho$ with function
$\widehat H$ satisfying~\eqref{eq:parah}, $m\in\mathbb{N}$, $\tau>0$
\ENSURE approximate eigenbasis $X^*\in\UU^{n,k}$ of the NEPv 
\STATE Evaluate $\rho_0=\rho(X_0)$
\FOR{$i=0,1,2,\dots$}
\STATE Evaluate $H_{i}=\widehat H(\rho_{i})$
\STATE Solve $H_iX_{i+1}=X_{i+1}\Lambda_{i+1}$ for the $k$ smallest eigenvalues in $\Lambda_{i+1}$
\STATE Evaluate $y_i=\rho(X_{i+1})$ and $r_i=y_i-\rho_i$
\STATE {\bf if} {$\|r_i\|\leq\tau$}, {\bf then return} $X^*=X_{i+1}$
\STATE Set $m_i=\min\{m,i\}$
\IF{$m_i=0$}
    \STATE Set $\rho_{i+1}=y_i$
\ELSE
    \STATE Form $\Delta P_i=[\Delta \rho_{i-m_i},\dots,\Delta \rho_{i-1}]$ and
    $\Delta R_i=[\Delta r_{i-m_i},\dots,\Delta r_{i-1}]$, where
    $\Delta \rho_\ell=\rho_{\ell+1}-\rho_\ell$ and
    $\Delta r_\ell=r_{\ell+1}-r_\ell$
    \STATE\label{alg:aarho:i:ls} Compute $\beta^{(i)}\in\RR^{m_i}$ by solving
    $\beta^{(i)}=\argmin_{\beta\in\RR^{m_i}}
    \|r_i-\Delta R_i\beta\|$
    \STATE Set $\rho_{i+1}= \widetilde \rho_{i+1} +\widetilde r_{i+1}$,
	where 
	$\widetilde \rho_{i+1} := \rho_i - \Delta P_i\beta^{(i)}$
	and 
	$\widetilde r_{i+1} := r_i - \Delta R_i\beta^{(i)}$
	\label{rho:alg:aarho}
\ENDIF
\ENDFOR
\end{algorithmic}    
\label{alg:aarho}
\end{algorithm}
We caution that the least squares problem in line~\ref{alg:aarho:i:ls}
of~\Cref{alg:aarho} may involve complex-valued $r_i$ and $\Delta R_i$.
In this case, since the coefficients $\beta^{(i)}$ are restricted to be real
(as specified in the algorithm), the complex least-squares problem can be
equivalently recast as a real least-squares problem by separating the real
and imaginary parts of both the residual and the difference matrix.
Specifically, minimizing
$\|r_i - \Delta R_i \beta\|_2$ over $\beta \in \RR^{m_i}$ is equivalent to
solving
\[
\beta^{(i)} = \mathop{\mbox{argmin}}_{\beta \in \RR^{m_i}}
\left\|
\begin{bmatrix}
\operatorname{Re}(r_i) \\
\operatorname{Im}(r_i)
\end{bmatrix}
-
\begin{bmatrix}
\operatorname{Re}(\Delta R_i) \\
\operatorname{Im}(\Delta R_i)
\end{bmatrix}
\beta
\right\|_2,
\]
which is a standard real least-squares problem of size $2p \times m_i$.

%

\subsection{Alignment-based Anderson Acceleration} \label{sec:aaaligned}
We next describe a direct approach to resolving basis
non-uniqueness in SCF: at each step~\eqref{eq:scf},
we select a representative eigenbasis by aligning it
with a fixed reference basis.
This alignment can be viewed as a specialized
ui-parameterization, but merits separate treatment
because it is explicit, inexpensive to compute when
$k\ll n$, and does not rely on any particular functional
structure of $H(X)$.

\subsubsection{Basis Alignment}
We begin with a definition of an alignment operation. 
\begin{definition}\label{def:alignment}
Let $M\in\UU^{n,k}$ be a given reference basis matrix.
The alignment of an orthogonal basis $X\in\UU^{n,k}$ 
with respect to $M$ is defined as 
the orthogonal basis of $\range(X)$ that is closest to $M$
in the Frobenius norm, namely,
\begin{equation}\label{eq:alignment}
    \ALIGNM{X}:=\argmin \lbrace\|Y-M\|_F:Y\in\UU^{n,k},
	~\range(Y)=\range(X)\rbrace, 
\end{equation}
where the minimizer is assumed to be unique; 
otherwise $\ALIGNM{X}$ denotes any particular minimizer.
\end{definition}

Basis alignment is a classical operation and 
is related to the orthogonal Procrustes problem; see,
e.g.,~\cite{GolV2013}. 
It is well known that the alignment $\ALIGNM{X}$ can be obtained 
by a unitary rotation, 
\begin{equation}
    \ALIGNM{X}=XQ,\quad \text{with $Q=UV^H$}, \label{eq:alignX}
\end{equation}
where $U,V\in\UU^{k,k}$ are obtained from 
the singular value decomposition $X^HM=U\Sigma V^H$.
The matrix $Q$ may also be interpreted as the unitary polar 
factor of the matrix $X^HM$.

The alignment $\ALIGNM{X}$ is
uniquely defined, unitarily invariant, and locally
differentiable in $X$, whenever $X^HM$ is nonsingular.
These properties follow from the explicit formula~\eqref{eq:alignX}
and standard properties of the polar
decomposition; see, e.g., \cite[Thm.~8.6]{Higham2008} and
\cite[Lem.~A.1]{LuL2024}.
We summarize them in the following lemma and
omit the proof.
\begin{lemma}\label{lem:alignment_smooth}
Let $M\in\UU^{n,k}$ be fixed.
For all $X\in\UU^{n,k}$ such that $X^HM$ is nonsingular,
the alignment $\ALIGNM{X}$ in~\eqref{eq:alignment}
is uniquely defined and given by~\eqref{eq:alignX}.
It is invariant under right unitary transformations,
i.e., $\ALIGNM{XQ}=\ALIGNM{X}$ for all $Q\in\UU^{k,k}$.
Moreover, the map $X\mapsto\ALIGNM{X}$ is well defined
and differentiable in a neighborhood of $X$ on $\UU^{n,k}$,
with respect to the real and imaginary parts of its entries.
\end{lemma}

Observe that basis alignment induces a particular
ui-parameterization according to~\Cref{def:uiparam}.  
Indeed, by the unitary
invariance of $H(X)$, we obtain the representation~\eqref{eq:parah}
using alignment as
\begin{equation}\label{eq:alignasuip}
	H(X) \equiv H(\ALIGNM{X}) =
	\widehat H(\rho(X)),
\end{equation}
where $\rho(X):=\vect(\ALIGNM{X})\in\CC^{nk}$, and $\widehat H(\rho)$ is
defined by applying $H$ to the matricization of the vector $\rho$.  
By~\Cref{lem:alignment_smooth},
the parameterization $\rho(X)$ is well-defined and unitarily invariant 
whenever $\det(X^HM)\neq 0$.


\subsubsection{SCF with Alignment}
Applying the alignment~\eqref{eq:alignment} at each
SCF step selects a representative basis for the computed
eigenspace and resolves the ambiguity in the basis choice.
This gives the following aligned SCF iteration: 
starting from $\widehat X_0=\ALIGNM{X_0}$, compute
\begin{equation}\label{eq:alignedscf}
    \widehat X_{i+1}=\ALIGNM{X_{i+1}}
    \quad\text{with}\quad
    H(\widehat X_i)X_{i+1}=X_{i+1}\Lambda_{i+1},
\end{equation}
for $i=0,1,\dots$, where $X_{i+1}\in\UU^{n,k}$
is an eigenbasis for the $k$ smallest eigenvalues of $H(\widehat X_i)$.
In fixed-point form, this becomes
\begin{equation}
    \widehat X_{i+1}=g(\widehat X_i),
\end{equation}
where
\begin{equation}\label{eq:fpalignedscf}
    g(X):=\left\{
    \ALIGNM{Y}:Y\in\UU^{n,k}
	\text{is an eigenbasis of the $k$ smallest eigenvalues
		of $H(X)$}
    \right\}.
\end{equation}
The following lemma establishes the well-definedness
and local differentiability of $g$ under suitable conditions.
The result essentially follows from~\Cref{lem:param_scf_smooth}
applied to the alignment-induced
ui-parameterization~\eqref{eq:alignasuip}, together
with the regularity of the alignment map
in~\Cref{lem:alignment_smooth}.

\begin{lemma}\label{lem:aligned_scf_smooth}
Let $M\in\UU^{n,k}$ be fixed.
Suppose that, at $X\in\UU^{n,k}$, the $k$ smallest
eigenvalues of $H(X)$ are separated from the rest
of the spectrum, i.e.,
$\lambda_k(H(X))<\lambda_{k+1}(H(X))$,
and let $Y\in\UU^{n,k}$ be a corresponding eigenbasis.
Assume further that $Y^HM$ is nonsingular.
Then,
\begin{enumerate}[(i)]
\item\label{lem:aligned_scf_smooth:i}
the fixed-point mapping $g$ in~\eqref{eq:fpalignedscf}
is single-valued at $X$ (i.e., $g(X)$ is uniquely defined);
\item\label{lem:aligned_scf_smooth:ii}
the mapping $g$ is continuous (and differentiable) at $X$
if $H$ is continuous (and differentiable) at $X$,
in terms of the real and imaginary parts of its argument.
\end{enumerate}
\end{lemma}

\begin{proof}
The proof follows the same argument as that
of~\Cref{lem:param_scf_smooth}, with $\rho$ and
$\widehat H$ replaced by the alignment map
$\ALIGNM{\cdot}$ and $H$, respectively.
The required unitary invariance and local regularity
of the alignment map follow from~\Cref{lem:alignment_smooth},
since $Y^HM$ is nonsingular.
\end{proof}

\subsubsection{Anderson Acceleration with Aligned Bases}
Applying \AAm to the aligned SCF
mapping~\eqref{eq:fpalignedscf} gives
\Cref{alg:alignedaa}.
This is essentially the parameterized scheme in~\Cref{alg:aarho}
with alignment as the parameterization $\rho(X)=\ALIGNM{X}\in\CC^{n,k}$, 
and is presented directly with basis matrices $\ALIGNM{X}$ rather than their vectorizations. 

\begin{algorithm}[!htbp]
	\caption{AA($m$) for aligned SCF in difference form}
    \label{alg:alignedaa}
	\begin{algorithmic}[1]
        \REQUIRE Initial $X_0\in\UU^{n,k}$, 
		matrix function $H$, 
		reference matrix $M\in\UU^{n,k}$, 
		$m\in\mathbb{N}$, $\tau>0$
        \ENSURE approximate eigenbasis $X^*\in\UU^{n,k}$ of the NEPv~\eqref{eq:nepv}
        \STATE Compute the SVD $X_0^HM=U_0\Sigma_0V_0^H$ and set
        $\widehat X_0=X_0U_0V_0^H$
        \FOR{$i=0,1,2,\dots$}
        \STATE Evaluate $H_{i}=H(\widehat X_{i})$
        \STATE Solve $H_iX_{i+1}=X_{i+1}\Lambda_{i+1}$ for the $k$ smallest eigenvalues in $\Lambda_{i+1}$
        \STATE Compute the SVD $X_{i+1}^HM=U_{i+1}\Sigma_{i+1}V_{i+1}^H$ and set
        $Y_{i}=X_{i+1}U_{i+1}V_{i+1}^H$ and $R_i=Y_{i}-\widehat X_i$
        \STATE {\bf if} {$\|R_i\|_F\leq\tau$}, {\bf then return} $X^*=Y_{i}$
		\STATE Set $m_i=\min\{m,i\}$
        \IF{$m_i=0$}
            \STATE Set $\widehat X_{i+1}=Y_{i}$
        \ELSE
            \STATE Form $\Delta \widehat\X_i=[\Delta \widehat X_{i-m_i},\dots,\Delta \widehat X_{i-1}]$ and
            $\Delta \mathcal{R}_i=[\Delta R_{i-m_i},\dots,\Delta R_{i-1}]$, where
            $\Delta \widehat X_\ell=\widehat X_{\ell+1}-\widehat X_\ell$ and
            $\Delta R_\ell=R_{\ell+1}-R_\ell$
            \STATE Compute $\beta^{(i)}\in\RR^{m_i}$ by solving
            $\beta^{(i)}=\argmin_{\beta\in\RR^{m_i}}
            \|R_i-\Delta \mathcal{R}_i\odot\beta\|_F$
    		\label{ls:alg:alignedaa}
			\STATE Set $\widehat X_{i+1}=\widetilde X_{i+1} +\widetilde \R_{i+1}$,
			where 
			$\widetilde X_{i+1} := \widehat X_i - \Delta \widehat\X_i\odot\beta^{(i)}$
			and 
			$\widetilde R_{i+1} := R_i - \Delta R_i\odot\beta^{(i)}$
			\label{xi:alg:alignedaa}
        \ENDIF
        \ENDFOR
	\end{algorithmic}
\end{algorithm}

In~\Cref{alg:alignedaa},  the differences $\Delta \widehat \X_i$ and $\Delta \R_i$ are block matrices of the form 
$\Delta \equiv [\Delta_1,\dots,\Delta_m]$ with each block column
$\Delta_j\in\CC^{n,k}$,
and the block linear combination is $\Delta \odot\beta :=
\sum_{j=1}^m \beta_j \Delta_j$.
To solve the least-squares problem in line~\ref{ls:alg:alignedaa},
one can reformulate it in the standard matrix form
\begin{equation}\label{eq:lsblock}
    \min_{\beta\in\RR^{m_i}}
    \left\|
    \begin{bmatrix}
    \operatorname{Re}(\vect(R_i))\\
    \operatorname{Im}(\vect(R_i))
    \end{bmatrix}
    -
    \begin{bmatrix}
    \operatorname{Re}(\mathcal V_i)\\
    \operatorname{Im}(\mathcal V_i)
    \end{bmatrix}
    \beta
    \right\|_2,
\end{equation}
where
$\mathcal V_i=[\vect(\Delta R_{i-m_i}),\dots,\vect(\Delta
R_{i-1})]\in\CC^{nk,m_i}$,
and the real and imaginary parts of the complex variables have been
separated.
Alternatively, to avoid explicitly forming and operating on vectorized
matrices, we can regard $\min_\beta \|R - \sum_{j=1}^m \Delta_j\beta_j\|_F$ 
directly as a least-squares problem in the vector space $\CC^{n,k}$,
equipped with the Frobenius inner product 
$\langle A ,B\rangle=\Re(\TR(A^HB))$.


\subsection{Riemannian Anderson Acceleration} \label{sec:aamanifold}
We now consider another technique to address the unitary invariance 
of the SCF~\eqref{eq:scf}
by regarding it as a fixed-point iteration on the Grassmann manifold
$\Gr(n,k)$.
In this formulation, each SCF iterate is viewed as a $k$-dimensional subspace,
so the ambiguity caused by the choice of orthogonal basis is removed at
the level of the iteration.
The Anderson mixing is then carried out, using the local coordinates of
the iterates, in a tangent space of Grassmannian,
which provides a local linearization of the manifold. 
This approach is intrinsic to unitary invariance of the NEPv and SCF, 
but also requires additional geometric tools and operations for its
derivation and analysis.
We describe the essential ingredients and a self-consistent treatment below;
related operations and Anderson mixing on matrix manifold 
can also be found in~\cite{AbsMS2008,LiB2025}.

\subsubsection{Grassmannian Notation and Basics}
We briefly introduce the notation and concepts needed in our
derivation; see~\cite{AbsMS2008,EdeAS1998} for details.
\begin{itemize}
\setlength{\itemsep}{1pt}
\setlength{\parskip}{2pt}
\item \emph{Tangent space and metric.}
The tangent space $T_{\X}$ at $\X\in\Gr(n,k)$ is the real
vector space of velocities of smooth curves through $\X$.
The canonical (Euclidean) inner products on these tangent spaces define
the canonical Riemannian metric,
and $\|\cdot\|_F$ is the induced tangent norm,
which agrees with the Frobenius norm in the matrix
representations used below.

\item \emph{Geodesics and exponential/logarithmic maps.}
Geodesics generalize straight lines. 
The geodesic distance $\operatorname{dist}(\X,\Y)$ is the
length of a shortest geodesic connecting $\X$ and $\Y$.
The smooth exponential map
$\Exp_{\X}\colon T_{\X}\to\Gr(n,k)$ maps $\xi$ to the point
reached at time $1$ along the geodesic starting at $\X$
with initial velocity $\xi$.
Its smooth local inverse $\Log_{\X}\colon\Gr(n,k)\to T_{\X}$ 
maps a nearby point $\Y$ back to $\xi\in T_{\X}$.
These maps serve as local analogues of addition
(symbolically, $\X+\xi$) and subtraction
(symbolically, $\Y-\X$), respectively.
In particular, $\operatorname{dist}(\X,\Y)=\|\Log_{\X}(\Y)\|_F$,
whenever the logarithmic map is defined.

\item \emph{Parallel transport.}
For sufficiently close $\X$ and $\Y$, parallel transport
$\T_{\X}^{\Y}\colon T_{\X}\to T_{\Y}$
moves tangent vectors along the unique shortest geodesic
connecting them.
This map is linear and preserves inner products, allowing
vectors associated with different iterates to be combined
in a common tangent space.
\end{itemize}

\subsubsection{Fixed-Point Iteration on the Grassmannian}
Due to the unitary invariance of the NEPv~\eqref{eq:nepv},
each SCF iterate $X_i$ in~\eqref{eq:scf} can be regarded as 
a $k$-dimensional subspace $\mbox{range}(X_i)$, i.e., 
a point in the Grassmannian.
Therefore,  SCF iteration~\eqref{eq:scf} can be interpreted as 
the fixed-point iteration
\begin{equation}\label{eq:fpg}
	\X_{i+1}=g(\X_i),
\end{equation}
where $\X_i\in \Gr(n,k)$ is a $k$-dimensional subspace of $\CC^{n}$,
and $g:\Gr(n,k)\to \Gr(n,k)$ is given by 
\begin{equation}\label{eq:scff}
g(\X) := \lbrace \Y: \text{$\Y\in\Gr(n,k)$ is an eigenspace for the
	$k$ smallest eigenvalues of $H(X)$}\rbrace,
\end{equation}
and $X\in\UU^{n,k}$ denotes an orthogonal basis of $\X$.
The fixed-point mapping $g$ is well defined 
and regular under the eigengap assumption of $H(X)$,
as precisely described in the following lemma.

\begin{lemma}\label{lem:grass_scf_smooth}
Let $\X\in\Gr(n,k)$, and let $X\in\UU^{n,k}$ be any orthonormal basis of
$\X$. Suppose that the $k$ smallest eigenvalues of $H(X)$ are separated
from the rest of the spectrum, i.e.,
$\lambda_k(H(X))<\lambda_{k+1}(H(X))$.
Then the fixed-point mapping $g$ in~\eqref{eq:scff} is single-valued at
$\X$. If, in addition, $H$ is continuous (and differentiable) in a
neighborhood of $X$, then $g$ is locally continuous (and differentiable)
at $\X$.\footnote{
	The continuity and differentiability on $\Gr(n,k)$ is understood in the standard manifold sense; see, e.g.,~\cite{AbsMS2008}.
} 
\end{lemma}

\begin{proof}
By representing the Grassmannian as 
$$\Gr(n,k)=\{P\in\CC^{n,n}\colon P=P^H,\, P^2=P,\, \mbox{rank}(P)=k\},$$
each $\X_i\in\Gr(n,k)$ can be identified with $P_i=X_iX_i^H$;
see, e.g.,~\cite{BenZA2024}.
Hence, the fixed-point mapping $g$ in~\eqref{eq:scff}
is obtained from~\eqref{eq:fprho} by using the ui-parameterization
$\rho(X):=XX^H \in\CC^{n,n}$ (up to vectorization),
with the associated map $\widehat H(P):=H(X(P))$.
Here, $X(P):=PM(M^HPM)^{-1/2}$, where $M\in\CC^{n,k}$ is a fixed
rank-$k$ matrix such that $M^HPM$ is nonsingular.
The function $X(P)$ selects an orthonormal basis of the range of $P$,
and it is locally continuous and differentiable in $P$ whenever
$M^HPM$ is nonsingular.
Therefore, the result follows directly from~\Cref{lem:param_scf_smooth}.
\end{proof}

\subsubsection{Anderson Acceleration on the Grassmannian}
We next describe the AA for the fixed-point
iteration~\eqref{eq:fpg} on the Grassmannian.
Since the iteration subspaces $\{\X_i\}$
cannot be directly added or subtracted in the Euclidean sense,
the Riemannian Anderson acceleration performs the mixing
operations on $\{\X_i\}$ using their local coordinates
in a tangent space (a Euclidean vector space) of the manifold.

The basic idea is to represent nearby subspaces by tangent vectors at a
reference point, allowing the linear combinations and residual
evaluation required by Anderson acceleration to be defined in the usual
Euclidean manner.
Particularly, using $\X_{i}$ as a reference point, 
we can formulate a mixing scheme analogous to the \AAm in~\eqref{eq:aa} as
\begin{subequations}\label{eq:raa0}
\begin{equation}\label{eq:raa0step}
	\X_{i+1} = 
	\Exp_{\Y_i}
	\left( \sum_{\ell=0}^{m-1} \alpha_\ell \Log_{\Y_i} (\Y_{i-\ell}) \right),
\end{equation}
where $\Y_\ell=g(\X_\ell)$, and the coefficients $\{\alpha_\ell\}$ are obtained by solving
\begin{equation}\label{eq:raa0coeff}
\min_{\sum_{\ell=0}^{m-1}\alpha_\ell =1} \left\|\sum_{\ell=0}^{m-1}
\alpha_\ell \cdot \T_{\X_{i-\ell}}^{\X_i}r(\X_{i-\ell})\right\|_F
\qquad \text{with} \quad r(\X):=\Log_{\X}(g(\X)).
\end{equation}
\end{subequations}
Recall that the exponential map $\Exp_{\X}(\xi)\in\Gr(n,k)$,
with $\X\in\Gr(n,k)$ and $\xi\in T_{\X}$, is 
a Riemannian analogue of the addition (informally, $\X+\xi$),
the logarithm map $\Log_{\X}(\Y)\in T_{\X}$, with $\X,\Y\in\Gr(n,k)$,
a Riemannian analogue of the difference (informally, $\Y-\X$),
and $\T_{\X}^{\Y}$ denotes the parallel transport along the geodesic
from $\X$ to $\Y$.
Consequently, the mixing step~\eqref{eq:raa0step}
is a direct Riemannian analogue of 
$x_{i+1} = g(x_i) + \sum_{\ell=0}^{m-1}\alpha_\ell (g(x_{i-\ell})-g(x_i))$ from~\eqref{eq:aastep},
while the mapping $r(\X)$ in~\eqref{eq:raa0coeff},
measuring the differences between $g(\X)$ and $\X$ on the tangent space $T_{\X}$,
is a Riemannian analogue of the residual $g(x)-x$
in~\eqref{eq:aacoeff}.

The following result establishes the Riemannian analogue of the expansion
relation~\eqref{eq:resd} for standard Anderson acceleration,
and it partially justifies the coefficient choice in~\eqref{eq:raa0coeff}.
It shows that, after transporting residuals to the common tangent space
$T_{\X_i}$, the residual of the mixed iterate is approximated by the
same affine combination used in~\eqref{eq:raa0coeff}.

\begin{lemma}\label{lem:grass_residual_expansion}
Let $r(\X)=\Log_{\X}(g(\X))$ be locally differentiable in a neighborhood
of $\X_i\in\Gr(n,k)$.
Let $\X_{i-m+1},\dots,\X_{i-1}\in\Gr(n,k)$ be sufficiently close to
$\X_i$, and let the coefficients $\{\alpha_\ell\}_{\ell=0}^{m-1}$
satisfy $\sum_{\ell=0}^{m-1}\alpha_\ell=1$ and be uniformly bounded
(independent of $\X_\ell$).
Then the mixed iterate
\begin{equation}\label{eq:raa0mixedstep}
	\widetilde \X_{i+1} =
	\Exp_{\X_i}
	\left( \sum_{\ell=0}^{m-1} \alpha_\ell \Log_{\X_i}(\X_{i-\ell})
	\right)
\end{equation}
satisfies
\begin{equation}\label{eq:exprtx}
\T_{\widetilde \X_{i+1}}^{\X_i} r(\widetilde \X_{i+1}) =
\sum_{\ell=0}^{m-1}
\alpha_\ell \T_{\X_{i-\ell}}^{\X_i}r(\X_{i-\ell}) +
\lo\left(\max_{\ell=1,\dots,m-1}\mathrm{dist}(\X_i,\X_{i-\ell})\right).
\end{equation}
Here, $\widetilde \X_{i+1}$ differs from the update
$\X_{i+1}$ in~\eqref{eq:raa0step} by replacing
$\Y_\ell = g(\X_{\ell})$ with $\X_{\ell}$.
\end{lemma}

\begin{proof}
Define the local representation of the residual vector field at $\X_i$ by
\[
    \widehat r_i(\xi)
    :=
    \T_{\Exp_{\X_i}(\xi)}^{\X_i}
    r(\Exp_{\X_i}(\xi)),
    \qquad \xi\in T_{\X_i}.
\]
Since $r$ is locally differentiable, so is $\widehat r_i$, and hence
$\widehat r_i(\xi) = \widehat r_i(0)+D\widehat r_i(0)[\xi]+\lo(\|\xi\|_F), $
where $D\widehat r_i(0)$ denotes the differential of $\widehat r_i$ at
the origin of $T_{\X_i}$.
Applying this expansion to
$\xi_\alpha=\sum_{\ell=0}^{m-1}\alpha_\ell\xi_\ell$ and to each
$\xi_\ell=\Log_{\X_i}(\X_{i-\ell})$, and using
$\sum_{\ell=0}^{m-1}\alpha_\ell=1$ together with the uniform boundedness
of $\{\alpha_\ell\}$, gives
\[
    \widehat r_i(\xi_\alpha)
    =
    \sum_{\ell=0}^{m-1}\alpha_\ell \widehat r_i(\xi_\ell)
    +
    \lo\left(\max_{\ell=1,\dots,m-1}\|\xi_\ell\|_F\right).
\]
Hence,~\eqref{eq:exprtx} follows from
$\|\xi_\ell\|_F=\mbox{dist}(\X_i,\X_{i-\ell})$ and
$\widehat r_i(\xi_\ell)=\T_{\X_{i-\ell}}^{\X_i}r(\X_{i-\ell})$.
\end{proof} 

The extension of AA to the Grassmann manifold through
tangent-space operations, as in~\eqref{eq:raa0}, is conceptually
straightforward, although technically involved.
Similar ideas are explored in the general framework for
Anderson mixing on the Riemannian manifold developed in~\cite{LiB2025}, 
primarily from a root-finding perspective.
In contrast, our derivation of~\eqref{eq:raa0} and its
justification in~\Cref{lem:grass_residual_expansion} are based
directly on fixed-point iterations, 
following closely the rationale of classical Anderson acceleration 
using the local expansion~\eqref{eq:resd}.

\subsubsection{Difference Form and Implementation}
For numerical implementation, it is more convenient to use 
a difference form of \AAm on the Grassmannian, 
in direct analogy with~\eqref{eq:aa:inc}. 
Following the same derivation as for~\eqref{eq:raa0},
we define the Riemannian analogues of the differences 
for the sequences $\{\X_\ell\}$ and $\{r(\X_\ell)\}$
as 
\begin{equation}\label{eq:diff}
\Delta \X_\ell^{(i)} := \T_{\X_{\ell+1}}^{\X_{i}} \left( \Log_{\X_{\ell+1}}(\X_{\ell})\right)
\quad\text{and}\quad
\Delta r_\ell^{(i)} := 
\T_{\X_{\ell+1}}^{\X_{i}} \left(
\T_{\X_{\ell}}^{\X_{\ell+1}} \left(r(\X_\ell) \right)
- r(\X_{\ell+1})\right),
\end{equation}
where the parallel transports place all difference vectors in a common
tangent space $T_{\X_i}$. 
We can then recast the 
difference form~\eqref{eq:aa:inc} of \AAm
into the Riemannian version\footnote{
The difference form \eqref{eq:aa:inc} is equivalent to 
$x_{i+1} = x_i+r(x_i)+\sum_{\ell=0}^m
\beta_\ell\left[(x_{\ell} - x_{\ell+1})+ (r(x_{\ell})-r(x_{\ell+1})) \right]$.
}
\begin{subequations}\label{eq:grassaa:inc}
\begin{equation}\label{eq:grassaa:incstep}
    \X_{i+1}
    =
    \Exp_{\X_{i}}\left(r(\X_i)+\sum_{\ell=i-m}^{i-1} \beta_\ell
	\cdot (\Delta \X_\ell^{(i)} + \Delta r_\ell^{(i)}) \right),
\end{equation}
where the coefficients $\{\beta_\ell\}$ are obtained from the least squares
\begin{equation}\label{eq:grassaa:inccoeff}
    \min_{\beta}
    \left\|
    r(\X_i)+\sum_{\ell=i-m}^{i-1}
    \beta_\ell \Delta r_\ell^{(i)}
    \right\|_F.
\end{equation}
\end{subequations}
We caution that, unlike the Euclidean difference form~\eqref{eq:aa:inc},
which is equivalent to the basic form~\eqref{eq:aa},
the formulation~\eqref{eq:grassaa:inc} is not exactly equivalent to~\eqref{eq:raa0},
since different transports are applied to the individual difference vectors.
Nevertheless, when all iterates are sufficiently close, the two
formulations agree up to higher-order terms in their distances,
since the composition and direct parallel transports differ 
only by second-order terms due to manifold curvature; see
also~\cite{LiB2025}.

\begin{algorithm}[!htbp]
\caption{Riemannian AA($m$) for SCF on the Grassmannian in difference form}
\label{alg:grassaa_diff}
\begin{algorithmic}[1]
\REQUIRE Initial $X_0\in\UU^{n,k}$, matrix function $H$,
$m\in\mathbb{N}$, $\tau>0$
\ENSURE approximate eigenbasis $X^*\in\UU^{n,k}$ of the NEPv~\eqref{eq:nepv}
\FOR{$i=0,1,2,\dots$}
\STATE Evaluate $H_i=H(X_i)$
\STATE Solve $H_iY_i=Y_i\Lambda_i$ for the $k$ smallest eigenvalues in $\Lambda_i$
\STATE Compute the Riemannian residual $r_i=\Log_{X_i}(Y_i)\in\CC^{n,k}$
\STATE  {\bf if} {$\|r_i\|_F\leq\tau$}, {\bf then return} $X^*=Y_i$
\STATE Set $m_i=\min\{m,i\}$
\IF{$m_i=0$}
    \STATE Set $X_{i+1}=Y_i$ and compute $\T_{X_i}^{X_{i+1}}$ by SVD $r_i=U_i\Theta_i V_i^H$ (see \eqref{eq:paratransport}) 
\ELSE
	\STATE 
	Compute $\Delta X_{i-1}^{(i)}=\Log_{X_i}(X_{i-1})\in\CC^{n,k}$ and 
	$\Delta r_{i-1}^{(i)}= \T_{X_{i-1}}^{X_{i}} r_{i-1}-r_i \in\CC^{n,k}$

	\STATE 
	Transport 
	$\Delta X_{i-\ell}^{(i)}=\T_{X_{i-1}}^{X_{i}} \Delta X_{i-\ell}^{(i-1)} \in\CC^{n,k}$ and 
	$\Delta r_{i-\ell}^{(i)}= \T_{X_{i-1}}^{X_{i}} \Delta r_{i-\ell}^{(i-1)} \in\CC^{n,k}$,
	for $\ell=2,\dots,m_i$
	\label{alg:grassaa_diff:trans}

    \STATE Form
    $\Delta\mathcal X_i=[\Delta X_{i-m_i}^{(i)},\dots,\Delta X_{i-1}^{(i)}]\in\CC^{n,km_i}$
    and
    $\Delta\mathcal R_i=[\Delta r_{i-m_i}^{(i)},\dots,\Delta r_{i-1}^{(i)}]\in\CC^{n,km_i}$
    \STATE Compute $\beta^{(i)}\in\RR^{m_i}$ by solving
        $\beta^{(i)} = \argmin_{\beta\in\RR^{m_i}} \|r_i+\Delta\mathcal R_i\odot\beta\|_F $
	\label{alg:grassaa_diff:ls}
    \STATE Set up $\widetilde{r}_i=r_i+ (\Delta\mathcal R_i+\Delta\mathcal X_i)\odot\beta^{(i)}$ 
    \STATE Compute $X_{i+1}=\Exp_{X_i}(\widetilde{r}_i)\in\UU^{n,k}$ and $\T_{X_i}^{X_{i+1}}$,
	using SVD $\widetilde{r}_i=U_i\Theta_iV_i^H$ (see \eqref{eq:expmapping}, \eqref{eq:paratransport}) 
\ENDIF
\ENDFOR
\end{algorithmic}
\end{algorithm}
Algorithm~\ref{alg:grassaa_diff} outlines the Riemannian \AAm.
We highlight three implementation details.
\begin{enumerate}[(i)]
\item 
As is standard in computations on the Grassmannian
(see, e.g.,~\cite{Absil:2004,AbsMS2008,BenZA2024}),
we use Stiefel representatives $X\in\UU^{n,k}$ for
$\X=\mbox{range}(X)\in\Gr(n,k)$ and represent tangent vectors 
by their unique horizontal lifts $\Delta\in\CC^{n,k}$ satisfying $X^H\Delta=0$.
In particular, with these representations, the exponential map takes the closed-form
\begin{equation}\label{eq:expmapping}
\Exp_X(\Delta)=(XV\cos\Theta+U\sin\Theta)V^H\in\UU^{n,k}
\quad \text{with SVD $\Delta=U\Theta V^H$},
\end{equation}
the logarithm map, for nonsingular $X^HY$, takes the form 
\begin{equation}\label{eq:logmapping}
\Log_X(Y)=U\arctan(\Sigma)V^H\in\CC^{n,k}
\quad \text{with SVD $(I-XX^H)Y(X^HY)^{-1}=U\Sigma V^H$},
\end{equation}
and parallel transport takes the form 
\begin{equation}\label{eq:paratransport}
\T_X^Y(\Delta)=
\bigl[I+(U(\cos\Theta-I)-XV\sin\Theta)U^H\bigr]\Delta 
\in\CC^{n,k}
\quad \text{with SVD $\Log_{X}(Y) = U\Theta V^H$}.
\end{equation}
We refer to~\cite{BenZA2024} for details and 
the treatment of logarithms when $X^HY$ is singular.

\item
For computational efficiency, line~\ref{alg:grassaa_diff:trans}
reuses the differences
$\Delta X_\ell^{(i-1)},\Delta r_\ell^{(i-1)}
\in T_{\X_{i-1}}$
from iteration $i-1$ by parallel transport to $T_{\X_i}$.
This amounts to approximating the direct transport in~\eqref{eq:diff}
by the composition
$\T_{\X_\ell}^{\X_i}
\approx
\T_{\X_{i-1}}^{\X_i}\circ\T_{\X_\ell}^{\X_{i-1}}$
-- For nearby iterates, these transports agree to first order,
and hence would not affect the local behavior of the algorithm.
Since the same transport $\T_{\X_{i-1}}^{\X_i}$ is applied to all
stored differences, 
the computation via~\eqref{eq:paratransport} requires only matrix multiplications
after computing $\Log_{X_{i-1}}(X_i)$ and its SVD once.

\item
The least-squares problem in line~\ref{alg:grassaa_diff:ls}
involves a potentially complex block matrix $\Delta\R_i$ but real
coefficients $\beta^{(i)}$.
It can be solved using the same techniques discussed
for~\Cref{alg:alignedaa}, either by converting it to the real
standard form~\eqref{eq:lsblock} or by computing a block QR decomposition.
\end{enumerate}


\section{Basis-Aligned Newton Acceleration} \label{sec:newton}
In this section, we present a Newton-type acceleration scheme.
Newton's method for NEPv has been considered in~\cite{GaoYM2007,Werner2025}.
One major issue is that the unitary invariance of NEPv leads to a singular
derivative operator at a solution as noted
in~\cite{Werner2025}.
This singularity may cause numerical difficulties when solving
the Newton correction equation.
We address this issue through a correction scheme 
in which the updated basis aligns with the current basis $X$. 
This alignment is incorporated into the correction equation through the condition
$X^H\Delta X=0$ to remove the ambiguity associated with
basis rotations and yields a nonsingular correction equation under mild assumptions. 
We then combine the Newton correction with SCF to accelerate
convergence towards desired eigenvalues and discuss its implementations.
These developments complement previous work and 
provide further theoretical and algorithmic support 
for applying Newton's method to NEPv.

\subsection{Newton Correction for Aligned Basis of NEPv}\label{sec:grassntderive}
Throughout this section, we assume that $H$ is Fr\'echet differentiable with a
locally Lipschitz continuous derivative in an open neighborhood of $X\in\UU^{n,k}$, 
that is, there exists a linear derivative
operator $\DD{H}(X)[\cdot]\colon\CC^{n,k}\to\CC^{n,n}$ that is locally
Lipschitz in $X$ and satisfies 
\begin{equation}\label{eq:expandH}
    H(X+E)=H(X)+\DD{H}(X)[E]+\mathcal O(\|E\|^2).
\end{equation}
Here, the complex matrix space $\CC^{n,k}$ is viewed as a vector space
over $\RR$ (i.e., treat real and imaginary parts as individual variables),
so the linearity really means $\RR$-linearity:
\begin{equation}\label{eq:rlinear}
    \DD{H}(X)[\alpha_1E_1+\alpha_2 E_2]
    =\alpha_1\DD{H}(X)[E_1]+\alpha_2\DD{H}(X)[E_2],
    \qquad E_1,E_2\in\CC^{n,k},\quad \alpha_1,\alpha_2\in\RR.
\end{equation}
We assume the derivative operator $\DD{H}(X)[\cdot]$ is available
in the sense that its action $\DD{H}(X)[E]$ can be evaluated 
for any given direction $E$.
Observe that this action coincides with the directional derivative 
\begin{equation}\label{eq:directional}
    \DD{H}(X)[E]
    =\left.\frac{d}{dt}H(X+tE)\right|_{t=0},
    \qquad t\in\RR,
\end{equation}
which provides a convenient expression for both derivation and
computation.

\subsubsection{Aligned Basis Correction}
We consider the following correction scheme for the NEPv.
Given a current approximate solution
$(X,\Lambda)\in\UU^{n,k}\times\CC^{k,k}$ of the NEPv~\eqref{eq:nepv},
we seek the next approximation in the form 
\begin{subequations}\label{eq:ntcorrections}
\begin{equation}
    X_{\mathrm{next}}=\mathrm{orth}(X+\dX),
    \qquad
    \Lambda_{\mathrm{next}}=\Lambda+\dL,
\end{equation}
where the corrections
$(\dX,\dL)\in\CC^{n,k}\times\CC^{k,k}$
satisfy
\begin{equation}\label{eq:orthxdx}
    X^H\dX=0,
\end{equation}
\end{subequations}
and $\mathrm{orth}(\cdot)$ denotes the orthonormal polar factor of a
full-column-rank matrix $A$, given by
\begin{equation}\label{eq:polarorth}
	\mathrm{orth}(A):=A(A^HA)^{-1/2}.
\end{equation}
The motivation behind the update~\eqref{eq:ntcorrections} is to seek 
an orthonormal basis $X_{\mathrm{next}}$ aligned with the current basis $X$, 
satisfying $X^HX_{\mathrm{next}}\succ0$; see, e.g.,~\Cref{sec:aaaligned}.
This removes the ambiguity in choosing a basis for the updated
subspace. Moreover, each such aligned basis corresponds to a unique
correction $E$ satisfying $X^HE=0$.

\begin{lemma}\label{prop:aligned_correction}
Let $X,Y\in\UU^{n,k}$. Then $X^HY\succ0$ if and only if
there exists a unique $E\in\CC^{n,k}$ satisfying $X^HE=0$
such that $Y=\mathrm{orth}(X+E)$.
\end{lemma}
\begin{proof}
If $Y=\mathrm{orth}(X+E)$ with $X^HE=0$, then
\begin{equation}\label{eq:propal}
Y\equiv (X+E)(I+E^HE)^{-1/2}\quad\text{and}\quad X^HY\equiv (I+E^HE)^{-1/2}\succ0.
\end{equation}
Conversely, let $C:=X^HY\succ0$ and $E:=YC^{-1}-X$.
Then, $X^HE=0$ and 
$\mathrm{orth}(X+E)\equiv\mathrm{orth}(YC^{-1})=Y$.
Uniqueness follows as any such $E$ must satisfy~\eqref{eq:propal}, 
so $Y=(X+E)X^HY$, i.e., $E=Y(X^HY)^{-1}-X$.
\end{proof}

\subsubsection{Newton Correction Equation}
To solve the NEPv,
an {\em exact} correction $(\dX,\dL)$
in~\eqref{eq:ntcorrections} would satisfy the nonlinear systems of equations 
\begin{equation}\label{eq:correctexact}
    \begin{bmatrix}
        H(X_{\mathrm{next}})X_{\mathrm{next}}
        -X_{\mathrm{next}}\Lambda_{\mathrm{next}}\\
        -X^H\dX
    \end{bmatrix}=0,
\end{equation}
where the upper block represents the residual of the NEPv and the 
lower comes from the constraints~\eqref{eq:orthxdx} 
(we negate the equation for convenience of later discussions).
Rather than solving this equation exactly, we can obtain
an approximate solution using a Newton correction by linearization.
Since 
\[
    X_{\mathrm{next}}
    =(X+\dX)(I+\dX^H\dX)^{-1/2}
    =X+\dX+\mathcal O(\|\dX\|_F^2),
\]
it follows from the local expansion~\eqref{eq:expandH} of $H(X_{\mathrm{next}})$ at $X$ that 
\[
\begin{aligned}
    H(X_{\mathrm{next}})X_{\mathrm{next}}
    -X_{\mathrm{next}}\Lambda_{\mathrm{next}}
    = H(X)X&-X\Lambda
       +\DD{H}(X)[\dX]X+H(X)\dX
       -\dX\Lambda-X\dL\\
       &+\mathcal O(\|\dX\|^2+\|\dL\|^2).
\end{aligned}
\]
Plugging into~\eqref{eq:correctexact}, and dropping the higher-order terms,
yields the linear system in $(\dX,\dL)$: 
\begin{subequations}\label{eq:grassnt}
\begin{equation}\label{eq:grassnt:system}
    \LL{(X,\Lambda)}[\dX,\dL]
    = -R{(X,\Lambda)},
\end{equation}
where
\begin{equation}\label{eq:grassnt:operator}
    \LL{(X,\Lambda)}[\dX,\dL]:
    =
    \begin{bmatrix}
        \DD{H}(X)[\dX]X+H(X)\dX
        -\dX\Lambda-X\dL\\
        -X^H\dX
    \end{bmatrix},\quad
    R{(X,\Lambda)}
    :=\begin{bmatrix}H(X)X-X\Lambda\\0\end{bmatrix}.
\end{equation}
\end{subequations}
The correction operator
$\LL(X,\Lambda)[\cdot]\colon
\CC^{n,k}\times\CC^{k,k}\to\CC^{n,k}\times\CC^{k,k}$
is $\RR$-linear and is nonsingular at a solution under mild
assumptions, as discussed in~\Cref{sec:nttheory}.

Let us highlight three major benefits of the basis-aligned Newton correction equation~\eqref{eq:grassnt} 
over the standard Newton approach in~\cite{GaoYM2007,Werner2025}
(also described in~\Cref{rmk:stdnt}).
First, basis alignment yields a correction operator that is
nonsingular under suitable assumptions, as proved
in~\Cref{sec:nttheory}, simplifying both the theoretical
analysis and the numerical solution.
Second, for NEPv arising from certain optimization problems,
the correction operator is self-adjoint, allowing efficient
iterative solvers with short recurrences (see~\Cref{rmk:symnt});
This property generally does not hold for the standard
Newton approach. 
Finally, the elementary derivation and subsequent analysis
apply to both real and complex NEPv by treating derivatives
and correction operators as $\RR$-linear maps.
This preserves the original matrix form of the correction
equation, without explicitly separating complex variables
into real and imaginary parts, as used in related works~\cite{JarKM2014,JarU2022}.

\begin{remark}[Standard Newton method]
\label{rmk:stdnt}
\rm
The correction equation~\eqref{eq:grassnt} is intrinsically
different from that obtained by applying the standard Newton
method directly to the root-finding formulation of the NEPv:
\begin{equation}\label{eq:nepvroot}
    F(X,\Lambda)=0
    \quad\text{with}\quad
    F(X,\Lambda):=
    \begin{bmatrix}
        H(X)X-X\Lambda\\
        I_k-X^HX
    \end{bmatrix},
\end{equation}
as considered in the real setting
in~\cite{GaoYM2007,Werner2025}.
Starting from $X\in\UU^{n,k}$ and $\Lambda\in\CC^{k,k}$,
a standard Newton step takes the form
$
\begin{bmatrix}X_{\mathrm{next}}\\\Lambda_{\mathrm{next}}\end{bmatrix}
    =
    \begin{bmatrix}X\\\Lambda\end{bmatrix}
    +\begin{bmatrix}\dX\\\dL\end{bmatrix},
$
where the corrections $\dX$ and $\dL$ satisfy
\begin{subequations}\label{eq:stnt}
\begin{equation}\label{eq:stnt:system}
    \widehat\LL(X,\Lambda)[\dX,\dL]
    =-R(X,\Lambda),
\end{equation}
with the $\RR$-linear operator
$\widehat\LL(X,\Lambda)[\cdot]\colon
\CC^{n,k}\times\CC^{k,k}\to\CC^{n,k}\times\CC^{k,k}$
denoting the derivatives of $F$ as given by
\begin{equation}\label{eq:stnt:operator}
    \widehat\LL(X,\Lambda)[\dX,\dL]
    :=
    \begin{bmatrix}
        \DD{H}(X)[\dX]X+H(X)\dX
        -\dX\Lambda-X\dL\\
        -(X^H\dX+\dX^HX)
    \end{bmatrix}.
\end{equation}
\end{subequations}
Here we used $F(X,\Lambda)=R(X,\Lambda)$ since $X^HX=I_k$.
Compared with our correction equation~\eqref{eq:grassnt},
the standard Newton equation~\eqref{eq:stnt} uses
$X^H\dX+\dX^HX$ in the lower block instead of $X^H\dX$.
This seemingly minor difference has an important consequence:
the standard Newton operator $\widehat\LL(X,\Lambda)$ is
singular at a solution (see, e.g.,~\cite{Werner2025}),
whereas our operator $\LL(X,\Lambda)$ is nonsingular
under mild assumptions.
The singularity of the standard Newton operator
$\widehat\LL(X,\Lambda)$ arises from unitary
invariance of the NEPv: if $(X,\Lambda)$ is a solution
of~\eqref{eq:nepvroot}, then $(XQ,Q^H\Lambda Q)$ 
remains a solution for every unitary $Q$, so the solution is not
locally unique and the derivative operator of $F$ must be singular.
\end{remark}

\begin{remark}[Self-adjointness of the correction operator]
\label{rmk:symnt}
\rm
For a useful property of the basis-aligned Newton method, 
the correction equation~\eqref{eq:grassnt} involves
a self-adjoint $\RR$-linear operator $\LL(X,\Lambda)[\cdot]$,
if $H(X)X=\nabla f(X)$
for some real-valued twice-differentiable function
$f\colon\CC^{n,k}\to\RR$ and $\Lambda=\Lambda^H$.
Here, $\nabla$ denotes the gradient 
(w.r.t. the inner product $\langle X,Y\rangle=\Re(\TR(X^HY))$). 
The condition $H(X)X=\nabla f(X)$ holds for many practical NEPv arising
from constrained optimization of the form
$\min f(X)$ s.t. $X^HX=I$, where the NEPv 
serves as the first-order optimality condition; see,
e.g.,~\cite{LuL2024} and references therein.
Specifically, under these assumptions,
\begin{equation}\label{eq:selfadjl}
    \Big\langle \LL(X,\Lambda)[E_1,D_1],\, (E_2,D_2) \Big\rangle
    =
    \Big\langle (E_1,D_1),\, \LL(X,\Lambda)[E_2,D_2] \Big\rangle
\end{equation}
for all $(E_i,D_i)\in\CC^{n,k}\times\CC^{k,k}$, where
$\langle(E_1,D_1),(E_2,D_2)\rangle :=\Re\TR(E_1^HE_2+D_1^HD_2)$.
The identity~\eqref{eq:selfadjl} follows by direct verification
using the self-adjointness property of the Hessian operator $\nabla^2 f(X)[E]=\DD{H}(X)[E]X+H(X)E$.
We omit the details of the derivation.
As a benefit, the self-adjointness allows the correction equation~\eqref{eq:grassnt}
to be solved by, e.g., an $\RR$-linear MINRES method, instead of GMRES,
through short recurrences that reduce storage and orthogonalization costs. 
We emphasize that this self-adjointness is not shared
by the standard Newton correction~\eqref{eq:stnt},
whose operator $\widehat\LL(X,\Lambda)$ is generally
not self-adjoint on $\CC^{n,k}\times\CC^{k,k}$.
\end{remark}

\subsection{Local Convergence Analysis}\label{sec:nttheory}
In this section, we establish the local quadratic convergence of the update
scheme~\eqref{eq:ntcorrections} with corrections determined
by~\eqref{eq:grassnt}.
Following the framework for local NEPv analysis
in~\cite{BaiLL2022}, consider a solution
$(X^*,\Lambda^*)\in\UU^{n,k}\times\CC^{k,k}$
of NEPv~\eqref{eq:nepv}, with the eigenvalue decomposition
\begin{equation}\label{eq:fulleig}
    H(X^*)=
    \begin{bmatrix}X^*&X_\perp^*\end{bmatrix}
    \begin{bmatrix}
        \Lambda^*&\\
        &\Lambda_\perp^*
    \end{bmatrix}
    \begin{bmatrix}X^*&X_\perp^*\end{bmatrix}^H,
\end{equation}
where the columns of $X^*$ are orthonormal eigenvectors
associated with the $k$ smallest eigenvalues of $H(X^*)$,
and the columns of $X_\perp^*\in\UU^{n,n-k}$ form an
orthonormal basis for the complementary invariant subspace.
The diagonal matrices $\Lambda^*\in\CC^{k,k}$ and
$\Lambda_\perp^*\in\CC^{n-k,n-k}$ contain the corresponding
eigenvalues. We impose the following assumptions. 

\begin{assumption}\label{ass:qspd}
\begin{enumerate}[(i)] 
    \item\label{ass:qspd:i}
        There is a positive eigengap:
        $\lambda_k(H(X^*))<\lambda_{k+1}(H(X^*))$.
    \item\label{ass:qspd:ii}
		$H$ is Fr\'echet differentiable with a locally Lipschitz 
		continuous derivative in a neighborhood of $X^*$.
    \item\label{ass:qspd:iii}
        The $\RR$-linear operator
        $\mathcal Q^*\colon\CC^{n-k,k}\to\CC^{n-k,k}$,
        defined by
        \begin{equation}\label{eq:rlinearq}
            \mathcal Q^*(Z)
            ={X_\perp^*}^H\DD{H}(X^*)[X_\perp^*Z]X^*
            +\Lambda_\perp^*Z-Z\Lambda^*,
        \end{equation}
        is self-adjoint and positive definite 
		w.r.t. the Frobenius inner product
        $\langle T,Z\rangle=\Re(\TR(T^HZ))$;
		i.e.,
        \begin{equation}\label{eq:qspd}
            \langle T,\mathcal Q^*(Z)\rangle
            =\langle\mathcal Q^*(T),Z\rangle
            \qquad
            \langle Z,\mathcal Q^*(Z)\rangle>0,
        \end{equation}
        for all $T,Z\in\CC^{n-k,k}$, with $Z\ne0$
        in the second inequality.
\end{enumerate}
\end{assumption}

Note that the eigengap and differentiability conditions in
\Cref{ass:qspd}~\eqref{ass:qspd:i} and~\eqref{ass:qspd:ii} 
are standard in the local analyses of NEPv and SCF iteration 
(e.g.,~\cite{BaiLL2022,CanKL2021,Liu:2014,Upadhyaya:2021a}).
The positive-definiteness condition in \Cref{ass:qspd}\,\eqref{ass:qspd:iii}
typically holds for practical NEPv arising in optimization 
on the Grassmann manifold,
where it is related to second-order sufficient conditions
for optimality and also implies
that $\mathcal Q^*$ is invertible; see~\cite{BaiLL2022}
for further discussion.
Under these assumptions, we obtain the following result.

\begin{theorem}[Local quadratic convergence]\label{thm:convergence}
Let $(X^*,\Lambda^*)\in\UU^{n,k}\times\CC^{k,k}$
be a solution of NEPv~\eqref{eq:nepv} satisfying
\Cref{ass:qspd}. Then the following statements hold.
\begin{enumerate}[(a)]
    \item\label{thm:convergence:i}
        The $\RR$-linear operator
        $\LL(X^*,\Lambda^*)[\,\cdot\,]$
        defined in~\eqref{eq:grassnt:operator}
        is invertible on $\CC^{n,k}\times\CC^{k,k}$.

    \item\label{thm:convergence:ii}
        Let $X\in\UU^{n,k}$ have a range space
        sufficiently close to $\range(X^*)$, and let
        $\Lambda\in\CC^{k,k}$ satisfy
        $\Lambda=X^HH(X)X+\mathcal O(\|\Theta(X,X^*)\|)$.
        The update~\eqref{eq:ntcorrections}, with corrections
        determined by~\eqref{eq:grassnt}, satisfies
        \begin{equation}
            \|\Theta(X_{\mathrm{next}},X^*)\|
            =\mathcal O\bigl(\|\Theta(X,X^*)\|^2\bigr)
            \quad \mbox{as} \quad 
            \|\Theta(X,X^*)\|\to0. 
        \end{equation}
        Here, $\Theta(X,Y)$ is the diagonal matrix
        of canonical angles between the subspaces $\range(X)$ and $\range(Y)$,
        i.e.,  $\Theta(X,Y)=\mathrm{Diag}(\arccos\sigma_1,\dots,\arccos\sigma_k)$ with $\sigma_1\geq\cdots\geq\sigma_k$ being the singular values of $X^HY$.
\end{enumerate}
\end{theorem}
\begin{proof}
We first prove item~\eqref{thm:convergence:i}.
It suffices to show that $\LL(X^*,\Lambda^*)$ has a trivial
null space. Let $(\dX,\dL)$ be a null vector,
that is, a solution of $\LL(X^*,\Lambda^*)[\dX,\dL]=0$.
By the definition of $\LL$ in~\eqref{eq:grassnt:operator},
this is equivalent to
\begin{equation}\label{eq:YDnull}
    \begin{bmatrix}
        H(X^*)\dX+\DD{H}(X^*)[\dX]X^*
        -\dX\Lambda^*-X^*\dL\\
        -{X^*}^H\dX
    \end{bmatrix}=0.
\end{equation}
The lower block gives ${X^*}^H\dX=0$, so
$\dX=X_\perp^*Z$ for some $Z\in\CC^{n-k,k}$.
Multiplying the upper block by ${X_\perp^*}^H$ from the
left and substituting $\dX=X_\perp^*Z$, we obtain
\begin{equation}\label{eq:topeqnull}
    0={X_\perp^*}^HH(X^*)X_\perp^*Z
    +{X_\perp^*}^H\DD{H}(X^*)[X_\perp^*Z]X^*
    -Z\Lambda^*
    =\mathcal Q^*(Z),
\end{equation}
where the last equality follows from
${X_\perp^*}^HH(X^*)X_\perp^*=\Lambda_\perp^*$
and the definition of $\mathcal Q^*$ in~\eqref{eq:rlinearq}.
By Assumption~\ref{ass:qspd}\,\eqref{ass:qspd:iii},
$\mathcal Q^*$ is positive definite, so
\eqref{eq:topeqnull} implies $Z=0$.
Substituting $\dX=X_\perp^*Z=0$ into the upper block
of~\eqref{eq:YDnull} gives $X^*\dL=0$, hence $\dL=0$.
Thus $(\dX,\dL)=(0,0)$, proving that
$\LL(X^*,\Lambda^*)$ is invertible.

For item~\eqref{thm:convergence:ii}, we consider the exact
corrections $\dX_*\in\CC^{n,k}$ and $\dL_*\in\CC^{k,k}$
that take $(X,\Lambda)$ to the solution $(X^*,\Lambda^*)$
through~\eqref{eq:ntcorrections}, that is,
\begin{equation}\label{eq:ntupdate2}
    X^*=\mathrm{orth}(X+\dX_*),
    \qquad
    \Lambda^*=\Lambda+\dL_*.
\end{equation}
Here, we assume without loss of generality that the target solution $X^*$
is aligned with the given basis $X$ so that
$X^HX^*\succ0$; see~\Cref{sec:aaaligned}.
Then a direct calculation gives
\[
    \dX_*\equiv X^*(X^HX^*)^{-1}-X
	=(I-XX^H)X^* (X^HX^*)^{-1}
    =\mathcal O(\|\Theta(X,X^*)\|),
\]
where the bound follows from $\|\dX_*\|=\|\tan\Theta(X,X^*)\|$;
see, e.g.,~\cite{Zhu:2013}.
This also implies 
\[
    \dL_* \equiv {X^*}^HH(X^*)X^* - X^HH(X)X +  \mathcal O(\|\Theta(X,X^*)\|) 
	= \mathcal O(\|\Theta(X,X^*)\|),
\]
where the 
second equality is due to 
the smoothness of the Rayleigh quotient map $X\mapsto \Lambda:=X^HH(X)X$,
which implies the difference in Rayleigh quotient is $\mathcal O(\|X^*-X\|)= \mathcal O(E^*\|)$ by~\eqref{eq:ntupdate2}.
Now, following the derivation from~\eqref{eq:correctexact}
to~\eqref{eq:grassnt}, with
$(X_{\mathrm{next}},\Lambda_{\mathrm{next}})$ replaced
by $(X^*,\Lambda^*)$ and retaining the higher-order terms,
we obtain
\[
    \LL(X,\Lambda)[\dX_*,\dL_*]
    + R(X,\Lambda)
	=\mathcal O(\|\dX_*\|^2+\|\dL_*\|^2)
    =\mathcal O(\|\Theta(X,X^*)\|^2).
\]
Recalling $R(X,\Lambda) = - \LL(X,\Lambda)[\dX,\dL]$
from the equation~\eqref{eq:grassnt:system}
for the correction $(\dX,\dL)$,
we obtain
\begin{equation}\label{eq:ntexact}
    \LL(X,\Lambda)[\dX_*,\dL_*] - \LL(X,\Lambda)[\dX,\dL] = \mathcal O(\|\Theta(X,X^*)\|^2).
\end{equation}
Observe that the linear operator $\LL(X,\Lambda)$ 
is Lipschitz continuous in a neighborhood of $(X^*,\Lambda^*)$.
Its invertibility at $(X^*,\Lambda^*)$ in
item~\eqref{thm:convergence:i} therefore implies that it
remains invertible, with a uniformly bounded inverse,
in a sufficiently small neighborhood of the solution.
Applying this inverse to~\eqref{eq:ntexact} gives
\[
    (\dX_*,\dL_*)-(\dX,\dL)
	= \mathcal O(\|\Theta(X,X^*)\|^2).
\]
Hence, it follows from the update formulas~\eqref{eq:ntcorrections}
and~\eqref{eq:ntupdate2}
\[
    X_{\mathrm{next}} - X^*
    =\mathrm{orth}(X+\dX) - \mathrm{orth}(X+\dX_*)
	=\mathcal O(\dX-\dX_*) 
	=\mathcal O(\|\Theta(X,X^*)\|^2),
\]
where the second equality is due to the continuity of the polar-factor
$\mbox{orth}(\cdot)$.
Consequently, we must also have the geodesic distance 
$\|\Theta(X_{\mathrm{next}},X^*)\| =\mathcal O(\|\Theta(X,X^*)\|^2)$, 
which completes the proof.
\end{proof}

\subsection{SCF Iteration with Basis-Aligned Newton Acceleration}

\Cref{alg:scf-newton} outlines an SCF iteration combined with
the Newton acceleration developed in~\Cref{sec:grassntderive}.
Each iteration first computes an SCF iterate and then refines it
by a Newton correction.
The SCF step serves to target the desired eigenvalues:
the governing equation~\eqref{eq:correctexact} does not distinguish
the $k$ smallest eigenvalues from other sets of $k$ eigenvalues,
so Newton's method alone may converge to an undesired solution.
By selecting the $k$ smallest eigenvalues at each iteration,
the SCF step helps direct the iteration toward the desired
invariant subspace.
For simplicity, the algorithm applies a Newton correction at
every iteration; in practice, one may activate Newton acceleration
only after the residual norm falls below a prescribed threshold,
since its convergence guarantee is local.

\begin{algorithm}[!htbp]
\caption{SCF with basis-aligned Newton acceleration}
\label{alg:scf-newton}
\begin{algorithmic}[1]
\REQUIRE Initial $X_0\in\UU^{n,k}$, matrix function $H$
and its derivative $\DD{H}$, residual tolerance $\tau>0$. 
\ENSURE approximate eigenbasis $X^*\in\UU^{n,k}$ of the NEPv~\eqref{eq:nepv}
\FOR{$i=0,1,2,\dots$}
    \STATE Evaluate $H_i=H(X_i)$
    \STATE Solve $H_iY_i=Y_i\Lambda_i$ for the $k$ smallest
        eigenvalues in $\Lambda_i$
	\STATE Evaluate $\widehat H_i=H(Y_i)$ and the residual $R_i=\widehat H_iY_i-Y_i\Lambda_i$
    \STATE {\bf if} {$\|R_i\|_F\leq\tau$}, {\bf then return} $X^*=Y_i$
    \STATE Compute $(\dX_i,\dL_i)$ by solving
        $\LL(Y_i,\Lambda_i)[\dX_i,\dL_i]
        =-\begin{bmatrix}R_i\\0\end{bmatrix}$ 
		\label{alg:scf-newton:ls}
	\STATE Compute $X_{i+1}=\mathrm{orth}(Y_i+\dX_i)$
		\label{alg:scf-newton:retraction}
\ENDFOR
\end{algorithmic}
\end{algorithm}

\begin{remark}\label{rmk:ntls}
\rm
The linear system in line~\ref{alg:scf-newton:ls} of~\Cref{alg:scf-newton}
can be solved inexactly by an $\RR$-linear global GMRES
method with the Frobenius inner product
$\langle A,B\rangle=\Re(\TR(A^HB))$
\cite{EirHV2003,HuhP2013,Jbilou:1999}.
This is essentially the GMRES method, but with
the solution $(\dX,\dL)$ represented as a block matrix
$\begin{bmatrix}\dX^H,\dL^H\end{bmatrix}^H\in\CC^{n+k,k}$
instead of its vectorization as a long column vector.
In particular, starting with the initial basis matrix $V_1=-[R_i^H,0_{k,k}^H]^H/\|R_i\|_F$,
the Arnoldi process generates
$\mathcal V_j=[V_1,\dots,V_j]$, whose column blocks
$V_\ell\in\CC^{n+k,k}$ are orthonormal with respect to
$\langle\cdot,\cdot\rangle$ and satisfy
\begin{equation}\label{eq:rlineararnoldi}
    \LL{(Y_i,\Lambda_i)}[\mathcal V_j]
    =\mathcal V_{j+1}(\underline H_{j}\otimes I_k),
\end{equation}
where the operator on the left acts on each column block
and $\underline H_{j}\in\RR^{j+1,j}$ is an upper
Hessenberg matrix.
The approximate solution is given by $\mathcal V_j(y_j\otimes I_k)$, with $y_j\in\RR^j$
solving $\min_y\left\|\|R_i\|_Fe_1-\underline H_{j}y\right\|_2$.
Restarting may be applied to limit storage and orthogonalization costs,
and the iteration stops once the residual of the $\RR$-linear system 
is below a tolerance chosen according to standard techniques for 
inexact Newton methods, such as the Eisenstat-Walker adaptive forcing strategy; see, e.g.,~\cite{EisW1994,EisW1996,Saad:2007}.
Note that GMRES requires only operator actions,
so there is no need to construct a matrix representation of
$\LL(Y_i,\Lambda_i)$ explicitly.
When these operator actions are expensive due to the derivative
$\DD{H}$, preconditioning may reduce the number of derivative evaluations.
Using the splitting
\begin{equation}\label{eq:split}
    \LL{(X,\Lambda)}[\dX,\dL]
    =
    \begin{bmatrix}
        H(X) & -X\\
        -X^H & 0
    \end{bmatrix}
    \begin{bmatrix}\dX\\ \dL\end{bmatrix}
    -
    \begin{bmatrix}\dX\\ \dL\end{bmatrix}\Lambda
    +
    \begin{bmatrix}
        \DD{H}(X)[\dX]X\\
        \dL\Lambda
    \end{bmatrix},
\end{equation}
we approximate
\begin{equation}\label{eq:precondll}
    \LL{(Y_i,\Lambda_i)}[\dX,\dL]
    \approx \mathcal S_i[Z]
    :=
    \begin{bmatrix}
        \widehat H_i & -Y_i\\
        -Y_i^H & 0
    \end{bmatrix}Z-Z\Lambda_i,
    \qquad
    Z=\begin{bmatrix}\dX\\\dL\end{bmatrix} \in\CC^{n+k,k}.
\end{equation}
Since $\mathcal S_i$ is a Hermitian Sylvester operator,
a few steps of global MINRES applied to
$\mathcal S_i[Z]=V$ provide an approximate inverse
preconditioner without $\DD{H}$.
In this case, the preconditioner is implicit and may change between iterations,
making flexible right-preconditioned GMRES a suitable choice;
see, e.g.,~\cite{Saad1993,Saad:2007}.
Note that even when $\LL(X,\Lambda)$ is self-adjoint and MINRES
can be applied directly, preconditioned GMRES may still
be beneficial by reducing the number of expensive
$\DD{H}(X)$ evaluations, since applying the preconditioner
$\mathcal S_i^{-1}$ does not require this derivative.
\end{remark}

\section{Numerical Experiments} \label{sec:experiments}
We demonstrate the convergence behavior and the effectiveness of the
acceleration techniques proposed in previous sections.  
For a set of NEPv, we compare the plain SCF~\eqref{eq:scf}
with its accelerated variants. 
These include several versions of \AAm, with a fixed window size $m\in\mathbb{N}$,
as well as the basis-aligned Newton acceleration:
\begin{enumerate}[(i)]
\item  
	\AAm--$X$: \AAm applied directly to $\{X_i\}$ by SCF,
	i.e., \Cref{alg:aag} with $g$ given by~\eqref{eq:scffpmap1};
\item 
	\AAm--$\rho$: \AAm for a problem-specific ui-parameterization
	$\rho(X)$ in~\Cref{alg:aarho}; 
\item 
	\AAm--$\ALIGN{X}$: \AAm for aligned basis
	$\{\ALIGNM{X_i}\}$ with a common $M\in\UU^{n,k}$
	in~\Cref{alg:alignedaa};%
	\footnote{We use a random orthonormal reference matrix $M$,
		generated in MATLAB by \texttt{M = orth(randn(n,k))}.}
\item 
	\AAm-Gr: the Riemannian version of \AAm in~\Cref{alg:grassaa_diff};
\item 
	AlignNT: the basis-aligned Newton acceleration in~\Cref{alg:scf-newton}, with the adaptive forcing~\cite{EisW1996}.
\end{enumerate}
For the SCF step in each algorithm,
we use \textsc{MATLAB}'s built-in function \texttt{eig} or \texttt{eigs} 
to compute the desired eigenpairs of the Hermitian matrix $H(X_i)$. 
We measure the accuracy of an approximate solution
$\widehat X \in\UU^{n,k}$ using the residual norm
\begin{equation}\label{eq:nepvresidual}
    \mathrm{Resd}(\widehat X):=\|H(\widehat X)\widehat X- \widehat X \widehat\Lambda\|_F
    \quad\text{with}\quad
    \widehat \Lambda=\widehat X^HH(\widehat X) \widehat X.
\end{equation}
All methods use the common stopping
criterion $\mathrm{Resd}(X_i)\leq 10^{-14}$,
unless otherwise specified.

All experiments are performed using \textsc{MATLAB 2024b} on a
laptop with AMD Ryzen 7 PRO 4750U CPU and 32GB RAM.
The code used to reproduce the numerical experiments in this
paper is available at
\url{https://github.com/towerner/Local-acceleration-for-uiNEPv}.

\subsection{Model Kohn-Sham Problem}
We consider a model Kohn--Sham NEPv
\begin{equation}\label{eq:1dks}
H(X)X=X\Lambda
\quad\text{with}\quad
H(X):=L+\gamma\mathrm{Diag}(L^{-1}\rho(X)),
\quad \rho(X):=\mathrm{diag}(XX^T),
\end{equation}
where $L=\mathrm{tridiag}(-1,2,-1)\in\RR^{n,n}$ is a discrete 1D
Laplacian, $\gamma\geq0$ is a given parameter,
and the eigenvalues of $\Lambda\in\RR^{k,k}$
are the $k$ smallest eigenvalues of $H(X)$.
This simplified Kohn--Sham model omits the external potential
in the Hamiltonian $H(X)$ (see~\Cref{rmk:parameterization})
while retaining its nonlinear dependence on the charge-density $\rho(X)$.
The parameter $\gamma$ controls the strength of the nonlinearity of $H(X)$.
By varying $\gamma$, the NEPv~\eqref{eq:1dks} provides a convenient test 
problem for studying SCF convergence and
has been widely used as a benchmark in previous work; 
see, e.g.,~\cite{BaiLL2022,CaiZBL2018,Liu:2014,Yang:2009}.

The NEPv~\eqref{eq:1dks} admits a natural ui-parameterization through the density 
$\rho(X)= \mathrm{diag}(XX^T)\in\RR^n$.
Moreover, the Hamiltonian $H$ is Fr\'echet differentiable 
with the derivative operator given by 
\begin{equation}
\DD H[X](E)=2\gamma \mathrm{Diag}\left(L^{-1}\mathrm{diag}(EX^T)\right),
\end{equation}
which allows convenient application of the basis-aligned Newton acceleration
in~\Cref{alg:scf-newton}.

\subsubsection{Effectiveness of Acceleration}
\begin{figure}[!htbp]
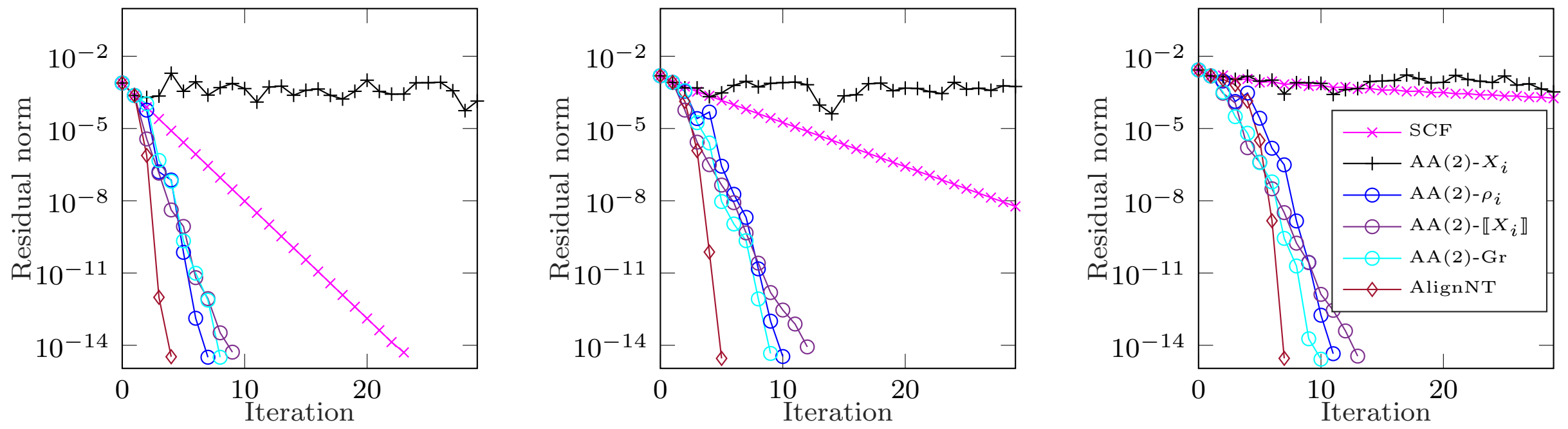

    \centering
    \begin{minipage}{\threeplotscale}
    \include{Figs/tikz/ks1d_beta2_iter}
    \end{minipage}
    \begin{minipage}{\threeplotscale}
    \include{Figs/tikz/ks1d_beta4_iter}
    \end{minipage}
    \begin{minipage}{\threeplotscale}
    \include{Figs/tikz/ks1d_beta7_iter}
    \end{minipage}
    \vspace{-1cm}
    \caption{Convergence histories of the algorithms for 
		the 1D Kohn-Sham NEPv~\eqref{eq:1dks} with $n=200$, $k=2$, 
		and $\gamma_0=2,4,7$ (left, middle, and right).}
    \label{fig:1dks}
\end{figure}

We first demonstrate the effectiveness of the proposed acceleration techniques.
For experiments, we consider test NEPv~\eqref{eq:1dks} of size $n=200$, $k=2$, and
parameters $\gamma=\frac{\gamma_0}{n^2}$ for ${\gamma_0}\in\lbrace 2,4,7\rbrace$.
We apply SCF iteration and the accelerated variants with a common initial guess 
$X_0\in\RR^{n,k}$, chosen as the eigenbasis for the $k$ smallest
eigenvalues of $L$. 
For AA, we use a fixed window size $m=2$.
Because the NEPv is relatively small, we solve the eigenvalue problem 
in each SCF step using the dense solver \texttt{eig}. 
We also solve the Newton equation in line~\ref{alg:scf-newton:ls} of
\Cref{alg:scf-newton} to numerical precision 
using unrestarted GMRES with a maximum of $(n+k)k$ iterations.
Convergence histories of the algorithms are reported in~\Cref{fig:1dks}.

We observe that the convergence of SCF becomes
increasingly slow as $\gamma$ increases, highlighting the need
for acceleration. Direct application of \AAm to the iterates
$\{X_i\}$ fails to converge in all three cases, illustrating
the difficulty caused by the non-uniqueness of the basis under
rotations. In contrast, the three proposed variants of \AAm,
based on the ui-parameterization $\rho(X)$, the aligned basis
$\ALIGNM{X}$, and the Grassmannian formulation, all address this
ambiguity and significantly accelerate SCF.
These variants exhibit comparable convergence behavior,
requiring similar numbers of iterations.
Basis-aligned Newton acceleration exhibits quadratic convergence in all three cases,
reaching the prescribed accuracy in at most seven iterations.

\subsubsection{Nonsingularity and Preconditioning of the Basis-Aligned Newton Equation}
Next we consider the nonsingularity of the correction operator
in the basis-aligned Newton equation~\eqref{eq:grassnt}, 
with its advantages over the standard Newton formulation for NEPv~\cite{GaoYM2007,Werner2025}
and the effectiveness of the preconditioning strategy described in~\Cref{rmk:ntls}.
We use the model problem~\eqref{eq:1dks} with $\gamma_0=7$.
Let $\widehat X\in\RR^{n,k}$ be a computed solution 
($\|{\widehat X}^T \widehat X-I_k\|\leq 10^{-15}$, residual normal satisfying $\mbox{Resd}(\widehat X)\leq 10^{-14}$) 
and set $\widehat \Lambda={\widehat X}^T H(\widehat X) \widehat X$.
The condition numbers of the correction operators
$\mathcal L(\widehat X, \widehat \Lambda)$ and
$\widehat{\mathcal L}(\widehat X, \widehat \Lambda)$,
corresponding to the basis-aligned Newton~\eqref{eq:grassnt} and
standard Newton~\eqref{eq:stnt}, respectively, are%
\footnote{For a linear operator $\mathcal L$ on
$\RR^{n,k}\times\RR^{k,k}$, we construct its matrix representation
$\mathbf L\in\RR^{N,N}$, where $N=(n+k)k$, using an orthonormal
basis with respect to the Frobenius inner product.
We then compute
$\kappa(\mathcal L)=\sigma_1(\mathbf L)/\sigma_N(\mathbf L)$,
where $\sigma_1$ and $\sigma_N$ denote the largest and smallest
singular values, respectively.}
\begin{equation}
    \kappa(\mathcal L(\widehat X, \widehat \Lambda))=2.25\times10^3
    \qquad\text{and}\qquad
    \kappa(\widehat{\mathcal L}(\widehat X,  \widehat \Lambda))=1.68\times10^{16}.
\end{equation}
These values indicate that $\mathcal L(\widehat X, \widehat\Lambda)$ is
numerically nonsingular, consistent with the analysis
in~\Cref{thm:convergence}, while
$\widehat{\mathcal L}(\widehat X, \widehat \Lambda)$ is numerically singular,
as previously observed in~\cite{Werner2025}.

\begin{figure}[!htbp]
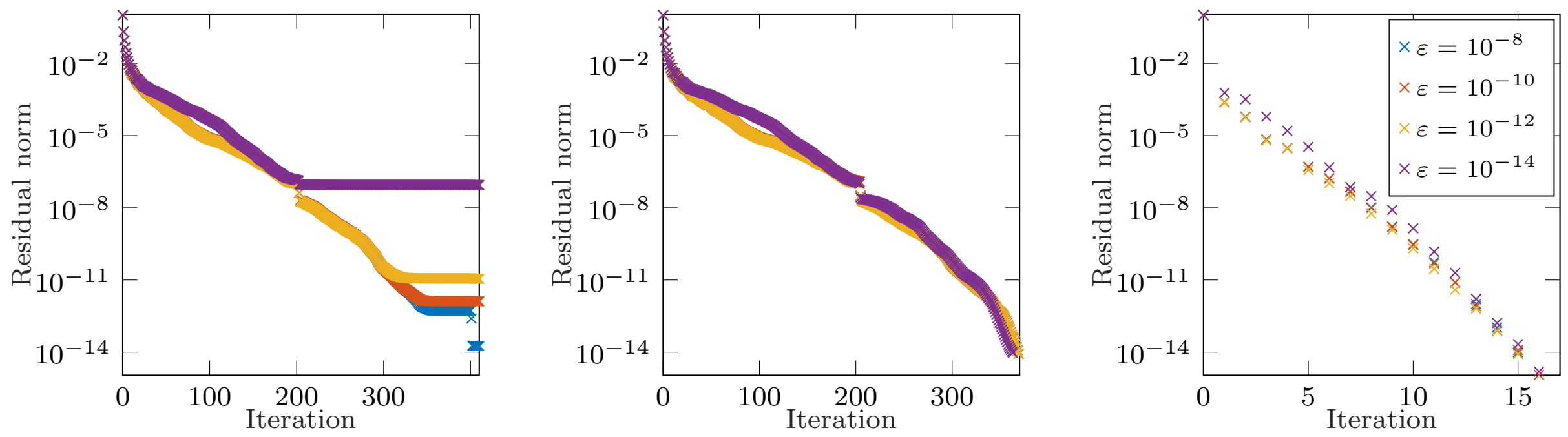

    \centering
  \begin{minipage}{\threeplotscale}
       \include{Figs/tikz/GMRES_residual_NT}
   \end{minipage}
   \begin{minipage}{\threeplotscale}
        \include{Figs/tikz/GMRES_residual_NTAcc}
   \end{minipage}
   \begin{minipage}{\threeplotscale}
       \include{Figs/tikz/GMRES_residual_NTAcc_MR}
   \end{minipage}
    \vspace{-1cm}
	\caption{Residual histories for the Newton correction equations near
		a solution: GMRES for standard Newton~\eqref{eq:stnt} (left),
		GMRES for basis-aligned Newton~\eqref{eq:grassnt} (middle), 
		and flexible GMRES for basis-aligned Newton with MINRES preconditioning (right).  
		The curves correspond to different perturbation sizes $\varepsilon$.}
    \label{fig:jacsingular}
\end{figure}

To demonstrate how this singularity affects the solution of the Newton
correction equations near convergence, we consider the perturbed
solutions of a perturbation size $\varepsilon$:
\begin{equation}
    X_\varepsilon:=\mathrm{orth}(\widehat X+\varepsilon P),
    \qquad
    \varepsilon\in\{10^{-8},10^{-10},10^{-12},10^{-14}\},
\end{equation}
where $P\in\RR^{n,k}$ is a Gaussian random matrix
(generated using \textsc{MATLAB}'s \texttt{randn(n,k)}).
For each $X_\varepsilon$, we set
$\Lambda_\varepsilon=X_\varepsilon^T H(X_\varepsilon)X_\varepsilon$
and apply non-restarted GMRES to the corresponding correction
equations~\eqref{eq:grassnt} and~\eqref{eq:stnt}, respectively, where
the relative convergence tolerance is set to $\tau=10^{-14}$. The
convergence histories are shown in the left and middle panels
of~\Cref{fig:jacsingular}. Note that, in theory, an exact solution
should be contained in the Krylov subspace of size $(n+k)k=404$.
For the standard Newton equation~\eqref{eq:stnt}, GMRES stagnates,
consistent with the ill-conditioning of the correction operator
near a solution. In particular, the condition numbers of the systems are between $1.7\times10^{16}$ to $2.1\times10^{17}$ across all tests.
In contrast, GMRES converges without stagnation for the basis-aligned Newton 
equation~\eqref{eq:grassnt}, with little variation in convergence
behavior across the four perturbation sizes.  
The condition numbers are almost constant at around $2.25\times10^3$ across all different perturbations.

Finally, the right panel of~\Cref{fig:jacsingular} shows the
convergence histories of flexible GMRES for the basis-aligned Newton 
equation~\eqref{eq:grassnt}, using the preconditioning strategy
described in~\Cref{rmk:ntls}.
Each application of the preconditioner approximately solves
$\mathcal S_i^{-1}$ by performing $40$ steps of MINRES
on the Sylvester equation~\eqref{eq:precondll}.
Compared with the unpreconditioned solves as shown in the middle panel,
preconditioning reduces the number of outer iterations by
approximately a factor of ten across all four values of
$\varepsilon$. The corresponding computation time decreases
from approximately $0.22$ to $0.045$ seconds.
These results suggest that the proposed preconditioner
can substantially improve the efficiency of solving the basis-aligned Newton equation.

\subsection{Gross-Pitaevskii NEPv} \label{sec:synthc}
Consider a discrete Gross--Pitaevskii NEPv in the form of 
\begin{equation}\label{eq:gpe}
    H(x)x=\lambda x
    \quad\text{with}\quad
    H(x):=A_f+\gamma\mathrm{Diag}(|x|^2),
\end{equation}
where $x\in\CC^n$ satisfies $\|x\|_2=1$,
$|x|$ denotes entry-wise absolute value, 
$\gamma\geq0$ is a given parameter, and $\lambda$ is
the smallest eigenvalue of $H(x)$.
The constant Hermitian matrix $A_f\in\CC^{n,n}$ is given by
\begin{equation} \label{eq:gpedis}
    A_f:=-\frac{1}{2}L_{\mathrm{2D}}
    -\imath\omega L_C+\mathrm{Diag}(v_f),
\end{equation}
where $L_{\mathrm{2D}},L_C\in\RR^{n,n}$ are central-difference
discretizations of the two-dimensional Laplacian and the angular
derivative, respectively, on $\Omega=[-1,1]^2$, using
$N=\sqrt{n}$ interior grid points in each coordinate direction,
$\imath$ denotes the imaginary unit, $\omega\geq0$ is
the rotation parameter, and $v_f\in\RR^n$ contains the values
of an external potential $f:\RR^2\to\RR$ at the grid points.
The Gross-Pitaevskii NEPv arises in the study of rotating Bose--Einstein
condensates; see, e.g.,~\cite{BaoD2004,Bao:2003,Henning:2025}.
The parameters $\gamma$ and $\omega$ control the strength of
the nonlinear interaction and the rotation, respectively.
For details of the discretization~\eqref{eq:gpedis}, 
we refer to~\cite{BaiLL2022,JarKM2014}.

The solution of the NEPv~\eqref{eq:gpe} 
is a single eigenvector $x$,
which is generally complex, $x\in\CC^n$,
since $H(x)$ is Hermitian.
Nevertheless, the non-uniqueness issue for the 
iterates $\{x_i\}$ of the SCF iteration~\eqref{eq:scf}
persists, as each eigenvector $x_i$ is only determined up to a phase rotation 
$x_i\alpha$ with $\alpha\in\CC$ and $|\alpha|=1$.
Hence, the techniques introduced in~\Cref{sec:aa,sec:newton} are still needed 
to effectively accelerate the SCF iteration. 
By the definition~\eqref{eq:gpe},
the Hamiltonian $H(x)$ admits a straightforward ui-parameterization
\begin{equation}\label{eq:gperho}
\rho(x)=|x|^2\in\RR^n.
\end{equation}
Moreover, $H$ is Fr\'echet differentiable over the underlying
real vector space, with derivative
\begin{equation}\label{eq:gpeder}
    \DD H(x)[e]
    =2\gamma\mathrm{Diag}
    \left(\Re\left(\overline{x}\circledast e\right)\right),
\end{equation}
where $\circledast$ denotes the Hadamard product (entry-wise product)
and $\overline x$ is the conjugate vector.

For the experiments, we consider a relatively large-size 
NEPv~\eqref{eq:gpe} with $n=10,000$, $\omega=0.85$,
and the harmonic potential $f(\xi,\eta)=(\xi^2+\eta^2)/2$,
following a setup in~\cite{BaiLL2022,JarKM2014}.
Since $n$ is relatively large, we use MATLAB's iterative
solver \texttt{eigs} to solve the Hermitian eigenvalue problem
$H(x_i)x=\lambda x$ in each SCF step.
For computational efficiency, we warm-start \texttt{eigs}
with the current iterate $x_i$. For the Newton-type acceleration in Algorithm~5, we use an
inexact Newton implementation: the correction equation is
solved approximately by flexible GMRES with a maximum Krylov
subspace dimension of $100$, using $100$ MINRES steps for
each preconditioner application; see~\Cref{rmk:ntls}. 


\begin{figure}
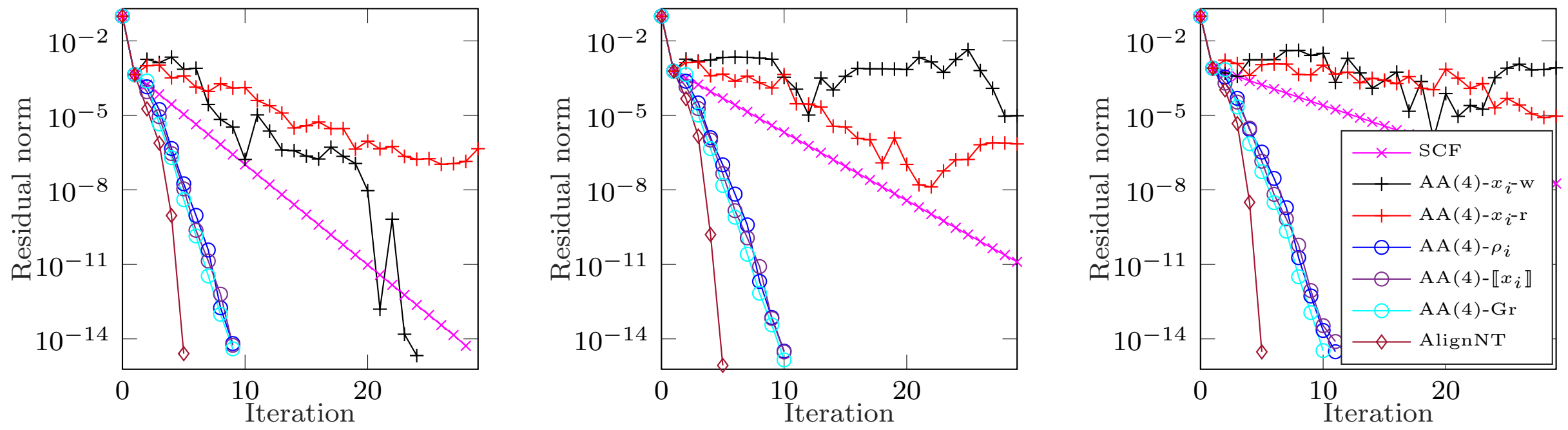

    \centering
    \begin{minipage}{\threeplotscale}
       \include{Figs/tikz/gpe_N100_beta40}
   \end{minipage}
   \begin{minipage}{\threeplotscale}
       \include{Figs/tikz/gpe_N100_beta55}
   \end{minipage}
   \begin{minipage}{\threeplotscale}
       \include{Figs/tikz/gpe_N100_beta72}
   \end{minipage}
   \vspace{-1cm}
   \caption{Convergence histories of the algorithms for the
	   Gross--Pitaevskii NEPv~\eqref{eq:gpe} with $n=10,000$ 
	   and $\gamma=4,5.5,7.2$ (left, middle, and right). 
	   We use window size $m=4$ for \AAm.
	   For AA-$x_i$, results are shown with both
	   warm-started (labeled $w$) and randomly initialized (labeled $r$)
   		\texttt{eigs} runs for SCF steps.
	}
  	\label{fig:gpeN100warmstart}
\end{figure}

\begin{table}
{
\footnotesize 
	\caption{Computation times (seconds) for the Gross--Pitaevskii
	NEPv~\eqref{eq:gpe} with 
    sizes  $n\in\{10,000,40,000\}$ 
    and parameters $\gamma\in\{4,5.5,7.2\}$.
	All methods use warm-started \texttt{eigs} for the SCF steps.  
	A dash indicates failure to converge within the prescribed iteration limit.}
    \label{tab:gpe}
    \begin{center}
    \begin{tabular}{cc||c|c|c|c|c|c}
    \toprule 
     $n$&$\gamma$ & ~~~SCF~~~ & AA($4$)-$x_i$  & AA($4$)-$\rho_i$  & AA($4$)-$\ALIGN{x_i}$ & AA($4$)-Gr& AlignNT \\
     \midrule
     \multirow{3}{*}{$10,000$} & 4.0 & 3.72 & 4.35  & 1.55 & 1.52 & 1.51 & 1.35\\
     \cline{2-8}
     & 5.5 & 3.72 & 4.76  & 1.61 & 1.66 & 1.60 & 1.40\\
     \cline{2-8} 
     &7.2 & 5.52 & --  & 1.85 & 1.78 & 1.74 & 1.33 \\
     \midrule
     \multirow{3}{*}{$40,000$} & 4.0 & 21.50 & 58.97  & 9.18 & 9.36 & 8.88 & 6.63\\
     \cline{2-8}
     & 5.5 & 30.10 & -- & 10.00 & 10.19 & 9.52 & 6.50\\
     \cline{2-8} 
     &7.2 & 48.74 & --  & 10.91 & 11.06 & 10.54 & 6.47 \\
    \bottomrule 
    \end{tabular}
    \end{center}
}
\end{table}

\Cref{fig:gpeN100warmstart} reports the convergence histories
for the NEPv~\eqref{eq:gpe} with $\gamma\in\{4,\, 5.5,\, 7.2\}$, 
where all methods start from the common initial vector
$x_0=(\mathbf{1}_n+\imath\mathbf{1}_n)/\sqrt{2n}$. 
The corresponding computation times are listed in~\Cref{tab:gpe}, along with timings for a larger scale problem of size $n=200^2=40,000$ using the same setup.
We observe that the plain SCF converges in all three cases,
with slower convergence as $\gamma$ increases, similarly
to the preceding Kohn--Sham example.
Direct application of \AAm to the SCF iterates $x_i$
(\AAm-$x_i$) does not effectively accelerate convergence,
due to the phase rotation $x_i\alpha$ of computed eigenvectors.
Indeed, its acceleration is substantially delayed in the first
two cases of $\gamma$, with convergence requiring $24$ and
$62$ iterations, respectively;%
\footnote{
The eventual convergence may be explained by local approximate 
phase alignment induced by warm-starting \texttt{eigs} with $x_i$,
which near a solution is already an approximate eigenvector of
$H(x_i)$ and tends to yield $x_i^Hx_{i+1}\approx 1$.
Such implicit alignment, however, depends on the internal implementation 
of \texttt{eigs} and is not guaranteed in general.}
in the last case, convergence is not achieved within the prescribed
limit of $80$ iteration.
We note that the behavior of \AAm-$x_i$ also depends on the 
eigensolver used in the SCF step, 
which can lead to different phase choices
for the computed eigenvectors.
For comparison, we also report the convergence histories
of \AAm-$x_i$ when \texttt{eigs} uses a random 
starting vector\footnote{
We set the \texttt{StartVector} option of \texttt{eigs}
to a normalized vector generated by
\texttt{randn(n,1)+1i*randn(n,1)}.}.
In this case, convergence is not achieved within the
prescribed iteration limit in all three tests.

In comparison, all four proposed acceleration techniques
significantly reduce the number of SCF iterations in all three
test cases.
The three variants of AA
exhibit comparable performance
in terms of both iteration counts and computation times.
The Newton acceleration exhibits 
superlinear convergence, reaching the prescribed accuracy in six iterations.
It also achieves the least computation times in all three cases.
This efficiency can be attributed to the inexpensive evaluation
of the derivative~\eqref{eq:gpeder}, together with effective
preconditioning of the Newton equation.

\subsection{Sparse PCA}
In this example, we consider NEPv arising from 
sparse PCA (principal component analysis)~\cite{JolTU2003,ZouHT2006}.
The optimization problem of interest takes the form
\begin{equation}\label{eq:21max}
    \min_{X\in\RR^{n,k},\,X^TX=I_k}
    -\mathrm{tr}(X^TAX)+\alpha\cdot \|X\|_{2,1},
\end{equation}
where $A\in\RR^{n,n}$ is symmetric positive semidefinite, 
$\alpha\geq0$ is a penalty parameter, 
and $\|X\|_{2,1}:=\sum_{j=1}^n\|e_j^TX\|_2$ denotes the $(2,1)$-norm 
(i.e., the sum of $2$-norms of the rows of $X$).
Without the penalty, i.e., when $\alpha=0$, the problem reduces
to standard PCA and can be solved by computing the eigenvectors
of the $k$ smallest eigenvalues of $-A$ (equivalently, the $k$ largest of $A$).
The penalty term promotes row sparsity and thus favors a sparse
solution $X$.
With a smooth approximation of the penalty,
the first-order optimality condition of~\eqref{eq:21max} leads to the NEPv
\begin{equation}\label{eq:spcanepv}
    H_\varepsilon(X)X=X\Lambda
    \quad\text{with}\quad
    H_\varepsilon(X):=
	-2A+\alpha\cdot\mathrm{Diag}
    \left((\mathrm{diag}(XX^T)+\varepsilon^2\mathbf{1}_n)^{-1/2}\right),
\end{equation}
where $\Lambda$ corresponds to the $k$ smallest eigenvalues of $H_\varepsilon(X)$
for the purpose of minimization, the power ${-1/2}$ is taken entrywise, 
and $\varepsilon>0$ (e.g., $\varepsilon=10^{-3}\sqrt{k/n}$) 
is a smoothing parameter introduced to improve numerical stability when entries of $\mathrm{diag}(XX^T)$ are small 
-- Specifically, each row norm $\|e_j^TX\|_2$ in the penalty
$\|X\|_{2,1}$~\eqref{eq:21max}
is replaced by $(\|e_j^TX\|_2^2+\varepsilon^2)^{1/2}$.
We refer to~\cite{Li:2026,WanWL2025} for further details on the
derivation of~\eqref{eq:spcanepv}.


The NEPv~\eqref{eq:spcanepv} admits a natural ui-parameterization 
\begin{equation} \label{eq:spcanepvuip}
\rho(X):=\left(\mathrm{diag}(XX^T)+\varepsilon^2\mathbf{1}_n\right)^{1/2}\in\RR^{n}.
\end{equation}
Moreover, $H_\varepsilon$ is Fr\'echet differentiable, with derivative
\begin{equation}\label{eq:spcanepvder}
    \DD H_\varepsilon(X)[E]
    =-\alpha\,\mathrm{Diag}\left(
    (\mathrm{diag}(X)+\varepsilon^2\mathbf{1}_n)^{-3/2}
    \circledast\mathrm{diag}(EX^T)\right),
\end{equation}
where the power $-3/2$ is taken entrywise and $\circledast$ denotes
the Hadamard product (entrywise product).

For the experiments, we consider two real-world data sets,
USPS and COIL20~\cite{Hull1994,NenNM1996}, also used
in~\cite{WanWL2025}, with preprocessed data available
from~\cite{LiCWMTTL2017}.
The USPS data set contains $p=9,298$ samples of handwritten
digits with $n=256$ features, and COIL20 contains $p=1{,}440$
grayscale images with $n=1,024$ features.
We collect each dataset into a matrix $Y\in\RR^{p\times n}$,
with each row containing the feature vector of one sample.
The sample covariance matrix is constructed as 
$A=\frac{1}{p-1}\widehat Y^T\widehat Y$,
with 
$\widehat Y=Y-\mathbf{1}_p\mu^T$
and 
the sample mean $\mu=\frac{1}{p}Y^T\mathbf{1}_p$.
We then extract the sparse principal components by solving
the NEPv~\eqref{eq:spcanepv} using SCF iteration.
We aim to find $k=50$ eigenpairs for the USPS dataset 
and $k=30$ eigenpairs for the COIL20 dataset. 
The regularization parameters in both NEPv are set to 
$\varepsilon=10^{-3}\sqrt{k/n}$ and $\alpha=0.1$.

\begin{figure}[!htbp]
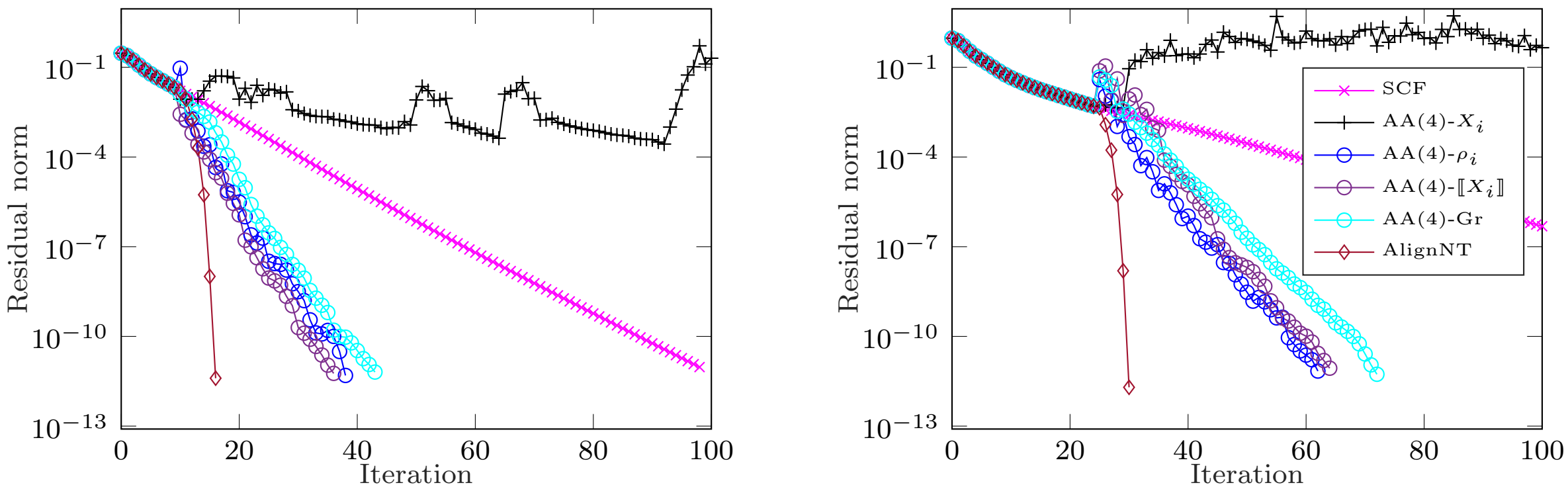

    \centering
	\begin{minipage}{\twoplotscale}
       \centering
       \include{Figs/tikz/residual_iter_usps}
   \end{minipage}
   \begin{minipage}{\twoplotscale}
       \centering
       \include{Figs/tikz/residual_iter_coil}
   \end{minipage}
   \vspace{-1cm}
	\caption{Convergence histories of the algorithms 
		for NEPv~\eqref{eq:spcanepv} 
		on the data sets
		USPS (left, $n=256$, $k=50$,  $10$ preliminary SCF steps) 
		and 
		COIL20 (right, $n=1,024$, $k=30$, $25$ preliminary SCF steps). 	}
    \label{fig:pca21}
\end{figure}

\Cref{fig:pca21} shows the convergence histories for the two
test NEPv, and \Cref{tab:pca21data} summarizes the parameters
and computation times.
The SCF iteration starts from the standard PCA solution $X_0\in\RR^{n,k}$, consisting of orthonormal eigenvectors
associated with the $k$ largest eigenvalues of $A$.
To approach the local convergence regime, we perform
$10$ preliminary SCF steps for USPS
and $25$ steps for COIL20 before enabling acceleration.
Since the problem sizes are moderate, the symmetric
eigenvalue problem in each SCF step is solved using
the dense solver \texttt{eig}.
For basis-aligned Newton acceleration, the correction equations 
are solved 
using MINRES, which exploits the self-adjointness of $\LL$;
see~\Cref{rmk:symnt}.
All algorithms stop when the NEPv residual norm falls
below $\tau=10^{-11}$, and  \AAm variants use a fixed window size  $m=4$.

We observe that the convergence behaviors of the algorithms are
largely consistent with the previous experiments.
Direct application of \AAm to the SCF basis iterates again
leads to divergence.
All other \AAm variants converge and accelerate SCF,
with similar convergence behavior and comparable computational costs.
In comparison, Newton acceleration exhibits superlinear convergence
and requires the least outer iterations.
However, its computation time shown in~\Cref{tab:pca21data}
is even higher than that of the plain SCF.
This is mainly due to the high cost of iterative solver 
for the Newton correction equation due to 
large condition number
partly attributed to the small row norms of the sparse PCA solution
that leads to large variations in diagonals of $H_\varepsilon$ and 
its derivative~\eqref{eq:spcanepvder}. The preconditioner in~\Cref{rmk:ntls}
cannot effectively accelerate the solution of these equations.

\begin{table}[!htbp]
{\footnotesize
    \caption{
		Computation times (seconds)
		for the sparse PCA NEPv~\eqref{eq:spcanepv} on the 
		data sets USPS and COIL20 as shown in~\Cref{fig:pca21}.
		A dash indicates failure to converge within the given iteration limit.
	}
    \label{tab:pca21data}
    \begin{center}
    	\begin{tabular}{c|c|c|c|c|c|c}
            \toprule
            Data set & ~~~SCF~~~ & AA($4$)-$X_i$ & AA($4$)-$\rho_i$
            & AA($4$)-$\ALIGN{X_i}$ & AA($4$)-Gr & AlignNT \\
            \midrule
            USPS & 0.77 & -- & 0.30 & 0.37 & 0.51 & 11.17 \\
            \hline
            COIL20 & 44.22 & -- & 15.41 & 16.16 & 17.93 & 96.26 \\
            \bottomrule
        \end{tabular}
    \end{center}
}
\end{table}

\subsection{Uniform Multidimensional Scaling}
As our final example, we consider the uniform multidimensional
scaling (UMDS) problem from multi-view learning introduced
in~\cite{ZhaZL2017}. 
Consider observations of the same $n$ objects from $s$ different views,
collected in matrices $Y_j\in\RR^{d_j,n}$, $j=1,\dots,s$,
where each column of $Y_j$ contains the $d_j$-dimensional observation
of one object in the $j$th view.
UMDS seeks a common orthogonal projection $X\in\RR^{n,k}$
that represents the $n$ objects in a latent space while preserving
information from the different views.
After relaxing the constraints, the projection can be obtained
from the optimization problem\footnote{
    We negate the original objective to obtain an equivalent
    minimization problem.
}
\begin{equation}\label{eq:umds}
    \min_{X\in\RR^{n,k},\,X^TX=I_k}
    -\sum_{j=1}^s\|X^TA_jX\|_F^2,
\end{equation}
where $A_j\in\RR^{n,n}$ is the symmetric centered squared
distance matrix for the $j$th view $Y_j$; see~\cite{ZhaZL2017} for details.
The first-order optimality condition for~\eqref{eq:umds}
leads directly to the NEPv
\begin{equation}\label{eq:umdsnepv}
    H(X)X=X\Lambda
    \quad\text{with}\quad
    H(X):=-\sum_{j=1}^s A_jXX^TA_j,
\end{equation}
where we target the $k$ smallest eigenvalues of $H(X)$ for
minimization; see, e.g.,~\cite{Li:2024}.
The matrix valued function $H(X)$ in~\eqref{eq:umdsnepv} is unitarily invariant.
However unlike the previous examples, a compact ui-parameterization is not readily
available for the NEPv~\eqref{eq:umdsnepv}, so we skip the AA-$\rho$ scheme 
for this test.
The function $H(X)$ is differentiable, with a Fr\'echet derivative
\begin{equation}\label{eq:umdsdirectional}
    \DD H(X)[E] = -\sum_{j=1}^s \left(A_jEX^TA_j+A_jXE^TA_j\right).
\end{equation}
Iterative solvers used in the SCF and Newton steps require
only products involving $H(X)z$ and $\DD H(X)[E]X$,
and they can be evaluated  using the factored expressions in~\eqref{eq:umdsnepv} and~\eqref{eq:umdsdirectional},
without explicitly forming $H(X)$ or the derivative operator.
With the factors $A_jX$ precomputed at each outer iteration,
the computational complexity for $H(X)z$ is about $\mathcal{O}(nsk)$,
while for $\DD H(X)[E]X$ it is about  $\mathcal{O}(n^2sk)$.
Thus, derivative applications are substantially more expensive for large $n$.

\begin{figure}[!htbp]
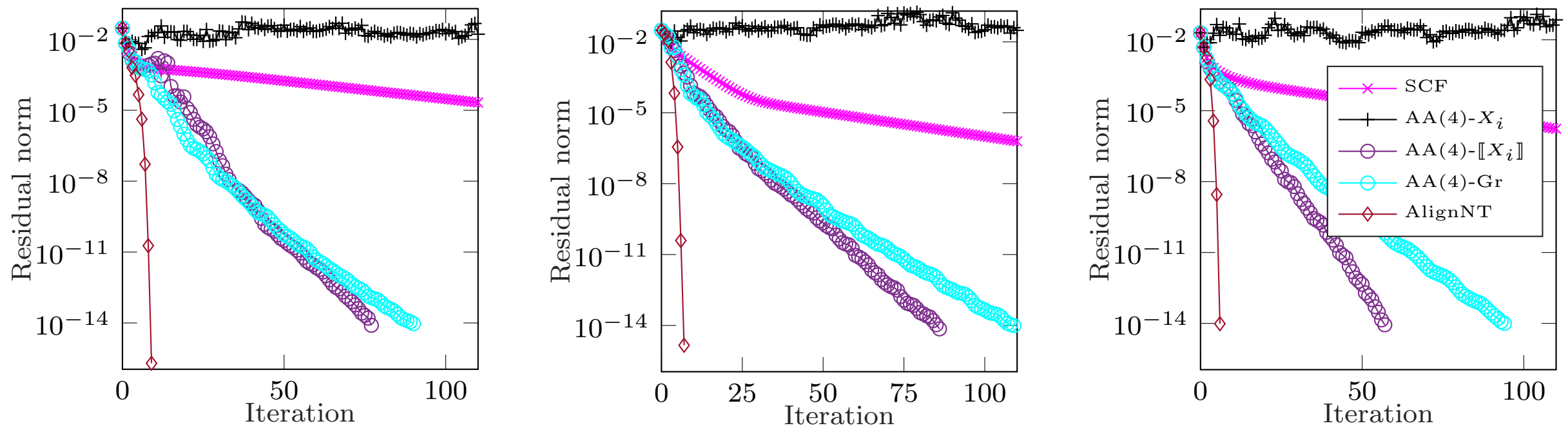

    \centering
    \begin{minipage}{\threeplotscale}
       \include{Figs/tikz/citeseer_iter}
   \end{minipage}
   \begin{minipage}{\threeplotscale}
       \include{Figs/tikz/cora_iter}
   \end{minipage}
   \begin{minipage}{\threeplotscale}
       \include{Figs/tikz/reuters_iter}
   \end{minipage}
   \vspace{-1cm}
	\caption{Convergence histories of the algorithms
    for the UMDS NEPv~\eqref{eq:umdsnepv}
    on the data sets
    CiteSeer (left, $n=3,312$, $k=6$),
    Cora (middle, $n=2,708$, $k=7$), and
    Reuters1500 (right, $n=1,500$, $k=6$).}
    \label{fig:umds}
\end{figure}

For the experiments, we consider three benchmark data sets:
CiteSeer~\cite{GilBL:1998,Sen:2008},
Cora~\cite{McCallum:2000,Sen:2008}, and
Reuters1500~\cite{AmiUG:2009,WanC:2017},
with problem sizes listed in~\Cref{tab:umds}.
We compare the plain SCF, 
the AA variants \AAm-$X$, \AAm-$\ALIGN{X}$, and \AAm-Gr with a fixed window size $m=4$,
and basis-aligned Newton acceleration.
Given the relatively large problem sizes,
we use a setup similar to that in~\Cref{sec:synthc}: warm-started
\texttt{eigs} for the SCF steps and an inexact implementation
of basis-aligned Newton, with the correction equations solved
by flexible GMRES with a maximum Krylov subspace dimension
of $25$, using $40$ MINRES steps for each preconditioner
application.
Since applying $\DD H(X)[E]X$ is substantially more expensive
than $H(X)z$, each MINRES step is much cheaper
than a GMRES step because it avoids derivative evaluations.
All methods use a common initial guess $X_0\in\RR^{n,k}$,
formed from the $k$ eigenvectors corresponding to the smallest
eigenvalues of $-\sum_{j=1}^s A_j$.
The convergence histories and computation times are reported
in~\Cref{fig:umds,tab:umds}, respectively.

In all three experiments, SCF converges slowly, and AA applied
directly to $X_i$ stagnates at the initial residual. 
Both AA-$\ALIGN{X_i}$ and Grassmannian AA substantially accelerate convergence, 
and AA-$\ALIGN{X_i}$ slightly outperforms Grassmannian. 
Newton acceleration also improves the slow convergence of SCF in all of
our tests. 
It requires the least number of iterations among the algorithms,
converged within nine iterations as compared with $55$--$105$ 
for the two convergent AA variants.
In terms of computation time, it is not the fastest, but it is still competitive due to the efficient preconditioner for solving the correction equation. 
We observed that without preconditioning, the plain MINRES can 
require up to four times as much computation time as the 
preconditioned solver to solve a Newton equation 
to the same residual tolerance. 

\begin{table}[!htbp]
{\footnotesize
    \caption{
        Problem sizes and computation times (seconds) for
        the UMDS NEPv~\eqref{eq:umdsnepv} on the data sets
        CiteSeer, Cora, and Reuters1500
        as shown in~\Cref{fig:umds}.
        Here, $n$, $s$, and $k$ denote the numbers of
        samples, views, and clusters, respectively.
        A dash indicates failure to converge within the
        given iteration limit.
    }
    \label{tab:umds}
    \begin{center}
    \begin{tabular}{c|ccc||c|c|c|c|c}
        \toprule
        Data set &  $s$ & $n$ & $k$ & SCF
        & AA($4$)-$X_i$ & AA($4$)-$\ALIGN{X_i}$
        & AA($4$)-Gr & AlignNT \\
        \midrule
        CiteSeer& 2 & 3,312 & 6 & 10.64 & -- & 1.31 & 1.68 & 1.63 \\
        \hline
        Cora& 4 & 2,708 & 7 & 7.93 & -- & 1.59 & 2.13 & 2.18 \\
        \hline
        Reuters1500& 5 & 1,500 & 6 & 3.97 & -- & 0.51 & 0.87 & 0.77  \\
        \bottomrule
    \end{tabular}
    \end{center}
}
\end{table}

\newpage
\section{Concluding Remarks}\label{sec:conclusion}
We studied acceleration techniques of the SCF iteration for unitarily invariant NEPv based on the AA and Newton methods. To address the non-uniqueness of the eigenvector bases,
we presented three AA variants based on invariant
parameterization, basis alignment, and Grassmannian geometry.
We also developed a basis-aligned Newton correction that removes
the singularity caused by unitary invariance and established
its local quadratic convergence under suitable assumptions.
Numerical experiments demonstrated the effectiveness of these
approaches in accelerating the convergence of SCF.

In our experiments, the three AA variants exhibited a similar
convergence behavior and comparable computational costs.
The choice among them can therefore be guided by ease of implementation.
Parameterization-based AA is convenient when a suitable invariant
parameterization is readily available.
Otherwise, aligned AA can be implemented using basis alignment,
while Grassmannian AA is a natural choice when the required
manifold operations are already available.
Newton acceleration is attractive when high accuracy is required
and derivative actions can be evaluated efficiently.
Its local quadratic convergence can reduce the number of outer
iterations, although the additional cost of solving the correction
equation must also be considered.
%

This work  provides a starting point for several directions for further research.
A systematic convergence study of the proposed AA variants
remains open, particularly concerning the influence of 
different choices of ui-parameterization.
The basis-aligned Newton formulation also provides
a theoretical foundation for developing more efficient solvers
through inexact correction and preconditioning.
One could also explore the implementation and
benefits of combining accelerations with stabilization
techniques, such as the level-shifted SCF,
whose global convergence may come at the cost of slow local convergence.
Finally, beyond NEPv, we expect the ideas and methods developed in 
this work to extend naturally to other iterative methods that operate on subspaces, 
offering new ways to accelerate these methods.

\section*{Acknowledgments}
The work of Ding Lu was supported in part by NSF under grant DMS-2110731.\\
Tom Werner acknowledges the hospitality of the Depatment of Mathematics at University of Kentucky during his visits in Spring 2024 and Fall 2025.
\section*{CRediT Author Statement} The authors have contributed to this
work in the following way.\\
{\bf Ding Lu:} Conceptualization, Investigation, Methodology, Software, Supervision, Writing - Original Draft;\\
{\bf Tom Werner:} Conceptualization, Investigation, Methodology, Software, Writing - Original Draft;\\
{\bf Zhaojun Bai} Conceptualization, Supervision, Writing - Review \& Editing.

\newpage 
\bibliographystyle{siam}
\bibliography{references.bib}
\end{document}

%% file: Figs/tikz/ks1d_beta2_iter.tex
%
%
\definecolor{mycolor1}{rgb}{1.00000,0.00000,1.00000}%
\definecolor{mycolor2}{rgb}{0.49400,0.18400,0.55600}%
\definecolor{mycolor3}{rgb}{0.00000,1.00000,1.00000}%
\definecolor{mycolor4}{rgb}{0.63500,0.07800,0.18400}%
\begin{tikzpicture}

\begin{axis}[%
width=3.65cm,
height=3.7cm,
scale only axis,
tick label style={font=\footnotesize},
xmin=0,
xmax=29,
xlabel style={font=\footnotesize\color{white!15!black},yshift=.2cm},
xlabel={Iteration},
ymode=log,
ymin=1.1e-15,
ymax=.95,
yminorticks=true,
ytickten = {-14, -11, -8, -5, -2},
ylabel style={font=\footnotesize\color{white!15!black},yshift=-.2cm},
ylabel={Residual norm},
axis background/.style={fill=white},
legend style={at={(axis cs:24.5,4.5e-1)}, anchor=north east, legend cell align=left, align=left, font = \tiny, draw=white!15!black}
]
\addplot [color=myscf, mark=x, mark options={solid, myscf}]
  table[row sep=crcr]{%
0	0.000770323859631843\\
1	0.000235100471963844\\
2	7.66070530580152e-05\\
3	2.48827820347243e-05\\
4	8.10784923856096e-06\\
5	2.63979300693333e-06\\
6	8.59727542934734e-07\\
7	2.79970684122632e-07\\
8	9.11753272041363e-08\\
9	2.96918932401349e-08\\
10	9.66940627478019e-09\\
11	3.14891818416597e-09\\
12	1.02547004155316e-09\\
13	3.33952367238609e-10\\
14	1.087542268328e-10\\
15	3.54168210202613e-11\\
16	1.15336653137875e-11\\
17	3.75562941308181e-12\\
18	1.22285227053702e-12\\
19	3.98372195543074e-13\\
20	1.29760072912854e-13\\
21	4.21204114887317e-14\\
22	1.36011813875066e-14\\
23	5.06150860972305e-15\\
};
\addlegendentry{SCF}

\addplot [color=black, mark=+, mark options={solid, black}]
  table[row sep=crcr]{%
0	0.000770323859631843\\
1	0.000235100471963844\\
2	0.000194530329626275\\
3	0.000222975655827858\\
4	0.00196980238000191\\
5	0.000350188860646446\\
6	0.000864733133340715\\
7	0.000243935473836303\\
8	0.000498269366599108\\
9	0.000737683033856946\\
10	0.000453764939471481\\
11	0.000125148669740175\\
12	0.000528099457585981\\
13	0.000573864382272188\\
14	0.0002402397320883\\
15	0.000386159728828657\\
16	0.000436790066583778\\
17	0.000238240228876691\\
18	0.000167899418445888\\
19	0.00035002975897477\\
20	0.00101005568470392\\
21	0.000347796615447007\\
22	0.000267062452885241\\
23	0.000269574278769082\\
24	0.000788420502613227\\
25	0.000792977016129419\\
26	0.000829732859522169\\
27	0.000368914575020376\\
28	5.50423593144669e-05\\
29	0.000137728333861786\\
30	4.87961052140253e-05\\
31	2.97848406760406e-05\\
32	2.7662661268313e-07\\
33	1.07956350937213e-05\\
34	4.78561065437638e-06\\
35	3.60645312118917e-06\\
36	1.40120389778017e-05\\
37	6.07968351278431e-05\\
38	0.00194741249680512\\
39	0.000840154515417659\\
40	0.000758536933984174\\
41	0.000816492544599302\\
42	0.000419714820324326\\
43	4.07257989104341e-05\\
44	4.85953717789434e-05\\
45	7.7186814653397e-05\\
46	1.94057756951186e-05\\
47	1.07859365538412e-05\\
48	0.000175992374422036\\
49	0.000210540357196576\\
50	0.000811069529599731\\
51	0.000573788144914607\\
52	0.000948484965204328\\
53	0.000252439579381867\\
54	0.000298340727865951\\
55	0.000378588680557038\\
56	0.00056942060517796\\
57	0.000739545533627655\\
58	0.000885727327399781\\
59	0.000179114192044857\\
60	0.000265248119337006\\
61	0.000159197920879442\\
62	0.00142557880923896\\
63	0.00162921001107669\\
64	0.000265723751716573\\
65	8.5547193985067e-05\\
66	8.5868158942147e-05\\
67	7.48713151138799e-05\\
68	5.22078435196659e-05\\
69	0.00175859274019546\\
70	0.000339621300219476\\
71	0.00034276437892162\\
72	0.000773214142897641\\
73	0.00077590105337104\\
74	0.000766668850443509\\
75	0.000785749326796619\\
76	0.000807104384006309\\
77	7.65739234520139e-05\\
78	0.000398770071537317\\
79	0.000199840586344561\\
80	0.000335151504546845\\
81	5.53027132530932e-05\\
82	4.50025240051925e-05\\
83	8.54422360140527e-05\\
84	0.000150921290693617\\
85	0.00138344925099991\\
86	0.000255768274469182\\
87	0.000165775343260163\\
88	0.000396485693235253\\
89	9.14278718588196e-05\\
90	0.000130566513271283\\
91	0.000122757696645472\\
92	3.64397959112594e-05\\
93	3.51170015607573e-06\\
94	8.12125923409805e-06\\
95	4.63753153749421e-06\\
96	1.63716739929414e-06\\
97	0.000111858498709282\\
98	0.00079466580776264\\
99	0.000827620373141157\\
100	0.000408584696626573\\
101	0.000131263783879601\\
102	0.000130917251497434\\
103	0.00010093372716612\\
104	0.000182483313720959\\
105	0.0012566331056946\\
106	0.000462689472450616\\
107	0.000580930234048726\\
108	0.00077156154497213\\
109	0.000854796506550701\\
110	0.000326605502781047\\
111	0.000766465617379209\\
112	0.000551183321670797\\
113	0.000284568296138046\\
114	0.000486441437126514\\
115	0.000152287698549492\\
116	0.000156911252098231\\
117	0.000159576709237445\\
118	8.27142546433418e-05\\
119	0.000166792402165288\\
120	0.000820838869846683\\
121	0.000234880465885209\\
122	0.000195081266999295\\
123	0.000787093270901294\\
124	0.000638715871937528\\
125	0.000785061507772257\\
126	0.000259233006838257\\
127	0.000762847295128826\\
128	0.000239294306620611\\
129	0.000675531898179355\\
130	0.000243785765759799\\
131	7.15406610465764e-05\\
132	3.15789605067343e-05\\
133	0.000192095405688146\\
134	0.00248538944514365\\
135	0.000206543818447572\\
136	0.000652650413395735\\
137	0.000615896698378768\\
138	0.000824640598805005\\
139	0.000466044997446678\\
140	0.000143769717587296\\
141	0.000153497196723659\\
142	0.000143513526699384\\
143	5.05169618255798e-05\\
144	0.000101781812683839\\
145	0.00049891500655622\\
146	9.12306022576162e-05\\
147	7.91414403016776e-05\\
148	7.8634345954807e-05\\
149	2.67431993669447e-05\\
150	0.000131203529463195\\
151	0.000855002626225324\\
152	0.00083095034812711\\
153	0.000465352252933499\\
154	0.000219534661065864\\
155	0.000306909335894619\\
156	0.000601187219424846\\
157	0.000882026023804741\\
158	0.000763639113559364\\
159	0.000833310902447221\\
160	0.00082383180616787\\
161	0.000285415625657341\\
162	9.21236982052392e-05\\
163	9.26190471726748e-05\\
164	8.95500135433408e-05\\
165	3.42831406740378e-05\\
166	6.92302734510758e-05\\
167	0.00163740182573512\\
168	0.00108726723924748\\
169	0.000820600555997474\\
170	0.00063721230089166\\
171	0.000595899637904754\\
172	0.00074526956461007\\
173	0.000515056449744708\\
174	0.000261547660261379\\
175	0.00070654380840641\\
176	0.000721547568836131\\
177	0.000261226074123833\\
178	0.000222935313986242\\
179	0.000239288708156169\\
180	2.8683356548748e-06\\
181	3.78276246661752e-05\\
182	8.74750146760299e-05\\
183	0.00151191800155411\\
184	0.000776812558869899\\
185	0.000718706953610494\\
186	0.000367355636618228\\
187	0.000104223030105537\\
188	1.91062336700676e-05\\
189	7.91746833696338e-05\\
190	0.000162945547656573\\
191	1.50278733303757e-05\\
192	1.45081210826131e-05\\
193	5.82069037004519e-06\\
194	6.51903561656576e-06\\
195	0.000126145325059279\\
196	0.00181874536109081\\
197	0.000630096838275669\\
198	0.000674411072607385\\
199	0.000758203655472706\\
200	0.000424041591264581\\
201	9.63827226697231e-05\\
202	0.000467498828368874\\
203	0.000800472626081979\\
204	0.000795290060174521\\
205	0.000689920909997268\\
206	0.000223085143857447\\
207	0.000129653702435809\\
208	1.03483331202478e-05\\
209	2.3796992290777e-05\\
210	7.19099099792938e-05\\
211	0.00201043874759372\\
212	5.8031895361236e-05\\
213	4.86850133556224e-05\\
214	0.000105436653159164\\
215	6.22231335101647e-05\\
216	0.000323681953715639\\
217	0.000938757440907065\\
218	0.000835116283098821\\
219	0.000907965905417375\\
220	0.000665089468807735\\
221	0.000812349685951445\\
222	0.000822577444017357\\
223	0.000403637161427629\\
224	0.000130356710564396\\
225	0.000136680420090661\\
226	0.000133812335823142\\
227	4.81968432997566e-05\\
228	1.4889872959659e-05\\
229	1.56822309693328e-05\\
230	1.40855080527353e-05\\
231	5.04588795749458e-06\\
232	0.000177396998004348\\
233	0.000800286929096501\\
234	0.000606590215685207\\
235	0.000191694057310808\\
236	0.00019080120742656\\
237	0.000189041248983271\\
238	6.04339982288138e-05\\
239	8.86789694244985e-05\\
240	0.000187417388453441\\
241	1.87625509470456e-05\\
242	0.000484493322074632\\
243	7.35893351469001e-05\\
244	0.000137158428321954\\
245	0.000159493923727858\\
246	0.000792512362417037\\
247	0.000822085510683446\\
248	0.00060412124541774\\
249	0.000989444643148488\\
250	0.000393703889763964\\
251	0.000785090695241321\\
252	0.000819305694626922\\
253	0.000279502018382544\\
254	0.000723482771899075\\
255	0.000692004066891189\\
256	0.000120143443400331\\
257	6.0749710437293e-05\\
258	0.000134441468186502\\
259	8.91669564954353e-05\\
260	0.000124656965765565\\
261	0.00262280487712785\\
262	0.000204948927184712\\
263	0.000195931901444878\\
264	0.000103652268375512\\
265	0.000533684285113811\\
266	0.000686061777647098\\
267	0.000179729352974629\\
268	0.000355493108701419\\
269	9.1296396743355e-05\\
270	0.000118279096334474\\
271	0.000110640235743937\\
272	3.51096389811436e-05\\
273	2.74179515039552e-06\\
274	7.08005310643051e-06\\
275	4.27081159632976e-06\\
276	5.51719162622023e-07\\
277	3.38436419660009e-05\\
278	2.17870331830677e-06\\
279	0.000200929860509066\\
280	0.000810232464707313\\
281	0.000690228981012138\\
282	0.000345989783234245\\
283	0.000830079671283899\\
284	0.000698133675322135\\
285	0.000636850924459558\\
286	0.000771713334991663\\
287	0.000156876311855983\\
288	0.000145617288550304\\
289	0.000347701220577992\\
290	0.0015893431967179\\
291	0.000389555982910792\\
292	0.000725824326314416\\
293	2.3986754138891e-05\\
294	0.000829468458203548\\
295	0.000827329406816712\\
296	0.0010075181841494\\
297	0.000281002926690972\\
298	0.00083267301617565\\
299	0.000602759940058518\\
300	0.000741143842781557\\
301	0.00083212543684175\\
302	0.000771867059158846\\
303	7.11725318136563e-05\\
304	0.000427152773717271\\
305	6.90490948460867e-05\\
306	0.000247627260697342\\
307	6.74257790504569e-05\\
308	8.46104320679726e-05\\
309	5.76560051862509e-05\\
310	6.31166128659441e-05\\
311	0.00258993141320829\\
312	2.66054574456141e-05\\
313	3.01849333447483e-05\\
314	8.83732204592202e-06\\
315	0.000109033156923515\\
316	0.000853920324814944\\
317	0.000274472507360436\\
318	0.000673052346158688\\
319	0.000532784691626252\\
320	0.000496507696903617\\
321	0.000840029198849619\\
322	0.000835330254533014\\
323	0.000365919574112853\\
324	0.000348863640123437\\
325	0.000822907933523636\\
326	0.00081997465574572\\
327	0.000286736574555794\\
328	0.000534239538144219\\
329	0.000591891080884562\\
330	0.000758114548003767\\
331	0.000757380173923147\\
332	0.000594738723581307\\
333	0.000389702830417924\\
334	0.000596649922179977\\
335	0.000289669249757615\\
336	0.00044120236070263\\
337	0.000830629802662687\\
338	0.000827677280561076\\
339	0.00013693027706095\\
340	0.00073652693590599\\
341	0.000604004415447547\\
342	0.000406395288915059\\
343	8.39238376830353e-05\\
344	3.63942269969065e-05\\
345	0.00011983569371266\\
346	0.000236759428714112\\
347	7.75650820684362e-05\\
348	7.78020048188411e-05\\
349	7.50054200139551e-05\\
350	0.00172291057035395\\
351	0.000835034028057148\\
352	0.000144675993374693\\
353	0.000851008196900371\\
354	0.000826478306284308\\
355	0.000835540381047852\\
356	4.81550565073308e-05\\
357	0.000427376639737698\\
358	0.000431630823001981\\
359	0.000846076580640006\\
360	0.000880189138418033\\
361	0.000714296307946351\\
362	0.000120148607646737\\
363	0.000203892416573462\\
364	0.000497607056254398\\
365	0.000825343188646153\\
366	0.000579904549278266\\
367	0.000577188996309\\
368	0.000108707207607644\\
369	0.000101139551316079\\
370	0.000406435546517407\\
371	0.000821617131504892\\
372	0.000446243999943513\\
373	0.000791836865364619\\
374	0.000802462990871314\\
375	0.000804535271361851\\
376	0.000159814609903155\\
377	0.000292625523728022\\
378	0.000260183272259731\\
379	0.00011811630043685\\
380	6.61962073514846e-05\\
381	7.65147550449076e-05\\
382	0.000550655440310322\\
383	0.000181228894647543\\
384	0.000182062338517451\\
385	0.000179435795613703\\
386	4.82529761784008e-05\\
387	0.0017636080527294\\
388	3.00140381180591e-05\\
389	2.54529836216118e-05\\
390	0.000466823379805108\\
391	8.65009218715362e-05\\
392	0.000133260547739415\\
393	0.000136675555968053\\
394	0.000804185742986188\\
395	0.000822952462281474\\
396	0.000506481227201403\\
397	0.000218135640010396\\
398	0.000138288698905706\\
399	0.000208294608925505\\
400	8.78085713723233e-05\\
401	0.000823757868285912\\
402	0.000344349392036054\\
403	0.000861663103286265\\
404	0.000742372727572005\\
405	0.000808525373541941\\
406	0.0001315669590437\\
407	0.000222899235545262\\
408	0.000472817816734794\\
409	0.000699139540027147\\
410	0.000846455767111214\\
411	0.000976596164428727\\
412	0.000245203907522165\\
413	0.000112275155517036\\
414	0.000627864288437464\\
415	0.000271887544153979\\
416	2.4876096128101e-05\\
417	9.83548046096154e-06\\
418	8.29387908090276e-05\\
419	0.000103111148243414\\
420	0.000422264553837\\
421	0.000226666993735976\\
422	0.000368241440967938\\
423	9.25315528873894e-05\\
424	0.00209145849926838\\
425	0.000111442122941167\\
426	0.000654829472092615\\
427	0.000455782475242119\\
428	0.000605984703114599\\
429	0.00019896739627721\\
430	0.000303735521023522\\
431	0.000172916380170785\\
432	0.000836293062670547\\
433	0.000613874473662871\\
434	0.000825547090047516\\
435	0.000130963192224006\\
436	0.000374449014525527\\
437	0.000788798186372868\\
438	0.000474682267664034\\
439	0.000287959653205176\\
440	7.20762746435861e-06\\
441	0.000137857185727095\\
442	0.000413801949845025\\
443	5.40358056121036e-05\\
444	4.48780062413208e-05\\
445	0.000115056348080695\\
446	2.27171483078016e-05\\
447	2.35341076880803e-05\\
448	3.75733462285231e-05\\
449	9.36165739442518e-05\\
450	0.000560194831891068\\
451	0.000182218282172424\\
452	0.000182327230759158\\
453	0.000178657332870624\\
454	5.24433368586866e-05\\
455	9.32818434215431e-05\\
456	0.000856122818432684\\
457	0.00159422422644655\\
458	0.000294779722419751\\
459	0.000756109488856596\\
460	0.000827140613204418\\
461	0.000820894511518715\\
462	0.000375375032355633\\
463	0.000200581659982786\\
464	0.000812267651040154\\
465	0.000744351358020888\\
466	0.000823155996273029\\
467	0.000745582523909725\\
468	0.000547117014385461\\
469	0.000828163310426938\\
470	0.000832406803266605\\
471	0.000159697032111788\\
472	0.000149243640506581\\
473	0.000388876174355238\\
474	0.000117502123608814\\
475	0.000203229734560536\\
476	7.42510438206632e-05\\
477	0.000510418769404473\\
478	0.000102345329580527\\
479	0.00016937816194881\\
480	0.00039778473153594\\
481	0.00084529306905091\\
482	0.000570133496436044\\
483	0.000212025596381132\\
484	0.000150349586637571\\
485	0.000224776747376289\\
486	6.52384449577424e-05\\
487	0.000843931101294823\\
488	0.000588381785139491\\
489	0.000193121691415837\\
490	0.000152298810698766\\
491	0.000224901758062524\\
492	8.79805529379069e-05\\
493	9.48514065532163e-05\\
494	0.000775456132690454\\
495	0.000258176254986764\\
496	0.000423354679269296\\
497	0.000778512033930737\\
498	0.000562317818130794\\
499	0.000220086648629031\\
500	0.000133466022892284\\
501	7.16772418773712e-05\\
502	0.00011448625191821\\
503	0.000276730616403787\\
504	0.000791649630882001\\
505	0.000820812255715876\\
506	0.00036009555604599\\
507	0.000251762683854879\\
508	0.000149503734362242\\
509	9.96861816017953e-05\\
510	0.000112697386523521\\
511	7.95757071756459e-05\\
512	0.000232770499181574\\
513	0.000355137436006836\\
514	0.000244782381048205\\
515	0.00245831669979899\\
516	0.000611103833992165\\
517	0.000970732612111886\\
518	0.0024584879736683\\
519	0.0001757825968747\\
520	0.000114766931484549\\
521	4.42624609810343e-05\\
522	0.000187797107724533\\
523	0.000847463892539389\\
524	0.000550370340408472\\
525	0.000170930114798641\\
526	0.000157869779416978\\
527	6.96221942968882e-05\\
528	0.000215504433152424\\
529	0.00113751637038249\\
530	0.000440922069900809\\
531	0.000592911201182956\\
532	0.00038930125440699\\
533	0.0001227285866259\\
534	0.00012262603638785\\
535	0.000126358970094546\\
536	6.80135653491399e-05\\
537	2.15425598275583e-05\\
538	2.20409351222341e-05\\
539	2.20320987105557e-05\\
540	6.91799805211162e-06\\
541	0.000104470353556897\\
542	0.000262291827580715\\
543	0.000131441291630379\\
544	0.000856090394507088\\
545	0.000195156141657958\\
546	6.36221940485376e-05\\
547	6.36331831786457e-05\\
548	6.09322972466873e-05\\
549	2.31331596568779e-05\\
550	7.33146454327213e-05\\
551	0.00202503663718384\\
552	0.000722781775560727\\
553	0.000497763466700695\\
554	0.000422974661776917\\
555	0.000833069038373145\\
556	0.000365620228610956\\
557	0.000800870198774442\\
558	0.000854703449005192\\
559	0.000127203495181082\\
560	4.09375777175452e-05\\
561	4.14674626805267e-05\\
562	4.16263836260521e-05\\
563	1.39436042362507e-05\\
564	4.74722355589223e-06\\
565	9.73514066621428e-07\\
566	9.56588145268698e-05\\
567	0.000772579931952093\\
568	0.00059713425313376\\
569	0.000942168050303861\\
570	0.000261328945974884\\
571	0.000639940304753864\\
572	0.000799074950591314\\
573	0.000814439308978295\\
574	0.000647820648916979\\
575	0.000218214419669507\\
576	0.000126185012560582\\
577	6.60337987382133e-06\\
578	2.37170048987133e-05\\
579	7.37699219009965e-05\\
580	0.00199501071070843\\
581	3.25765878181438e-05\\
582	4.62055577152041e-06\\
583	4.6218420419067e-06\\
584	4.60810017988432e-06\\
585	1.40886202276242e-06\\
586	7.88977578493806e-07\\
587	0.000102589807980174\\
588	5.59931165639819e-05\\
589	0.000756079992411781\\
590	0.000651553100877572\\
591	0.000120265192583375\\
592	0.000576943838624898\\
593	0.00076436947317016\\
594	0.000173913524319467\\
595	5.29969508904388e-05\\
596	5.48984815545992e-05\\
597	5.47465318433025e-05\\
598	1.80001649918045e-05\\
599	9.77503048380318e-05\\
};
\addlegendentry{AA(2)-$X_i$}

\addplot [color=blue, mark=o, mark options={solid, blue}]
  table[row sep=crcr]{%
0	0.000770323859631843\\
1	0.000235100471963844\\
2	5.78580095816898e-05\\
3	1.58559292803132e-07\\
4	7.25296700024834e-08\\
5	7.22029643071963e-11\\
6	1.32502621743247e-13\\
7	3.19341370262166e-15\\
};
\addlegendentry{AA(2)-$\rho_i$}
\addplot [color=myalign, mark=o, mark options={solid, myalign}]
  table[row sep=crcr]{%
0	0.000770323859631845\\
1	0.000235100471963843\\
2	3.67831527957807e-06\\
3	1.3996059120468e-07\\
4	4.16094702583693e-09\\
5	8.58866965116244e-10\\
6	6.42991394474939e-12\\
7	8.6424051770249e-13\\
8	3.29317175459099e-14\\
9	5.19240534817318e-15\\
};
\addlegendentry{AA(2)-$\ALIGN{X_i}$}

\addplot [color=mygrassaa, mark=o, mark options={solid, mygrassaa}]
  table[row sep=crcr]{%
0	0.000770323859631843\\
1	0.000235100471963844\\
2	0.000106307415592508\\
3	4.71451680975027e-07\\
4	6.4642781650735e-08\\
5	2.14870101981433e-10\\
6	9.77864548928054e-12\\
7	7.49597985766571e-13\\
8	3.18616947778094e-15\\
};
\addlegendentry{AA(2)-Gr}

\addplot [color=mygnt, mark=diamond, mark options={solid, mygnt}]
  table[row sep=crcr]{%
0	0.000770323859631843\\
1	0.000235100471963844\\
2	7.59208853612306e-07\\
3	9.79896767924028e-13\\
4	3.33599687548767e-15\\
};
\addlegendentry{AlignNT}
\legend{}
\end{axis}
\end{tikzpicture}%

%% file: Figs/tikz/ks1d_beta4_iter.tex
%
%
\definecolor{mycolor1}{rgb}{1.00000,0.00000,1.00000}%
\definecolor{mycolor2}{rgb}{0.49400,0.18400,0.55600}%
\definecolor{mycolor3}{rgb}{0.00000,1.00000,1.00000}%
\definecolor{mycolor4}{rgb}{0.63500,0.07800,0.18400}%
\begin{tikzpicture}

\begin{axis}[%
width=3.65cm,
height=3.7cm,
scale only axis,
tick label style={font=\footnotesize},
xmin=0,
xmax=29,
xlabel style={font=\footnotesize\color{white!15!black},yshift=.2cm},
xlabel={Iteration},
ymode=log,
ymin=1.1e-15,
ymax=.95,
yminorticks=true,
ytickten = {-14, -11, -8, -5, -2},
ylabel style={font=\footnotesize\color{white!15!black},yshift=-.2cm},
ylabel={Residual norm},
axis background/.style={fill=white},
legend style={at={(axis cs:24.5,9e-1)}, anchor=north east, legend cell align=left, align=left, font = \scriptsize, draw=white!15!black}
]
\addplot [color=myscf, mark=x, mark options={solid, myscf}]
  table[row sep=crcr]{%
0	0.00154064771926348\\
1	0.000820403789956607\\
2	0.000551184268782526\\
3	0.000338491097098383\\
4	0.000227627397730847\\
5	0.000145642344085184\\
6	9.67040298698267e-05\\
7	6.27787902353876e-05\\
8	4.13937718326311e-05\\
9	2.70253037719911e-05\\
10	1.77614118851276e-05\\
11	1.16232671130258e-05\\
12	7.6279356827437e-06\\
13	4.99672682095722e-06\\
14	3.27710868532608e-06\\
15	2.14759178415054e-06\\
16	1.40811699651372e-06\\
17	9.22948956769083e-07\\
18	6.05081337176121e-07\\
19	3.96630574959445e-07\\
20	2.60016169573584e-07\\
21	1.70446139035877e-07\\
22	1.11735690169666e-07\\
23	7.3246167116659e-08\\
24	4.80159585220936e-08\\
25	3.14761242957439e-08\\
26	2.06338469143644e-08\\
27	1.3526236792489e-08\\
28	8.86696942768711e-09\\
29	5.81262736561113e-09\\
30	3.81039844073698e-09\\
31	2.49785858558066e-09\\
32	1.63744114321969e-09\\
33	1.07340453489828e-09\\
34	7.03656399566968e-10\\
35	4.61272860012621e-10\\
36	3.0238155353021e-10\\
37	1.98222269183499e-10\\
38	1.29942245976218e-10\\
39	8.5182518330937e-11\\
40	5.58404818736134e-11\\
41	3.66059579989552e-11\\
42	2.39970979509094e-11\\
43	1.57305661887013e-11\\
44	1.03118009884564e-11\\
45	6.75957111111059e-12\\
46	4.43046952984354e-12\\
47	2.90433012978088e-12\\
48	1.90372460566709e-12\\
49	1.24811242307001e-12\\
50	8.18439075849863e-13\\
51	5.36849545094339e-13\\
52	3.51872289280142e-13\\
53	2.3078632208769e-13\\
54	1.51366719801472e-13\\
55	9.88742907070426e-14\\
56	6.4270451835797e-14\\
57	4.19503878322078e-14\\
58	2.73758718333315e-14\\
59	1.80663964564352e-14\\
60	1.16825389706596e-14\\
61	8.3122643039354e-15\\
};
\addlegendentry{SCF}

\addplot [color=black, mark=+, mark options={solid, black}]
  table[row sep=crcr]{%
0	0.00154064771926348\\
1	0.000820403789956607\\
2	0.000486143061109116\\
3	0.000482349029208183\\
4	0.000209218683320886\\
5	0.000300987227165901\\
6	0.00061604801274812\\
7	0.000902493646599572\\
8	0.000533254017092038\\
9	0.000749890775306338\\
10	0.000802938918700671\\
11	0.000844081982063105\\
12	0.000656121357030697\\
13	9.37780820237866e-05\\
14	4.08546985241304e-05\\
15	0.00022683819805547\\
16	0.000251069195852655\\
17	0.000722217254833257\\
18	0.000767292128902441\\
19	0.000369761344224301\\
20	0.000475029478284957\\
21	0.000458907108090548\\
22	0.000341129965193969\\
23	0.00028191693295966\\
24	0.000823360687843079\\
25	0.000428161353019195\\
26	0.000470632103155698\\
27	0.00038428532610692\\
28	0.000590001836512939\\
29	0.000554867369626976\\
30	0.000941109038995465\\
31	0.00092265026082716\\
32	0.000842429844552567\\
33	0.00042408611299224\\
34	0.000538195975710088\\
35	0.000545848332815762\\
36	0.000132692528571263\\
37	0.000819700395499147\\
38	0.00103079626828335\\
39	6.82522088205004e-05\\
40	8.52064287740001e-05\\
41	0.000544441219393617\\
42	0.000324916187640903\\
43	0.000379644535190722\\
44	0.000352515790824068\\
45	0.000783372112938071\\
46	0.000238007195114406\\
47	0.000210368617191935\\
48	0.000606583139456826\\
49	0.000912430141043577\\
50	0.000436754034067883\\
51	0.000419294644917477\\
52	0.000335383027664697\\
53	0.000364478545141952\\
54	0.000236034380391386\\
55	0.000239789362582087\\
56	6.2769151362395e-05\\
57	0.000155119890601252\\
58	0.00108769806593543\\
59	0.000751848101732391\\
60	0.000530787287749201\\
61	0.000675453598199989\\
62	0.000419377848662062\\
63	0.00041420693270624\\
64	0.000456801823625941\\
65	0.000145292416883002\\
66	0.00136566546616317\\
67	0.000596643381271662\\
68	0.000428533815065147\\
69	0.000416358746280488\\
70	0.000394521926505224\\
71	0.000209807142270518\\
72	0.000841056169316315\\
73	0.000502080852248089\\
74	0.000698805247402359\\
75	0.00074981550967733\\
76	0.000814731703662059\\
77	0.000615956500216147\\
78	0.000803394105267307\\
79	0.00049325925041709\\
80	0.000360805161269384\\
81	0.00107234311041138\\
82	0.000503648006932504\\
83	0.000629347526722618\\
84	0.00042608072056725\\
85	0.00014365949730293\\
86	0.000351353063328305\\
87	0.000326587208892804\\
88	0.000317700128705598\\
89	0.000171802640686214\\
90	3.67939311126727e-05\\
91	0.000216570985850826\\
92	0.000419227511641398\\
93	0.00135410583330962\\
94	0.000552989232144579\\
95	0.000761001359938921\\
96	0.000716147656199244\\
97	0.00070657098864875\\
98	0.000608936861445986\\
99	0.00044175957844099\\
100	0.000762003465508534\\
101	0.0010812600122221\\
102	0.000646799383078743\\
103	0.000662523154074209\\
104	0.000730178388963427\\
105	0.000846996238846246\\
106	0.00111640961215525\\
107	0.000928020218396406\\
108	0.000724671401269154\\
109	0.000559764889864906\\
110	0.000281968116458894\\
111	0.000792316999752912\\
112	0.0010504382110911\\
113	0.000156659553087245\\
114	0.00017808611478724\\
115	0.000517147529388146\\
116	0.00061623664681942\\
117	0.00106221291543064\\
118	0.000271303896371207\\
119	0.000754673228323326\\
120	0.000773696648007801\\
121	0.000977781630387073\\
122	0.000308053016363396\\
123	0.00188332967832245\\
124	0.00073931700057355\\
125	0.000194405510646944\\
126	0.000931992737640362\\
127	0.000202829454223956\\
128	0.00110495185113732\\
129	0.000439189131240559\\
130	0.000909026550438484\\
131	0.001152580666629\\
132	0.000705907262711277\\
133	0.00067926861454487\\
134	0.000754792758900378\\
135	0.00080748829463293\\
136	0.000377517180661276\\
137	0.000417544553752117\\
138	0.000170612007222504\\
139	0.000117005006458427\\
140	0.000236314535760237\\
141	0.000520308449474932\\
142	0.00105266392660768\\
143	0.00117426590526177\\
144	0.000837025981993508\\
145	0.000816440099738415\\
146	0.000523290866404942\\
147	0.00118977163163743\\
148	0.000297649163578956\\
149	0.000878411055252627\\
150	0.000944581487321883\\
151	0.00057581829618506\\
152	0.000693399625528461\\
153	0.000785788448076582\\
154	0.000298745063368959\\
155	0.000289723539405162\\
156	0.000287510166337886\\
157	0.000240856228444382\\
158	0.000290786682846173\\
159	0.000378279196861601\\
160	0.000308747722675115\\
161	0.000491013225247663\\
162	0.000298250037005744\\
163	0.000455184965201075\\
164	0.000478697550340134\\
165	0.000657751772688549\\
166	0.00072464348283727\\
167	0.00112746943770092\\
168	0.00096027449148677\\
169	0.000899073802674826\\
170	0.000882414882778269\\
171	0.00064391227829883\\
172	0.000467649205699226\\
173	0.000445332062848586\\
174	0.000428343786224323\\
175	0.000254030279389645\\
176	0.000225504796801221\\
177	0.000131396355468885\\
178	0.000413198344707804\\
179	0.000462450607351263\\
180	0.000320211116764978\\
181	0.000311667886979031\\
182	0.000320094160141736\\
183	3.5807334155453e-05\\
184	1.94441583674473e-05\\
185	1.93415464370229e-05\\
186	1.91084107551924e-05\\
187	8.38003023954557e-05\\
188	0.000304355837578013\\
189	0.000773742398023506\\
190	0.000861940630328887\\
191	0.000284432379973485\\
192	0.000403186099310432\\
193	0.000197683603167287\\
194	0.00021107959950679\\
195	9.65256561317523e-05\\
196	2.57922758632081e-05\\
197	0.000122443063848168\\
198	0.000187157641978047\\
199	0.000364465054004267\\
200	0.000420365889243946\\
201	0.00077505862689404\\
202	0.0010670579220941\\
203	0.000804836147949099\\
204	0.000855719064316138\\
205	0.000773456551705558\\
206	0.000537104137860609\\
207	0.000231707794277977\\
208	0.000790872386354282\\
209	0.000387220039078705\\
210	0.000726663786747779\\
211	0.00068859252470414\\
212	0.000504851921915523\\
213	0.000733930610475888\\
214	0.000488957682301413\\
215	0.00066483362296982\\
216	0.00101482283906066\\
217	0.00127168902432449\\
218	0.000836260657778637\\
219	0.000551253952945735\\
220	0.000297807708615018\\
221	0.000124153265322058\\
222	0.000391612426793453\\
223	0.000860488807920996\\
224	0.000646923131803084\\
225	0.000664017993225785\\
226	0.000514404909298336\\
227	0.000467158856664904\\
228	0.000982967708290161\\
229	0.000570418958832636\\
230	0.000695950766949056\\
231	0.000463272175670861\\
232	0.000276264103049835\\
233	0.000816129144238659\\
234	0.000287784606814961\\
235	0.000501310667501591\\
236	0.000546077013278619\\
237	0.000415274039186435\\
238	0.000452400413234065\\
239	0.000203901715652666\\
240	0.000100104142621546\\
241	9.93259632284371e-05\\
242	0.000110521423928271\\
243	5.75689235975144e-05\\
244	0.000262464606050751\\
245	0.000169066289777615\\
246	0.000164646875316217\\
247	5.9527774783432e-05\\
248	0.000462884819721461\\
249	0.000321880187937844\\
250	0.00135377721766531\\
251	0.000764142974238095\\
252	0.00057721445608748\\
253	0.00055105095795824\\
254	0.000514270598125435\\
255	0.00091299339761214\\
256	0.000945302065674576\\
257	0.000254393288808307\\
258	0.00019038837294333\\
259	0.000405327724816172\\
260	0.000471746112642788\\
261	0.000343435100805799\\
262	0.000352006838890399\\
263	0.000297070552776239\\
264	0.000953701702139023\\
265	0.000736332685634572\\
266	0.000748367713334908\\
267	0.000709571031776641\\
268	0.000859585090125051\\
269	0.000915528866557011\\
270	0.000657143537326207\\
271	0.000832679185179115\\
272	0.000964423435517114\\
273	0.000943983762921614\\
274	0.000999272461214142\\
275	0.000613648297070202\\
276	0.000429369450782952\\
277	0.000423009513196464\\
278	0.000418506550998067\\
279	0.000261886835588749\\
280	0.000183525407217219\\
281	0.000165788191821638\\
282	0.000157216473027348\\
283	0.00108452232349636\\
284	0.000660529791409244\\
285	0.000567250512955648\\
286	0.000585701078635652\\
287	0.000949320600768027\\
288	0.00115861344573187\\
289	0.000399928668864126\\
290	0.00036773344776817\\
291	0.000780830969493097\\
292	0.000562982866023474\\
293	0.000515934122986196\\
294	0.000577607785483949\\
295	0.000297411643111472\\
296	0.000657821795424488\\
297	0.000239402196560886\\
298	0.000228021981340326\\
299	0.000279153637703145\\
300	0.00140385585823034\\
301	0.000278286047222572\\
302	0.000142970444605229\\
303	0.000909386491785271\\
304	0.000297336959182782\\
305	0.000422184629697439\\
306	0.000383276373627631\\
307	0.000171193519502998\\
308	0.000758472361516984\\
309	0.000381960814628247\\
310	0.000471344211170592\\
311	0.000657515387099566\\
312	0.000194033803221666\\
313	0.000223066471466303\\
314	0.000295617452437303\\
315	0.000248938836431941\\
316	0.000177233371637479\\
317	0.000237700351760481\\
318	0.000377414526278804\\
319	0.000320395040829075\\
320	0.000239564191710325\\
321	0.000188260740351593\\
322	0.000985057129776908\\
323	0.000505260006271185\\
324	0.000106575140464199\\
325	0.000512628523306142\\
326	0.000177025411002482\\
327	0.000179046854775379\\
328	0.000106046477946148\\
329	1.90682374392236e-05\\
330	0.00100112139538692\\
331	0.000367806574265425\\
332	0.000944676065782322\\
333	0.000742550306426441\\
334	0.000548608600191704\\
335	0.000168647826596736\\
336	0.000775187699247723\\
337	5.34314903140204e-05\\
338	6.32598831805698e-05\\
339	0.000812276385932935\\
340	0.000836122918609097\\
341	0.00141350195175795\\
342	0.000692351887132718\\
343	0.000786846292503369\\
344	0.000291348762938213\\
345	0.000755164554465579\\
346	0.000942529784309207\\
347	0.000829356380322616\\
348	0.000783619755960869\\
349	0.000744186464373462\\
350	0.000611023783922404\\
351	0.000585084856829928\\
352	0.000513402026740164\\
353	0.000382949149595956\\
354	0.000327994083010578\\
355	0.000370740258920394\\
356	6.15412756255473e-05\\
357	1.14078427036448e-05\\
358	1.46456563152464e-05\\
359	9.39048468424648e-06\\
360	0.000657079507185096\\
361	0.00046774047379303\\
362	0.000365894223225993\\
363	0.000664871408464755\\
364	0.0010824605226545\\
365	0.000697233752869409\\
366	0.000685434467336224\\
367	0.000573647283625521\\
368	0.000555260577505467\\
369	0.000156438164661915\\
370	0.000247338053307304\\
371	0.000994990329187445\\
372	0.000787857470472215\\
373	0.000645925336866601\\
374	0.000167786359193886\\
375	0.000128995117710285\\
376	0.000991179259706629\\
377	0.00038690957666102\\
378	0.000720174191651556\\
379	0.00050932639542258\\
380	0.00100724920071669\\
381	0.00108931751770993\\
382	0.00052404231375288\\
383	0.000270525602522155\\
384	0.000194943428246237\\
385	0.00016921009077606\\
386	0.000176622852603834\\
387	0.000139137185867867\\
388	9.0806488949477e-05\\
389	8.01288615639012e-05\\
390	9.16451030161509e-05\\
391	0.000556284437520033\\
392	0.00028771626011227\\
393	0.000149176661238236\\
394	0.000511400276719866\\
395	0.000913971799135005\\
396	0.000497726245514161\\
397	0.000563783327909593\\
398	0.000307514011065361\\
399	0.0009865037162948\\
400	0.000571081529479937\\
401	0.000577026721135808\\
402	0.000542503187677696\\
403	0.000700529879532812\\
404	0.000406029969519024\\
405	0.000421775079498933\\
406	0.000423622768329309\\
407	0.000241483744804317\\
408	0.000175185801169815\\
409	0.000129466639105735\\
410	0.000148005147724521\\
411	9.99595481040723e-05\\
412	0.000701330339559901\\
413	0.000965949997852173\\
414	0.000611773206034724\\
415	0.000565835705634625\\
416	0.000584500630010328\\
417	0.000763872398484608\\
418	0.000594400430811385\\
419	0.000694570982777729\\
420	0.000748941888939025\\
421	0.000392410937172054\\
422	0.000453827156243693\\
423	0.000804417537908087\\
424	0.00103851084165674\\
425	0.000597895557470196\\
426	0.000266434938596441\\
427	0.000532276082782737\\
428	0.00108513588306164\\
429	0.000381211256405758\\
430	0.00248186492444745\\
431	0.000374193680223114\\
432	0.000453305449778348\\
433	0.000662686255662137\\
434	0.000400943568047324\\
435	0.000687032030785669\\
436	0.000454836428192821\\
437	0.000350344920356814\\
438	0.000225498380219722\\
439	0.000628542573427988\\
440	0.00115777761821327\\
441	0.000349207823844261\\
442	0.000403733393580014\\
443	0.000959791904086599\\
444	0.000235645465078131\\
445	0.00098687760130534\\
446	0.000599554366175568\\
447	0.00096360648488454\\
448	0.00113646002673868\\
449	0.000709231896393872\\
450	0.000682965352526202\\
451	0.000731342171248313\\
452	0.000961073222449253\\
453	0.000157778489315541\\
454	0.000818935255020524\\
455	0.000387593617879631\\
456	0.00135074104925458\\
457	0.000280203936885925\\
458	0.000335793805502808\\
459	0.000444252408410807\\
460	0.000514759126487112\\
461	0.00107947114248943\\
462	0.000599549944349568\\
463	0.000105007026583719\\
464	0.000701190695191385\\
465	0.000741795405232761\\
466	0.00050218957808776\\
467	0.000739744070844272\\
468	0.000878892712554314\\
469	0.000212542119847368\\
470	0.00108233217126958\\
471	0.000982752826749084\\
472	0.000442710466247246\\
473	0.000337891173299632\\
474	0.000586862572901079\\
475	0.000778680506682941\\
476	0.000473244540492621\\
477	0.000634032948657522\\
478	0.000579322852580393\\
479	0.000377722826148201\\
480	0.000556854497285447\\
481	0.00104146778339268\\
482	0.000695268911252563\\
483	0.000485754197941148\\
484	0.000480744825553011\\
485	0.000466915927368166\\
486	0.000278759897025827\\
487	0.00021562827803745\\
488	0.000155537326131966\\
489	0.000583310147685467\\
490	0.000773663092257931\\
491	0.000566611919870179\\
492	0.000565212915732038\\
493	0.000561824325759264\\
494	0.000517497204050745\\
495	0.000392007168158299\\
496	0.000317801696951362\\
497	0.000369385213223469\\
498	0.000196904393002907\\
499	0.000894117159166521\\
500	8.57736896496535e-05\\
501	0.00025369644098731\\
502	0.000826746264843412\\
503	0.000542597863697207\\
504	0.000551069459274358\\
505	0.000572291076626316\\
506	0.000331456790657376\\
507	0.000216674193044686\\
508	0.000145129552214639\\
509	2.20765814166168e-06\\
510	1.84994234922908e-05\\
511	1.05320746954822e-07\\
512	1.21764815304034e-05\\
513	1.06283091899844e-05\\
514	2.54301479086815e-06\\
515	4.26561548567946e-05\\
516	1.87466071825392e-05\\
517	7.81090177618271e-06\\
518	0.000190386022642943\\
519	0.000123804517038932\\
520	0.000123829097563375\\
521	0.000114813549265202\\
522	7.25248809301798e-05\\
523	0.00032366876124878\\
524	0.000868725220323963\\
525	0.000269052821816188\\
526	0.000168300150209129\\
527	0.000178302293095357\\
528	0.000152812464407364\\
529	0.000128225081183955\\
530	0.000150570057759368\\
531	0.000422505103921372\\
532	0.000793085206131684\\
533	0.00058019282647861\\
534	0.000603733746473643\\
535	0.00042306047621397\\
536	0.000485422472629484\\
537	0.000802800030038001\\
538	0.00025428988261452\\
539	0.000198393993371717\\
540	0.000234215543574067\\
541	0.000852675175886393\\
542	0.000897809921859865\\
543	0.000978055349040792\\
544	0.000709603305577267\\
545	0.000527639384918283\\
546	0.000621283809795674\\
547	0.000716427352610584\\
548	0.000335690612321631\\
549	0.00023887096956402\\
550	0.000198053665356921\\
551	0.000431130177698011\\
552	0.000515902081481922\\
553	0.000333718539787681\\
554	0.000897036403083783\\
555	0.000280345490321062\\
556	0.000427042309805051\\
557	0.000354776814387017\\
558	0.000345022030850635\\
559	0.000340757168766302\\
560	0.000673890054973975\\
561	0.000547510899051937\\
562	0.000164221194052799\\
563	0.000930722447169514\\
564	0.0008512349819385\\
565	0.00106728353486036\\
566	0.00071787806277842\\
567	0.000775815313394007\\
568	0.000863394798698225\\
569	0.000515523156675\\
570	0.000562296765344483\\
571	0.000928088662774012\\
572	0.000969078852208438\\
573	0.00093208009190343\\
574	0.000632747483922107\\
575	0.000222797847920331\\
576	0.000844086930887845\\
577	0.000579847239664662\\
578	9.52111828979155e-05\\
579	0.00060295461171462\\
580	0.000734670037165659\\
581	0.000771938050671018\\
582	0.000788386128228819\\
583	0.000277209206232552\\
584	0.00100989921040436\\
585	0.00101055858645671\\
586	0.000999601257197842\\
587	0.000617285368672632\\
588	0.000625948093387843\\
589	0.000333069491463755\\
590	0.000480747138322226\\
591	0.000268406242509476\\
592	0.000931364172659619\\
593	0.000385237485865742\\
594	0.000809636137850576\\
595	0.00086869737765063\\
596	0.000191836916385684\\
597	0.000203898440103536\\
598	0.000335861796050138\\
599	0.000539517152888265\\
};
\addlegendentry{AA(2)-$X_i$}

\addplot [color=blue, mark=o, mark options={solid, blue}]
  table[row sep=crcr]{%
0	0.00154064771926348\\
1	0.000820403789956607\\
2	0.000350128390879639\\
3	2.58483014899325e-05\\
4	4.82276369884318e-05\\
5	2.7178191860958e-07\\
6	1.85696295003182e-08\\
7	2.0424597059118e-09\\
8	1.48941176522815e-11\\
9	1.00328672832162e-13\\
10	3.41758310778565e-15\\
};
\addlegendentry{AA(2)-$\rho_i$}

\addplot [color=myalign, mark=o, mark options={solid, myalign}]
  table[row sep=crcr]{%
0	0.00154064771926347\\
1	0.000820403789956604\\
2	5.79903299802225e-05\\
3	2.71555496109217e-06\\
4	3.14916027590177e-07\\
5	4.47521517783969e-08\\
6	8.33406990910289e-09\\
7	4.62824856775606e-10\\
8	2.56856004897138e-11\\
9	1.54438311776379e-12\\
10	2.90482562501961e-13\\
11	7.63679174319031e-14\\
12	8.57404176422874e-15\\
};
\addlegendentry{AA(2)-$\ALIGN{X_i}$}

\addplot [color=mygrassaa, mark=o, mark options={solid, mygrassaa}]
  table[row sep=crcr]{%
0	0.00154064771926348\\
1	0.000820403789956607\\
2	0.000372967811175703\\
3	1.74439733359143e-05\\
4	2.50795875439732e-06\\
5	9.25739601964327e-09\\
6	1.08856305390067e-09\\
7	2.20068829806374e-10\\
8	8.40397728832697e-13\\
9	4.58345231201336e-15\\
};
\addlegendentry{AA(2)-Gr}

\addplot [color=mygnt, mark=diamond, mark options={solid, mygnt}]
  table[row sep=crcr]{%
0	0.00154064771926348\\
1	0.000820403789956607\\
2	0.00013325669020215\\
3	1.20514815449128e-06\\
4	7.43938010239175e-11\\
5	2.85068410225813e-15\\
};
\addlegendentry{NTAcc}
\legend{}
\end{axis}
\end{tikzpicture}%

%% file: Figs/tikz/ks1d_beta7_iter.tex
%
%
\definecolor{mycolor1}{rgb}{1.00000,0.00000,1.00000}%
\definecolor{mycolor2}{rgb}{0.49400,0.18400,0.55600}%
\definecolor{mycolor3}{rgb}{0.00000,1.00000,1.00000}%
\definecolor{mycolor4}{rgb}{0.63500,0.07800,0.18400}%
\begin{tikzpicture}

\begin{axis}[%
width=3.65cm,
height=3.7cm,
scale only axis,
tick label style={font=\footnotesize},
xmin=0,
xmax=29,
xlabel style={font=\footnotesize\color{white!15!black},yshift=.2cm},
xlabel={Iteration},
ymode=log,
ymin=1.1e-15,
ymax=.95,
yminorticks=true,
ytickten = {-14, -11, -8, -5, -2},
ylabel style={font=\footnotesize\color{white!15!black},yshift=-.2cm},
ylabel={Residual norm},
axis background/.style={fill=white},
legend style={at={(axis cs:28.5,6e-5)}, anchor=north east, legend cell align=left, align=left, font = \tiny, draw=white!15!black}
]
\addplot [color=myscf, mark=x, mark options={solid, myscf}]
  table[row sep=crcr]{%
0	0.0026961335087109\\
1	0.0015223818035677\\
2	0.00159004912005195\\
3	0.00104126971025474\\
4	0.00114247883577157\\
5	0.000819108653403263\\
6	0.000895664849234087\\
7	0.000679344812193931\\
8	0.000734599134237397\\
9	0.000579277824642882\\
10	0.000618941918326986\\
11	0.000502342024984223\\
12	0.000530727696916927\\
13	0.000440493043367893\\
14	0.000460644212561792\\
15	0.000389255997955168\\
16	0.000403330965969593\\
17	0.000345899740093857\\
18	0.000355455131730911\\
19	0.000308639688514258\\
20	0.000314818036404478\\
21	0.000276246427411465\\
22	0.000279899037685933\\
23	0.000247836114365633\\
24	0.000249604923231074\\
25	0.000222750484288754\\
26	0.00022312389368802\\
27	0.000200484401668715\\
28	0.000199836298344128\\
29	0.000180640189703453\\
30	0.00017925783893756\\
31	0.00016289778910391\\
32	0.00016100218743278\\
33	0.00014699465990986\\
34	0.000144755654572714\\
35	0.000132711886623609\\
36	0.000130259582089649\\
37	0.000119864348092269\\
38	0.000117297824191826\\
39	0.000108293617185584\\
40	0.000105687664506297\\
41	9.78627332259299e-05\\
42	9.52731020838431e-05\\
43	8.84522834785573e-05\\
44	8.59198024528823e-05\\
45	7.99574150806978e-05\\
46	7.75112392106389e-05\\
47	7.22855184745025e-05\\
48	6.99457004084945e-05\\
49	6.53544024352697e-05\\
50	6.31339323470077e-05\\
51	5.90908339021304e-05\\
52	5.69972597445663e-05\\
53	5.34293522139357e-05\\
54	5.14660666364965e-05\\
55	4.8311292486791e-05\\
56	4.64785539293066e-05\\
57	4.36839705778088e-05\\
58	4.19797117769032e-05\\
59	3.94999946017637e-05\\
60	3.79204607683104e-05\\
61	3.57166769310842e-05\\
62	3.4256927352621e-05\\
63	3.22955271009818e-05\\
64	3.09498272102747e-05\\
65	2.92018107292353e-05\\
66	2.7963936420331e-05\\
67	2.64041630421108e-05\\
68	2.52676347490859e-05\\
69	2.3874248101125e-05\\
70	2.28325088113718e-05\\
71	2.1586456767407e-05\\
72	2.06330053961513e-05\\
73	1.95176378190185e-05\\
74	1.86461271748868e-05\\
75	1.76468577327716e-05\\
76	1.68511645277354e-05\\
77	1.59551854495135e-05\\
78	1.52294583576045e-05\\
79	1.44254990726725e-05\\
80	1.37641896743054e-05\\
81	1.30423119401269e-05\\
82	1.24401924293016e-05\\
83	1.17916159041757e-05\\
84	1.12437866714777e-05\\
85	1.06607399427501e-05\\
86	1.01626295167161e-05\\
87	9.63822249885808e-06\\
88	9.18558181684227e-06\\
89	8.71369612858084e-06\\
90	8.30258870673393e-06\\
91	7.87778322391019e-06\\
92	7.50457244615299e-06\\
93	7.12200170571428e-06\\
94	6.78333621552386e-06\\
95	6.43867972010943e-06\\
96	6.13147763409149e-06\\
97	5.82087844601673e-06\\
98	5.54231099081659e-06\\
99	5.26232225317161e-06\\
100	5.00979724270626e-06\\
101	4.75733549791425e-06\\
102	4.52848099104266e-06\\
103	4.30078532185065e-06\\
104	4.09343368476956e-06\\
105	3.88802988015406e-06\\
106	3.70020244011383e-06\\
107	3.51487150055227e-06\\
108	3.34476389257362e-06\\
109	3.17751428482084e-06\\
110	3.0234825710038e-06\\
111	2.87252574512482e-06\\
112	2.73307335031504e-06\\
113	2.59680209299663e-06\\
114	2.47056756395324e-06\\
115	2.34753682634216e-06\\
116	2.23328242177311e-06\\
117	2.12219231048118e-06\\
118	2.0187934013337e-06\\
119	1.91847406604047e-06\\
120	1.82490932333662e-06\\
121	1.73430751417233e-06\\
122	1.64964983905257e-06\\
123	1.56781692728553e-06\\
124	1.49122510945638e-06\\
125	1.41730640085056e-06\\
126	1.34801743926734e-06\\
127	1.28124263662941e-06\\
128	1.21856469793883e-06\\
129	1.15823936663198e-06\\
130	1.10154533537547e-06\\
131	1.04704327400796e-06\\
132	9.95764849689948e-07\\
133	9.46521253758452e-07\\
134	9.00143558678764e-07\\
135	8.55648890418702e-07\\
136	8.137055515257e-07\\
137	7.73500046915683e-07\\
138	7.3556870608459e-07\\
139	6.9923743769647e-07\\
140	6.64935665843635e-07\\
141	6.32104112330203e-07\\
142	6.01085691871963e-07\\
143	5.71415756014337e-07\\
144	5.4336729239599e-07\\
145	5.16553717171158e-07\\
146	4.91191568657787e-07\\
147	4.66958716644956e-07\\
148	4.44026196761908e-07\\
149	4.22125145533114e-07\\
150	4.01389985809558e-07\\
151	3.81595922840192e-07\\
152	3.6284797162039e-07\\
153	3.44957837766432e-07\\
154	3.28006968916866e-07\\
155	3.11837339412294e-07\\
156	2.96511558228875e-07\\
157	2.81896734045287e-07\\
158	2.68040458384541e-07\\
159	2.54830739531739e-07\\
160	2.42303242220883e-07\\
161	2.30363380187373e-07\\
162	2.19037379788484e-07\\
163	2.08245171793001e-07\\
164	1.98005552026725e-07\\
165	1.88250580966046e-07\\
166	1.78993233237988e-07\\
167	1.70175723825018e-07\\
168	1.6180649835256e-07\\
169	1.53836287286591e-07\\
170	1.4627004550766e-07\\
171	1.39065656349024e-07\\
172	1.32225407132829e-07\\
173	1.25713208277709e-07\\
174	1.19529335977552e-07\\
175	1.13642783103834e-07\\
176	1.08052336768453e-07\\
177	1.02731290886479e-07\\
178	9.76773532221875e-08\\
179	9.28674620876939e-08\\
180	8.82985675966801e-08\\
181	8.39507060075832e-08\\
182	7.9820322464672e-08\\
183	7.58900925978772e-08\\
184	7.21561493588968e-08\\
185	6.8603420830943e-08\\
186	6.52278798058524e-08\\
187	6.2016382706264e-08\\
188	5.89648521960279e-08\\
189	5.60618016965631e-08\\
190	5.33031897648044e-08\\
191	5.06789539226004e-08\\
192	4.81851495008702e-08\\
193	4.58129438273554e-08\\
194	4.35585335053965e-08\\
195	4.1414148354672e-08\\
196	3.93761558301469e-08\\
197	3.74377078713035e-08\\
198	3.55953612998078e-08\\
199	3.38430694618734e-08\\
200	3.21775898663766e-08\\
201	3.05935738565604e-08\\
202	2.90879846616505e-08\\
203	2.76560824177437e-08\\
204	2.62950354007101e-08\\
205	2.50006376736305e-08\\
206	2.37702585711158e-08\\
207	2.26001594512779e-08\\
208	2.14879049142237e-08\\
209	2.04301672109237e-08\\
210	1.94246968002782e-08\\
211	1.84685297724816e-08\\
212	1.75595927263499e-08\\
213	1.66952423169467e-08\\
214	1.58735717847774e-08\\
215	1.50922194100804e-08\\
216	1.43494375357939e-08\\
217	1.36431139310912e-08\\
218	1.29716470844307e-08\\
219	1.23331469027155e-08\\
220	1.17261482147867e-08\\
221	1.1148957940598e-08\\
222	1.06002381491942e-08\\
223	1.00784709750022e-08\\
224	9.58243514816225e-09\\
225	9.11076885701083e-09\\
226	8.66235891780892e-09\\
227	8.23598190943251e-09\\
228	7.83062510930962e-09\\
229	7.44518934478101e-09\\
230	7.07875219714993e-09\\
231	6.73032631692996e-09\\
232	6.39907166360709e-09\\
233	6.08410140241735e-09\\
234	5.78465203242916e-09\\
235	5.49992480379262e-09\\
236	5.22922745746402e-09\\
237	4.97184011676369e-09\\
238	4.7271327039704e-09\\
239	4.49445904526973e-09\\
240	4.27324698466231e-09\\
241	4.06291470117245e-09\\
242	3.86294254876558e-09\\
243	3.67280605422887e-09\\
244	3.49203439414379e-09\\
245	3.32015431900905e-09\\
246	3.15673938374067e-09\\
247	3.00136265477968e-09\\
248	2.85363842943658e-09\\
249	2.71318080653541e-09\\
250	2.5796404382746e-09\\
251	2.45266964234219e-09\\
252	2.33195161822286e-09\\
253	2.21717153835998e-09\\
254	2.10804434914036e-09\\
255	2.00428569489954e-09\\
256	1.90563624406437e-09\\
257	1.81184077742546e-09\\
258	1.72266370115512e-09\\
259	1.63787474712382e-09\\
260	1.55725947438902e-09\\
261	1.48061078152608e-09\\
262	1.4077362260355e-09\\
263	1.33844690057207e-09\\
264	1.27256933650564e-09\\
265	1.20993340835576e-09\\
266	1.15038114641322e-09\\
267	1.09375930504909e-09\\
268	1.03992554097387e-09\\
269	9.8873995264498e-10\\
270	9.40074568390683e-10\\
271	8.93804078195699e-10\\
272	8.49811292806076e-10\\
273	8.07983505492346e-10\\
274	7.68214420137388e-10\\
275	7.30403003359584e-10\\
276	6.94452979679874e-10\\
277	6.60271969948867e-10\\
278	6.2777395330671e-10\\
279	5.96875232541639e-10\\
280	5.6749716430475e-10\\
281	5.39565053645039e-10\\
282	5.13008061405472e-10\\
283	4.87758026997055e-10\\
284	4.63750694862944e-10\\
285	4.4092520098353e-10\\
286	4.19222532883071e-10\\
287	3.98588236288529e-10\\
288	3.78970084819855e-10\\
289	3.60317245096853e-10\\
290	3.42582616134337e-10\\
291	3.25721176270854e-10\\
292	3.09689398628752e-10\\
293	2.94446174970781e-10\\
294	2.79953696886402e-10\\
295	2.66174347947247e-10\\
296	2.53073561113898e-10\\
297	2.40617118949444e-10\\
298	2.28773931442544e-10\\
299	2.17513898090771e-10\\
300	2.06807755931447e-10\\
301	1.96629014501389e-10\\
302	1.86950953976781e-10\\
303	1.77749726938208e-10\\
304	1.6900099982491e-10\\
305	1.60682933618961e-10\\
306	1.52774344224671e-10\\
307	1.4525474949951e-10\\
308	1.38105311327642e-10\\
309	1.31308077521459e-10\\
310	1.24845204464341e-10\\
311	1.18700156414248e-10\\
312	1.12857371271075e-10\\
313	1.07303266912916e-10\\
314	1.02021773036565e-10\\
315	9.70003703197325e-11\\
316	9.22258347356052e-11\\
317	8.76864588282134e-11\\
318	8.33704023666443e-11\\
319	7.92666327538924e-11\\
320	7.5365212796588e-11\\
321	7.16557664076967e-11\\
322	6.81291981417371e-11\\
323	6.47761910671415e-11\\
324	6.15880789489657e-11\\
325	5.85569725837261e-11\\
326	5.56751549314495e-11\\
327	5.29350807564706e-11\\
328	5.03300522648341e-11\\
329	4.78530197605007e-11\\
330	4.54978390455432e-11\\
331	4.3258703759667e-11\\
332	4.11294271098941e-11\\
333	3.91046908349746e-11\\
334	3.71798030990307e-11\\
335	3.53492335472884e-11\\
336	3.36093263643948e-11\\
337	3.19554041248344e-11\\
338	3.03821231177927e-11\\
339	2.88866493838839e-11\\
340	2.74645136928311e-11\\
341	2.61124012391644e-11\\
342	2.48266911018918e-11\\
343	2.36044897583257e-11\\
344	2.2442543184723e-11\\
345	2.1337609661151e-11\\
346	2.02870454863407e-11\\
347	1.92882316722317e-11\\
348	1.83388145108024e-11\\
349	1.74360851261057e-11\\
350	1.65780018511394e-11\\
351	1.57628094962022e-11\\
352	1.49873905450873e-11\\
353	1.42495459731061e-11\\
354	1.35478114908442e-11\\
355	1.28810688490794e-11\\
356	1.2246648886519e-11\\
357	1.16438171295646e-11\\
358	1.10709992944249e-11\\
359	1.05257123450585e-11\\
360	1.00076992178423e-11\\
361	9.51513398378699e-12\\
362	9.04677699379503e-12\\
363	8.60125507756907e-12\\
364	8.17775918849114e-12\\
365	7.77526090530918e-12\\
366	7.39314168721221e-12\\
367	7.02911881167728e-12\\
368	6.68271526418164e-12\\
369	6.35415630637007e-12\\
370	6.04137502650716e-12\\
371	5.74406390824431e-12\\
372	5.46157982738188e-12\\
373	5.19296977509108e-12\\
374	4.93752804198916e-12\\
375	4.69496912401629e-12\\
376	4.46386159943933e-12\\
377	4.24376096769034e-12\\
378	4.03506538768297e-12\\
379	3.83673994098887e-12\\
380	3.64755283394633e-12\\
381	3.46825650452735e-12\\
382	3.29731596410177e-12\\
383	3.13449847012503e-12\\
384	2.97965939757929e-12\\
385	2.83252758828481e-12\\
386	2.69361213354484e-12\\
387	2.56138269419482e-12\\
388	2.43498967052567e-12\\
389	2.31490584899215e-12\\
390	2.20126312049883e-12\\
391	2.09330389546358e-12\\
392	1.99064771407746e-12\\
393	1.89258966099806e-12\\
394	1.79954234081156e-12\\
395	1.71085917045861e-12\\
396	1.62685579137788e-12\\
397	1.54712520971815e-12\\
398	1.47047342745574e-12\\
399	1.39763408260627e-12\\
400	1.32877936985698e-12\\
401	1.26327983950172e-12\\
402	1.20089085740321e-12\\
403	1.14145979562024e-12\\
404	1.08552179906072e-12\\
405	1.03201754896516e-12\\
406	9.8091344957235e-13\\
407	9.32871394602183e-13\\
408	8.87482886110665e-13\\
409	8.4336094856199e-13\\
410	8.01589292255788e-13\\
411	7.61892603087515e-13\\
412	7.24061930636797e-13\\
413	6.88671816498267e-13\\
414	6.55155237681586e-13\\
415	6.23076253711387e-13\\
416	5.92653210140565e-13\\
417	5.63097560204676e-13\\
418	5.35453236158888e-13\\
419	5.09141680972804e-13\\
420	4.84077569020541e-13\\
421	4.60057741640973e-13\\
422	4.36740229215811e-13\\
423	4.16478079190939e-13\\
424	3.96482605859417e-13\\
425	3.7758058122847e-13\\
426	3.59286309194002e-13\\
427	3.41740246977498e-13\\
428	3.25139301057858e-13\\
429	3.09739439022379e-13\\
430	2.94704390729093e-13\\
431	2.80551757517911e-13\\
432	2.66554062581641e-13\\
433	2.53984411741746e-13\\
434	2.415273551764e-13\\
435	2.29846906168149e-13\\
436	2.18368389771599e-13\\
437	2.06823430920042e-13\\
438	1.96253429004385e-13\\
439	1.86611910601394e-13\\
440	1.77317932468375e-13\\
441	1.68520137098097e-13\\
442	1.60517940056714e-13\\
443	1.52294306900257e-13\\
444	1.44836162944924e-13\\
445	1.3813672734427e-13\\
446	1.31442879523338e-13\\
447	1.25559228971496e-13\\
448	1.19196581899261e-13\\
449	1.13853534255206e-13\\
450	1.08356438507474e-13\\
451	1.02680206802296e-13\\
452	9.74283361554505e-14\\
453	9.33213384686369e-14\\
454	8.88712805787637e-14\\
455	8.46186918071708e-14\\
456	8.08148012811573e-14\\
457	7.67806937371338e-14\\
458	7.23477684224399e-14\\
459	6.79774428425103e-14\\
460	6.42293154110749e-14\\
461	6.06538411349666e-14\\
462	5.71045435050632e-14\\
463	5.41964120927721e-14\\
464	5.16355959811505e-14\\
465	4.8730187268985e-14\\
466	4.6188538452399e-14\\
467	4.44487042861363e-14\\
468	4.24599358037343e-14\\
469	4.05391958050515e-14\\
470	3.84863901286484e-14\\
471	3.65173519003787e-14\\
472	3.46063688302762e-14\\
473	3.28950491163338e-14\\
474	3.16799652275968e-14\\
475	2.97423023230993e-14\\
476	2.85944526562348e-14\\
477	2.73925009086486e-14\\
478	2.55511022587177e-14\\
479	2.40543078539486e-14\\
480	2.29314834489677e-14\\
481	2.17873785029077e-14\\
482	2.0585671444609e-14\\
483	1.93141089160189e-14\\
484	1.83251012084052e-14\\
485	1.74104981875068e-14\\
486	1.59232186568161e-14\\
487	1.51720693228716e-14\\
488	1.40643059558765e-14\\
489	1.33977733708176e-14\\
490	1.29041027266037e-14\\
491	1.2092160729626e-14\\
492	1.16375226847025e-14\\
493	1.10233520002795e-14\\
494	1.03987997347425e-14\\
495	1.0127127675877e-14\\
496	1.03255444861202e-14\\
497	1.00622108079558e-14\\
498	9.63903611909184e-15\\
};
\addlegendentry{SCF}

\addplot [color=black, mark=+, mark options={solid, black}]
  table[row sep=crcr]{%
0	0.0026961335087109\\
1	0.0015223818035677\\
2	0.00141640095597572\\
3	0.00109286504516236\\
4	0.0014786186760199\\
5	0.000926093253576387\\
6	0.00101968993423823\\
7	0.000268509704522905\\
8	0.000797484260408955\\
9	0.000787795985158083\\
10	0.000747509400553606\\
11	0.000258769946563444\\
12	0.000411420034423197\\
13	0.000454619012567523\\
14	0.000886024839087124\\
15	0.000939433361422338\\
16	0.000983120199562017\\
17	0.0016653453540267\\
18	0.00116777082144813\\
19	0.000779616634596101\\
20	0.000826400320336522\\
21	0.0016078903058888\\
22	0.00108368991105286\\
23	0.000917158414963543\\
24	0.00082032913055474\\
25	0.0015110984820302\\
26	0.000634282166879696\\
27	0.000714886964957095\\
28	0.000438065501526716\\
29	0.000330878396875546\\
30	0.000330232235601248\\
31	0.00016801766516724\\
32	0.000811405284025549\\
33	0.00102217650622745\\
34	0.00103927748836023\\
35	0.000737145983911915\\
36	0.0017505024331626\\
37	0.0008207926187226\\
38	0.000413283907109703\\
39	0.00095907936735711\\
40	0.00164592502628273\\
41	0.00085982283074693\\
42	0.00020736481959552\\
43	0.000150322593717516\\
44	0.000305473192696252\\
45	0.000336645695233571\\
46	0.000852546158506241\\
47	0.00143695048120691\\
48	0.00105329680236032\\
49	0.000690314910635298\\
50	0.000645312688997971\\
51	0.000633266335120186\\
52	0.00174381858506205\\
53	0.00112790701514546\\
54	0.000666627126668255\\
55	0.00105078068118721\\
56	0.00206073337053024\\
57	0.00101342689017716\\
58	0.0011224031176533\\
59	0.00125494449665226\\
60	0.00103521953753759\\
61	0.000932148398372918\\
62	0.000321549247636451\\
63	0.000334477639776759\\
64	0.000213182351467992\\
65	0.000197313609070072\\
66	0.000308457899820062\\
67	0.000327276656122081\\
68	0.000342389814791049\\
69	0.000858961830959209\\
70	0.000861962319635153\\
71	0.00102371249035874\\
72	0.000853199380714707\\
73	0.000730801430649895\\
74	0.000632898860230608\\
75	0.000530461043771479\\
76	0.00101858710795604\\
77	0.000848287469115876\\
78	0.000660906377314258\\
79	0.000485930389356151\\
80	0.000143714101199871\\
81	0.000635943184802561\\
82	0.000701507245364646\\
83	0.00126072109035652\\
84	0.00105202268200203\\
85	0.000870297793848259\\
86	0.00121453472913561\\
87	0.000785132036103106\\
88	0.000855309822270247\\
89	0.000902467044155014\\
90	0.000893982501639475\\
91	0.000922778207327699\\
92	0.000393803028041685\\
93	0.00108792234068823\\
94	0.000555841554379687\\
95	0.00069026873267396\\
96	0.000573773424295378\\
97	0.000996116803960838\\
98	0.000477666107148906\\
99	0.00054652020046857\\
100	0.00126926623450658\\
101	0.000961319228690334\\
102	0.00136449524746234\\
103	0.0012109367278279\\
104	0.00082394368654848\\
105	0.000900335903795621\\
106	0.000411106266641657\\
107	0.0005143260569581\\
108	0.000342033923544215\\
109	0.000921756988794551\\
110	0.000708935547677974\\
111	0.00076734767627952\\
112	0.000538970286327822\\
113	0.000781537397632014\\
114	0.00102492244131095\\
115	0.000785502880243794\\
116	0.000552778505719158\\
117	0.00147823618310891\\
118	0.00121005603585691\\
119	0.00109170135487015\\
120	0.000932386345907578\\
121	0.00164321844775142\\
122	0.000743155155824293\\
123	0.000944901386055127\\
124	0.00114578890107564\\
125	0.000777209524204557\\
126	0.00059983624705194\\
127	0.00113879298466205\\
128	0.00165958486373608\\
129	0.00126804855078984\\
130	0.000510737352198505\\
131	0.000476973348279597\\
132	0.00083476378319542\\
133	0.00199721235623168\\
134	0.000929674578465367\\
135	0.000922652162580618\\
136	0.000666463952026048\\
137	0.000904570507023822\\
138	0.000719111507349214\\
139	0.000536110490047836\\
140	0.000937232678217564\\
141	0.000404731254449181\\
142	0.000177152117988627\\
143	0.00143337727088229\\
144	0.00115735494967539\\
145	0.000673387287272812\\
146	0.000888938540622511\\
147	0.00205141743639624\\
148	0.00116646231654199\\
149	0.000863623624228723\\
150	0.00030065189760627\\
151	0.000932936891343347\\
152	0.000990097393259898\\
153	0.000275296470351482\\
154	0.000838733299617324\\
155	0.00133270022091616\\
156	0.000450906699585127\\
157	0.00056089924201431\\
158	0.000815799004472922\\
159	0.000183451674602373\\
160	0.000755444782654711\\
161	0.00107406231057802\\
162	0.0012753762504751\\
163	0.000819260247847475\\
164	0.000911245936202765\\
165	0.000791626338223048\\
166	0.00151740561950154\\
167	0.000585467972739311\\
168	0.000144117555019336\\
169	0.000590794891872805\\
170	0.00083336818429529\\
171	0.00161095603555903\\
172	0.00116197074397267\\
173	0.00108932795739567\\
174	0.000674032477398983\\
175	0.000393005623470265\\
176	0.00092481239599272\\
177	0.000764076739904351\\
178	0.000469593988417191\\
179	0.00172124542081858\\
180	0.00100834157392934\\
181	0.000653951336232263\\
182	0.000736649185800347\\
183	0.000159250628697211\\
184	0.000104497007064945\\
185	0.000561638220148789\\
186	0.000677964222408894\\
187	0.00129912695917186\\
188	0.00111394004335482\\
189	0.0009969163871723\\
190	0.000926759517858914\\
191	0.000925150571927421\\
192	0.00111001651746259\\
193	0.0011242861249878\\
194	0.00094471268502548\\
195	0.000590162693299139\\
196	0.000735368023050234\\
197	0.000698927860256025\\
198	0.000637236944655795\\
199	0.000919094871077365\\
200	0.000816087561240682\\
201	0.000876785903781787\\
202	0.000630826183040881\\
203	0.000168677618997332\\
204	0.000207363389214271\\
205	0.000262948064539479\\
206	0.000108702053946212\\
207	0.000969019353295731\\
208	0.000899752579409608\\
209	0.000540227794598899\\
210	0.000351558940935539\\
211	0.00110392283172246\\
212	0.000962238412000445\\
213	0.00100748818351481\\
214	0.000992977975686098\\
215	0.0015917218030931\\
216	0.00106424971446187\\
217	0.00108161786397671\\
218	0.000744670744323163\\
219	0.000831540180121157\\
220	0.000684052087999205\\
221	0.00144630584688871\\
222	0.00103100522534865\\
223	0.000737822674653158\\
224	0.000419295000982358\\
225	0.000499908539612934\\
226	0.00049371789361523\\
227	0.000328338991257766\\
228	0.000475678532775363\\
229	0.000393184160026365\\
230	0.000848784134624481\\
231	0.00086912833977859\\
232	0.000834466317169479\\
233	0.00132495593262132\\
234	0.000384205428028824\\
235	0.00148277510027478\\
236	0.00107166316824303\\
237	0.000906408546864628\\
238	0.00129451600151036\\
239	0.000498407735126516\\
240	0.000461266456426547\\
241	0.000816860910006397\\
242	0.000850357404253056\\
243	0.000923244312474458\\
244	0.000934054773962291\\
245	0.000359335837671135\\
246	0.00158207591142828\\
247	0.000556978516172323\\
248	0.0010001268643078\\
249	0.000580696730125695\\
250	0.000139447520078925\\
251	0.000776534346987114\\
252	0.000857313799729345\\
253	0.000364636328817506\\
254	0.000830470678926745\\
255	0.000758684863897089\\
256	0.00098006502886672\\
257	0.000519414100516049\\
258	6.84147318998084e-05\\
259	4.31135237235942e-05\\
260	4.59948348381127e-05\\
261	4.15268494825656e-05\\
262	0.000153064730675411\\
263	0.000131629094615969\\
264	0.000152039938396248\\
265	0.000153313039432799\\
266	0.00095021203487891\\
267	0.000734528065435757\\
268	0.000601292339031072\\
269	0.000178405558065295\\
270	0.000150659920351612\\
271	0.000730898210750259\\
272	0.000792111911323419\\
273	0.000174713753943674\\
274	0.000766985735631496\\
275	0.00036326880833409\\
276	0.000457420553893832\\
277	0.00134693796983072\\
278	0.00094441785739483\\
279	0.000297857076015855\\
280	0.000292512269808777\\
281	0.000212218783002623\\
282	0.000725158037756242\\
283	0.000931015412390422\\
284	0.00124319671714132\\
285	0.000546707006934195\\
286	0.00096249202671636\\
287	0.00113345909711275\\
288	0.000747914429483043\\
289	0.0013983856001524\\
290	0.0012522356742656\\
291	0.00118947374512143\\
292	0.000853042233972703\\
293	0.00106333361080208\\
294	0.00133615627585225\\
295	0.000836394607368139\\
296	0.00170442670580129\\
297	0.00173531891978787\\
298	0.000817981902220732\\
299	0.00115119951133701\\
300	0.00114284207748029\\
301	0.00118945692031874\\
302	0.000900730526846862\\
303	0.000411102445231359\\
304	0.00042107562744582\\
305	0.000638559951190394\\
306	0.000265728899862384\\
307	0.00065177213400478\\
308	0.000494065294554424\\
309	0.000616612779457654\\
310	0.000414751427518158\\
311	0.00105067772457394\\
312	0.00121529065617555\\
313	0.00155827989141096\\
314	0.000780698099534939\\
315	0.000460155919202928\\
316	0.000422427869205182\\
317	0.000637434832814294\\
318	0.000189013749479665\\
319	0.000880522105789601\\
320	0.00150456587303048\\
321	0.000885676446934998\\
322	0.00105865387576764\\
323	0.0010101399317887\\
324	0.000753938663638302\\
325	0.000871792741745472\\
326	0.000487701565092592\\
327	0.00165829847705724\\
328	0.00100774264816277\\
329	0.000899761580389788\\
330	0.000956219865381418\\
331	0.00100813204599234\\
332	0.00152190100342114\\
333	0.00156860903730178\\
334	0.000993184874644815\\
335	0.000470212051546901\\
336	0.0012157159668168\\
337	0.00108041076890886\\
338	0.000881555457144456\\
339	0.000579579489592173\\
340	0.000292396491121177\\
341	0.000344188678347952\\
342	0.00053269645995325\\
343	0.000266060072577062\\
344	0.000905589636316842\\
345	0.000936647659386625\\
346	0.00119893684223908\\
347	0.000467144253775524\\
348	0.0016768297886231\\
349	0.00102498329577838\\
350	0.000965349550842499\\
351	0.000364450921913741\\
352	0.00140701881651831\\
353	0.0012503209637784\\
354	0.00088016530973881\\
355	0.00099574317550287\\
356	0.000825287854213789\\
357	0.000225211268839175\\
358	0.000923337788670036\\
359	0.000819244325088409\\
360	0.00068571949693754\\
361	0.000669424279495795\\
362	0.0008187622418697\\
363	0.000644243718290735\\
364	0.000654379794627145\\
365	0.000654752976410126\\
366	0.000118761617599983\\
367	0.000300503616586091\\
368	0.000310296780645602\\
369	0.000328740467969958\\
370	0.000308469723093318\\
371	0.000289106445265846\\
372	0.000979790383958037\\
373	0.00152994626013613\\
374	0.00158474797703066\\
375	0.00105487989358959\\
376	0.000299934708628951\\
377	0.00176241575898689\\
378	0.000608458394011066\\
379	0.00140625214809502\\
380	0.000752585373421155\\
381	0.00113012454317799\\
382	0.000392870579457206\\
383	0.000592209834223381\\
384	0.000396311546735991\\
385	0.000427951378466306\\
386	0.000454501604319681\\
387	0.000583274549272483\\
388	0.000530638162930511\\
389	0.00060436501546169\\
390	0.000816656247724636\\
391	0.000787927941390293\\
392	0.000589036068095808\\
393	0.000618028860799092\\
394	0.000628590511391833\\
395	0.000854910202671143\\
396	0.00064540030545142\\
397	0.000658924608108869\\
398	0.000654374083599302\\
399	0.000556772048046631\\
400	0.000649216600756314\\
401	0.000640736379856203\\
402	0.000899296924414132\\
403	0.000532158827212657\\
404	0.00073869820494178\\
405	0.000724627063500621\\
406	0.000725178090785599\\
407	0.00104281320673061\\
408	0.000943089533202591\\
409	0.000907421550904434\\
410	0.000531098744161906\\
411	0.000389901693586889\\
412	0.000363034317762709\\
413	0.000529932436037829\\
414	0.000498687593146165\\
415	0.000252040558614122\\
416	0.000295994760782888\\
417	0.000492236958085146\\
418	0.00113091412550957\\
419	0.000405405029910756\\
420	0.00127520253134758\\
421	0.000119058219467625\\
422	0.000300371105365287\\
423	0.000216791844593132\\
424	0.000211372422898076\\
425	0.000185922943017812\\
426	0.000182609138600587\\
427	0.000186734314773377\\
428	0.000163522370538367\\
429	0.000164405604232996\\
430	0.00013133860772847\\
431	0.00018394935028778\\
432	2.62834408010772e-05\\
433	0.00126575983373592\\
434	0.00115558059276531\\
435	0.000256084680911876\\
436	0.000995247036823167\\
437	0.0014828420957827\\
438	0.000409060004225201\\
439	0.00130842099207627\\
440	0.000891753389909663\\
441	0.000921644413389208\\
442	0.00033310089682816\\
443	0.000474183074017404\\
444	0.000631518861176926\\
445	0.00155874051781515\\
446	0.00085192515784851\\
447	0.000858372927449316\\
448	0.000280745692570019\\
449	0.000807403967957993\\
450	0.00022899459463434\\
451	0.000418705973345346\\
452	0.000322905871586895\\
453	0.000372376061490282\\
454	0.000531960152944542\\
455	0.00031562348478833\\
456	0.000804047795797802\\
457	0.000575102747932696\\
458	0.0001858800430097\\
459	0.000134654347094432\\
460	0.000537696404482919\\
461	0.00107291034829357\\
462	0.00105124844822198\\
463	0.000811782661187485\\
464	0.000928788643453112\\
465	0.000643440907304942\\
466	0.00100405536121249\\
467	0.000991619744320763\\
468	0.000714331532374051\\
469	0.000759881864965426\\
470	0.00175882179278959\\
471	0.000881688314429236\\
472	0.00101929409576683\\
473	0.000653287810401043\\
474	0.000444215824542638\\
475	0.000535058976880896\\
476	0.000583406506492275\\
477	0.000134676170848337\\
478	9.73468439068442e-05\\
479	0.000116803636323149\\
480	0.000130252963547324\\
481	8.58722224834143e-05\\
482	0.000581149386122982\\
483	0.00112262381360848\\
484	0.000281288814632974\\
485	0.000129884324915285\\
486	0.000666888672397715\\
487	0.00111187331983422\\
488	0.000502075510218345\\
489	0.00151459518511132\\
490	0.000502241976385768\\
491	0.000222847566595013\\
492	0.000258169516889453\\
493	0.00135145742964944\\
494	0.00102229190835639\\
495	0.000719599759671637\\
496	0.000819667090128689\\
497	0.000963877218778139\\
498	0.000402857428689396\\
499	0.000523377066640061\\
500	0.000915635753963754\\
501	0.000894673591707532\\
502	0.00107248706574882\\
503	0.000800339659275199\\
504	0.0012727233851944\\
505	0.00106793998482883\\
506	0.000700644992166021\\
507	0.000447826896660282\\
508	0.000375057959930893\\
509	0.0012668087964108\\
510	0.000781165664511098\\
511	0.000876455486296476\\
512	0.00079549534730879\\
513	0.00204915000557276\\
514	0.00108106731104648\\
515	0.00116198425170714\\
516	0.00117352942155801\\
517	0.00107786917869654\\
518	0.00048193629717803\\
519	0.00128477518836873\\
520	0.00158706230269123\\
521	0.00180870780634748\\
522	0.00185109077155849\\
523	0.00183312340614621\\
524	0.00123031353439886\\
525	0.000842330602326955\\
526	0.000883679903212091\\
527	0.00126136179746285\\
528	0.000662136118156379\\
529	0.00132476525315203\\
530	0.00122262105421565\\
531	0.00105253951620387\\
532	0.000281396022005387\\
533	0.000533216927680284\\
534	0.000686001489314654\\
535	0.000844129433378812\\
536	0.00026806978179987\\
537	0.000107770167193619\\
538	0.000100519300805749\\
539	8.38564889122113e-05\\
540	1.81894881869967e-05\\
541	0.00119621892892186\\
542	0.00117484496956299\\
543	0.00195574991706458\\
544	0.00123556006873418\\
545	0.00145444710203628\\
546	0.00123704395712449\\
547	0.000561477474021361\\
548	0.00119993863450395\\
549	0.00123946777051332\\
550	0.000848479532337094\\
551	0.00109952762735347\\
552	0.000753993845091813\\
553	0.000508912987527326\\
554	0.000716252488964944\\
555	0.00047574795054975\\
556	0.000816291658372476\\
557	0.000179476158359909\\
558	0.000926672369498456\\
559	0.000693044776011723\\
560	0.000901296097315365\\
561	0.000586300760617776\\
562	0.001150547106493\\
563	0.000567889238869289\\
564	0.000489972138526308\\
565	0.000960097218750054\\
566	0.00115429304737485\\
567	0.000678668484808287\\
568	0.00129228494344474\\
569	0.00102556856874394\\
570	0.000938296787384207\\
571	0.00115994808665221\\
572	0.000400605709863208\\
573	0.000660601948023536\\
574	0.000402106962521702\\
575	0.00136406423531343\\
576	0.00034517560481017\\
577	0.000659206245106451\\
578	0.000779539930003681\\
579	0.000505123552148095\\
580	0.000510719674960179\\
581	0.000367594016438464\\
582	0.000342755378165798\\
583	0.000286904030488869\\
584	0.00013339544207005\\
585	0.000145641731483886\\
586	0.000144841678972605\\
587	2.74983581539582e-05\\
588	1.24856657845784e-05\\
589	9.45418879209979e-06\\
590	1.78624968998848e-05\\
591	1.880224806219e-06\\
592	1.60425591411311e-05\\
593	1.97091428001944e-05\\
594	0.000124491623095252\\
595	0.000400134259569838\\
596	0.00121010400969579\\
597	0.00141352613094469\\
598	0.00169961810371167\\
599	0.000742588170381181\\
};
\addlegendentry{AA(2)-$X_i$}

\addplot [color=blue, mark=o, mark options={solid, blue}]
  table[row sep=crcr]{%
0	0.0026961335087109\\
1	0.0015223818035677\\
2	0.000934243008783277\\
3	0.000132210034879149\\
4	0.000296194347890923\\
5	2.66596141952223e-05\\
6	1.51531779981564e-06\\
7	3.10441159519056e-07\\
8	1.448610294864e-09\\
9	2.72786906848794e-11\\
10	1.7697547009389e-13\\
11	4.48404607021951e-15\\
};
\addlegendentry{AA(2)-$\rho_i$}

\addplot [color=myalign, mark=o, mark options={solid, myalign}]
  table[row sep=crcr]{%
0	0.0026961335087109\\
1	0.0015223818035677\\
2	0.000291955694115765\\
3	0.000113371893823403\\
4	1.59176991888612e-06\\
5	3.99431854780604e-07\\
6	3.13611807105384e-08\\
7	3.22148130886114e-09\\
8	1.73727687475775e-10\\
9	2.75200619215978e-11\\
10	1.30042829009528e-12\\
11	2.82116193564394e-13\\
12	3.94483674570764e-14\\
13	3.54934341967726e-15\\
};
\addlegendentry{AA(2)-$\ALIGN{X_i}$}

\addplot [color=mygrassaa, mark=o, mark options={solid, mygrassaa}]
  table[row sep=crcr]{%
0	0.0026961335087109\\
1	0.0015223818035677\\
2	0.0003232984578474\\
3	3.10667982697972e-05\\
4	6.36774912881799e-06\\
5	3.81898393825317e-07\\
6	6.0545021668805e-08\\
7	2.74331375910896e-10\\
8	1.95360982773954e-11\\
9	1.86186004093585e-14\\
10	2.64919209675976e-15\\
};
\addlegendentry{AA(2)-Gr}

\addplot [color=mygnt, mark=diamond, mark options={solid, mygnt}]
  table[row sep=crcr]{%
0	0.0026961335087109\\
1	0.0015223818035677\\
2	0.00109787232248698\\
3	0.000657211480666542\\
4	0.000137814268458899\\
5	3.19508723409014e-06\\
6	1.47555938201806e-09\\
7	2.87392459622677e-15\\
};
\addlegendentry{AlignNT}
\end{axis}
\end{tikzpicture}%

%% file: Figs/tikz/GMRES_residual_NT.tex
%
%
\definecolor{mycolor1}{rgb}{0.00000,0.44700,0.74100}%
\definecolor{mycolor2}{rgb}{0.85000,0.32500,0.09800}%
\definecolor{mycolor3}{rgb}{0.92900,0.69400,0.12500}%
\definecolor{mycolor4}{rgb}{0.49400,0.18400,0.55600}%
\begin{tikzpicture}

\begin{axis}[%
width=3.65cm,
height=3.7cm,
scale only axis,
tick label style={font=\footnotesize},
xmin=0,
xmax=410,
xtick = {0, 100, 200, 300, 400},
xticklabels = {0, 100, 200, 300, $ $},
xlabel style={font=\footnotesize\color{white!15!black},yshift=.2cm},
xlabel={Iteration},
ymode=log,
ymin=1.1e-15,
ymax=1.1,
yminorticks=true,
ytickten = {-14, -11, -8, -5, -2},
ylabel style={font=\footnotesize\color{white!15!black},yshift=-.2cm},
ylabel={Residual norm},
axis background/.style={fill=white},
legend style={at={(axis cs:505,9e-1)}, anchor=north east, legend cell align=left, align=left, font = \scriptsize, draw=white!15!black}
]
\addplot [color=mycolor1, only marks, mark=x, mark options={solid, mycolor1}]
  table[row sep=crcr]{%
0	1\\
1	0.204873076072875\\
2	0.0896638419929104\\
3	0.0463024790184499\\
4	0.0264444369981361\\
5	0.0175688806122771\\
6	0.0123059108473381\\
7	0.00880392571253182\\
8	0.0064191542106809\\
9	0.00496301915342828\\
10	0.00396141286514167\\
11	0.00338938840271267\\
12	0.00326668579129264\\
13	0.00315795230714888\\
14	0.00272425523285366\\
15	0.0023057934263151\\
16	0.00197552814860765\\
17	0.00169737040660245\\
18	0.00146477352776768\\
19	0.0012935811270209\\
20	0.00122341117701275\\
21	0.00122115761383236\\
22	0.00111246518490445\\
23	0.000979235725970141\\
24	0.000880080589838911\\
25	0.000794070011635746\\
26	0.000713559520978982\\
27	0.000636896193631892\\
28	0.000578048818906549\\
29	0.00057118743754868\\
30	0.00055096260926942\\
31	0.000498757687988901\\
32	0.000459471243028583\\
33	0.000421172962701658\\
34	0.000388645307922558\\
35	0.000360379426128009\\
36	0.000339403034667663\\
37	0.000320773245009384\\
38	0.000298347881528757\\
39	0.000276941674310424\\
40	0.000261010920133988\\
41	0.000243460279828632\\
42	0.000228316928285722\\
43	0.000212943874595316\\
44	0.000199344343085087\\
45	0.000186295182481351\\
46	0.000174122842687743\\
47	0.000162182197119404\\
48	0.000149209149270473\\
49	0.000137239202025032\\
50	0.000125860564031184\\
51	0.00011631657344569\\
52	0.000109047448405801\\
53	0.000102865523025157\\
54	9.69669514509603e-05\\
55	9.15405105978209e-05\\
56	8.72664629956135e-05\\
57	8.32470738660649e-05\\
58	7.89167877659295e-05\\
59	7.4571694315209e-05\\
60	6.97489520248647e-05\\
61	6.57338345414187e-05\\
62	6.17582066257447e-05\\
63	5.79868598610548e-05\\
64	5.45581385124384e-05\\
65	5.14425000932816e-05\\
66	4.83544319229528e-05\\
67	4.56275313725573e-05\\
68	4.34904564140037e-05\\
69	4.1399523910366e-05\\
70	3.90969299155383e-05\\
71	3.62039089956327e-05\\
72	3.30913774405313e-05\\
73	3.01158326623909e-05\\
74	2.78556206368964e-05\\
75	2.58548143368345e-05\\
76	2.3823191516756e-05\\
77	2.20818027556776e-05\\
78	2.05325451481531e-05\\
79	1.94379542930343e-05\\
80	1.83260082625741e-05\\
81	1.72054902953204e-05\\
82	1.6302258740955e-05\\
83	1.53723158940558e-05\\
84	1.46669458628024e-05\\
85	1.40557833837473e-05\\
86	1.35078931176676e-05\\
87	1.29403122336386e-05\\
88	1.22867751547664e-05\\
89	1.17429433801881e-05\\
90	1.11771449515877e-05\\
91	1.06109932081334e-05\\
92	1.01203691407635e-05\\
93	9.69697324458682e-06\\
94	9.25440040943257e-06\\
95	8.8376421612065e-06\\
96	8.51243879329683e-06\\
97	8.23279134705954e-06\\
98	7.98605838478697e-06\\
99	7.76272311090873e-06\\
100	7.57361520274234e-06\\
101	7.38444585303021e-06\\
102	7.20944948820887e-06\\
103	7.04226340281758e-06\\
104	6.84340991329695e-06\\
105	6.64790006798046e-06\\
106	6.44571639098787e-06\\
107	6.28091206308648e-06\\
108	6.11582261155173e-06\\
109	5.93992076857093e-06\\
110	5.77838105441249e-06\\
111	5.63920846927357e-06\\
112	5.48505820750688e-06\\
113	5.29366145291491e-06\\
114	5.11425938994836e-06\\
115	4.95011817799596e-06\\
116	4.75337097138781e-06\\
117	4.55532567471207e-06\\
118	4.38171559208682e-06\\
119	4.19372711912735e-06\\
120	4.03046875292634e-06\\
121	3.89293653297232e-06\\
122	3.74550722158903e-06\\
123	3.5923422437607e-06\\
124	3.43799815928443e-06\\
125	3.30563597943814e-06\\
126	3.16646153898207e-06\\
127	3.03700117496231e-06\\
128	2.90701825202156e-06\\
129	2.78780801619026e-06\\
130	2.6798557349191e-06\\
131	2.56654893807711e-06\\
132	2.46919757245083e-06\\
133	2.3745846984223e-06\\
134	2.26265070732837e-06\\
135	2.14983077024128e-06\\
136	2.06793907490551e-06\\
137	1.98591382494867e-06\\
138	1.91188238354938e-06\\
139	1.85050744811714e-06\\
140	1.77986924958279e-06\\
141	1.70196303558922e-06\\
142	1.62874702574299e-06\\
143	1.56364397751387e-06\\
144	1.50071954595125e-06\\
145	1.44647837848897e-06\\
146	1.38610052202972e-06\\
147	1.32878091781414e-06\\
148	1.28145464415402e-06\\
149	1.23007577954926e-06\\
150	1.17191723093195e-06\\
151	1.11804575515052e-06\\
152	1.05943554479439e-06\\
153	9.97285039927271e-07\\
154	9.41691993340488e-07\\
155	8.89037982656683e-07\\
156	8.30379913861785e-07\\
157	7.82710229526357e-07\\
158	7.47275968947695e-07\\
159	7.14077484326499e-07\\
160	6.88229433356305e-07\\
161	6.6460839335771e-07\\
162	6.47095146295994e-07\\
163	6.27356460003953e-07\\
164	6.0852966989102e-07\\
165	5.83306337793504e-07\\
166	5.57380448760256e-07\\
167	5.35266465462929e-07\\
168	5.07786723139162e-07\\
169	4.78263512664458e-07\\
170	4.54580521412834e-07\\
171	4.32759859452179e-07\\
172	4.12057090273505e-07\\
173	3.86707679753273e-07\\
174	3.65117209532375e-07\\
175	3.44299976562716e-07\\
176	3.26070292086222e-07\\
177	3.04017259123192e-07\\
178	2.84259605735418e-07\\
179	2.6982392355372e-07\\
180	2.58168178936337e-07\\
181	2.4569729622449e-07\\
182	2.32145182947659e-07\\
183	2.20388581103825e-07\\
184	2.08318967313009e-07\\
185	1.96899388150402e-07\\
186	1.87590301201616e-07\\
187	1.78524732085192e-07\\
188	1.71979541496205e-07\\
189	1.64755781874604e-07\\
190	1.59340966592163e-07\\
191	1.53472216215485e-07\\
192	1.48748803930748e-07\\
193	1.44719437351485e-07\\
194	1.41737759337885e-07\\
195	1.39548644925956e-07\\
196	1.37223483089579e-07\\
197	1.35159280995479e-07\\
198	1.32772102444901e-07\\
199	1.31188810358617e-07\\
200	1.29533721758046e-07\\
201	1.27970059027299e-07\\
202	1.25759632786564e-07\\
203	4.14906968537095e-08\\
204	1.9881745421546e-08\\
205	1.92619460076411e-08\\
206	1.87521693143487e-08\\
207	1.82276988801308e-08\\
208	1.76395650431475e-08\\
209	1.6994207562665e-08\\
210	1.62915282925908e-08\\
211	1.55843875254697e-08\\
212	1.48849130934549e-08\\
213	1.42995212447539e-08\\
214	1.37611210259403e-08\\
215	1.32713805891814e-08\\
216	1.2797440209977e-08\\
217	1.23112593803939e-08\\
218	1.18018397084389e-08\\
219	1.13229928537865e-08\\
220	1.08047174939189e-08\\
221	1.02781420767827e-08\\
222	9.78359945618668e-09\\
223	9.27383656353009e-09\\
224	8.71469997140983e-09\\
225	8.18676096766282e-09\\
226	7.69962397413163e-09\\
227	7.02481607332816e-09\\
228	6.40569553769402e-09\\
229	6.0000845717958e-09\\
230	5.61638982373194e-09\\
231	5.24057578355021e-09\\
232	4.89376551284137e-09\\
233	4.60870318264171e-09\\
234	4.41930242450294e-09\\
235	4.24711918866297e-09\\
236	4.03930130197183e-09\\
237	3.78919758633681e-09\\
238	3.5635250222718e-09\\
239	3.37392147715611e-09\\
240	3.20428523623382e-09\\
241	3.03487786842266e-09\\
242	2.87603556843378e-09\\
243	2.68958978972775e-09\\
244	2.52398257830025e-09\\
245	2.36045312165574e-09\\
246	2.19111225482645e-09\\
247	2.05897460555043e-09\\
248	1.95626542874641e-09\\
249	1.86387400457222e-09\\
250	1.79491874884493e-09\\
251	1.7372703993896e-09\\
252	1.65213233047158e-09\\
253	1.56027452370018e-09\\
254	1.47616285457597e-09\\
255	1.39316557870768e-09\\
256	1.31187836547329e-09\\
257	1.23701656160222e-09\\
258	1.17410392288035e-09\\
259	1.12045494622256e-09\\
260	1.06533903927469e-09\\
261	1.02121627248475e-09\\
262	9.74720901140051e-10\\
263	9.31696920455441e-10\\
264	8.80983088815149e-10\\
265	8.34163174187345e-10\\
266	7.82732888861698e-10\\
267	7.27366947328125e-10\\
268	6.72090884239252e-10\\
269	6.24152726175966e-10\\
270	5.75207282239386e-10\\
271	5.41075500643755e-10\\
272	5.15875138320664e-10\\
273	4.90530435344157e-10\\
274	4.61554560047655e-10\\
275	4.31843419473062e-10\\
276	4.07592369485524e-10\\
277	3.77030057251163e-10\\
278	3.40846947938766e-10\\
279	3.12864095327044e-10\\
280	2.81079506350811e-10\\
281	2.48404837714463e-10\\
282	2.26348084191293e-10\\
283	2.05899540751534e-10\\
284	1.84958882681239e-10\\
285	1.66357348163522e-10\\
286	1.52296677603898e-10\\
287	1.39752783759616e-10\\
288	1.26136414080341e-10\\
289	1.15217097036051e-10\\
290	1.01569160726577e-10\\
291	9.07037077574089e-11\\
292	7.91935901207014e-11\\
293	6.80764321472546e-11\\
294	5.83878728448902e-11\\
295	5.0681690542807e-11\\
296	4.37595795634032e-11\\
297	3.90582178712742e-11\\
298	3.54225877904084e-11\\
299	3.17974769654799e-11\\
300	2.98145349526913e-11\\
301	2.79759617540804e-11\\
302	2.59859537233385e-11\\
303	2.39773661064171e-11\\
304	2.21934023676977e-11\\
305	2.04963986236917e-11\\
306	1.89117058040525e-11\\
307	1.75315685467617e-11\\
308	1.59056903955982e-11\\
309	1.47716982978158e-11\\
310	1.36576529729807e-11\\
311	1.23901603585963e-11\\
312	1.12855792935037e-11\\
313	1.03817512412786e-11\\
314	9.51256235597742e-12\\
315	8.73367445542104e-12\\
316	8.09208070028221e-12\\
317	7.49560462688354e-12\\
318	6.8667825124034e-12\\
319	6.23490444046777e-12\\
320	5.82075945278638e-12\\
321	5.33288916487644e-12\\
322	4.85153206388474e-12\\
323	4.49169907996164e-12\\
324	4.20005509541812e-12\\
325	3.92095993820922e-12\\
326	3.69557292029218e-12\\
327	3.50960893478848e-12\\
328	3.2857128027634e-12\\
329	3.09917020734444e-12\\
330	2.9343256710742e-12\\
331	2.72324917323446e-12\\
332	2.48856066206075e-12\\
333	2.2201090292773e-12\\
334	2.01754901015386e-12\\
335	1.85373526009388e-12\\
336	1.70943585936379e-12\\
337	1.55774853227264e-12\\
338	1.43811122311592e-12\\
339	1.31012352793523e-12\\
340	1.15808987898888e-12\\
341	1.02413298889594e-12\\
342	9.3505284524617e-13\\
343	8.54897687202101e-13\\
344	7.9216194202605e-13\\
345	7.48100592075553e-13\\
346	7.24515953490572e-13\\
347	7.03241445104924e-13\\
348	6.80422647946023e-13\\
349	6.52669779529667e-13\\
350	6.27848201469921e-13\\
351	6.04878998732279e-13\\
352	5.85770417899984e-13\\
353	5.71354969981586e-13\\
354	5.59320084684387e-13\\
355	5.5144797003988e-13\\
356	5.46401164853595e-13\\
357	5.42423092389786e-13\\
358	5.3985560955161e-13\\
359	5.38058112273924e-13\\
360	5.36940204523145e-13\\
361	5.36455755581864e-13\\
362	5.36112747436155e-13\\
363	5.35736940887629e-13\\
364	5.35525389687971e-13\\
365	5.354074988518e-13\\
366	5.35349353986767e-13\\
367	5.35316079178016e-13\\
368	5.35296624062713e-13\\
369	5.35290283049652e-13\\
370	5.35287715930491e-13\\
371	5.35285605087823e-13\\
372	5.35283544593243e-13\\
373	5.3528243858128e-13\\
374	5.35281383893879e-13\\
375	5.35280423729163e-13\\
376	5.35279876751442e-13\\
377	5.35279531394129e-13\\
378	5.35279346342544e-13\\
379	5.35279263210201e-13\\
380	5.35279193692374e-13\\
381	5.3527914926656e-13\\
382	5.35279129592966e-13\\
383	5.35279119195483e-13\\
384	5.35279114291195e-13\\
385	5.35279111358897e-13\\
386	5.35279109824452e-13\\
387	5.35279108729194e-13\\
388	5.35279108182209e-13\\
389	5.35279107806554e-13\\
390	5.35279107712304e-13\\
391	5.35279107679432e-13\\
392	5.35279107672409e-13\\
393	5.35279107670694e-13\\
394	5.35279107670201e-13\\
395	5.35279107670074e-13\\
396	5.35279107670067e-13\\
397	5.35279107670048e-13\\
398	5.35279107669755e-13\\
399	5.35279107665155e-13\\
400	5.35278217336914e-13\\
401	2.50145971591972e-13\\
402	1.85707176766868e-14\\
403	1.8570651656525e-14\\
404	1.85309204817424e-14\\
405	1.79657492628122e-14\\
406	1.79644098864886e-14\\
407	1.79395067949614e-14\\
408	1.79360489259315e-14\\
409	1.79301838881997e-14\\
410	1.79301449370623e-14\\
411	1.79167495409676e-14\\
412	1.78983281282222e-14\\
413	1.78153971324305e-14\\
414	1.76573297776721e-14\\
415	1.75990104997549e-14\\
416	1.73509966247636e-14\\
417	1.73509781722511e-14\\
418	1.72960670864124e-14\\
419	1.72829773965749e-14\\
420	1.72613179795594e-14\\
421	1.70833515621026e-14\\
422	1.70442157773959e-14\\
423	1.7043999358355e-14\\
424	1.70374610229136e-14\\
425	1.70306580005498e-14\\
426	1.70285870622673e-14\\
427	1.6953346699412e-14\\
428	1.69002313642822e-14\\
429	1.68560664834387e-14\\
430	1.68527077914394e-14\\
431	1.68311729177635e-14\\
432	1.66929580087546e-14\\
433	1.66887262970055e-14\\
434	1.66687143378863e-14\\
435	1.66621058655121e-14\\
436	1.6661505233369e-14\\
437	1.66472071014098e-14\\
438	1.66336279718235e-14\\
439	1.66207514520612e-14\\
440	1.6534718107864e-14\\
441	1.65193629343217e-14\\
442	1.64923294941336e-14\\
443	1.64914558448677e-14\\
444	1.64865897251971e-14\\
445	1.64175038422923e-14\\
446	1.6332298208475e-14\\
447	1.63320377127502e-14\\
448	1.63204238986245e-14\\
449	1.63186994887327e-14\\
450	1.61757221431064e-14\\
};
\addlegendentry{$\varepsilon=10^{-8}$}

\addplot [color=mycolor2, only marks, mark=x, mark options={solid, mycolor2}]
  table[row sep=crcr]{%
0	1\\
1	0.204873076853183\\
2	0.0896638436885824\\
3	0.0463024822230975\\
4	0.026444441475166\\
5	0.0175688832550397\\
6	0.012305913433854\\
7	0.00880394317875808\\
8	0.00641928434821926\\
9	0.00496391766677556\\
10	0.0039676741428587\\
11	0.00342507847467089\\
12	0.00335515833075814\\
13	0.00318411403938736\\
14	0.00272699806309163\\
15	0.00230600249647438\\
16	0.00197554154309841\\
17	0.00169740618355902\\
18	0.00146513344323079\\
19	0.00129623468816415\\
20	0.00123430563601103\\
21	0.00122868248532083\\
22	0.00111338880309911\\
23	0.00097929595553097\\
24	0.000880071043711044\\
25	0.000794054297299536\\
26	0.00071355158390216\\
27	0.000636969118179518\\
28	0.00057869609260352\\
29	0.000572527561013378\\
30	0.000551278209418292\\
31	0.000498762732041113\\
32	0.000459451218183874\\
33	0.000421149962032199\\
34	0.000388621671831447\\
35	0.000360355505328162\\
36	0.000339379053412178\\
37	0.000320749423806627\\
38	0.000298324464306588\\
39	0.000276918760334415\\
40	0.000260988496392557\\
41	0.000243438538577909\\
42	0.00022829586347358\\
43	0.000212923563213669\\
44	0.000199324745611857\\
45	0.000186276331702862\\
46	0.000174104733277287\\
47	0.000162164852994784\\
48	0.000149192634876973\\
49	0.000137223393325162\\
50	0.000125845343204551\\
51	0.00011630177623496\\
52	0.000109032949647926\\
53	0.000102851285874385\\
54	9.69529876796453e-05\\
55	9.1526811539011e-05\\
56	8.72530160287387e-05\\
57	8.32339112762881e-05\\
58	7.89039759261875e-05\\
59	7.45592547703305e-05\\
60	6.97369121993351e-05\\
61	6.57221086341215e-05\\
62	6.17467792973021e-05\\
63	5.79757051387326e-05\\
64	5.45472264065342e-05\\
65	5.14318170174499e-05\\
66	4.83439814261071e-05\\
67	4.56172940570441e-05\\
68	4.34804141283036e-05\\
69	4.13897038821425e-05\\
70	3.908737992242e-05\\
71	3.61946889399263e-05\\
72	3.30824510999805e-05\\
73	3.01071077381775e-05\\
74	2.78470068324506e-05\\
75	2.58462763828652e-05\\
76	2.38146998799261e-05\\
77	2.2073321518083e-05\\
78	2.05240529296432e-05\\
79	1.94294547625201e-05\\
80	1.83175177815225e-05\\
81	1.71970062541814e-05\\
82	1.62937854346943e-05\\
83	1.53638616501595e-05\\
84	1.46585157367257e-05\\
85	1.40473940538537e-05\\
86	1.34995624028957e-05\\
87	1.29320709235861e-05\\
88	1.22786506811582e-05\\
89	1.17349240712535e-05\\
90	1.11692514517266e-05\\
91	1.06032307588696e-05\\
92	1.01127223636205e-05\\
93	9.68943568133372e-06\\
94	9.24698671645624e-06\\
95	8.83035105539514e-06\\
96	8.5052537442527e-06\\
97	8.22571410119725e-06\\
98	7.9790962261581e-06\\
99	7.75588376694868e-06\\
100	7.56690337744251e-06\\
101	7.37788514318447e-06\\
102	7.20304856746208e-06\\
103	7.03603804435831e-06\\
104	6.83741337414166e-06\\
105	6.64214276154945e-06\\
106	6.44021623201688e-06\\
107	6.27563251612641e-06\\
108	6.11077608769637e-06\\
109	5.93513141889584e-06\\
110	5.77383502509427e-06\\
111	5.63487926770061e-06\\
112	5.48097511844169e-06\\
113	5.28988089965023e-06\\
114	5.11075382247135e-06\\
115	4.94685832915017e-06\\
116	4.7503929807476e-06\\
117	4.55261230089858e-06\\
118	4.37921848851311e-06\\
119	4.19144765845995e-06\\
120	4.02836325664481e-06\\
121	3.89096845470994e-06\\
122	3.7436777736047e-06\\
123	3.59064575922727e-06\\
124	3.43642309936716e-06\\
125	3.30415551162931e-06\\
126	3.16507143583204e-06\\
127	3.03568721263614e-06\\
128	2.90577318678377e-06\\
129	2.78661965371995e-06\\
130	2.67871474404372e-06\\
131	2.56545403290276e-06\\
132	2.46813937376777e-06\\
133	2.37356119829938e-06\\
134	2.26166529827834e-06\\
135	2.14887913619319e-06\\
136	2.06701046823032e-06\\
137	1.98500897862511e-06\\
138	1.91099965021173e-06\\
139	1.8496444705885e-06\\
140	1.77903112615801e-06\\
141	1.70115235922314e-06\\
142	1.62796229091186e-06\\
143	1.56288267326219e-06\\
144	1.49998158206313e-06\\
145	1.44576206453688e-06\\
146	1.38540967047257e-06\\
147	1.32811465779569e-06\\
148	1.28080984335202e-06\\
149	1.22945622035104e-06\\
150	1.17132625273039e-06\\
151	1.11748064143995e-06\\
152	1.05889744407706e-06\\
153	9.96773288752948e-07\\
154	9.41201471704884e-07\\
155	8.88565310152969e-07\\
156	8.29923690287624e-07\\
157	7.82264808637795e-07\\
158	7.46837860743936e-07\\
159	7.13646429174323e-07\\
160	6.87804772688658e-07\\
161	6.64190931089628e-07\\
162	6.46684606433557e-07\\
163	6.26955695656546e-07\\
164	6.08139447136933e-07\\
165	5.82931536813227e-07\\
166	5.57021298545119e-07\\
167	5.3492063800707e-07\\
168	5.07456922866434e-07\\
169	4.77949125504114e-07\\
170	4.54277072549554e-07\\
171	4.32465649211947e-07\\
172	4.11771008903885e-07\\
173	3.86429979067957e-07\\
174	3.64844875503083e-07\\
175	3.44031328926182e-07\\
176	3.25803775488112e-07\\
177	3.0375135324205e-07\\
178	2.83991939660642e-07\\
179	2.6955383793106e-07\\
180	2.5789569094329e-07\\
181	2.45421997772868e-07\\
182	2.31866149085658e-07\\
183	2.20105609040202e-07\\
184	2.08031348496598e-07\\
185	1.96606616086342e-07\\
186	1.8729310270985e-07\\
187	1.78223117765667e-07\\
188	1.716750279309e-07\\
189	1.64448664850064e-07\\
190	1.59032466777549e-07\\
191	1.5316299232553e-07\\
192	1.48439693340142e-07\\
193	1.44411388962856e-07\\
194	1.41431428679706e-07\\
195	1.39244690487783e-07\\
196	1.36923286463988e-07\\
197	1.34863727671199e-07\\
198	1.32483208241226e-07\\
199	1.3090526973185e-07\\
200	1.29257131154706e-07\\
201	1.27701117014472e-07\\
202	1.25502497562117e-07\\
203	4.14743247617514e-08\\
204	1.98670962423369e-08\\
205	1.92470385751296e-08\\
206	1.87371373871325e-08\\
207	1.82126336711446e-08\\
208	1.76245572231624e-08\\
209	1.69793509509789e-08\\
210	1.62769121050162e-08\\
211	1.55700815422825e-08\\
212	1.4870969396168e-08\\
213	1.42859380962197e-08\\
214	1.37479332806777e-08\\
215	1.32586162892663e-08\\
216	1.27851476838997e-08\\
217	1.2299515861603e-08\\
218	1.17907061722622e-08\\
219	1.13124708840073e-08\\
220	1.07948879612951e-08\\
221	1.02690268677036e-08\\
222	9.77515879037933e-09\\
223	9.26609068576513e-09\\
224	8.7076941349272e-09\\
225	8.18042890954222e-09\\
226	7.69388913156395e-09\\
227	7.01984118452465e-09\\
228	6.40133329661035e-09\\
229	5.99608602020103e-09\\
230	5.61271529362731e-09\\
231	5.23719486803374e-09\\
232	4.8906329456306e-09\\
233	4.6057606611218e-09\\
234	4.416483794705e-09\\
235	4.24441608066964e-09\\
236	4.03673894507478e-09\\
237	3.78679476960684e-09\\
238	3.56125546134038e-09\\
239	3.37175805786301e-09\\
240	3.20221588256072e-09\\
241	3.03290242587725e-09\\
242	2.87414727622625e-09\\
243	2.68780113586205e-09\\
244	2.52227857039441e-09\\
245	2.35883005182582e-09\\
246	2.18956829033664e-09\\
247	2.05748929514911e-09\\
248	1.95482871163476e-09\\
249	1.86248411053114e-09\\
250	1.79356946730681e-09\\
251	1.73596210517052e-09\\
252	1.65089274767314e-09\\
253	1.55911124098918e-09\\
254	1.47507149509764e-09\\
255	1.39214735278378e-09\\
256	1.3109334126447e-09\\
257	1.23613991327702e-09\\
258	1.17328647800355e-09\\
259	1.11968971138241e-09\\
260	1.06463030437638e-09\\
261	1.02055414125436e-09\\
262	9.74109725849359e-10\\
263	9.31133553976917e-10\\
264	8.80476309884764e-10\\
265	8.33707824490671e-10\\
266	7.82332484309569e-10\\
267	7.2702248567203e-10\\
268	6.71798306813283e-10\\
269	6.23902012294564e-10\\
270	5.74996125425842e-10\\
271	5.40890060838265e-10\\
272	5.15708172849455e-10\\
273	4.90381961178333e-10\\
274	4.61426302720733e-10\\
275	4.31734532008657e-10\\
276	4.07498382277995e-10\\
277	3.76953675887438e-10\\
278	3.40789199732333e-10\\
279	3.12819239215706e-10\\
280	2.81047720742244e-10\\
281	2.48384729077341e-10\\
282	2.2633486900704e-10\\
283	2.05892045459762e-10\\
284	1.84956596943965e-10\\
285	1.66359186471914e-10\\
286	1.52301322167615e-10\\
287	1.39759706824219e-10\\
288	1.26145614118009e-10\\
289	1.15227973067976e-10\\
290	1.0158202504299e-10\\
291	9.07181990169275e-11\\
292	7.92098839299228e-11\\
293	6.80947749862134e-11\\
294	5.84084644554989e-11\\
295	5.07045286682495e-11\\
296	4.37850124952291e-11\\
297	3.90858138659577e-11\\
298	3.54521006651711e-11\\
299	3.18292601266416e-11\\
300	2.9847695489708e-11\\
301	2.80104190532776e-11\\
302	2.60219452287074e-11\\
303	2.40151690058418e-11\\
304	2.22331367099376e-11\\
305	2.05383258083178e-11\\
306	1.89561066705519e-11\\
307	1.75785424993266e-11\\
308	1.59563930278255e-11\\
309	1.48255678119021e-11\\
310	1.37151988518278e-11\\
311	1.24527884737776e-11\\
312	1.13536837967378e-11\\
313	1.04552693089185e-11\\
314	9.59230339178508e-12\\
315	8.82008258195373e-12\\
316	8.18496263441251e-12\\
317	7.59550239427828e-12\\
318	6.97539683279494e-12\\
319	6.35404722926143e-12\\
320	5.94802348496722e-12\\
321	5.47129570047574e-12\\
322	5.00306292419106e-12\\
323	4.65480015601751e-12\\
324	4.37390964322153e-12\\
325	4.10650490571518e-12\\
326	3.89175072914992e-12\\
327	3.71550219662574e-12\\
328	3.50463385995791e-12\\
329	3.3302477935278e-12\\
330	3.17729440942847e-12\\
331	2.98330271463612e-12\\
332	2.77056831657006e-12\\
333	2.53202683951134e-12\\
334	2.35632497617333e-12\\
335	2.21759666807851e-12\\
336	2.09839992444025e-12\\
337	1.97671946864586e-12\\
338	1.88382885961749e-12\\
339	1.78798494532165e-12\\
340	1.6797229745491e-12\\
341	1.59029436414697e-12\\
342	1.53442353446088e-12\\
343	1.48692273614575e-12\\
344	1.45175111440267e-12\\
345	1.42817939234638e-12\\
346	1.41596412055132e-12\\
347	1.40519316243761e-12\\
348	1.39390860594585e-12\\
349	1.38056835534268e-12\\
350	1.36900438316246e-12\\
351	1.35862035415477e-12\\
352	1.35021911190272e-12\\
353	1.34402659518656e-12\\
354	1.33895391615576e-12\\
355	1.33568418653721e-12\\
356	1.33360823358491e-12\\
357	1.33198311156006e-12\\
358	1.33093951213863e-12\\
359	1.33021135712343e-12\\
360	1.32975952839661e-12\\
361	1.32956396831743e-12\\
362	1.32942559218126e-12\\
363	1.32927407424091e-12\\
364	1.32918882326202e-12\\
365	1.32914132878467e-12\\
366	1.32911790800234e-12\\
367	1.32910450614163e-12\\
368	1.3290966709306e-12\\
369	1.32909411734174e-12\\
370	1.32909308363234e-12\\
371	1.32909223371972e-12\\
372	1.32909140410574e-12\\
373	1.32909095879766e-12\\
374	1.32909053415792e-12\\
375	1.32909014754742e-12\\
376	1.32908992728283e-12\\
377	1.32908978820225e-12\\
378	1.32908971368217e-12\\
379	1.32908968020952e-12\\
380	1.32908965221872e-12\\
381	1.32908963432854e-12\\
382	1.32908962640539e-12\\
383	1.32908962221909e-12\\
384	1.32908962024491e-12\\
385	1.32908961906449e-12\\
386	1.32908961844645e-12\\
387	1.32908961800515e-12\\
388	1.32908961778489e-12\\
389	1.32908961763381e-12\\
390	1.32908961759583e-12\\
391	1.3290896175825e-12\\
392	1.32908961757966e-12\\
393	1.32908961757897e-12\\
394	1.32908961757878e-12\\
395	1.32908961757872e-12\\
396	1.32908961757871e-12\\
397	1.3290896175787e-12\\
398	1.3290896175787e-12\\
399	1.3290896175787e-12\\
400	1.32908961742205e-12\\
401	1.32903748847442e-12\\
402	1.32767104416144e-12\\
403	1.32767104371833e-12\\
404	1.32488632025446e-12\\
405	1.32403755325971e-12\\
406	1.32018040910288e-12\\
407	1.31835346025723e-12\\
408	1.30659036239123e-12\\
409	1.30656790312642e-12\\
410	1.30516915802198e-12\\
411	1.30502243736324e-12\\
412	1.30440714755362e-12\\
413	1.30333646921339e-12\\
414	1.30291641210692e-12\\
415	1.29742093272487e-12\\
416	1.29618335738378e-12\\
417	1.29098930932238e-12\\
418	1.29017348988048e-12\\
419	1.28623491643437e-12\\
420	1.28550237407896e-12\\
421	1.26247704324215e-12\\
422	1.25991467514233e-12\\
423	1.25988013877365e-12\\
424	1.25975426501373e-12\\
425	1.2593707730084e-12\\
426	1.25933732346775e-12\\
427	1.24896051035594e-12\\
428	1.24460193305558e-12\\
429	1.24388223327933e-12\\
430	1.24317807852082e-12\\
431	1.24081295037811e-12\\
432	1.24061574412816e-12\\
433	1.23319521001008e-12\\
434	1.22426165599379e-12\\
435	1.2241097515881e-12\\
436	1.21479895362731e-12\\
437	1.21448744742964e-12\\
438	1.21433242969283e-12\\
439	1.21425706168687e-12\\
440	1.21378153659143e-12\\
441	1.21376679567666e-12\\
442	1.20988168864482e-12\\
443	1.20810555477956e-12\\
444	1.20746575660818e-12\\
445	1.19542889695642e-12\\
446	1.18961659342585e-12\\
447	1.18948353809083e-12\\
448	1.18628128865227e-12\\
449	1.18579977814486e-12\\
450	1.17946680452647e-12\\
};
\addlegendentry{$\varepsilon=10^{-10}$}

\addplot [color=mycolor3, only marks, mark=x, mark options={solid, mycolor3}]
  table[row sep=crcr]{%
0	1\\
1	0.20487249005328\\
2	0.0896637584651119\\
3	0.0463025178414623\\
4	0.0264447810631194\\
5	0.0175690797389075\\
6	0.012305915317045\\
7	0.00880386500852525\\
8	0.0064191165442984\\
9	0.00496370982774242\\
10	0.00396805958184122\\
11	0.00342939196965014\\
12	0.00336516571875689\\
13	0.00318610746200696\\
14	0.00272637965506713\\
15	0.00230495797845032\\
16	0.00197436029569957\\
17	0.00169615084534183\\
18	0.00146386684403284\\
19	0.00129527140150326\\
20	0.00123449884033008\\
21	0.00122828655841353\\
22	0.00111207713379842\\
23	0.00097780653707637\\
24	0.000878499512832847\\
25	0.000792434648526043\\
26	0.00071223856968024\\
27	0.000639006196768241\\
28	0.000600318535700703\\
29	0.00059906499997\\
30	0.000554449940129844\\
31	0.000497088293452811\\
32	0.000457358372261883\\
33	0.000418966392360043\\
34	0.000386392412709735\\
35	0.000358098590112339\\
36	0.000337115714417641\\
37	0.000318501579890074\\
38	0.000296116623096657\\
39	0.000274761709990763\\
40	0.00025888123107209\\
41	0.000241400607296547\\
42	0.000226326754506182\\
43	0.000211031141903022\\
44	0.000197504901899617\\
45	0.000184531872618664\\
46	0.000172434772618516\\
47	0.000160571428435107\\
48	0.000147682110638983\\
49	0.000135783589806731\\
50	0.000124464768570094\\
51	0.000114964068869078\\
52	0.000107725168143183\\
53	0.000101569182540387\\
54	9.56970805353291e-05\\
55	9.02958684796947e-05\\
56	8.60454521899982e-05\\
57	8.20524232769656e-05\\
58	7.77544151063129e-05\\
59	7.34434034428172e-05\\
60	6.86571287035053e-05\\
61	6.46705671694908e-05\\
62	6.07220850741039e-05\\
63	5.69756471778079e-05\\
64	5.35692866386733e-05\\
65	5.04750215614417e-05\\
66	4.74089535327912e-05\\
67	4.47024903301199e-05\\
68	4.25842332980516e-05\\
69	4.05147983614309e-05\\
70	3.8238406341534e-05\\
71	3.53777850309478e-05\\
72	3.22948529921645e-05\\
73	2.93406690592629e-05\\
74	2.70930888134718e-05\\
75	2.51015276415001e-05\\
76	2.30765761797464e-05\\
77	2.13382755700414e-05\\
78	1.97899144571728e-05\\
79	1.86959382062148e-05\\
80	1.75859935620716e-05\\
81	1.64671389015557e-05\\
82	1.55656300776233e-05\\
83	1.46380486593935e-05\\
84	1.39352095103364e-05\\
85	1.33278887944584e-05\\
86	1.27852972715993e-05\\
87	1.22256384446736e-05\\
88	1.15823407105869e-05\\
89	1.10476219943632e-05\\
90	1.04926220635207e-05\\
91	9.93760417845707e-06\\
92	9.45669152279122e-06\\
93	9.04239219395398e-06\\
94	8.61006909246188e-06\\
95	8.20338644262205e-06\\
96	7.88688020865718e-06\\
97	7.61610121002357e-06\\
98	7.37887678106877e-06\\
99	7.1657188868904e-06\\
100	6.98722768485253e-06\\
101	6.81068584211396e-06\\
102	6.64907765042936e-06\\
103	6.49663915545895e-06\\
104	6.31706117351515e-06\\
105	6.14177473127373e-06\\
106	5.96141532142068e-06\\
107	5.81540375790784e-06\\
108	5.6702309328842e-06\\
109	5.51639961694291e-06\\
110	5.37583240714828e-06\\
111	5.25542004865629e-06\\
112	5.12264783787333e-06\\
113	4.95769288547699e-06\\
114	4.80250222182081e-06\\
115	4.6601301846294e-06\\
116	4.4885310104264e-06\\
117	4.31428775277415e-06\\
118	4.16025998124729e-06\\
119	3.99208383383768e-06\\
120	3.84473144481805e-06\\
121	3.71977526567562e-06\\
122	3.58501645899648e-06\\
123	3.44397750158119e-06\\
124	3.3006660628308e-06\\
125	3.17683733968396e-06\\
126	3.04572680603787e-06\\
127	2.92295684794508e-06\\
128	2.79890424926104e-06\\
129	2.68444943879467e-06\\
130	2.58035069321923e-06\\
131	2.47066376829406e-06\\
132	2.37608902251341e-06\\
133	2.28402366149651e-06\\
134	2.17477769996837e-06\\
135	2.06419512317314e-06\\
136	1.98375890709441e-06\\
137	1.90323319154152e-06\\
138	1.83060632842762e-06\\
139	1.7705248803566e-06\\
140	1.70157756904952e-06\\
141	1.62555586395173e-06\\
142	1.55414730594735e-06\\
143	1.490710800978e-06\\
144	1.42948963948054e-06\\
145	1.3768871986724e-06\\
146	1.31849453742484e-06\\
147	1.26313394699193e-06\\
148	1.2175663730733e-06\\
149	1.16832020889959e-06\\
150	1.11262344087221e-06\\
151	1.06101847340772e-06\\
152	1.00481629297087e-06\\
153	9.45054771926282e-07\\
154	8.91420391175357e-07\\
155	8.40442479811654e-07\\
156	7.83356940580865e-07\\
157	7.367419450609e-07\\
158	7.02036718853825e-07\\
159	6.69553087417925e-07\\
160	6.44359970893351e-07\\
161	6.21480174991287e-07\\
162	6.0467980282361e-07\\
163	5.85947023025894e-07\\
164	5.68207461759245e-07\\
165	5.4458132053181e-07\\
166	5.20293204607077e-07\\
167	4.99586750975217e-07\\
168	4.73821296639363e-07\\
169	4.45981185315401e-07\\
170	4.23520217700726e-07\\
171	4.02752012948162e-07\\
172	3.82992152813198e-07\\
173	3.58648235871436e-07\\
174	3.3774092721447e-07\\
175	3.17433121187901e-07\\
176	2.99536107657923e-07\\
177	2.77678472302837e-07\\
178	2.57848181854104e-07\\
179	2.43229289500415e-07\\
180	2.31366698718027e-07\\
181	2.18632995776312e-07\\
182	2.04703229863318e-07\\
183	1.92524344821713e-07\\
184	1.7993053645367e-07\\
185	1.6790392898312e-07\\
186	1.58045471824884e-07\\
187	1.48400964784644e-07\\
188	1.41443841838477e-07\\
189	1.33801355296622e-07\\
190	1.28112521087752e-07\\
191	1.22007206800921e-07\\
192	1.1715001569249e-07\\
193	1.13095201409101e-07\\
194	1.10184127805837e-07\\
195	1.0815647825945e-07\\
196	1.06122765222223e-07\\
197	1.04442078174825e-07\\
198	1.02619288447619e-07\\
199	1.01496223539242e-07\\
200	1.00441340162731e-07\\
201	9.95372685922883e-08\\
202	9.8345879694681e-08\\
203	3.94422253022261e-08\\
204	1.84128694773651e-08\\
205	1.77716017506525e-08\\
206	1.72522346034486e-08\\
207	1.67266273723087e-08\\
208	1.6145844373102e-08\\
209	1.55164669479391e-08\\
210	1.48378350716202e-08\\
211	1.41608521683816e-08\\
212	1.34960000230137e-08\\
213	1.29446925343677e-08\\
214	1.24434606771431e-08\\
215	1.19933955186393e-08\\
216	1.15636096798213e-08\\
217	1.11288173546544e-08\\
218	1.06765981743278e-08\\
219	1.02552480449699e-08\\
220	9.80230719132811e-09\\
221	9.34346534771054e-09\\
222	8.91325054183267e-09\\
223	8.47019359099628e-09\\
224	7.98269484426424e-09\\
225	7.5205706688057e-09\\
226	7.09231389419021e-09\\
227	6.49328044382872e-09\\
228	5.93604696414512e-09\\
229	5.56757345880515e-09\\
230	5.21729457164992e-09\\
231	4.87210162738298e-09\\
232	4.55151981623963e-09\\
233	4.28679613283934e-09\\
234	4.11082354577508e-09\\
235	3.95130333583085e-09\\
236	3.7591001202258e-09\\
237	3.52695700071465e-09\\
238	3.31652618766673e-09\\
239	3.13920972197456e-09\\
240	2.98053438511899e-09\\
241	2.82211003482723e-09\\
242	2.67349368562547e-09\\
243	2.49878178437888e-09\\
244	2.34317222539753e-09\\
245	2.18919643863144e-09\\
246	2.02919848501544e-09\\
247	1.9039689386595e-09\\
248	1.80687229158971e-09\\
249	1.71978436657112e-09\\
250	1.65530513142587e-09\\
251	1.60205664264907e-09\\
252	1.52416250240424e-09\\
253	1.44027616254212e-09\\
254	1.36360796993903e-09\\
255	1.28811193986492e-09\\
256	1.21428141587114e-09\\
257	1.14632784911374e-09\\
258	1.08936849309768e-09\\
259	1.04094828145138e-09\\
260	9.91472548958402e-10\\
261	9.51997862473921e-10\\
262	9.10583845007774e-10\\
263	8.72335008613745e-10\\
264	8.272879043302e-10\\
265	7.8564190196139e-10\\
266	7.39769929063204e-10\\
267	6.90102099516123e-10\\
268	6.40150685163467e-10\\
269	5.96543346529097e-10\\
270	5.5172532434992e-10\\
271	5.20300855324888e-10\\
272	4.97063708642319e-10\\
273	4.73705315106844e-10\\
274	4.46923068980253e-10\\
275	4.19335451828312e-10\\
276	3.9673300959085e-10\\
277	3.68135007746603e-10\\
278	3.34044397396987e-10\\
279	3.07515223271198e-10\\
280	2.77211447527393e-10\\
281	2.45862212428205e-10\\
282	2.24589140542136e-10\\
283	2.04791510346384e-10\\
284	1.84441830040755e-10\\
285	1.6630533219015e-10\\
286	1.52559286615098e-10\\
287	1.40269780703663e-10\\
288	1.26905785896734e-10\\
289	1.16169865369449e-10\\
290	1.02736029847336e-10\\
291	9.20396142048962e-11\\
292	8.07104983998507e-11\\
293	6.97911130571595e-11\\
294	6.03113281136305e-11\\
295	5.28074784547217e-11\\
296	4.61121809089203e-11\\
297	4.1594249723605e-11\\
298	3.81162047563433e-11\\
299	3.4672477684998e-11\\
300	3.27959947631008e-11\\
301	3.10531746835868e-11\\
302	2.91726680893539e-11\\
303	2.72902159956925e-11\\
304	2.56381947495002e-11\\
305	2.40876091270989e-11\\
306	2.26643503343798e-11\\
307	2.14476475066015e-11\\
308	2.00518248905789e-11\\
309	1.91062692581577e-11\\
310	1.82026380769658e-11\\
311	1.72120903542407e-11\\
312	1.6388948394942e-11\\
313	1.57460546253642e-11\\
314	1.51565730401332e-11\\
315	1.46551734847857e-11\\
316	1.42627328306369e-11\\
317	1.39156502330322e-11\\
318	1.35700359142159e-11\\
319	1.32459481280565e-11\\
320	1.3046628572579e-11\\
321	1.28256391517466e-11\\
322	1.26234266410304e-11\\
323	1.24827094830984e-11\\
324	1.23752550307351e-11\\
325	1.22779115738145e-11\\
326	1.22032034909799e-11\\
327	1.21441820695488e-11\\
328	1.20764942724882e-11\\
329	1.20230938407914e-11\\
330	1.1978101116993e-11\\
331	1.19238197534653e-11\\
332	1.18682139146921e-11\\
333	1.18108812584945e-11\\
334	1.17720692164337e-11\\
335	1.1743401273867e-11\\
336	1.17201439292563e-11\\
337	1.16977327275087e-11\\
338	1.16815331312008e-11\\
339	1.16656598604617e-11\\
340	1.16488173541718e-11\\
341	1.1635786412245e-11\\
342	1.16280426836251e-11\\
343	1.16216889108887e-11\\
344	1.16171209087894e-11\\
345	1.16141193217387e-11\\
346	1.16125785124212e-11\\
347	1.16112221765414e-11\\
348	1.16098092521865e-11\\
349	1.16081553382884e-11\\
350	1.16067410022715e-11\\
351	1.16054864483313e-11\\
352	1.16044832674129e-11\\
353	1.16037504509923e-11\\
354	1.16031550943494e-11\\
355	1.16027734099163e-11\\
356	1.16025318258223e-11\\
357	1.16023430999853e-11\\
358	1.16022220396032e-11\\
359	1.16021376796511e-11\\
360	1.1602085393631e-11\\
361	1.16020627251017e-11\\
362	1.16020466579354e-11\\
363	1.16020291183979e-11\\
364	1.1602019290418e-11\\
365	1.16020138268649e-11\\
366	1.16020111413059e-11\\
367	1.16020096070629e-11\\
368	1.16020087121893e-11\\
369	1.16020084211068e-11\\
370	1.16020083038155e-11\\
371	1.16020082077777e-11\\
372	1.16020081142042e-11\\
373	1.16020080639639e-11\\
374	1.16020080160059e-11\\
375	1.16020079721766e-11\\
376	1.16020079470744e-11\\
377	1.16020079311746e-11\\
378	1.16020079226631e-11\\
379	1.16020079188619e-11\\
380	1.16020079156962e-11\\
381	1.16020079136735e-11\\
382	1.16020079127784e-11\\
383	1.16020079123089e-11\\
384	1.16020079120881e-11\\
385	1.1602007911956e-11\\
386	1.16020079118864e-11\\
387	1.16020079118361e-11\\
388	1.16020079118112e-11\\
389	1.16020079117943e-11\\
390	1.16020079117899e-11\\
391	1.16020079117882e-11\\
392	1.16020079117879e-11\\
393	1.16020079117879e-11\\
394	1.16020079117878e-11\\
395	1.16020079117878e-11\\
396	1.16020079117878e-11\\
397	1.16020079117878e-11\\
398	1.16020079117878e-11\\
399	1.16020079117876e-11\\
400	1.1602007794178e-11\\
401	1.15576629814924e-11\\
402	1.15423206372341e-11\\
403	1.15423206372337e-11\\
404	1.15382502090297e-11\\
405	1.14423549835269e-11\\
406	1.1340937511259e-11\\
407	1.13225590102094e-11\\
408	1.13008823718759e-11\\
409	1.12699489239231e-11\\
410	1.11498794852929e-11\\
411	1.11486591842296e-11\\
412	1.11458853252225e-11\\
413	1.11369594477046e-11\\
414	1.1128468868061e-11\\
415	1.10974115283353e-11\\
416	1.10738433264735e-11\\
417	1.09682088975578e-11\\
418	1.07983568683633e-11\\
419	1.07584312619696e-11\\
420	1.07576462063732e-11\\
421	1.07571449519337e-11\\
422	1.06301411781446e-11\\
423	1.06275108301172e-11\\
424	1.04990742971598e-11\\
425	1.04060837086503e-11\\
426	1.03244790226688e-11\\
427	1.02193732057722e-11\\
428	1.02160868752749e-11\\
429	1.01344789757166e-11\\
430	1.00865977253391e-11\\
431	1.00763963195619e-11\\
432	1.00639110119682e-11\\
433	9.87494089963187e-12\\
434	9.83607475043113e-12\\
435	9.64303303051068e-12\\
436	9.6403901289454e-12\\
437	9.63605037354701e-12\\
438	9.63469850795442e-12\\
439	9.63411493340951e-12\\
440	9.60641217751416e-12\\
441	9.60445356915996e-12\\
442	9.60357873250702e-12\\
443	9.60121055599458e-12\\
444	9.56796695560878e-12\\
445	9.52490226942555e-12\\
446	9.48712413648138e-12\\
447	9.41285009309626e-12\\
448	9.39146507541696e-12\\
449	9.31783798192359e-12\\
450	9.28389060035639e-12\\
};
\addlegendentry{$\varepsilon=10^{-12}$}

\addplot [color=mycolor4, only marks, mark=x, mark options={solid, mycolor4}]
  table[row sep=crcr]{%
0	1\\
1	0.204894082574935\\
2	0.0896841470305176\\
3	0.046363942611835\\
4	0.0265353175255953\\
5	0.0176531672342712\\
6	0.0123816540342985\\
7	0.00889949414586309\\
8	0.00661905891161701\\
9	0.00567771580045686\\
10	0.00567770555778162\\
11	0.00491271568415579\\
12	0.00400359352794484\\
13	0.00337014833273293\\
14	0.00284898695631678\\
15	0.00244800440098295\\
16	0.00223950163919241\\
17	0.00223879400524219\\
18	0.00207831226282461\\
19	0.00183813882202298\\
20	0.00163089429570918\\
21	0.00148238171064782\\
22	0.00134118043362558\\
23	0.0012220425016681\\
24	0.00113734213194291\\
25	0.00107187543405298\\
26	0.00105338581310693\\
27	0.00104256401183392\\
28	0.000984473036000102\\
29	0.000915985573785414\\
30	0.000857513859663388\\
31	0.000808102114462644\\
32	0.00077525762099701\\
33	0.000742815662134326\\
34	0.000714498485720318\\
35	0.000688818947427443\\
36	0.000668457956254586\\
37	0.000648715970275314\\
38	0.00062375029175377\\
39	0.000599384868606252\\
40	0.000580315210460185\\
41	0.000557919340303471\\
42	0.000537653681879763\\
43	0.000516124330831342\\
44	0.000496433748772787\\
45	0.000476886911785362\\
46	0.000458235267342108\\
47	0.000439662721428334\\
48	0.000419570564956493\\
49	0.000401313895770483\\
50	0.000383893163662847\\
51	0.000368802392132705\\
52	0.000356166748009116\\
53	0.000343919647368209\\
54	0.000330634542896946\\
55	0.000317375494202567\\
56	0.000305545542530896\\
57	0.000293139463156935\\
58	0.000278761342900579\\
59	0.000263812643761727\\
60	0.000247279190450497\\
61	0.000233775734112621\\
62	0.000220549397837644\\
63	0.00020840959245916\\
64	0.000197873830050784\\
65	0.000188617791100538\\
66	0.000179878285120293\\
67	0.000172507594931817\\
68	0.000166688300187972\\
69	0.000160840074816899\\
70	0.000154322004810118\\
71	0.000146564449070156\\
72	0.000139142331469799\\
73	0.000132946262516558\\
74	0.000128544591898058\\
75	0.000124654197591142\\
76	0.000120677139537221\\
77	0.00011718826040479\\
78	0.000113885429887769\\
79	0.000111223281325143\\
80	0.000107961077092455\\
81	0.000104362008307177\\
82	0.000101055090807758\\
83	9.71979574202781e-05\\
84	9.39156982318133e-05\\
85	9.05767952473625e-05\\
86	8.71203294198837e-05\\
87	8.30971668664495e-05\\
88	7.85223714352986e-05\\
89	7.48322371519197e-05\\
90	7.09603983027954e-05\\
91	6.7282447775403e-05\\
92	6.42146196808703e-05\\
93	6.14902984765503e-05\\
94	5.85730029620537e-05\\
95	5.5862071453895e-05\\
96	5.3686244470428e-05\\
97	5.16496308268724e-05\\
98	4.96585629315117e-05\\
99	4.76979659590891e-05\\
100	4.58254353600654e-05\\
101	4.3818716040418e-05\\
102	4.19118344140723e-05\\
103	4.00239269231413e-05\\
104	3.78652616084014e-05\\
105	3.59213295969807e-05\\
106	3.4108184943265e-05\\
107	3.26751396793086e-05\\
108	3.12225287809296e-05\\
109	2.9676040366258e-05\\
110	2.82339701443502e-05\\
111	2.69152028784765e-05\\
112	2.53791555580807e-05\\
113	2.3540584505527e-05\\
114	2.19253332762471e-05\\
115	2.04645914182571e-05\\
116	1.88077964238099e-05\\
117	1.72919475255986e-05\\
118	1.60350872995493e-05\\
119	1.47328411293946e-05\\
120	1.36675239073688e-05\\
121	1.27849999152401e-05\\
122	1.18410745951484e-05\\
123	1.09080423396801e-05\\
124	1.00419417272809e-05\\
125	9.33825506063706e-06\\
126	8.62832954413417e-06\\
127	7.99722269254892e-06\\
128	7.39529855206007e-06\\
129	6.8756393724439e-06\\
130	6.4270977974118e-06\\
131	5.98194086720525e-06\\
132	5.62378840133714e-06\\
133	5.29308580985141e-06\\
134	4.93788586859207e-06\\
135	4.61688468333784e-06\\
136	4.3946008223677e-06\\
137	4.17246137445542e-06\\
138	3.97103731396159e-06\\
139	3.7985758142101e-06\\
140	3.59402735176868e-06\\
141	3.37587454331535e-06\\
142	3.17762354063254e-06\\
143	3.00537700279495e-06\\
144	2.84071039509093e-06\\
145	2.69788697206778e-06\\
146	2.53950814471387e-06\\
147	2.3930743237615e-06\\
148	2.27380238154057e-06\\
149	2.14431171529868e-06\\
150	2.00307610007971e-06\\
151	1.87688942990071e-06\\
152	1.74478901215902e-06\\
153	1.61041751758181e-06\\
154	1.49184754731337e-06\\
155	1.37971457692962e-06\\
156	1.25820176008242e-06\\
157	1.16094316734182e-06\\
158	1.08882946686995e-06\\
159	1.02182387931616e-06\\
160	9.69295950504576e-07\\
161	9.20298471654457e-07\\
162	8.82775119381805e-07\\
163	8.39013989922196e-07\\
164	7.96631627570694e-07\\
165	7.39729615552595e-07\\
166	6.84631138467281e-07\\
167	6.39690709570124e-07\\
168	5.87158961045142e-07\\
169	5.35922412564122e-07\\
170	4.9876170834867e-07\\
171	4.67448956764445e-07\\
172	4.40380266583135e-07\\
173	4.10475762745289e-07\\
174	3.87095782888033e-07\\
175	3.66143684233072e-07\\
176	3.48735545982636e-07\\
177	3.28022449741563e-07\\
178	3.09195000846733e-07\\
179	2.94898447416746e-07\\
180	2.82909853251921e-07\\
181	2.69314995714798e-07\\
182	2.53763330599975e-07\\
183	2.3973800258662e-07\\
184	2.25319198899522e-07\\
185	2.12104301330649e-07\\
186	2.01725241906862e-07\\
187	1.92190040876729e-07\\
188	1.85707427279356e-07\\
189	1.79048315696936e-07\\
190	1.74368744876825e-07\\
191	1.69463839623414e-07\\
192	1.65510170495595e-07\\
193	1.61966531273947e-07\\
194	1.5915408362395e-07\\
195	1.56855888778633e-07\\
196	1.54141048035522e-07\\
197	1.51456190525647e-07\\
198	1.4824957915042e-07\\
199	1.46050566146759e-07\\
200	1.43654361116366e-07\\
201	1.41328203584126e-07\\
202	1.38151812090441e-07\\
203	9.47754854431273e-08\\
204	9.25925136073299e-08\\
205	9.25477659372598e-08\\
206	9.24995710159634e-08\\
207	9.24372011580099e-08\\
208	9.23551801123625e-08\\
209	9.22569835632237e-08\\
210	9.21481872594475e-08\\
211	9.20406791827015e-08\\
212	9.19384730469836e-08\\
213	9.18536570218349e-08\\
214	9.17739860634804e-08\\
215	9.16966414831702e-08\\
216	9.16155797784343e-08\\
217	9.15247702404215e-08\\
218	9.14266002327001e-08\\
219	9.13307898253474e-08\\
220	9.1225641102849e-08\\
221	9.11219661450637e-08\\
222	9.10275523543785e-08\\
223	9.09319390813903e-08\\
224	9.08313838288021e-08\\
225	9.07400549164453e-08\\
226	9.06581148988734e-08\\
227	9.0553369810818e-08\\
228	9.04667957099337e-08\\
229	9.04114547332647e-08\\
230	9.03556750270064e-08\\
231	9.02986499560361e-08\\
232	9.02436824878623e-08\\
233	9.01954807387383e-08\\
234	9.01584326033679e-08\\
235	9.01185671856781e-08\\
236	9.00656519485896e-08\\
237	9.00060660615395e-08\\
238	8.99586442054504e-08\\
239	8.99238749365542e-08\\
240	8.98962476772153e-08\\
241	8.98722993971855e-08\\
242	8.98528789827995e-08\\
243	8.9833413049329e-08\\
244	8.98187865622666e-08\\
245	8.98061728968924e-08\\
246	8.97947564408798e-08\\
247	8.97866712234178e-08\\
248	8.97802813789697e-08\\
249	8.97744191580656e-08\\
250	8.9769682250015e-08\\
251	8.97652670101301e-08\\
252	8.97584457228727e-08\\
253	8.97515779799854e-08\\
254	8.97456650659398e-08\\
255	8.97401603961487e-08\\
256	8.97351864499665e-08\\
257	8.97309381783463e-08\\
258	8.97274413997616e-08\\
259	8.97244756054931e-08\\
260	8.97213723239274e-08\\
261	8.97188982524157e-08\\
262	8.97162572752407e-08\\
263	8.97138861018231e-08\\
264	8.97112467638566e-08\\
265	8.97090818269658e-08\\
266	8.97070638024165e-08\\
267	8.97053384639332e-08\\
268	8.97039944768796e-08\\
269	8.97030431437628e-08\\
270	8.97022504810236e-08\\
271	8.97017861923081e-08\\
272	8.97014779542857e-08\\
273	8.97011937227993e-08\\
274	8.97009064743253e-08\\
275	8.97006565759882e-08\\
276	8.97004851981123e-08\\
277	8.97003083503526e-08\\
278	8.97001422240122e-08\\
279	8.970003810856e-08\\
280	8.9699938337542e-08\\
281	8.96998517153054e-08\\
282	8.9699798838975e-08\\
283	8.96997534973692e-08\\
284	8.9699710678025e-08\\
285	8.96996760843048e-08\\
286	8.96996518598772e-08\\
287	8.96996318892239e-08\\
288	8.96996117026528e-08\\
289	8.96995961275503e-08\\
290	8.96995772485831e-08\\
291	8.96995628741032e-08\\
292	8.96995476895227e-08\\
293	8.96995336198827e-08\\
294	8.96995218573471e-08\\
295	8.96995126611927e-08\\
296	8.96995045819167e-08\\
297	8.96994989838959e-08\\
298	8.96994944066845e-08\\
299	8.96994900726012e-08\\
300	8.96994877263944e-08\\
301	8.96994855268597e-08\\
302	8.96994834510166e-08\\
303	8.96994817637449e-08\\
304	8.9699480552355e-08\\
305	8.969947960877e-08\\
306	8.96994788823125e-08\\
307	8.96994783412425e-08\\
308	8.96994778102612e-08\\
309	8.96994775073428e-08\\
310	8.96994772510927e-08\\
311	8.96994770033677e-08\\
312	8.96994768229261e-08\\
313	8.96994766931347e-08\\
314	8.96994765751352e-08\\
315	8.96994764724997e-08\\
316	8.96994763869482e-08\\
317	8.96994763047239e-08\\
318	8.96994762149332e-08\\
319	8.96994761268481e-08\\
320	8.96994760697553e-08\\
321	8.96994760026931e-08\\
322	8.96994759382888e-08\\
323	8.96994758895813e-08\\
324	8.96994758487883e-08\\
325	8.96994758084186e-08\\
326	8.96994757747911e-08\\
327	8.96994757475216e-08\\
328	8.96994757187448e-08\\
329	8.96994756979655e-08\\
330	8.96994756817398e-08\\
331	8.96994756649242e-08\\
332	8.96994756501994e-08\\
333	8.96994756370942e-08\\
334	8.96994756293356e-08\\
335	8.96994756238669e-08\\
336	8.96994756196374e-08\\
337	8.96994756158212e-08\\
338	8.96994756133797e-08\\
339	8.96994756112108e-08\\
340	8.96994756092629e-08\\
341	8.96994756079691e-08\\
342	8.96994756072891e-08\\
343	8.96994756067811e-08\\
344	8.96994756064608e-08\\
345	8.96994756062695e-08\\
346	8.96994756061912e-08\\
347	8.9699475606143e-08\\
348	8.96994756061101e-08\\
349	8.96994756060839e-08\\
350	8.96994756060681e-08\\
351	8.96994756060573e-08\\
352	8.96994756060508e-08\\
353	8.96994756060471e-08\\
354	8.9699475606045e-08\\
355	8.96994756060436e-08\\
356	8.96994756060428e-08\\
357	8.96994756060422e-08\\
358	8.9699475606042e-08\\
359	8.96994756060418e-08\\
360	8.96994756060417e-08\\
361	8.96994756060416e-08\\
362	8.96994756060415e-08\\
363	8.96994756060415e-08\\
364	8.96994756060414e-08\\
365	8.96994756060414e-08\\
366	8.96994756060413e-08\\
367	8.96994756060412e-08\\
368	8.96994756060412e-08\\
369	8.96994756060412e-08\\
370	8.96994756060412e-08\\
371	8.96994756060411e-08\\
372	8.96994756060411e-08\\
373	8.96994756060411e-08\\
374	8.96994756060411e-08\\
375	8.96994756060411e-08\\
376	8.96994756060411e-08\\
377	8.96994756060411e-08\\
378	8.96994756060411e-08\\
379	8.96994756060411e-08\\
380	8.96994756060411e-08\\
381	8.96994756060411e-08\\
382	8.96994756060411e-08\\
383	8.96994756060411e-08\\
384	8.96994756060411e-08\\
385	8.9699475606041e-08\\
386	8.9699475606041e-08\\
387	8.9699475606041e-08\\
388	8.9699475606041e-08\\
389	8.96994756060408e-08\\
390	8.96994756060407e-08\\
391	8.96994756060406e-08\\
392	8.96994756060395e-08\\
393	8.96994756060348e-08\\
394	8.96994756060208e-08\\
395	8.96994756059779e-08\\
396	8.9699475605415e-08\\
397	8.96994756003236e-08\\
398	8.96994755641853e-08\\
399	8.96994750622862e-08\\
400	8.96989389450974e-08\\
401	8.96977388946383e-08\\
402	8.9697738886784e-08\\
403	8.96977388353903e-08\\
404	8.96789183029953e-08\\
405	8.96629339591504e-08\\
406	8.96183528848391e-08\\
407	8.92743608633068e-08\\
408	8.92626551826336e-08\\
409	8.85833603571725e-08\\
410	8.77307063963498e-08\\
411	8.73004004187541e-08\\
412	8.7192164491393e-08\\
413	8.71650348938956e-08\\
414	8.68631983358312e-08\\
415	8.67347476032123e-08\\
416	8.62498832127352e-08\\
417	8.44909745050243e-08\\
418	8.43997878808873e-08\\
419	8.42196965148942e-08\\
420	8.40945302348031e-08\\
421	8.40894470458075e-08\\
422	8.40870976266358e-08\\
423	8.36498568294806e-08\\
424	8.36433547308405e-08\\
425	8.34782057941045e-08\\
426	8.34769796470152e-08\\
427	8.33668724734173e-08\\
428	8.3251576596549e-08\\
429	8.29588277840748e-08\\
430	8.21908851170157e-08\\
431	8.21659973569796e-08\\
432	8.20408807181241e-08\\
433	8.04230330349953e-08\\
434	7.97220040328105e-08\\
435	7.93583364224785e-08\\
436	7.83524154638357e-08\\
437	7.71349970306269e-08\\
438	7.67411537706049e-08\\
439	7.62975749243406e-08\\
440	7.60254966543097e-08\\
441	7.59618823998389e-08\\
442	7.59370320980663e-08\\
443	7.56817609899482e-08\\
444	7.55718825648833e-08\\
445	7.5523288675925e-08\\
446	7.55221150805584e-08\\
447	7.35541581330744e-08\\
448	7.33292405257705e-08\\
449	7.33278014587025e-08\\
450	7.31490987074055e-08\\
};
\addlegendentry{$\varepsilon=10^{-14}$}
\legend{}
\end{axis}

\end{tikzpicture}%

%% file: Figs/tikz/GMRES_residual_NTAcc.tex
\definecolor{mycolor1}{rgb}{0.00000,0.44700,0.74100}%
\definecolor{mycolor2}{rgb}{0.85000,0.32500,0.09800}%
\definecolor{mycolor3}{rgb}{0.92900,0.69400,0.12500}%
\definecolor{mycolor4}{rgb}{0.49400,0.18400,0.55600}%
\begin{tikzpicture}

\begin{axis}[%
width=3.65cm,
height=3.7cm,
scale only axis,
tick label style={font=\footnotesize},
xmin=0,
xmax=370,
xlabel style={font=\footnotesize\color{white!15!black},yshift=.2cm},
xlabel={Iteration},
ymode=log,
ymin=1.1e-15,
ymax=1.1,
yminorticks=true,
ytickten = {-14, -11, -8, -5, -2},
ylabel style={font=\footnotesize\color{white!15!black},yshift=-.2cm},
ylabel={Residual norm},
axis background/.style={fill=white},
legend style={at={(axis cs:348,9e-1)}, anchor=north east, legend cell align=left, align=left, font = \scriptsize, draw=white!15!black}
]
\addplot [color=mycolor1, only marks, mark=x, mark options={solid, mycolor1}]
  table[row sep=crcr]{%
0	1\\
1	0.204873076072875\\
2	0.089663841992867\\
3	0.0463024790152402\\
4	0.0264444369400659\\
5	0.0175688800798046\\
6	0.0123059061169049\\
7	0.00880388652289965\\
8	0.00641885221533357\\
9	0.00496085452459367\\
10	0.00394571178603344\\
11	0.00328994118237524\\
12	0.00280591037970442\\
13	0.00250323977954526\\
14	0.00249164417546957\\
15	0.00234742213632333\\
16	0.00205268993960589\\
17	0.00176852910216773\\
18	0.00152181638307428\\
19	0.00132513904337307\\
20	0.00117674654736241\\
21	0.00103117345247605\\
22	0.000936940591246331\\
23	0.000922006906939215\\
24	0.000906112714244648\\
25	0.000836613315622882\\
26	0.000758286342980023\\
27	0.000677832640462804\\
28	0.000598801072363435\\
29	0.000535059355821351\\
30	0.000482663530073714\\
31	0.000448609460881055\\
32	0.000428200983425794\\
33	0.000428185265847336\\
34	0.000411416659089777\\
35	0.00038201095246543\\
36	0.000356975801254782\\
37	0.000337819006019057\\
38	0.000316607536341122\\
39	0.0002933668524078\\
40	0.000275301560678758\\
41	0.000257378684601973\\
42	0.000241454268812185\\
43	0.000225592206092689\\
44	0.000210943502755423\\
45	0.000198031056964654\\
46	0.000189083979913555\\
47	0.000188085871360903\\
48	0.000184579213858316\\
49	0.000174980154906968\\
50	0.000162845210851581\\
51	0.000150167507408532\\
52	0.000137356125866578\\
53	0.000126014981212848\\
54	0.000116773450295495\\
55	0.000109977415189548\\
56	0.000103043002966853\\
57	9.71556676792766e-05\\
58	9.20980018135593e-05\\
59	8.77434736338755e-05\\
60	8.35211427115721e-05\\
61	7.93292172273112e-05\\
62	7.46218036554663e-05\\
63	7.01607118079285e-05\\
64	6.60411550829421e-05\\
65	6.1927431900102e-05\\
66	5.82916058722293e-05\\
67	5.49605680375269e-05\\
68	5.14974937604214e-05\\
69	4.84655658372462e-05\\
70	4.58896887055266e-05\\
71	4.36644160523638e-05\\
72	4.16161345840121e-05\\
73	3.93012296629485e-05\\
74	3.62084246545712e-05\\
75	3.30290547179265e-05\\
76	3.03928443601671e-05\\
77	2.81444852724115e-05\\
78	2.59561957714694e-05\\
79	2.3942876133406e-05\\
80	2.19949398118819e-05\\
81	2.07545450661716e-05\\
82	1.95318675659842e-05\\
83	1.83246242145728e-05\\
84	1.73161044977849e-05\\
85	1.63071186405783e-05\\
86	1.54661465323688e-05\\
87	1.47013511749883e-05\\
88	1.41255877821137e-05\\
89	1.358468751092e-05\\
90	1.28998165867682e-05\\
91	1.23295384720041e-05\\
92	1.17963115659728e-05\\
93	1.12046886727758e-05\\
94	1.06482215612972e-05\\
95	1.01767138694265e-05\\
96	9.71684551851301e-06\\
97	9.25940305943014e-06\\
98	8.86250430428572e-06\\
99	8.54792716271384e-06\\
100	8.24510934188213e-06\\
101	7.99468120458535e-06\\
102	7.79121447382585e-06\\
103	7.58261374025624e-06\\
104	7.39904695415671e-06\\
105	7.22991899698292e-06\\
106	7.05618597842475e-06\\
107	6.85713867396531e-06\\
108	6.65101843920387e-06\\
109	6.46804768015697e-06\\
110	6.29196733477189e-06\\
111	6.11693224227491e-06\\
112	5.9484910527736e-06\\
113	5.80271568442552e-06\\
114	5.66013624764397e-06\\
115	5.47941706396393e-06\\
116	5.30289614201969e-06\\
117	5.14133796770422e-06\\
118	4.95876384320712e-06\\
119	4.75323080335275e-06\\
120	4.5749823181146e-06\\
121	4.38924003771629e-06\\
122	4.20390770467137e-06\\
123	4.04805851712927e-06\\
124	3.90489104089011e-06\\
125	3.75739376345258e-06\\
126	3.60132527029087e-06\\
127	3.45127831646641e-06\\
128	3.30875756722382e-06\\
129	3.17798074055742e-06\\
130	3.04263535666199e-06\\
131	2.91246433470985e-06\\
132	2.80199690462785e-06\\
133	2.68280927264608e-06\\
134	2.57991353355636e-06\\
135	2.48336870040715e-06\\
136	2.37503870378238e-06\\
137	2.25450213612128e-06\\
138	2.15904178047336e-06\\
139	2.07460038864824e-06\\
140	1.99070975712928e-06\\
141	1.92532214825555e-06\\
142	1.85375671670549e-06\\
143	1.7806173827498e-06\\
144	1.70644306006917e-06\\
145	1.63276460090691e-06\\
146	1.56562705445834e-06\\
147	1.50942608426608e-06\\
148	1.45000848891897e-06\\
149	1.38551775105313e-06\\
150	1.33604851431971e-06\\
151	1.2874002803512e-06\\
152	1.23278282063383e-06\\
153	1.17632257078058e-06\\
154	1.1174611418695e-06\\
155	1.05904763264844e-06\\
156	1.00374971156743e-06\\
157	9.468291405244e-07\\
158	8.88632218468085e-07\\
159	8.33438051714798e-07\\
160	7.87017673692929e-07\\
161	7.4778640455819e-07\\
162	7.16120809763435e-07\\
163	6.88451091346161e-07\\
164	6.69204343888255e-07\\
165	6.47787242079783e-07\\
166	6.30265961724982e-07\\
167	6.08189685198363e-07\\
168	5.83243020489562e-07\\
169	5.60695238491654e-07\\
170	5.36605539137967e-07\\
171	5.07991962647239e-07\\
172	4.80301339571148e-07\\
173	4.56814535293527e-07\\
174	4.36259979812489e-07\\
175	4.12234985625706e-07\\
176	3.8746236902851e-07\\
177	3.65938676533111e-07\\
178	3.47398071104609e-07\\
179	3.25216174746676e-07\\
180	3.04606329488208e-07\\
181	2.85452413461847e-07\\
182	2.71448904795932e-07\\
183	2.59234713919754e-07\\
184	2.45862321638897e-07\\
185	2.34109993214903e-07\\
186	2.20758067769945e-07\\
187	2.08781484922698e-07\\
188	1.97853857011595e-07\\
189	1.87401415716105e-07\\
190	1.79575533736566e-07\\
191	1.71989804361842e-07\\
192	1.65783714820375e-07\\
193	1.59257599599389e-07\\
194	1.53861348608029e-07\\
195	1.48807751053719e-07\\
196	1.4502918308458e-07\\
197	1.42032081387579e-07\\
198	1.39615412128236e-07\\
199	1.37505469207273e-07\\
200	1.34883597284949e-07\\
201	1.33035716230585e-07\\
202	1.31486355805577e-07\\
203	1.28725345449268e-07\\
204	1.22928730754556e-07\\
205	3.16248109012343e-08\\
206	2.16877688308452e-08\\
207	2.00878825809631e-08\\
208	1.93591411452608e-08\\
209	1.87996824869669e-08\\
210	1.82605283151215e-08\\
211	1.76644367374745e-08\\
212	1.70395665347557e-08\\
213	1.63379756030525e-08\\
214	1.56071778424964e-08\\
215	1.49375124016342e-08\\
216	1.43401072946723e-08\\
217	1.37736995757149e-08\\
218	1.33366383714798e-08\\
219	1.28338657286752e-08\\
220	1.23183347511254e-08\\
221	1.18503468990211e-08\\
222	1.13656128613939e-08\\
223	1.08362094197247e-08\\
224	1.03464436922952e-08\\
225	9.80993229114711e-09\\
226	9.26667482517716e-09\\
227	8.75709779655533e-09\\
228	8.27611585663277e-09\\
229	7.67406755441442e-09\\
230	6.97467756590485e-09\\
231	6.4768319474835e-09\\
232	6.05290143163395e-09\\
233	5.63750197961191e-09\\
234	5.24175400960185e-09\\
235	4.89387028434504e-09\\
236	4.63328846644792e-09\\
237	4.45363375772427e-09\\
238	4.26104191527424e-09\\
239	4.03903357858867e-09\\
240	3.79049035550226e-09\\
241	3.57374207465436e-09\\
242	3.39672820494029e-09\\
243	3.21884763812597e-09\\
244	3.05634766940262e-09\\
245	2.87469930588435e-09\\
246	2.70517322350079e-09\\
247	2.54714620062436e-09\\
248	2.34650381719955e-09\\
249	2.20033519886503e-09\\
250	2.06644591449397e-09\\
251	1.9552604258954e-09\\
252	1.86859955716331e-09\\
253	1.8118441529167e-09\\
254	1.73984379612388e-09\\
255	1.65272677502259e-09\\
256	1.56472648376334e-09\\
257	1.48231560972794e-09\\
258	1.39368458649831e-09\\
259	1.3165242481338e-09\\
260	1.23643139477143e-09\\
261	1.18035772545718e-09\\
262	1.12136152170222e-09\\
263	1.07461648755287e-09\\
264	1.02425112867385e-09\\
265	9.78645719702967e-10\\
266	9.30572359123639e-10\\
267	8.86148112726646e-10\\
268	8.39290266543233e-10\\
269	7.83523913826782e-10\\
270	7.28066655933824e-10\\
271	6.80024271944824e-10\\
272	6.23598927273233e-10\\
273	5.72168444485498e-10\\
274	5.43840554335443e-10\\
275	5.20082669279937e-10\\
276	4.92113801990836e-10\\
277	4.60256002515735e-10\\
278	4.3764784834729e-10\\
279	4.09238323636783e-10\\
280	3.73762325244279e-10\\
281	3.46429510309756e-10\\
282	3.12335063684311e-10\\
283	2.7863470484518e-10\\
284	2.52285351301263e-10\\
285	2.28538827990074e-10\\
286	2.05994650719112e-10\\
287	1.84391876698092e-10\\
288	1.68662269709607e-10\\
289	1.53084522026728e-10\\
290	1.39620925155684e-10\\
291	1.28504017307803e-10\\
292	1.14583299312241e-10\\
293	1.02659038999379e-10\\
294	9.1669016503719e-11\\
295	7.92962346377495e-11\\
296	6.8753971730959e-11\\
297	5.87245095264252e-11\\
298	5.01508977775509e-11\\
299	4.43603637191559e-11\\
300	3.95857890902672e-11\\
301	3.4891820873926e-11\\
302	3.21385361141965e-11\\
303	3.00387044691086e-11\\
304	2.81514538756236e-11\\
305	2.60232739133592e-11\\
306	2.41282496710126e-11\\
307	2.22598490739826e-11\\
308	2.06586151486065e-11\\
309	1.90894971085619e-11\\
310	1.74257139753483e-11\\
311	1.5971468907246e-11\\
312	1.49379981754492e-11\\
313	1.36652077623906e-11\\
314	1.23802190660871e-11\\
315	1.13772423839126e-11\\
316	1.04441586846634e-11\\
317	9.52000145572905e-12\\
318	8.76104309281868e-12\\
319	8.14547226312639e-12\\
320	7.50282774534581e-12\\
321	6.81353534417184e-12\\
322	6.31665639538531e-12\\
323	5.79937117330571e-12\\
324	5.31683101917611e-12\\
325	4.84874637703461e-12\\
326	4.46993035190409e-12\\
327	4.16911852586223e-12\\
328	3.90241478879282e-12\\
329	3.69831501500569e-12\\
330	3.46245225728695e-12\\
331	3.25594290037033e-12\\
332	3.08716676482505e-12\\
333	2.88880309517841e-12\\
334	2.67811249083305e-12\\
335	2.41811625018654e-12\\
336	2.17254195931661e-12\\
337	1.9673546093255e-12\\
338	1.78031660239407e-12\\
339	1.61921531960316e-12\\
340	1.48295242943012e-12\\
341	1.35178765193246e-12\\
342	1.19791072933208e-12\\
343	1.02158651073377e-12\\
344	8.87393280749196e-13\\
345	7.68725501040973e-13\\
346	6.65442207323173e-13\\
347	5.7629319534705e-13\\
348	5.30296568650072e-13\\
349	4.94023712059391e-13\\
350	4.61750195964158e-13\\
351	4.20137347114532e-13\\
352	3.74259955124984e-13\\
353	3.28405893406634e-13\\
354	2.83269272397732e-13\\
355	2.41733674001192e-13\\
356	1.99033285141687e-13\\
357	1.64689804338387e-13\\
358	1.30813342798107e-13\\
359	1.09480329793809e-13\\
360	8.95235875657235e-14\\
361	6.93899053045033e-14\\
362	5.21992954251789e-14\\
363	4.3180212103199e-14\\
364	3.72251713474438e-14\\
365	2.8527728086204e-14\\
366	2.21921399913753e-14\\
367	1.60759732420342e-14\\
368	1.17550008453403e-14\\
369	8.63189181516715e-15\\
};
\addlegendentry{$\varepsilon=10^{-8}$}

\addplot [color=mycolor2, only marks, mark=x, mark options={solid, mycolor2}]
  table[row sep=crcr]{%
0	1\\
1	0.204873076853183\\
2	0.089663843688434\\
3	0.0463024822168889\\
4	0.0264444413831622\\
5	0.017568882472138\\
6	0.0123059067086805\\
7	0.00880388829777766\\
8	0.00641886473146619\\
9	0.00496092471172096\\
10	0.00394608664268371\\
11	0.00329175533527272\\
12	0.002815558750453\\
13	0.00254465297944286\\
14	0.00254302476983752\\
15	0.00236050427606625\\
16	0.00205443452108744\\
17	0.00176873500569283\\
18	0.00152183332025855\\
19	0.00132515152462149\\
20	0.00117685374133905\\
21	0.0010319195434168\\
22	0.000940574642169501\\
23	0.00092908266707697\\
24	0.000909038793635485\\
25	0.000837059173941081\\
26	0.000758338389484468\\
27	0.0006778236742757\\
28	0.00059878359337731\\
29	0.000535046163565212\\
30	0.0004826851881711\\
31	0.000448805191895324\\
32	0.000429057547146357\\
33	0.000429057178305145\\
34	0.000411604639897006\\
35	0.000382013706484234\\
36	0.000356954851144028\\
37	0.000337795493375212\\
38	0.000316583863183072\\
39	0.000293343541093821\\
40	0.000275278687962496\\
41	0.000257356461127915\\
42	0.000241433015966685\\
43	0.00022557400433411\\
44	0.000210940318543205\\
45	0.000198116253716646\\
46	0.000189558156893963\\
47	0.00018888525275229\\
48	0.000184826320407112\\
49	0.000175011145550657\\
50	0.000162835225859255\\
51	0.000150151992361552\\
52	0.000137340447278494\\
53	0.000125999757382993\\
54	0.000116758619102872\\
55	0.000109962866898296\\
56	0.000103028756419794\\
57	9.71416796948132e-05\\
58	9.20842616517624e-05\\
59	8.77299936686441e-05\\
60	8.35079642829462e-05\\
61	7.93163835481008e-05\\
62	7.46093644813302e-05\\
63	7.01486370903115e-05\\
64	6.60294040525017e-05\\
65	6.19159889726196e-05\\
66	5.82804225261541e-05\\
67	5.49496249671105e-05\\
68	5.14868043380889e-05\\
69	4.84550952655246e-05\\
70	4.58794215995624e-05\\
71	4.36543551512558e-05\\
72	4.16062994038043e-05\\
73	3.92916695495691e-05\\
74	3.61992114568546e-05\\
75	3.30201266916173e-05\\
76	3.03840946681046e-05\\
77	2.81358565375891e-05\\
78	2.59476522449243e-05\\
79	2.39343791879531e-05\\
80	2.19864474601668e-05\\
81	2.07460447486696e-05\\
82	1.95233702113871e-05\\
83	1.83161290547922e-05\\
84	1.73076174241101e-05\\
85	1.62986417829292e-05\\
86	1.5457683034091e-05\\
87	1.46929135714635e-05\\
88	1.41171878372759e-05\\
89	1.35763514997063e-05\\
90	1.28915827657883e-05\\
91	1.2321403213266e-05\\
92	1.1788283722953e-05\\
93	1.11967900011317e-05\\
94	1.06404465493229e-05\\
95	1.01690516516169e-05\\
96	9.70930357086231e-06\\
97	9.25198525491612e-06\\
98	8.85520091738741e-06\\
99	8.54072568807145e-06\\
100	8.2380248208387e-06\\
101	7.98771029638312e-06\\
102	7.78435621932792e-06\\
103	7.57589509557401e-06\\
104	7.39247154143617e-06\\
105	7.2234999283089e-06\\
106	7.04994794715779e-06\\
107	6.85112862001357e-06\\
108	6.64525644936847e-06\\
109	6.46251724223939e-06\\
110	6.28667373806503e-06\\
111	6.11188516525839e-06\\
112	5.94368839345312e-06\\
113	5.79813235026598e-06\\
114	5.65577705220754e-06\\
115	5.47534576480627e-06\\
116	5.29910126211437e-06\\
117	5.13779276707267e-06\\
118	4.95549363846905e-06\\
119	4.75025352687149e-06\\
120	4.57224405651582e-06\\
121	4.38673469756767e-06\\
122	4.20161638120283e-06\\
123	4.04593456558604e-06\\
124	3.90291175488138e-06\\
125	3.75555422573875e-06\\
126	3.5996213227636e-06\\
127	3.44969321309386e-06\\
128	3.30727482722285e-06\\
129	3.1765835240374e-06\\
130	3.04131824339767e-06\\
131	2.91121596594915e-06\\
132	2.80080188501218e-06\\
133	2.68166692601653e-06\\
134	2.57881272013618e-06\\
135	2.48230540245951e-06\\
136	2.37401569764075e-06\\
137	2.25351870322687e-06\\
138	2.15808673404565e-06\\
139	2.07366975841304e-06\\
140	1.98980293054902e-06\\
141	1.9244348878388e-06\\
142	1.85289323023147e-06\\
143	1.77977920188122e-06\\
144	1.70563075210132e-06\\
145	1.63197832720392e-06\\
146	1.56486454966046e-06\\
147	1.5086847606242e-06\\
148	1.44929110768522e-06\\
149	1.38482676933508e-06\\
150	1.33537868835585e-06\\
151	1.28675323237398e-06\\
152	1.23216225730526e-06\\
153	1.17572962690927e-06\\
154	1.11689667879243e-06\\
155	1.05850994087377e-06\\
156	1.00323540710361e-06\\
157	9.46337046690787e-07\\
158	8.88159768814891e-07\\
159	8.32980663597998e-07\\
160	7.86570938196939e-07\\
161	7.47347736768766e-07\\
162	7.15689000921522e-07\\
163	6.88025901569022e-07\\
164	6.68785084124618e-07\\
165	6.47376457807521e-07\\
166	6.29863842219966e-07\\
167	6.07800218996451e-07\\
168	5.8286828475466e-07\\
169	5.60334210377067e-07\\
170	5.36259303685732e-07\\
171	5.07662176875501e-07\\
172	4.7998586162546e-07\\
173	4.56509973956196e-07\\
174	4.35964492586933e-07\\
175	4.11949005143683e-07\\
176	3.87184441637813e-07\\
177	3.65666060959885e-07\\
178	3.47129052425696e-07\\
179	3.24949875917986e-07\\
180	3.04340330367012e-07\\
181	2.85184684929471e-07\\
182	2.7117895699341e-07\\
183	2.58962499332675e-07\\
184	2.45587068460306e-07\\
185	2.33831493470526e-07\\
186	2.20475220239886e-07\\
187	2.08493882497135e-07\\
188	1.97561438137282e-07\\
189	1.87103875931034e-07\\
190	1.79274269666235e-07\\
191	1.71685151735224e-07\\
192	1.65476825085895e-07\\
193	1.58948996551439e-07\\
194	1.5355196470613e-07\\
195	1.48498430738842e-07\\
196	1.44720643045793e-07\\
197	1.41725346691164e-07\\
198	1.39311223175623e-07\\
199	1.3720482330832e-07\\
200	1.34588631348974e-07\\
201	1.32745822884561e-07\\
202	1.3120165862673e-07\\
203	1.28452384951819e-07\\
204	1.22690037040618e-07\\
205	3.16116607939486e-08\\
206	2.16737351915869e-08\\
207	2.00732395248787e-08\\
208	1.93442387417565e-08\\
209	1.87846557716963e-08\\
210	1.82454687761612e-08\\
211	1.76494341738085e-08\\
212	1.70247029565502e-08\\
213	1.63233475891244e-08\\
214	1.55928564149562e-08\\
215	1.49235311620689e-08\\
216	1.43264904672091e-08\\
217	1.37604975540592e-08\\
218	1.33238119794601e-08\\
219	1.28215444352064e-08\\
220	1.23065830514782e-08\\
221	1.18391566795929e-08\\
222	1.13550431836575e-08\\
223	1.08263415467799e-08\\
224	1.0337237584411e-08\\
225	9.80146282371683e-09\\
226	9.25894302882221e-09\\
227	8.75003763674899e-09\\
228	8.26967944184425e-09\\
229	7.66837320796277e-09\\
230	6.9697543610617e-09\\
231	6.47240015996619e-09\\
232	6.04885767094906e-09\\
233	5.63381146290017e-09\\
234	5.23837216190697e-09\\
235	4.89073505593583e-09\\
236	4.63032670346722e-09\\
237	4.45079167550764e-09\\
238	4.25833206662583e-09\\
239	4.0364727807463e-09\\
240	3.78808606476004e-09\\
241	3.57146489111926e-09\\
242	3.39455126873971e-09\\
243	3.21676978188546e-09\\
244	3.05436066358701e-09\\
245	2.87281199652249e-09\\
246	2.70337585804447e-09\\
247	2.54543082736219e-09\\
248	2.34488685257442e-09\\
249	2.19878461469162e-09\\
250	2.06495599476008e-09\\
251	1.95382203899545e-09\\
252	1.86720442551883e-09\\
253	1.81048355599869e-09\\
254	1.73853542846013e-09\\
255	1.65148684628267e-09\\
256	1.56355905883503e-09\\
257	1.48121864594506e-09\\
258	1.39266571824189e-09\\
259	1.31557430772577e-09\\
260	1.23555426560505e-09\\
261	1.1795334112115e-09\\
262	1.12059484159682e-09\\
263	1.0738973802109e-09\\
264	1.02358521599526e-09\\
265	9.78029790099506e-10\\
266	9.30009867909065e-10\\
267	8.85635112646678e-10\\
268	8.38829013268236e-10\\
269	7.83122257576198e-10\\
270	7.27720694682845e-10\\
271	6.79723695222074e-10\\
272	6.23347649242779e-10\\
273	5.71958107757183e-10\\
274	5.43651537382806e-10\\
275	5.19911372951701e-10\\
276	4.91962899105378e-10\\
277	4.60127029554132e-10\\
278	4.37533602739075e-10\\
279	4.09141752729564e-10\\
280	3.73685806167864e-10\\
281	3.46366972130651e-10\\
282	3.12288199250218e-10\\
283	2.78601306325677e-10\\
284	2.52261067986768e-10\\
285	2.28521781869061e-10\\
286	2.05983579636511e-10\\
287	1.84385758142864e-10\\
288	1.68659271259278e-10\\
289	1.53084198936044e-10\\
290	1.3962260846228e-10\\
291	1.28507088972177e-10\\
292	1.14587814449862e-10\\
293	1.02664570371545e-10\\
294	9.1675291865847e-11\\
295	7.93031996311092e-11\\
296	6.87614626373481e-11\\
297	5.87324470190547e-11\\
298	5.01592288668107e-11\\
299	4.43689444614268e-11\\
300	3.95944533184216e-11\\
301	3.49005099275097e-11\\
302	3.21471471688335e-11\\
303	3.00470956198225e-11\\
304	2.81595007766674e-11\\
305	2.60308137639482e-11\\
306	2.41352676669406e-11\\
307	2.22662927941929e-11\\
308	2.06645447418693e-11\\
309	1.90948713074531e-11\\
310	1.74305073249279e-11\\
311	1.59757635956846e-11\\
312	1.4941928913001e-11\\
313	1.36686878734536e-11\\
314	1.23832866710149e-11\\
315	1.13800171961167e-11\\
316	1.04466815890631e-11\\
317	9.52229752174313e-12\\
318	8.76316825100904e-12\\
319	8.1474618675328e-12\\
320	7.50467497344678e-12\\
321	6.8152484435346e-12\\
322	6.3182805018865e-12\\
323	5.80090341876774e-12\\
324	5.31828674468161e-12\\
325	4.85013505583168e-12\\
326	4.47126380452388e-12\\
327	4.17040169872453e-12\\
328	3.90364456570489e-12\\
329	3.69948978809458e-12\\
330	3.46355082128335e-12\\
331	3.25696725755279e-12\\
332	3.08812018245304e-12\\
333	2.88966374532755e-12\\
334	2.67887188848689e-12\\
335	2.41875918403703e-12\\
336	2.17308835989012e-12\\
337	1.96782802593438e-12\\
338	1.78072528275228e-12\\
339	1.6195698716795e-12\\
340	1.48326102560472e-12\\
341	1.35205058808139e-12\\
342	1.19812307314535e-12\\
343	1.02175136837631e-12\\
344	8.87529393479382e-13\\
345	7.68840767342135e-13\\
346	6.6554274203871e-13\\
347	5.76383426105005e-13\\
348	5.30381236431005e-13\\
349	4.94102137097174e-13\\
350	4.61821070560624e-13\\
351	4.20197136876694e-13\\
352	3.7430767096873e-13\\
353	3.28442529511697e-13\\
354	2.83296458944994e-13\\
355	2.4175342703588e-13\\
356	1.99047015481757e-13\\
357	1.64699842531902e-13\\
358	1.30820420207796e-13\\
359	1.09485938401079e-13\\
360	8.95278180295623e-14\\
361	6.93930330129158e-14\\
362	5.22016453609069e-14\\
363	4.31821675146889e-14\\
364	3.72266880219581e-14\\
365	2.85286843680526e-14\\
366	2.21928424993123e-14\\
367	1.60766045360433e-14\\
368	1.17557799561885e-14\\
369	8.63299777949471e-15\\
};
\addlegendentry{$\varepsilon=10^{-10}$}

\addplot [color=mycolor3, only marks, mark=x, mark options={solid, mycolor3}]
  table[row sep=crcr]{%
0	1\\
1	0.204872490053279\\
2	0.0896637584648851\\
3	0.0463025178341983\\
4	0.0264447809626908\\
5	0.0175690789050974\\
6	0.0123059082360428\\
7	0.008803807535579\\
8	0.00641867836371924\\
9	0.00496058987040449\\
10	0.0039455877269407\\
11	0.00329113579940436\\
12	0.0028150674297701\\
13	0.00254521329204753\\
14	0.00254376930862128\\
15	0.0023598618833135\\
16	0.00205332838545437\\
17	0.00176749772760799\\
18	0.00152052513212604\\
19	0.0013237928897614\\
20	0.00117545350925901\\
21	0.00103044044940428\\
22	0.000938952651127277\\
23	0.000927363577381165\\
24	0.00090741461117064\\
25	0.000835425689602765\\
26	0.0007566281864556\\
27	0.000676018752895694\\
28	0.00059689381823772\\
29	0.00053319091382374\\
30	0.000481482491411324\\
31	0.000450953575159185\\
32	0.000440961716623942\\
33	0.000437076650170842\\
34	0.000411454734342885\\
35	0.000380023183148613\\
36	0.000354720513863212\\
37	0.000335536571189085\\
38	0.000314344184802905\\
39	0.000291145550303479\\
40	0.000273124962508136\\
41	0.000255262959921172\\
42	0.000239403803570657\\
43	0.000223617541100673\\
44	0.000209053144233582\\
45	0.000196275559745693\\
46	0.000187665743813867\\
47	0.000186919851993501\\
48	0.000183025593130865\\
49	0.00017332455075279\\
50	0.000161237470824693\\
51	0.000148636430875758\\
52	0.000135900324745422\\
53	0.000124617573915735\\
54	0.000115417558140095\\
55	0.000108650272976305\\
56	0.000101745730318766\\
57	9.58835225547713e-05\\
58	9.08495017996924e-05\\
59	8.65193721932781e-05\\
60	8.23250057259927e-05\\
61	7.816480088709e-05\\
62	7.34935295823587e-05\\
63	6.90656736646178e-05\\
64	6.49755431446319e-05\\
65	6.08898371179299e-05\\
66	5.72777116228924e-05\\
67	5.39688057915167e-05\\
68	5.0529394603784e-05\\
69	4.7518186821315e-05\\
70	4.49617997729632e-05\\
71	4.27564006695021e-05\\
72	4.0729890298784e-05\\
73	3.84416591355332e-05\\
74	3.53829471372345e-05\\
75	3.22324427965845e-05\\
76	2.96151564258442e-05\\
77	2.73802899332425e-05\\
78	2.52023030774399e-05\\
79	2.31956570539998e-05\\
80	2.12505406636708e-05\\
81	2.0010941516173e-05\\
82	1.878994651015e-05\\
83	1.75842011783097e-05\\
84	1.65773848491238e-05\\
85	1.55701791808978e-05\\
86	1.47310080049547e-05\\
87	1.39689358314211e-05\\
88	1.33967284765883e-05\\
89	1.28615997775662e-05\\
90	1.21858030885628e-05\\
91	1.16241621530594e-05\\
92	1.1100251849715e-05\\
93	1.05197272121485e-05\\
94	9.97376106743396e-06\\
95	9.51173450288526e-06\\
96	9.06189583173121e-06\\
97	8.61471478722152e-06\\
98	8.22721197383227e-06\\
99	7.92099625384344e-06\\
100	7.62780445019808e-06\\
101	7.38674872672779e-06\\
102	7.1926186485849e-06\\
103	6.99564078735397e-06\\
104	6.82403660801607e-06\\
105	6.66800670059412e-06\\
106	6.5094785957878e-06\\
107	6.3296346057938e-06\\
108	6.14448507562423e-06\\
109	5.98116191968931e-06\\
110	5.82524725326081e-06\\
111	5.6712815655947e-06\\
112	5.52381786953374e-06\\
113	5.39693941203473e-06\\
114	5.27373398449822e-06\\
115	5.11802269088975e-06\\
116	4.96567244173315e-06\\
117	4.82607727950342e-06\\
118	4.6678499019437e-06\\
119	4.48846150826808e-06\\
120	4.33170293364142e-06\\
121	4.16703860919222e-06\\
122	4.00119635944851e-06\\
123	3.86064289502479e-06\\
124	3.73070585262458e-06\\
125	3.59597909755621e-06\\
126	3.45228617083525e-06\\
127	3.31304098898925e-06\\
128	3.17975871327572e-06\\
129	3.05661675781632e-06\\
130	2.92831585463195e-06\\
131	2.80407015779236e-06\\
132	2.6980854580004e-06\\
133	2.58319095136144e-06\\
134	2.48356174317972e-06\\
135	2.38987514036849e-06\\
136	2.28451756282346e-06\\
137	2.16674709461786e-06\\
138	2.07316146787771e-06\\
139	1.9902833948254e-06\\
140	1.90788317720324e-06\\
141	1.84373937601782e-06\\
142	1.77375026969605e-06\\
143	1.70232196313442e-06\\
144	1.62991672496886e-06\\
145	1.5580498331801e-06\\
146	1.49259261841934e-06\\
147	1.43794312903139e-06\\
148	1.38033909797083e-06\\
149	1.31789090821401e-06\\
150	1.27010424014686e-06\\
151	1.22333017879745e-06\\
152	1.17094340533013e-06\\
153	1.11685760336612e-06\\
154	1.06049310431403e-06\\
155	1.00446587101402e-06\\
156	9.51285412309593e-07\\
157	8.96412685214522e-07\\
158	8.40058459118501e-07\\
159	7.86301537201062e-07\\
160	7.40919422746349e-07\\
161	7.02490754811468e-07\\
162	6.71519161562403e-07\\
163	6.44527304149728e-07\\
164	6.25890329371426e-07\\
165	6.05345319524497e-07\\
166	5.88715298458637e-07\\
167	5.67942402975487e-07\\
168	5.44524340556985e-07\\
169	5.23409608723501e-07\\
170	5.00879988257761e-07\\
171	4.74027471463447e-07\\
172	4.47899892442831e-07\\
173	4.2563104514298e-07\\
174	4.06103979543628e-07\\
175	3.83179927102649e-07\\
176	3.59376379915939e-07\\
177	3.38529804852854e-07\\
178	3.20474324774114e-07\\
179	2.98710488572434e-07\\
180	2.78256581086041e-07\\
181	2.59030391431666e-07\\
182	2.44863173964723e-07\\
183	2.32458325770265e-07\\
184	2.18802646551345e-07\\
185	2.06723915799403e-07\\
186	1.92907479740072e-07\\
187	1.80397679877098e-07\\
188	1.68902663493957e-07\\
189	1.5781899191778e-07\\
190	1.49503347005863e-07\\
191	1.41438908848126e-07\\
192	1.34874664080858e-07\\
193	1.280150270959e-07\\
194	1.2238976123682e-07\\
195	1.17187395393291e-07\\
196	1.13362846755554e-07\\
197	1.1044488620669e-07\\
198	1.08198778863779e-07\\
199	1.06366149020627e-07\\
200	1.04210930384134e-07\\
201	1.02792430372767e-07\\
202	1.01689281524998e-07\\
203	9.99427490687831e-08\\
204	9.70911001679267e-08\\
205	3.0141035645963e-08\\
206	2.02648133846874e-08\\
207	1.86182514538085e-08\\
208	1.78687195042254e-08\\
209	1.72999745669971e-08\\
210	1.675981705809e-08\\
211	1.61710876008678e-08\\
212	1.55610058234083e-08\\
213	1.48830383461088e-08\\
214	1.41821098485142e-08\\
215	1.35449772018401e-08\\
216	1.2982045381456e-08\\
217	1.24546326224533e-08\\
218	1.2052770303807e-08\\
219	1.15972723838333e-08\\
220	1.11350704271507e-08\\
221	1.07197286797015e-08\\
222	1.02933244853834e-08\\
223	9.83012781345029e-09\\
224	9.40307976312328e-09\\
225	8.93678452162104e-09\\
226	8.46435382210847e-09\\
227	8.01978088692023e-09\\
228	7.59966203841699e-09\\
229	7.07074406776242e-09\\
230	6.44831202937358e-09\\
231	6.00009136599944e-09\\
232	5.61573165536182e-09\\
233	5.23672880370278e-09\\
234	4.87316901013542e-09\\
235	4.5513200701777e-09\\
236	4.30929087678975e-09\\
237	4.14259427084897e-09\\
238	3.96448439855741e-09\\
239	3.75901232843179e-09\\
240	3.52810125969051e-09\\
241	3.3258907798275e-09\\
242	3.16046528368727e-09\\
243	2.99411833271918e-09\\
244	2.84223239483481e-09\\
245	2.67228387300696e-09\\
246	2.51335874247578e-09\\
247	2.36500063096852e-09\\
248	2.17599155529475e-09\\
249	2.03769346918438e-09\\
250	1.91093098224975e-09\\
251	1.80571356259445e-09\\
252	1.72396145310414e-09\\
253	1.6710016497996e-09\\
254	1.60460042489876e-09\\
255	1.52471910651137e-09\\
256	1.44430001623288e-09\\
257	1.36918316987453e-09\\
258	1.28856397472425e-09\\
259	1.2184228327039e-09\\
260	1.145694721907e-09\\
261	1.09493176714738e-09\\
262	1.04170816741277e-09\\
263	9.99713564940887e-10\\
264	9.54653885168651e-10\\
265	9.1403226996704e-10\\
266	8.71295548570192e-10\\
267	8.31828109813262e-10\\
268	7.90174863445197e-10\\
269	7.40435476111968e-10\\
270	6.90651496100739e-10\\
271	6.47265565744698e-10\\
272	5.95937872173613e-10\\
273	5.48780837373949e-10\\
274	5.22699994710593e-10\\
275	5.00825920991764e-10\\
276	4.75048107336888e-10\\
277	4.45570963986899e-10\\
278	4.24578407169063e-10\\
279	3.98126152205862e-10\\
280	3.64890609729246e-10\\
281	3.39129287541805e-10\\
282	3.06808253980453e-10\\
283	2.74639185918423e-10\\
284	2.49331650047917e-10\\
285	2.26415702681628e-10\\
286	2.04560368967362e-10\\
287	1.83528730102303e-10\\
288	1.68160113646945e-10\\
289	1.52893270287689e-10\\
290	1.39664007382411e-10\\
291	1.28711487592622e-10\\
292	1.14964208117353e-10\\
293	1.03163191804407e-10\\
294	9.22649274229728e-11\\
295	7.99777221976641e-11\\
296	6.95002422795483e-11\\
297	5.9525148439899e-11\\
298	5.09984288416941e-11\\
299	4.52373468658698e-11\\
300	4.04745290623814e-11\\
301	3.57864714786528e-11\\
302	3.30274570608955e-11\\
303	3.09071943814244e-11\\
304	2.89867838744936e-11\\
305	2.68090404107457e-11\\
306	2.48626623784662e-11\\
307	2.2937326445581e-11\\
308	2.12847821624339e-11\\
309	1.96597724955347e-11\\
310	1.79373351644474e-11\\
311	1.64325221913181e-11\\
312	1.5361976363749e-11\\
313	1.40431684162475e-11\\
314	1.27160238201043e-11\\
315	1.16829899907804e-11\\
316	1.07238842938042e-11\\
317	9.77599166041363e-12\\
318	8.99881264537685e-12\\
319	8.36851537048993e-12\\
320	7.71010708038431e-12\\
321	7.00563416970607e-12\\
322	6.4984321202626e-12\\
323	5.97027710716079e-12\\
324	5.47845546727974e-12\\
325	5.00204981697799e-12\\
326	4.61638027323243e-12\\
327	4.3094817467132e-12\\
328	4.03647563407454e-12\\
329	3.82608194172498e-12\\
330	3.5816337115913e-12\\
331	3.36682618619569e-12\\
332	3.19015669981115e-12\\
333	2.98150109103178e-12\\
334	2.75956769485706e-12\\
335	2.48659015592277e-12\\
336	2.23019801477104e-12\\
337	2.01681602699649e-12\\
338	1.82253903780154e-12\\
339	1.65545222313225e-12\\
340	1.51419407634191e-12\\
341	1.37816523039532e-12\\
342	1.21897469962517e-12\\
343	1.03771970100642e-12\\
344	9.00566426310079e-13\\
345	7.79760916336392e-13\\
346	6.74949454185528e-13\\
347	5.8471272387402e-13\\
348	5.38134266667435e-13\\
349	5.01233567592246e-13\\
350	4.68223891756555e-13\\
351	4.25555055483586e-13\\
352	3.78555601503815e-13\\
353	3.31695242574981e-13\\
354	2.85716943871085e-13\\
355	2.43529386404094e-13\\
356	2.00310682707003e-13\\
357	1.65655045124541e-13\\
358	1.31521806890919e-13\\
359	1.10060313530239e-13\\
360	8.99743867300859e-14\\
361	6.97299935029388e-14\\
362	5.24563753983256e-14\\
363	4.33918003171587e-14\\
364	3.73832222996927e-14\\
365	2.86095334083895e-14\\
366	2.222989633308e-14\\
367	1.60807884648118e-14\\
368	1.17445123479179e-14\\
369	8.61442034537403e-15\\
};
\addlegendentry{$\varepsilon=10^{-12}$}

\addplot [color=mycolor4, only marks, mark=x, mark options={solid, mycolor4}]
  table[row sep=crcr]{%
0	1\\
1	0.204894082544804\\
2	0.0896841464772763\\
3	0.0463639375554963\\
4	0.0265352783931417\\
5	0.0176529161830302\\
6	0.0123797881306357\\
7	0.00888540828155691\\
8	0.00651802709364896\\
9	0.00510096679384908\\
10	0.00425482203303825\\
11	0.0040760336580915\\
12	0.00397853819088769\\
13	0.00350080502939888\\
14	0.0029791168768418\\
15	0.00253057427214152\\
16	0.00219456399189656\\
17	0.00193259501075645\\
18	0.00176117845298482\\
19	0.00173649280033894\\
20	0.00169033563574947\\
21	0.00155983681613078\\
22	0.00141377586022867\\
23	0.0012825071090863\\
24	0.00118358086521409\\
25	0.00110140117336952\\
26	0.00102949383996783\\
27	0.000959045523392433\\
28	0.000895043447433894\\
29	0.000859197785708192\\
30	0.000857016340035311\\
31	0.000837383818916395\\
32	0.000803580776338997\\
33	0.000769667709223703\\
34	0.000738240983077079\\
35	0.00070922776365607\\
36	0.000685876748816096\\
37	0.000666767692920695\\
38	0.000643995157450699\\
39	0.000618267551925323\\
40	0.000597655307592304\\
41	0.000575974679876154\\
42	0.000556833906395553\\
43	0.00054288965245846\\
44	0.000541834764528029\\
45	0.000534698706856379\\
46	0.000517469459637508\\
47	0.00049776725945935\\
48	0.000478172173297532\\
49	0.000460117211177696\\
50	0.000440731184984285\\
51	0.000421159428881053\\
52	0.000401707873878925\\
53	0.000384560513307558\\
54	0.000369926356978079\\
55	0.000358010964220182\\
56	0.000344139821813465\\
57	0.00033125382704704\\
58	0.00031890936310481\\
59	0.000306797352752794\\
60	0.000293839317893102\\
61	0.000279859156687438\\
62	0.000263863640018464\\
63	0.000248683043124333\\
64	0.000234707855873809\\
65	0.000221018794037624\\
66	0.000209385473831444\\
67	0.000199040103333828\\
68	0.000188751104312684\\
69	0.00018029773578646\\
70	0.000173309827276894\\
71	0.000167186495606674\\
72	0.000161337122253728\\
73	0.000154730316655038\\
74	0.0001464835709909\\
75	0.00013909303322047\\
76	0.000133616994092746\\
77	0.000129151162520774\\
78	0.000124893277263787\\
79	0.000120964026925064\\
80	0.000117138061497718\\
81	0.000114484184138109\\
82	0.000111395189786887\\
83	0.000108002354173291\\
84	0.000104731853079724\\
85	0.000101073083469082\\
86	9.7676804592732e-05\\
87	9.41234543524887e-05\\
88	9.09751434337126e-05\\
89	8.74478907855223e-05\\
90	8.27611292148636e-05\\
91	7.88434416633293e-05\\
92	7.51217052636948e-05\\
93	7.11325145823618e-05\\
94	6.75942056852962e-05\\
95	6.45888452062023e-05\\
96	6.15770355480599e-05\\
97	5.86300929699795e-05\\
98	5.60787400211254e-05\\
99	5.39624895924602e-05\\
100	5.17479263570852e-05\\
101	4.97612440111329e-05\\
102	4.79421002671117e-05\\
103	4.58877325788442e-05\\
104	4.39807953553887e-05\\
105	4.20905945396398e-05\\
106	4.01317961393452e-05\\
107	3.79743222395341e-05\\
108	3.59794655983605e-05\\
109	3.43334626876217e-05\\
110	3.27714967258882e-05\\
111	3.12335621459746e-05\\
112	2.97754947566844e-05\\
113	2.84664008342331e-05\\
114	2.70725577563734e-05\\
115	2.5292637416039e-05\\
116	2.36458182205233e-05\\
117	2.21495539818712e-05\\
118	2.05084666905788e-05\\
119	1.88171852854265e-05\\
120	1.74377590162482e-05\\
121	1.60757203734385e-05\\
122	1.48122426018377e-05\\
123	1.37881546782715e-05\\
124	1.28477295147183e-05\\
125	1.18969529423916e-05\\
126	1.09615661450555e-05\\
127	1.01144026182487e-05\\
128	9.3525786420119e-06\\
129	8.68140534849801e-06\\
130	8.01916221304529e-06\\
131	7.41916827312543e-06\\
132	6.93227646431113e-06\\
133	6.43584747027605e-06\\
134	6.0325909148636e-06\\
135	5.66894118954468e-06\\
136	5.28991889194199e-06\\
137	4.91738507255563e-06\\
138	4.64486960721333e-06\\
139	4.41259684731129e-06\\
140	4.1873104030707e-06\\
141	4.00796521468642e-06\\
142	3.80290458373702e-06\\
143	3.59331784367759e-06\\
144	3.38713841035408e-06\\
145	3.18752026868297e-06\\
146	3.01041887040361e-06\\
147	2.86173898061541e-06\\
148	2.70414848866605e-06\\
149	2.53669203179352e-06\\
150	2.41035557734892e-06\\
151	2.28558838929377e-06\\
152	2.14879928179254e-06\\
153	2.01123934917036e-06\\
154	1.87299916428648e-06\\
155	1.74253055256635e-06\\
156	1.62242376064032e-06\\
157	1.4990332872108e-06\\
158	1.37584925786543e-06\\
159	1.26149165181388e-06\\
160	1.16593901632536e-06\\
161	1.0863762101981e-06\\
162	1.02224742245613e-06\\
163	9.66016812993452e-07\\
164	9.25717711829092e-07\\
165	8.79614169507876e-07\\
166	8.40402262319654e-07\\
167	7.89672373883237e-07\\
168	7.34001930763413e-07\\
169	6.85347903609077e-07\\
170	6.35302195819976e-07\\
171	5.80219491616609e-07\\
172	5.31513797636122e-07\\
173	4.93762713403201e-07\\
174	4.63483292833e-07\\
175	4.31533931333895e-07\\
176	4.0140369792375e-07\\
177	3.77474434926006e-07\\
178	3.58338300180396e-07\\
179	3.36490365469967e-07\\
180	3.16317780076541e-07\\
181	2.97211997595729e-07\\
182	2.82763167578212e-07\\
183	2.69490676617126e-07\\
184	2.54113924216284e-07\\
185	2.39664184167856e-07\\
186	2.2256126222611e-07\\
187	2.07201296429744e-07\\
188	1.93255658262718e-07\\
189	1.80394655805901e-07\\
190	1.71093568726596e-07\\
191	1.62595085973794e-07\\
192	1.56049574404004e-07\\
193	1.49519491420875e-07\\
194	1.44241122895086e-07\\
195	1.39218682756672e-07\\
196	1.35279587130158e-07\\
197	1.3185142065138e-07\\
198	1.28780970006492e-07\\
199	1.25696318861983e-07\\
200	1.21587390805273e-07\\
201	1.18492907242234e-07\\
202	1.15792323123419e-07\\
203	1.10650050027272e-07\\
204	9.9719950052993e-08\\
205	2.66892452684202e-08\\
206	2.35056592219729e-08\\
207	2.30560063583018e-08\\
208	2.28252158024286e-08\\
209	2.26030987482798e-08\\
210	2.23383278129424e-08\\
211	2.19917554009451e-08\\
212	2.15937859480449e-08\\
213	2.11278517285405e-08\\
214	2.06455214694554e-08\\
215	2.02060890934643e-08\\
216	1.9811271050796e-08\\
217	1.94170086222592e-08\\
218	1.90876221089072e-08\\
219	1.86646478718237e-08\\
220	1.81947304859278e-08\\
221	1.77332034204134e-08\\
222	1.72223276709497e-08\\
223	1.66498138746953e-08\\
224	1.61163433490545e-08\\
225	1.5519547636893e-08\\
226	1.49128365534445e-08\\
227	1.43467515754416e-08\\
228	1.37971244970245e-08\\
229	1.31137969139803e-08\\
230	1.23671931950776e-08\\
231	1.18500243742921e-08\\
232	1.13870932844443e-08\\
233	1.08976235507991e-08\\
234	1.03901967695871e-08\\
235	9.91251030186482e-09\\
236	9.50252944181211e-09\\
237	9.15384369800372e-09\\
238	8.68737945029389e-09\\
239	8.09208369873676e-09\\
240	7.41696525203377e-09\\
241	6.84758020985465e-09\\
242	6.39917265968306e-09\\
243	5.97597842154543e-09\\
244	5.61388777218552e-09\\
245	5.2436372890328e-09\\
246	4.93196827732681e-09\\
247	4.66477084982205e-09\\
248	4.35860193117479e-09\\
249	4.15461883162373e-09\\
250	3.96850552877549e-09\\
251	3.81024120331355e-09\\
252	3.67879842863673e-09\\
253	3.57915335636596e-09\\
254	3.43501144000027e-09\\
255	3.25331423908865e-09\\
256	3.06732847869666e-09\\
257	2.8906766165504e-09\\
258	2.7029085045964e-09\\
259	2.543295639606e-09\\
260	2.37660246406069e-09\\
261	2.25485726537834e-09\\
262	2.1199593289294e-09\\
263	2.0061186513816e-09\\
264	1.87429830362123e-09\\
265	1.74526729649075e-09\\
266	1.60201917109472e-09\\
267	1.46634778170633e-09\\
268	1.32409335338572e-09\\
269	1.16670717185498e-09\\
270	1.02736287880257e-09\\
271	9.1679778559981e-10\\
272	7.98352122165368e-10\\
273	7.00561680337871e-10\\
274	6.49000651900849e-10\\
275	6.06537288908724e-10\\
276	5.57126481874537e-10\\
277	5.04082742364496e-10\\
278	4.68463135036837e-10\\
279	4.27691492509509e-10\\
280	3.82517698290479e-10\\
281	3.51935677211354e-10\\
282	3.1736954175227e-10\\
283	2.86529797758866e-10\\
284	2.63254625336498e-10\\
285	2.42766200753421e-10\\
286	2.23348687767781e-10\\
287	2.04957643655079e-10\\
288	1.91566194508593e-10\\
289	1.78507224622602e-10\\
290	1.67496105013526e-10\\
291	1.5831062780994e-10\\
292	1.46560152514617e-10\\
293	1.36203554674381e-10\\
294	1.26048757164588e-10\\
295	1.13921785840242e-10\\
296	1.02951839742283e-10\\
297	9.15852209876315e-11\\
298	8.10983468896245e-11\\
299	7.32248820065505e-11\\
300	6.54823855127549e-11\\
301	5.71832895702227e-11\\
302	5.17562246178338e-11\\
303	4.69766048282593e-11\\
304	4.24420556351048e-11\\
305	3.75430553137218e-11\\
306	3.35203528045504e-11\\
307	2.98942184443967e-11\\
308	2.70970803058235e-11\\
309	2.45032400234644e-11\\
310	2.19781956028426e-11\\
311	1.99853499546535e-11\\
312	1.8683051293457e-11\\
313	1.71830298097213e-11\\
314	1.58427635574642e-11\\
315	1.48899799391066e-11\\
316	1.40482871528689e-11\\
317	1.32207288407361e-11\\
318	1.25225523467755e-11\\
319	1.19163799195173e-11\\
320	1.12179391948424e-11\\
321	1.0416830066731e-11\\
322	9.80437112010873e-12\\
323	9.1215305109587e-12\\
324	8.46145184559725e-12\\
325	7.77196754896944e-12\\
326	7.15658129561705e-12\\
327	6.60936840894688e-12\\
328	6.0640533296398e-12\\
329	5.58136185326802e-12\\
330	5.01783502958191e-12\\
331	4.52878831972264e-12\\
332	4.12948229076321e-12\\
333	3.68879262685427e-12\\
334	3.26460366815859e-12\\
335	2.7912335066032e-12\\
336	2.39431460608572e-12\\
337	2.07844752402533e-12\\
338	1.79250597238101e-12\\
339	1.5503409006306e-12\\
340	1.35334810854154e-12\\
341	1.16635350858785e-12\\
342	9.62237608634033e-13\\
343	7.52449241288571e-13\\
344	6.025589437415e-13\\
345	4.73606285876135e-13\\
346	3.62957532409122e-13\\
347	2.60620961133045e-13\\
348	2.06415971264402e-13\\
349	1.67093690092719e-13\\
350	1.39324764358601e-13\\
351	1.11101711230557e-13\\
352	8.74796829529601e-14\\
353	6.8839978370861e-14\\
354	5.39426573062977e-14\\
355	4.27126926444384e-14\\
356	3.3268583843046e-14\\
357	2.61346785420445e-14\\
358	2.01982719068712e-14\\
359	1.69252096880442e-14\\
360	1.44094160221518e-14\\
361	1.23818451294737e-14\\
362	1.05132403575099e-14\\
363	9.52975584530999e-15\\
};
\addlegendentry{$\varepsilon=10^{-14}$}
\legend{}
\end{axis}
\end{tikzpicture}

%% file: Figs/tikz/GMRES_residual_NTAcc_MR.tex
\definecolor{mycolor1}{rgb}{0.00000,0.44700,0.74100}%
\definecolor{mycolor2}{rgb}{0.85000,0.32500,0.09800}%
\definecolor{mycolor3}{rgb}{0.92900,0.69400,0.12500}%
\definecolor{mycolor4}{rgb}{0.49400,0.18400,0.55600}%
\begin{tikzpicture}

\begin{axis}[%
width=3.65cm,
height=3.7cm,
scale only axis,
tick label style={font=\footnotesize},
xmin=0,
xmax=17,
xlabel style={font=\footnotesize\color{white!15!black},yshift=.2cm},
xlabel={Iteration},
ymode=log,
ymin=1.1e-15,
ymax=1.1,
yminorticks=true,
ytickten = {-14, -11, -8, -5, -2},
ylabel style={font=\footnotesize\color{white!15!black},yshift=-.2cm},
ylabel={Residual norm},
axis background/.style={fill=white},
legend style={at={(axis cs:16.75,7e-1)}, anchor=north east, legend cell align=left, align=left, font = \scriptsize, draw=white!15!black}
]
\addplot [color=mycolor1, only marks, mark=x, mark options={solid, mycolor1}]
  table[row sep=crcr]{%
0	1\\
1	0.000246415344988842\\
2	6.00027760743915e-05\\
3	6.85397357753923e-06\\
4	3.03533771484673e-06\\
5	4.98355504325039e-07\\
6	1.68240496268028e-07\\
7	4.84238882445531e-08\\
8	1.02662754773091e-08\\
9	1.65339574576442e-09\\
10	2.89788361283949e-10\\
11	5.51342801000868e-11\\
12	7.82927544870987e-12\\
13	1.00749817700086e-12\\
14	1.05317737036711e-13\\
15	9.77714186705537e-15\\
};
\addlegendentry{$\varepsilon=10^{-8}$}

\addplot [color=mycolor2, only marks, mark=x, mark options={solid, mycolor2}]
  table[row sep=crcr]{%
0	1\\
1	0.000240967758946996\\
2	5.88323281334302e-05\\
3	6.74677295719915e-06\\
4	2.96241057575515e-06\\
5	5.10932707152857e-07\\
6	1.63669065872518e-07\\
7	4.61328127751707e-08\\
8	9.79942969657817e-09\\
9	1.57907951986143e-09\\
10	3.09178447274716e-10\\
11	4.86968242952493e-11\\
12	7.88971658097128e-12\\
13	8.14022002508799e-13\\
14	8.01711179853938e-14\\
15	1.22344139235495e-14\\
16	1.15227338688276e-15\\
};
\addlegendentry{$\varepsilon=10^{-10}$}

\addplot [color=mycolor3, only marks, mark=x, mark options={solid, mycolor3}]
  table[row sep=crcr]{%
0	1\\
1	0.000238402316880132\\
2	5.67353675788495e-05\\
3	6.28653793387165e-06\\
4	2.91073368071028e-06\\
5	3.73266491238539e-07\\
6	1.05080219740138e-07\\
7	3.09539589293456e-08\\
8	5.74685103400481e-09\\
9	1.2023160537885e-09\\
10	2.04630841213289e-10\\
11	2.91548566567338e-11\\
12	3.90123179660983e-12\\
13	6.42876676528489e-13\\
14	7.58263744308684e-14\\
15	8.44350187881746e-15\\
};
\addlegendentry{$\varepsilon=10^{-12}$}

\addplot [color=mycolor4, only marks, mark=x, mark options={solid, mycolor4}]
  table[row sep=crcr]{%
0	1\\
1	0.000602951428246474\\
2	0.000319952568869379\\
3	6.10464641027837e-05\\
4	1.57459398589476e-05\\
5	3.40139394978971e-06\\
6	4.78972340014141e-07\\
7	7.22912191338113e-08\\
8	3.01430527578213e-08\\
9	8.15785690040249e-09\\
10	1.39833047846908e-09\\
11	1.50674244638642e-10\\
12	1.98705193753001e-11\\
13	1.62040870247687e-12\\
14	1.6363065287968e-13\\
15	2.15792211059975e-14\\
16	1.5729767750935e-15\\
};
\addlegendentry{$\varepsilon=10^{-14}$}

\end{axis}
\end{tikzpicture}

%% file: Figs/tikz/gpe_N100_beta40.tex
%
%
\definecolor{mycolor1}{rgb}{1.00000,0.00000,1.00000}%
\definecolor{mycolor2}{rgb}{0.49400,0.18400,0.55600}%
\definecolor{mycolor3}{rgb}{0.00000,1.00000,1.00000}%
\definecolor{mycolor4}{rgb}{0.63500,0.07800,0.18400}%
\begin{tikzpicture}

\begin{axis}[%
width=3.65cm,
height=3.7cm,
scale only axis,
tick label style={font=\footnotesize},
xmin=0,
xmax=29,
xlabel style={font=\footnotesize\color{white!15!black},yshift=.2cm},
xlabel={Iteration},
ymode=log,
ymin=6e-16,
ymax=2e-1,
yminorticks=true,
ytickten = {-14, -11, -8, -5, -2},
ylabel style={font=\footnotesize\color{white!15!black},yshift=-.2cm},
ylabel={Residual norm},
axis background/.style={fill=white},
legend columns=4, 
legend style={at={(axis cs:0,10)}, anchor=north west, legend cell align=left, align=left, font = \tiny, draw=white!15!black}
]
\addplot [color=myscf, mark=x, mark options={solid, myscf}]
  table[row sep=crcr]{%
0	0.0990244932978807\\
1	0.000440734420970345\\
2	0.000174563106559469\\
3	7.09346060075529e-05\\
4	2.8037430572767e-05\\
5	1.10993947108036e-05\\
6	4.37639480570924e-06\\
7	1.72645870811211e-06\\
8	6.80707167377759e-07\\
9	2.68417522056109e-07\\
10	1.05834779375595e-07\\
11	4.17305697084085e-08\\
12	1.6454149800777e-08\\
13	6.48780804378831e-09\\
14	2.55811362689898e-09\\
15	1.00865326352418e-09\\
16	3.9770756649277e-10\\
17	1.56814372265001e-10\\
18	6.18312241889036e-11\\
19	2.43797826979197e-11\\
20	9.61284107397333e-12\\
21	3.79030318093901e-12\\
22	1.49450236815096e-12\\
23	5.89280153159493e-13\\
24	2.32344640429958e-13\\
25	9.16439117151027e-14\\
26	3.61259008272575e-14\\
27	1.43812062091751e-14\\
28	5.25801965638311e-15\\
};
\addlegendentry{SCF}

\addplot [color=black, mark=+, mark options={solid, black}]
  table[row sep=crcr]{%
0	0.0990244932978807\\
1	0.000440734420970345\\
2	0.00181275422616667\\
3	0.00132789086877625\\
4	0.00222017856897675\\
5	0.00072247859985214\\
6	0.000789967762143567\\
7	2.76915207900852e-05\\
8	6.90052632672672e-06\\
9	3.38156379829322e-06\\
10	1.68128229761561e-07\\
11	1.04544725091721e-05\\
12	2.33531323841824e-06\\
13	4.08282510581799e-07\\
14	3.80261268223759e-07\\
15	2.2919800222929e-07\\
16	1.75634411633883e-07\\
17	5.23069009374108e-07\\
18	2.24405271313019e-07\\
19	1.15935618187646e-07\\
20	9.32267583292391e-09\\
21	1.57731022350224e-13\\
22	6.61226936421242e-10\\
23	1.54382235217669e-14\\
24	2.12876735645714e-15\\
};
\addlegendentry{AA(4)-$X_i$ w/$X_{i-1}$}

\addplot [color=red, mark=+, mark options={solid, red}]
  table[row sep=crcr]{%
0	0.0990244932978807\\
1	0.000440734420970392\\
2	0.000997026126379697\\
3	0.00106869658465397\\
4	0.000325583144665157\\
5	0.000392070187749878\\
6	0.000138576002435364\\
7	9.18913137478888e-05\\
8	0.0001869487714958\\
9	0.000126499366161238\\
10	0.000131389937769752\\
11	3.98921907468188e-05\\
12	2.45368770835783e-05\\
13	1.2417564961616e-05\\
14	3.12999444370661e-06\\
15	3.82913805674838e-06\\
16	5.41583495706658e-06\\
17	2.95182364264356e-06\\
18	2.95809311512399e-06\\
19	4.39934361176124e-07\\
20	9.18992344794201e-07\\
21	4.49797013016734e-07\\
22	5.58253985536428e-07\\
23	2.23294527583161e-07\\
24	1.71650839240699e-07\\
25	1.80692998770112e-07\\
26	1.07646910487885e-07\\
27	1.14148324829546e-07\\
28	1.40331684550893e-07\\
29	4.55732314589051e-07\\
30	4.81378950731165e-08\\
31	6.05159323202198e-08\\
32	1.2464642106038e-08\\
33	1.41247433415331e-09\\
34	1.32854246188368e-09\\
35	4.52997463157511e-09\\
36	2.43145871090576e-09\\
37	2.48469004348923e-09\\
38	1.0331802096817e-09\\
39	6.56458868123349e-10\\
40	3.43936017866296e-10\\
41	3.18327447323757e-10\\
42	2.42599641113218e-09\\
43	1.71632311589715e-09\\
44	1.27609303342357e-09\\
45	4.03786486718899e-10\\
46	1.10598693887427e-10\\
47	8.62227903086653e-11\\
48	1.15661893814852e-10\\
49	2.13212025066518e-10\\
50	2.48117000354166e-11\\
51	6.54927510558164e-11\\
52	1.73520699320702e-11\\
53	3.25509104463501e-12\\
54	1.37313780342496e-12\\
55	1.71293860639004e-12\\
56	2.52778541466818e-12\\
57	6.04652674736837e-13\\
58	8.12787026976996e-13\\
59	6.37494682909077e-13\\
60	1.58369601113931e-12\\
61	1.03391969317253e-12\\
62	1.00943318910091e-12\\
63	8.6999518992277e-13\\
64	9.45971757962521e-13\\
65	3.47516574767585e-13\\
66	1.59171630077764e-13\\
67	9.76711310790219e-13\\
68	5.11917182123042e-13\\
69	2.32525086795486e-13\\
70	1.00315786806681e-13\\
71	2.57083709034334e-13\\
72	1.11475634627976e-13\\
73	3.07048005325437e-13\\
74	4.60159932932142e-13\\
75	7.81845137609191e-13\\
76	6.03398283503396e-13\\
77	4.18658522462746e-13\\
78	5.39428275719443e-13\\
79	4.373148567998e-13\\
};
\addlegendentry{AA(4)-$X_i$ w/\texttt{randn}}

\addplot [color=blue, mark=o, mark options={solid, blue}]
  table[row sep=crcr]{%
0	0.0990244932978807\\
1	0.000440734420970345\\
2	0.00014197985538763\\
3	1.72987592230172e-05\\
4	4.64275247185612e-07\\
5	1.76534434090331e-08\\
6	9.40195648611261e-10\\
7	3.80320678660745e-11\\
8	1.79100273327716e-13\\
9	6.50882519321499e-15\\
};
\addlegendentry{AA(4)-$\rho_i$}

\addplot [color=myalign, mark=o, mark options={solid, myalign}]
  table[row sep=crcr]{%
0	0.0990244932978807\\
1	0.00044073442097042\\
2	9.00568222049277e-05\\
3	8.95575576765357e-06\\
4	2.89714195731618e-07\\
5	1.1014616686652e-08\\
6	2.35812939446479e-10\\
7	1.37042684292111e-11\\
8	6.31738645631626e-13\\
9	5.52202313022991e-15\\
};
\addlegendentry{AA(4)-$\ALIGN{X_i}$}

\addplot [color=mygrassaa, mark=o, mark options={solid, mygrassaa}]
  table[row sep=crcr]{%
0	0.0990244932978807\\
1	0.000440734420970345\\
2	0.000243850920009095\\
3	4.55675325090239e-06\\
4	1.98214360106827e-07\\
5	4.05459335398761e-09\\
6	1.39505017019356e-10\\
7	3.42248128256877e-12\\
8	9.81102769720254e-14\\
9	3.98226187515187e-15\\
};
\addlegendentry{AA(4)-Gr}

\addplot [color=mygnt, mark=diamond, mark options={solid, mygnt}]
  table[row sep=crcr]{%
0	0.0990244932978807\\
1	0.000440734420970345\\
2	1.8295007090761e-05\\
3	7.6483052874305e-07\\
4	9.27814325587921e-10\\
5	2.55258079830719e-15\\
};
\addlegendentry{NTAcc}
\legend{}
\end{axis}
\end{tikzpicture}%

%% file: Figs/tikz/gpe_N100_beta55.tex
%
%
\definecolor{mycolor1}{rgb}{1.00000,0.00000,1.00000}%
\definecolor{mycolor2}{rgb}{0.49400,0.18400,0.55600}%
\definecolor{mycolor3}{rgb}{0.00000,1.00000,1.00000}%
\definecolor{mycolor4}{rgb}{0.63500,0.07800,0.18400}%
\begin{tikzpicture}

\begin{axis}[%
width=3.65cm,
height=3.7cm,
scale only axis,
tick label style={font=\footnotesize},
xmin=0,
xmax=29,
xlabel style={font=\footnotesize\color{white!15!black},yshift=.2cm},
xlabel={Iteration},
ymode=log,
ymin=6e-16,
ymax=2e-1,
yminorticks=true,
ytickten = {-14, -11, -8, -5, -2},
ylabel style={font=\footnotesize\color{white!15!black},yshift=-.2cm},
ylabel={Residual norm},
axis background/.style={fill=white},
legend style={at={(axis cs:28.3,2e-15)}, anchor=south east, legend cell align=left, align=left, font = \tiny, draw=white!15!black}
]
\addplot [color=myscf, mark=x, mark options={solid, myscf}]
  table[row sep=crcr]{%
0	0.0990244932978807\\
1	0.0006060098288344\\
2	0.00031429398156418\\
3	0.000175334008457136\\
4	9.2009561215979e-05\\
5	4.94810472055626e-05\\
6	2.61683805893961e-05\\
7	1.39479033087535e-05\\
8	7.40075687609335e-06\\
9	3.9358375989077e-06\\
10	2.09051626200927e-06\\
11	1.11110224553148e-06\\
12	5.90339442860846e-07\\
13	3.1371129285195e-07\\
14	1.66692286475793e-07\\
15	8.85775459579621e-08\\
16	4.70673400063777e-08\\
17	2.5010484485209e-08\\
18	1.32898821783481e-08\\
19	7.06190669302343e-09\\
20	3.75250993326428e-09\\
21	1.99398660917594e-09\\
22	1.05955218017315e-09\\
23	5.63018412468733e-10\\
24	2.99173254688067e-10\\
25	1.58972871503606e-10\\
26	8.44740332053772e-11\\
27	4.48873062276033e-11\\
28	2.3851956793315e-11\\
29	1.26743193003061e-11\\
30	6.73481501992411e-12\\
31	3.57870349624052e-12\\
32	1.90162139956805e-12\\
33	1.01047151078827e-12\\
34	5.36936401641668e-13\\
35	2.85310379622689e-13\\
36	1.51641930643608e-13\\
37	8.05658345327162e-14\\
38	4.27021071332506e-14\\
39	2.27520475044775e-14\\
40	1.18790725253109e-14\\
41	6.43754600297808e-15\\
};
\addlegendentry{SCF}

\addplot [color=black, mark=+, mark options={solid, black}]
  table[row sep=crcr]{%
0	0.0990244932978807\\
1	0.0006060098288344\\
2	0.0018430716029794\\
3	0.0015237937792179\\
4	0.00169512787543837\\
5	0.00216576296884576\\
6	0.00222101339794175\\
7	0.00218145233219489\\
8	0.00204833813515135\\
9	0.00183291498361826\\
10	0.000345729875479521\\
11	0.00011372023549249\\
12	1.03364448109919e-05\\
13	0.000319987589990385\\
14	0.000106197290672857\\
15	0.000374597408663545\\
16	0.000787373298026649\\
17	0.000753483602194397\\
18	0.000748856366076792\\
19	0.000734909042491998\\
20	0.00070577072758237\\
21	0.00214997736374217\\
22	0.00142129515799262\\
23	0.000547179963758792\\
24	0.00179012385705461\\
25	0.00444174455831632\\
26	0.00064160856404378\\
27	0.000123311934074098\\
28	9.26987092519294e-06\\
29	9.72166658990461e-06\\
30	5.53421194536339e-07\\
31	1.4329262258546e-08\\
32	9.30762448931315e-07\\
33	2.09863474887567e-06\\
34	7.70762250656015e-07\\
35	2.11210217478562e-06\\
36	4.79693150446818e-06\\
37	5.08851465220216e-08\\
38	1.65172664058801e-10\\
39	1.12213802939643e-07\\
40	4.715337685438e-07\\
41	5.50207738806656e-07\\
42	8.96434958105917e-08\\
43	1.07322016961829e-06\\
44	4.45854724316795e-06\\
45	5.41815719882755e-06\\
46	3.10061747194092e-06\\
47	2.62245487583359e-08\\
48	5.72697839132284e-11\\
49	4.32862095968659e-08\\
50	4.32803091794176e-08\\
51	6.20709510178677e-08\\
52	3.48164984013235e-08\\
53	1.79657342425356e-08\\
54	2.00613388001442e-08\\
55	1.23477905338064e-10\\
56	2.80561407760292e-10\\
57	1.58025209269877e-13\\
58	8.98732519564699e-14\\
59	2.0718710586068e-14\\
60	2.10822934345247e-14\\
61	2.57697719385862e-14\\
62	4.03827012044764e-15\\
};
\addlegendentry{AA(4)-$x_i$-w}

\addplot [color=red, mark=+, mark options={solid, red}]
  table[row sep=crcr]{%
0	0.0990244932978807\\
1	0.000606009828834274\\
2	0.0013373185440472\\
3	0.00140726725289529\\
4	0.000395287257977606\\
5	0.000462172922726009\\
6	0.000241123525092409\\
7	0.000382035880178876\\
8	0.00020800849020486\\
9	0.000127332086795987\\
10	0.000443663812604766\\
11	2.94213674350908e-05\\
12	2.85155864878381e-05\\
13	2.16516196473261e-05\\
14	3.66708392512307e-06\\
15	3.41993345490493e-06\\
16	1.18439136863009e-06\\
17	1.05073659740272e-06\\
18	1.23278409815387e-07\\
19	1.22827277532335e-06\\
20	1.06545165544133e-07\\
21	1.60643564737931e-08\\
22	1.32753975110904e-08\\
23	5.79587190295993e-08\\
24	1.62705628986145e-07\\
25	1.69259289587548e-07\\
26	6.59314727888355e-07\\
27	7.93370166086577e-07\\
28	7.7426894686272e-07\\
29	7.11349661522875e-07\\
30	7.47708618813009e-07\\
31	4.07173877392586e-07\\
32	9.03667541882152e-06\\
33	2.05557675870746e-07\\
34	3.39236134883291e-07\\
35	4.89974835110889e-07\\
36	4.39287588781191e-07\\
37	1.67884829071137e-07\\
38	3.23098401278303e-07\\
39	8.46049413553474e-07\\
40	2.2851581121629e-07\\
41	2.31219708811809e-07\\
42	1.73342475140317e-06\\
43	1.09007160307067e-07\\
44	6.03888355531467e-08\\
45	6.42703744480758e-08\\
46	7.93289684502069e-08\\
47	4.09198530675102e-08\\
48	9.74708256357993e-08\\
49	3.05672486673793e-08\\
50	5.21805379520064e-09\\
51	2.14590796802601e-08\\
52	9.08823371217546e-08\\
53	1.92874900406313e-08\\
54	1.51444850702369e-08\\
55	2.61987667630151e-08\\
56	4.63215066153847e-08\\
57	9.72985865838972e-09\\
58	1.0151877367668e-08\\
59	9.67451267622477e-09\\
60	5.42386654129399e-09\\
61	6.03379423415088e-10\\
62	1.28111470901526e-09\\
63	8.7809075781789e-10\\
64	5.98991445838151e-10\\
65	1.01666833669991e-09\\
66	1.27823821065196e-09\\
67	1.73974576616195e-09\\
68	1.58180431230435e-09\\
69	8.69566718299288e-10\\
70	1.17039476300615e-09\\
71	3.23657406221344e-09\\
72	1.76939861227166e-09\\
73	1.3379850666932e-09\\
74	2.50193746469139e-09\\
75	3.24983516367793e-09\\
76	4.99218846042892e-09\\
77	7.9447391333733e-10\\
78	6.91740459564882e-10\\
79	5.46596264880605e-10\\
};
\addlegendentry{AA(4)-$x_i$-r}

\addplot [color=blue, mark=o, mark options={solid, blue}]
  table[row sep=crcr]{%
0	0.0990244932978807\\
1	0.0006060098288344\\
2	0.000237132862322743\\
3	3.17486629257929e-05\\
4	1.26986727482607e-06\\
5	1.0185302897451e-07\\
6	6.7685299602606e-09\\
7	3.89345014027176e-10\\
8	2.06924112405787e-12\\
9	7.27449755535174e-14\\
10	3.21895924122797e-15\\
};
\addlegendentry{AA(4)-$\rho_i$}

\addplot [color=myalign, mark=o, mark options={solid, myalign}]
  table[row sep=crcr]{%
0	0.0990244932978807\\
1	0.000606009828834374\\
2	0.000142497228616539\\
3	1.8956692656041e-05\\
4	1.02157012269216e-06\\
5	4.55957417576265e-08\\
6	1.45533809359527e-09\\
7	1.17223629884909e-10\\
8	8.19233741348825e-12\\
9	6.66558173009353e-14\\
10	2.98104559979405e-15\\
};
\addlegendentry{AA(4)-$\ALIGN{x_i}$}

\addplot [color=mygrassaa, mark=o, mark options={solid, mygrassaa}]
  table[row sep=crcr]{%
0	0.0990244932978807\\
1	0.0006060098288344\\
2	0.000440470678764396\\
3	1.03340795799567e-05\\
4	4.61536233184273e-07\\
5	1.51303950908162e-08\\
6	7.9211380843134e-10\\
7	2.57651783948421e-11\\
8	6.89514827505797e-13\\
9	3.66507456972503e-14\\
10	1.4278210652714e-15\\
};
\addlegendentry{AA(4)-Gr}

\addplot [color=mygnt, mark=diamond, mark options={solid, mygnt}]
  table[row sep=crcr]{%
0	0.0990244932978807\\
1	0.0006060098288344\\
2	4.70912700039277e-05\\
3	1.4711328369189e-06\\
4	1.58752520542954e-10\\
5	8.5726796707918e-16\\
};
\addlegendentry{NTAcc}
\legend{}
\end{axis}
\end{tikzpicture}%

%% file: Figs/tikz/gpe_N100_beta72.tex
%
%
\definecolor{mycolor1}{rgb}{1.00000,0.00000,1.00000}%
\definecolor{mycolor2}{rgb}{0.49400,0.18400,0.55600}%
\definecolor{mycolor3}{rgb}{0.00000,1.00000,1.00000}%
\definecolor{mycolor4}{rgb}{0.63500,0.07800,0.18400}%
\begin{tikzpicture}

\begin{axis}[%
width=3.65cm,
height=3.7cm,
scale only axis,
tick label style={font=\footnotesize},
xmin=0,
xmax=29,
xlabel style={font=\footnotesize\color{white!15!black},yshift=.2cm},
xlabel={Iteration},
ymode=log,
ymin=6e-16,
ymax=2e-1,
yminorticks=true,
ytickten = {-14, -11, -8, -5, -2},
ylabel style={font=\footnotesize\color{white!15!black},yshift=-.2cm},
ylabel={Residual norm},
axis background/.style={fill=white},
legend style={at={(axis cs:28.75,9e-16)}, anchor=south east, legend cell align=left, align=left, font=\tiny, draw=white!15!black}
]
\addplot [color=myscf, mark=x, mark options={solid, myscf}]
  table[row sep=crcr]{%
0	0.0990244932978807\\
1	0.000793321957746615\\
2	0.000505902288737185\\
3	0.000375236883587034\\
4	0.000243614286204862\\
5	0.000172778777935657\\
6	0.000114996968715432\\
7	7.99157257039438e-05\\
8	5.39193708579956e-05\\
9	3.71234268651337e-05\\
10	2.52095886554883e-05\\
11	1.72815438734583e-05\\
12	1.17708040511018e-05\\
13	8.05267831106825e-06\\
14	5.4925161825353e-06\\
15	3.75398968280503e-06\\
16	2.56216372671796e-06\\
17	1.75039478249554e-06\\
18	1.195037664495e-06\\
19	8.16245597051288e-07\\
20	5.57349821152624e-07\\
21	3.80649443356402e-07\\
22	2.59932663454417e-07\\
23	1.775164625244e-07\\
24	1.21223731877733e-07\\
25	8.27858871958723e-08\\
26	5.65342390552258e-08\\
27	3.86078800532041e-08\\
28	2.63653841677971e-08\\
29	1.80051409054083e-08\\
30	1.22957784046173e-08\\
31	8.39687181200001e-09\\
32	5.7342636736869e-09\\
33	3.91596430219169e-09\\
34	2.67423244149754e-09\\
35	1.82624909892692e-09\\
36	1.24715542983769e-09\\
37	8.51689618908281e-10\\
38	5.81623557157423e-10\\
39	3.97194033793644e-10\\
40	2.71246023450529e-10\\
41	1.85235440736555e-10\\
42	1.26498323021828e-10\\
43	8.63864171420062e-11\\
44	5.89937790527401e-11\\
45	4.02872008765324e-11\\
46	2.75123605824108e-11\\
47	1.87883517344696e-11\\
48	1.28306789071979e-11\\
49	8.76215854175845e-12\\
50	5.98373966281567e-12\\
51	4.08633491683417e-12\\
52	2.79058782213866e-12\\
53	1.90571492877448e-12\\
54	1.30144369765159e-12\\
55	8.88790683036007e-13\\
56	6.07002332056574e-13\\
57	4.14579328487932e-13\\
58	2.83164913556519e-13\\
59	1.93433469929323e-13\\
60	1.32218286728763e-13\\
61	9.03560161650898e-14\\
62	6.16814427077997e-14\\
63	4.22027561500122e-14\\
64	2.86782137333834e-14\\
65	1.72652708190374e-14\\
66	1.12842897677703e-14\\
67	4.28832330534473e-15\\
};
\addlegendentry{SCF}

\addplot [color=black, mark=+, mark options={solid, black}]
  table[row sep=crcr]{%
0	0.0990244932978807\\
1	0.000793321957746615\\
2	0.000505778884151731\\
3	0.00037467100402109\\
4	0.00172921540194969\\
5	0.00173753047198162\\
6	0.00176853341075523\\
7	0.00409019564146395\\
8	0.00415779581829346\\
9	0.00260704218651468\\
10	0.00315343519554293\\
11	0.000213345674986728\\
12	0.00198804966156806\\
13	0.000511840802576621\\
14	0.00012470439278988\\
15	0.000254092538357016\\
16	0.000554551818772542\\
17	1.46809768308443e-05\\
18	0.000231630213851129\\
19	1.15332613238812e-06\\
20	7.63845287788961e-05\\
21	9.16734676365396e-06\\
22	2.5001321326603e-05\\
23	1.78031773959058e-05\\
24	0.000330970016634234\\
25	0.000799318161993253\\
26	0.00115471465170157\\
27	0.000680213478418156\\
28	0.000703141118218964\\
29	0.000825862776147391\\
30	0.00126206520869214\\
31	0.00189208517270387\\
32	0.00324164655231657\\
33	0.00331954647901626\\
34	0.00264424887864881\\
35	0.00215386692320341\\
36	0.00217386878532268\\
37	0.0021581302657965\\
38	0.00214782976542551\\
39	0.000506545316704227\\
40	0.00389384715155838\\
41	0.00119517953642473\\
42	0.000684966642757394\\
43	0.000382516008668111\\
44	0.00152525379110968\\
45	0.000131965582376987\\
46	0.00012233562397619\\
47	4.98093500973896e-06\\
48	4.86879234627999e-05\\
49	4.6751703473815e-05\\
50	3.20763857175826e-05\\
51	0.000109223329723787\\
52	5.93506656925714e-06\\
53	9.06406913824871e-09\\
54	2.0686091887513e-06\\
55	2.36425005580497e-06\\
56	6.45632123092064e-09\\
57	6.56976899097489e-10\\
58	9.40977342650717e-11\\
59	1.65707623112808e-11\\
60	6.56981399833745e-11\\
61	1.68779989655085e-11\\
62	5.09695170019906e-12\\
63	1.70018656395729e-11\\
64	3.25258130046094e-11\\
65	3.34880992082025e-11\\
66	1.28532840931731e-10\\
67	4.98417347837963e-11\\
68	7.65452734320245e-12\\
69	3.01778495116622e-13\\
70	3.01779899975863e-13\\
71	1.02699083490972e-12\\
72	3.3950151389903e-13\\
73	3.89264787001999e-12\\
74	3.43341099771412e-12\\
75	2.96676146248394e-12\\
76	1.56205076974749e-13\\
77	1.15528024762251e-12\\
78	2.7565561204512e-12\\
79	1.34546601987646e-11\\
};
\addlegendentry{AA(4)-$x_i$\!-w}

\addplot [color=red, mark=+, mark options={solid, red}]
  table[row sep=crcr]{%
0	0.0990244932978807\\
1	0.000793321957746732\\
2	0.00167606429095005\\
3	0.00123189202513303\\
4	0.000407631600408379\\
5	0.00110963855957887\\
6	0.00118522168001589\\
7	0.00116342829658255\\
8	0.000444287327621313\\
9	0.000429631906160129\\
10	0.00105101505191495\\
11	0.000463968612327928\\
12	0.000553760719911351\\
13	0.000207570784040663\\
14	0.000311749880003851\\
15	0.000261150289655203\\
16	0.000191589194230664\\
17	0.000366335473318428\\
18	0.000125059650226373\\
19	0.000108516633542774\\
20	0.000720691167555819\\
21	0.000316440802227516\\
22	0.000123242182126973\\
23	0.000190230394401797\\
24	2.01580175537546e-05\\
25	4.80884156335729e-05\\
26	2.63809352389251e-05\\
27	1.17000305961577e-05\\
28	8.2506859407969e-06\\
29	9.36198073291547e-06\\
30	1.00831723633216e-05\\
31	2.37772142472549e-05\\
32	1.9214507756088e-05\\
33	1.65659659643152e-05\\
34	1.31840203735104e-05\\
35	0.000183611890783038\\
36	3.27229703125766e-05\\
37	8.22427781545993e-06\\
38	4.48600727533577e-05\\
39	7.38861358437118e-06\\
40	6.38457179142422e-06\\
41	3.10641473855297e-06\\
42	2.11031786889796e-06\\
43	4.65885238566391e-07\\
44	5.35556881645664e-08\\
45	5.79075595234955e-07\\
46	2.07390425048103e-07\\
47	1.6324444369023e-07\\
48	6.21066224234887e-07\\
49	8.133216405091e-07\\
50	2.69884615023796e-07\\
51	3.54849255466696e-07\\
52	1.80101011326073e-07\\
53	2.18159482212555e-07\\
54	1.78329660160235e-07\\
55	2.07436862524463e-07\\
56	1.16623173918629e-07\\
57	1.55590272034121e-07\\
58	2.55583982052775e-07\\
59	9.23336015715331e-08\\
60	8.97424584153906e-08\\
61	1.94308690755789e-07\\
62	7.45746772256629e-08\\
63	5.80015998306352e-08\\
64	4.45871008004501e-08\\
65	4.15108154011323e-08\\
66	2.10391073995274e-08\\
67	5.06325658878502e-08\\
68	4.44350397454593e-08\\
69	5.73023842135311e-08\\
70	3.57419432837921e-08\\
71	3.05090319978185e-08\\
72	3.02868546615455e-08\\
73	1.60172394157052e-08\\
74	2.25144402352638e-08\\
75	4.64027143045567e-09\\
76	1.14640468639775e-08\\
77	1.06108190305045e-08\\
78	1.4681123640061e-08\\
79	1.62940145130529e-08\\
};
\addlegendentry{AA(4)-$x_i$\!-r}

\addplot [color=blue, mark=o, mark options={solid, blue}]
  table[row sep=crcr]{%
0	0.0990244932978807\\
1	0.000793321957746615\\
2	0.000352705395786514\\
3	4.86528988489987e-05\\
4	2.72496583243581e-06\\
5	3.3261699394906e-07\\
6	2.90484396230364e-08\\
7	1.93103493728127e-09\\
8	1.8679625635354e-11\\
9	5.34832724711875e-13\\
10	2.18647041269052e-14\\
11	3.04732610141663e-15\\
};
\addlegendentry{AA(4)-$\rho_i$}

\addplot [color=myalign, mark=o, mark options={solid, myalign}]
  table[row sep=crcr]{%
0	0.0990244932978807\\
1	0.000793321957747072\\
2	0.000199367246257309\\
3	3.44331282110784e-05\\
4	3.06692137842371e-06\\
5	1.4305245513408e-07\\
6	6.61649046498924e-09\\
7	6.81951959035306e-10\\
8	5.95413603516613e-11\\
9	8.8042362020064e-13\\
10	3.41728703806538e-14\\
11	7.83663958637459e-15\\
};
\addlegendentry{AA(4)-$\ALIGN{x_i}$}

\addplot [color=mygrassaa, mark=o, mark options={solid, mygrassaa}]
  table[row sep=crcr]{%
0	0.0990244932978807\\
1	0.000793321957746615\\
2	0.000693272919837388\\
3	2.20242529694477e-05\\
4	7.29056244921836e-07\\
5	5.44333207312882e-08\\
6	3.06636424832027e-09\\
7	2.18868275118832e-10\\
8	3.17142200073024e-12\\
9	1.14948628445261e-13\\
10	3.37407772188932e-15\\
};
\addlegendentry{AA(4)-Gr}

\addplot [color=mygnt, mark=diamond, mark options={solid, mygnt}]
  table[row sep=crcr]{%
0	0.0990244932978807\\
1	0.000793321957746615\\
2	0.000106727893368305\\
3	4.70327057647392e-06\\
4	3.17719081741979e-09\\
5	3.00066411767392e-15\\
};
\addlegendentry{AlignNT}
\end{axis}
\end{tikzpicture}%

%% file: Figs/tikz/residual_iter_usps.tex
%
%
\definecolor{mycolor1}{rgb}{0.49400,0.18400,0.55600}%
\definecolor{mycolor2}{rgb}{0.00000,1.00000,1.00000}%
\definecolor{mycolor3}{rgb}{1.00000,0.00000,1.00000}%
\definecolor{mycolor4}{rgb}{0.63500,0.07800,0.18400}%
\begin{tikzpicture}

\begin{axis}[%
width=5.9cm,
height=4.2cm,
scale only axis,
tick label style={font=\footnotesize},
xmin=0,
xmax=100,
xlabel style={font=\footnotesize\color{white!15!black},yshift=.2cm},
xlabel={Iteration},
ymode=log,
ymin=8e-14,
ymax=9,
yminorticks=true,
ytickten = {-1, -4, -7,-10, -13},
ylabel style={font=\footnotesize\color{white!15!black},yshift=-.2cm},
ylabel={Residual norm},
axis background/.style={fill=white},
legend style={at={(axis cs:73,5)}, anchor=north east, legend cell align=left, align=left, font = \tiny, draw=white!15!black}
]
\addplot [color=myscf, mark=x, mark options={solid, myscf}]
  table[row sep=crcr]{%
0	0.295562954793329\\
1	0.24465741042241\\
2	0.173968699995449\\
3	0.118708517109271\\
4	0.0827637379041677\\
5	0.0603464862017537\\
6	0.0458080616176206\\
7	0.0354435527939858\\
8	0.0274125179348112\\
9	0.0211070051898434\\
10	0.016325119352158\\
11	0.0127782516799314\\
12	0.0101065101876658\\
13	0.00802723219163032\\
14	0.00636904014683456\\
15	0.00503393202900897\\
16	0.0039597599823655\\
17	0.00310036348939746\\
18	0.00241764618821004\\
19	0.00187899866038079\\
20	0.0014565567007467\\
21	0.00112687766328714\\
22	0.00087059080935847\\
23	0.000671949980323964\\
24	0.000518328596075856\\
25	0.000399711070097462\\
26	0.000308219846166789\\
27	0.000237699008052335\\
28	0.000183361476230417\\
29	0.000141498133301119\\
30	0.00010924269607228\\
31	8.43843714079862e-05\\
32	6.52201580813617e-05\\
33	5.04393095129828e-05\\
34	3.90334711589532e-05\\
35	3.02270797400518e-05\\
36	2.34236189289496e-05\\
37	1.81642085636958e-05\\
38	1.40957445742493e-05\\
39	1.09464109974961e-05\\
40	8.50686954337353e-06\\
41	6.61581461437948e-06\\
42	5.14888153790982e-06\\
43	4.01012870311886e-06\\
44	3.12549476722608e-06\\
45	2.43777123621387e-06\\
46	1.90273783902237e-06\\
47	1.48619028792579e-06\\
48	1.16165316844137e-06\\
49	9.08619004511081e-07\\
50	7.11191655546351e-07\\
51	5.57040508644205e-07\\
52	4.36593815719955e-07\\
53	3.42415990199722e-07\\
54	2.68726666877587e-07\\
55	2.11028945780339e-07\\
56	1.65821856093474e-07\\
57	1.30377794288128e-07\\
58	1.02570153805369e-07\\
59	8.07396821779482e-08\\
60	6.35908554914779e-08\\
61	5.01114159672674e-08\\
62	3.95098843613201e-08\\
63	3.1166957486173e-08\\
64	2.45977377527011e-08\\
65	1.94222856625709e-08\\
66	1.5342726367746e-08\\
67	1.21253452202895e-08\\
68	9.58666631641909e-09\\
69	7.58255795257784e-09\\
70	5.99971541651097e-09\\
71	4.74903908552565e-09\\
72	3.76039435106404e-09\\
73	2.97856610323729e-09\\
74	2.36004747652534e-09\\
75	1.87053538342897e-09\\
76	1.48299102753219e-09\\
77	1.17606190175738e-09\\
78	9.32901674458973e-10\\
79	7.40202279265978e-10\\
80	5.87446103473748e-10\\
81	4.66319854237941e-10\\
82	3.70245821785707e-10\\
83	2.94025476052122e-10\\
84	2.33542369921184e-10\\
85	1.85534417027152e-10\\
86	1.47416674975706e-10\\
87	1.17149198057793e-10\\
88	9.31087601240146e-11\\
89	7.40168472246063e-11\\
90	5.88483218326766e-11\\
91	4.67918512095521e-11\\
92	3.7209412196138e-11\\
93	2.95915104552263e-11\\
94	2.35374579899515e-11\\
95	1.87186829804527e-11\\
96	1.48934895193868e-11\\
97	1.18513808074344e-11\\
98	9.42738161740039e-12\\
};
\addlegendentry{SCF}

\addplot [color=black, mark=+, mark options={solid, black}]
  table[row sep=crcr]{%
0	0.295562954793329\\
1	0.24465741042241\\
2	0.173968699995449\\
3	0.118708517109271\\
4	0.0827637379041677\\
5	0.0603464862017537\\
6	0.0458080616176206\\
7	0.0354435527939858\\
8	0.0274125179348112\\
9	0.0211070051898434\\
10	0.00855976207773502\\
11	0.0085597739173349\\
12	0.00855961750340875\\
13	0.00855941808998451\\
14	0.0163250468847019\\
15	0.0343394921023706\\
16	0.0505552937329942\\
17	0.0512625852432959\\
18	0.0496320484396301\\
19	0.0440824676763179\\
20	0.00864184100181224\\
21	0.0198178691972482\\
22	0.00674238032927729\\
23	0.0250665137956577\\
24	0.0115509542003584\\
25	0.0186163224671195\\
26	0.0172640928274419\\
27	0.0137403388570465\\
28	0.0150429294372243\\
29	0.00381159561444881\\
30	0.00348203663389662\\
31	0.00283306149928192\\
32	0.00246179829994905\\
33	0.00229196756447262\\
34	0.00224748154004919\\
35	0.00221771517517079\\
36	0.00182756617984319\\
37	0.00174244811483649\\
38	0.00155878970401246\\
39	0.00147550725549642\\
40	0.00129223839936561\\
41	0.00121792267719686\\
42	0.00119247002023916\\
43	0.00115477556436341\\
44	0.000959166972421304\\
45	0.000883708756609552\\
46	0.000980836197273391\\
47	0.000944765071843617\\
48	0.00152563345315019\\
49	0.00117834859138782\\
50	0.00810740405696419\\
51	0.022872487591082\\
52	0.0164788072702629\\
53	0.00791134495422414\\
54	0.00806522952581417\\
55	0.00914150100166165\\
56	0.00141693997278858\\
57	0.00127172793448889\\
58	0.00107971932886615\\
59	0.000946836234928625\\
60	0.000840132168509053\\
61	0.000625858751380784\\
62	0.000581505390755904\\
63	0.000501273211383955\\
64	0.000424347925005537\\
65	0.0126734135315145\\
66	0.0155182958970878\\
67	0.0171320807524868\\
68	0.0305909530825353\\
69	0.00889601222518806\\
70	0.00913718778037092\\
71	0.00172727968779461\\
72	0.00176403334303202\\
73	0.00197421642318618\\
74	0.00143863118241541\\
75	0.00126596180515891\\
76	0.00109846975740863\\
77	0.000998332066079247\\
78	0.000889065161268811\\
79	0.000822163428898116\\
80	0.000775412968556423\\
81	0.000721091143887928\\
82	0.000647791392127797\\
83	0.00058205178717934\\
84	0.000544601549846215\\
85	0.000527887558393696\\
86	0.00049774553278172\\
87	0.000429838809518314\\
88	0.000401017854416035\\
89	0.000391293412280214\\
90	0.000378762739404164\\
91	0.000316682250574872\\
92	0.000270664073730986\\
93	0.0010011567488057\\
94	0.00398337250956617\\
95	0.0176148783191029\\
96	0.0550199308945358\\
97	0.104391573721066\\
98	0.521271712124632\\
99	0.129244920656899\\
100	0.201348764516199\\
101	0.126676582394572\\
102	0.141044732753788\\
103	0.039308936311785\\
104	0.0272977389772196\\
105	0.0192846457747963\\
106	0.0166308280978068\\
107	0.0134330425708526\\
108	0.0106147735771933\\
109	0.00951420637855884\\
110	0.00832049501568571\\
111	0.00732915970513874\\
112	0.0156327897158619\\
113	0.0120269209197333\\
114	0.00798967060974341\\
115	0.00626492350560606\\
116	0.0602108429667348\\
117	0.0198076273260321\\
118	0.0196016839974696\\
119	0.0185909667214227\\
};
\addlegendentry{AA(4)-$X_i$}

\addplot [color=blue, mark=o, mark options={solid, blue}]
  table[row sep=crcr]{%
0	0.295562954793329\\
1	0.24465741042241\\
2	0.173968699995449\\
3	0.118708517109271\\
4	0.0827637379041677\\
5	0.0603464862017537\\
6	0.0458080616176206\\
7	0.0354435527939858\\
8	0.0274125179348112\\
9	0.0211070051898434\\
10	0.0930096026496691\\
11	0.00175694851799568\\
12	0.00188405698314851\\
13	0.000739535936425526\\
14	0.000228914608439815\\
15	0.000255849052240036\\
16	4.46444329862445e-05\\
17	5.86383337982521e-05\\
18	7.42342499136609e-06\\
19	6.55569717486529e-06\\
20	3.1069213402354e-06\\
21	1.00655063011306e-06\\
22	2.4885482763758e-07\\
23	1.33213773690279e-07\\
24	1.93115594476504e-07\\
25	3.24813964703355e-08\\
26	2.72722338209201e-08\\
27	2.5663963897398e-08\\
28	1.65493549479853e-08\\
29	5.75263115819153e-09\\
30	3.14645356203016e-09\\
31	1.57427931471864e-09\\
32	3.57747908290138e-10\\
33	1.30871578291978e-10\\
34	1.21158454235688e-10\\
35	1.58821801033961e-10\\
36	1.00241757334632e-10\\
37	3.13201027936409e-11\\
38	4.95592329315471e-12\\
};
\addlegendentry{AA(4)-$\rho_i$}

\addplot [color=myalign, mark=o, mark options={solid, myalign}]
  table[row sep=crcr]{%
0	0.29556295479333\\
1	0.24465741042241\\
2	0.173968699995448\\
3	0.118708517109271\\
4	0.0827637379041685\\
5	0.0603464862017528\\
6	0.0458080616176201\\
7	0.035443552793984\\
8	0.0274125179348053\\
9	0.0211070051898374\\
10	0.00267843904464036\\
11	0.00169063295460186\\
12	0.000629286210127776\\
13	0.000263478894988807\\
14	0.000147678507419571\\
15	8.34124437349898e-05\\
16	2.99233767305901e-05\\
17	1.95539331984573e-05\\
18	6.14008147971697e-06\\
19	2.80733901743055e-06\\
20	1.16114134602877e-06\\
21	1.64983078612376e-07\\
22	9.93034412163772e-08\\
23	4.27628740486342e-08\\
24	1.82789795422355e-08\\
25	8.94466589898372e-09\\
26	7.19975798934941e-09\\
27	5.15611960321397e-09\\
28	2.19584816233736e-09\\
29	1.06291321885309e-09\\
30	1.980554790212e-10\\
31	1.28054374553062e-10\\
32	8.0411412350299e-11\\
33	4.61895538251761e-11\\
34	2.38377570519006e-11\\
35	1.09375313047344e-11\\
36	5.81220679530728e-12\\
};
\addlegendentry{AA(4)-$\ALIGN{X_i}$}

\addplot [color=mygrassaa, mark=o, mark options={solid, mygrassaa}]
  table[row sep=crcr]{%
0	0.295562954793329\\
1	0.24465741042241\\
2	0.173968699995449\\
3	0.118708517109271\\
4	0.0827637379041677\\
5	0.0603464862017537\\
6	0.0458080616176206\\
7	0.0354435527939858\\
8	0.0274125179348112\\
9	0.0211070051898434\\
10	0.01236792447723\\
11	0.00872492557360838\\
12	0.00393118883444084\\
13	0.00312916146689434\\
14	0.00227402529091961\\
15	0.00147890365540176\\
16	0.000644636561178727\\
17	0.000317051019230808\\
18	0.000114491323347224\\
19	5.78979447704772e-05\\
20	1.78947070110534e-05\\
21	9.37646720977421e-06\\
22	2.57404202230085e-06\\
23	1.0495962915944e-06\\
24	5.41022768417011e-07\\
25	2.90573703730006e-07\\
26	1.8893681555136e-07\\
27	9.47690712500981e-08\\
28	5.50253114967488e-08\\
29	2.55838859570613e-08\\
30	1.63078983783032e-08\\
31	8.85435940338906e-09\\
32	3.45282825339758e-09\\
33	1.90141200047973e-09\\
34	1.10373337453147e-09\\
35	6.28259510814534e-10\\
36	1.61019238760418e-10\\
37	9.97345431178926e-11\\
38	9.42047811006771e-11\\
39	6.09402554841871e-11\\
40	3.37349992785759e-11\\
41	1.8077463937901e-11\\
42	1.15143047239491e-11\\
43	6.43241725680484e-12\\
};
\addlegendentry{AA(2)-Gr}

\addplot [color=mygnt, mark=diamond, mark options={solid, mygnt}]
  table[row sep=crcr]{%
0	0.295562954793329\\
1	0.24465741042241\\
2	0.173968699995449\\
3	0.118708517109271\\
4	0.0827637379041677\\
5	0.0603464862017537\\
6	0.0458080616176206\\
7	0.0354435527939858\\
8	0.0274125179348112\\
9	0.0211070051898434\\
10	0.016325119352158\\
11	0.0084865986737871\\
12	0.00182569762328413\\
13	0.000202942287291261\\
14	5.45464364417573e-06\\
15	1.01412722925684e-08\\
16	4.006458571754e-12\\
};
\addlegendentry{AlignNT}
\legend{}

\end{axis}
\end{tikzpicture}%

%% file: Figs/tikz/residual_iter_coil.tex
%
%
\definecolor{mycolor1}{rgb}{0.49400,0.18400,0.55600}%
\definecolor{mycolor2}{rgb}{0.00000,1.00000,1.00000}%
\definecolor{mycolor3}{rgb}{1.00000,0.00000,1.00000}%
\definecolor{mycolor4}{rgb}{0.63500,0.07800,0.18400}%
\begin{tikzpicture}

\begin{axis}[%
width=5.9cm,
height=4.2cm,
scale only axis,
tick label style={font=\footnotesize},
xmin=0,
xmax=100,
xlabel style={font=\footnotesize\color{white!15!black},yshift=.2cm},
xlabel={Iteration},
ymode=log,
ymin=8e-14,
ymax=9,
yminorticks=true,
ytickten = {-1, -4, -7,-10, -13},
ylabel style={font=\footnotesize\color{white!15!black},yshift=-.2cm},
ylabel={Residual norm},
axis background/.style={fill=white},
legend style={at={(axis cs:97,10^-1)}, anchor=north east, legend cell align=left, align=left, font = \tiny, draw=white!15!black}
]
\addplot [color=myscf, mark=x, mark options={solid, myscf}]
  table[row sep=crcr]{%
0	0.939332707222954\\
1	0.775878018730737\\
2	0.52218978005057\\
3	0.34576181166399\\
4	0.239386810460419\\
5	0.179760329904193\\
6	0.135399905167854\\
7	0.100078293211898\\
8	0.0741415318293827\\
9	0.0567516444081244\\
10	0.0450400639797682\\
11	0.0365625764352043\\
12	0.0300653734711954\\
13	0.0249790487920784\\
14	0.0209894519894769\\
15	0.0178549525466922\\
16	0.0153538041075849\\
17	0.0132928735556167\\
18	0.0115332051279933\\
19	0.00999914733071855\\
20	0.00866150582882681\\
21	0.007509527051816\\
22	0.00653086085275066\\
23	0.00570589009210822\\
24	0.00501101509860637\\
25	0.00442343001493456\\
26	0.00392372331832096\\
27	0.00349622585348633\\
28	0.00312830254372508\\
29	0.00280954307359472\\
30	0.00253122731565901\\
31	0.00228606225317127\\
32	0.00206805961641786\\
33	0.00187243912387163\\
34	0.00169550098132477\\
35	0.00153446058612101\\
36	0.00138726344980023\\
37	0.00125240348993447\\
38	0.00112876241234893\\
39	0.00101547948065158\\
40	0.000911853661424721\\
41	0.000817275329362454\\
42	0.000731182319618284\\
43	0.000653034471094868\\
44	0.000582301216326472\\
45	0.000518457687592867\\
46	0.000460985863793282\\
47	0.000409378270192735\\
48	0.000363142569650261\\
49	0.000321806021987624\\
50	0.000284919250175491\\
51	0.000252059066784836\\
52	0.00022283031404673\\
53	0.000196866787798279\\
54	0.000173831374785366\\
55	0.000153415554571515\\
56	0.000135338416303661\\
57	0.000119345326783217\\
58	0.000105206367332852\\
59	9.27146351531231e-05\\
60	8.16844856247638e-05\\
61	7.19497730824557e-05\\
62	6.33621328303764e-05\\
63	5.57893341287615e-05\\
64	4.91137239715195e-05\\
65	4.32307733699647e-05\\
66	3.80477320994379e-05\\
67	3.34823931809035e-05\\
68	2.94619653429446e-05\\
69	2.59220493914297e-05\\
70	2.28057130877755e-05\\
71	2.00626581203273e-05\\
72	1.76484724502422e-05\\
73	1.55239610237558e-05\\
74	1.36545481194324e-05\\
75	1.20097446092799e-05\\
76	1.05626740905853e-05\\
77	9.28965179805785e-06\\
78	8.1698110431841e-06\\
79	7.18477211969402e-06\\
80	6.31834918376592e-06\\
81	5.55629092314082e-06\\
82	4.88605130818849e-06\\
83	4.29658706532345e-06\\
84	3.77817889498775e-06\\
85	3.32227364642635e-06\\
86	2.92134513849708e-06\\
87	2.5687713909681e-06\\
88	2.25872646768091e-06\\
89	1.98608506743829e-06\\
90	1.74633863395944e-06\\
91	1.53552122630611e-06\\
92	1.3501444766375e-06\\
93	1.18714014991149e-06\\
94	1.04380969347419e-06\\
95	9.17779734964644e-07\\
96	8.06962984289938e-07\\
97	7.0952378738195e-07\\
98	6.23847825142463e-07\\
99	5.48515472879212e-07\\
100	4.82278337986686e-07\\
101	4.24038608111815e-07\\
102	3.72830936803652e-07\\
103	3.27806447973931e-07\\
104	2.88218688935952e-07\\
105	2.53411291008006e-07\\
106	2.22807124596695e-07\\
107	1.95898691730294e-07\\
108	1.72239758855607e-07\\
109	1.51437960834289e-07\\
110	1.33148308221105e-07\\
111	1.17067455821891e-07\\
112	1.02928658562951e-07\\
113	9.04974180671595e-08\\
114	7.95675268906688e-08\\
115	6.99576587976547e-08\\
116	6.15084204707151e-08\\
117	5.40796303915709e-08\\
118	4.75480519098125e-08\\
119	4.18053457151647e-08\\
120	3.67562259315067e-08\\
121	3.2316916377689e-08\\
122	2.84137808429821e-08\\
123	2.49820547434548e-08\\
124	2.19648079977673e-08\\
125	1.93119748257489e-08\\
126	1.69795546812541e-08\\
127	1.49288352181676e-08\\
128	1.31258007113781e-08\\
129	1.1540534285403e-08\\
130	1.01467383485625e-08\\
131	8.92127671020634e-09\\
132	7.84382877904438e-09\\
133	6.89649824897689e-09\\
134	6.06359717777006e-09\\
135	5.33129070085453e-09\\
136	4.68742075576454e-09\\
137	4.12131973186419e-09\\
138	3.62358994059546e-09\\
139	3.18597365391432e-09\\
140	2.80121481557103e-09\\
141	2.46291855271076e-09\\
142	2.1654910587109e-09\\
143	1.90397280844564e-09\\
144	1.67404421421789e-09\\
145	1.47187996840228e-09\\
146	1.29413090215927e-09\\
147	1.13785011357563e-09\\
148	1.00043708163107e-09\\
149	8.79630139071352e-10\\
150	7.73407489312811e-10\\
151	6.80017744685258e-10\\
152	5.97897610432772e-10\\
153	5.25694154147222e-10\\
154	4.62218716968076e-10\\
155	4.06411335478449e-10\\
156	3.57337380637582e-10\\
157	3.14190113077125e-10\\
158	2.7625033069369e-10\\
159	2.42897739853282e-10\\
160	2.13570302744002e-10\\
161	1.87786163331221e-10\\
162	1.65112900667111e-10\\
163	1.45180425428328e-10\\
164	1.27652876671118e-10\\
165	1.12241925683343e-10\\
166	9.86958158946539e-11\\
167	8.6787809178939e-11\\
168	7.63076384965325e-11\\
169	6.71100710851292e-11\\
170	5.90068880471557e-11\\
171	5.18950879590438e-11\\
172	4.56390689081158e-11\\
173	4.01466995739716e-11\\
174	3.53026128116797e-11\\
175	3.10498140623592e-11\\
176	2.73313116663568e-11\\
177	2.40429285634068e-11\\
178	2.11688087721465e-11\\
179	1.86378513720432e-11\\
180	1.64108660000363e-11\\
181	1.44530283252278e-11\\
182	1.27448053943695e-11\\
183	1.12387071319953e-11\\
184	9.93331792629204e-12\\
};
\addlegendentry{SCF}

\addplot [color=black, mark=+, mark options={solid, black}]
  table[row sep=crcr]{%
0	0.939332707222954\\
1	0.775878018730737\\
2	0.52218978005057\\
3	0.34576181166399\\
4	0.239386810460419\\
5	0.179760329904193\\
6	0.135399905167853\\
7	0.100078293211898\\
8	0.0741415318293827\\
9	0.0567516444081243\\
10	0.0450400639797682\\
11	0.0365625764352043\\
12	0.0300653734711954\\
13	0.0249790487920783\\
14	0.0209894519894769\\
15	0.0178549525466922\\
16	0.0153538041075848\\
17	0.0132928735556167\\
18	0.0115332051279932\\
19	0.00999914733071849\\
20	0.00866150582882674\\
21	0.00750952705181597\\
22	0.00653086085275068\\
23	0.00570589009210827\\
24	0.00501101509860637\\
25	0.005392715353774\\
26	0.00555276245460875\\
27	0.00457946442462683\\
28	0.004423536334734\\
29	0.00442331468227278\\
30	0.0896223201172382\\
31	0.171256362133139\\
32	0.155046327269517\\
33	0.386352420511159\\
34	0.19987642249496\\
35	0.309464241706098\\
36	0.240213601989802\\
37	0.799199765362519\\
38	0.244793318967285\\
39	0.271433908301588\\
40	0.279961698443324\\
41	0.208736869669284\\
42	0.277499323112282\\
43	0.584834931207173\\
44	0.814533251698574\\
45	0.324636149021072\\
46	1.48321988507187\\
47	1.07060496117205\\
48	0.678512139625632\\
49	0.913467695666978\\
50	0.866404722805916\\
51	0.759331786029727\\
52	0.672475224282201\\
53	0.540804912025371\\
54	0.373003396443314\\
55	5.01253111910613\\
56	1.01804029745337\\
57	0.818009674658102\\
58	0.649401719856003\\
59	0.68914107768031\\
60	1.63082985588961\\
61	0.91089988060305\\
62	0.777545503590511\\
63	0.762208390435597\\
64	0.900708523591058\\
65	0.558650780255469\\
66	1.0620344190593\\
67	0.621757621299971\\
68	1.04225284006697\\
69	1.61610592790674\\
70	1.81939515551298\\
71	1.88171670812383\\
72	0.519580259728681\\
73	2.16084205278732\\
74	0.684928632227486\\
75	1.1319532168012\\
76	1.11284567531502\\
77	2.99422539173993\\
78	1.22893724190043\\
79	1.46907121675086\\
80	0.918390863782727\\
81	0.966330221644761\\
82	0.652143551027816\\
83	1.84877566890281\\
84	1.06202196424065\\
85	5.24448706484702\\
86	1.83760549016138\\
87	1.89406108023453\\
88	1.31059614522203\\
89	1.86978632449317\\
90	0.934054402798856\\
91	1.18726126060099\\
92	0.616192716010624\\
93	0.916131869571697\\
94	0.76776034827252\\
95	0.553126455961605\\
96	0.503411627941204\\
97	1.14572302139194\\
98	0.389200417029974\\
99	0.539152692076064\\
100	0.451772030234439\\
101	0.804396553089747\\
102	0.487050562888232\\
103	0.516354663348327\\
104	0.59459384829\\
105	0.382072586364811\\
106	1.54581802632702\\
107	0.330867362037909\\
108	0.342068280392464\\
109	0.383762777428811\\
110	1.21567210753253\\
111	0.742375209074862\\
112	1.37492550379923\\
113	1.88203682477371\\
114	1.32276420036481\\
115	0.684065629023518\\
116	0.628501334567543\\
117	0.527752951512266\\
118	0.633601199475225\\
119	0.410419774894354\\
120	0.557401380231019\\
121	0.367797499063891\\
122	0.510220869124661\\
123	0.789072998552272\\
124	0.401207387434035\\
125	0.378230606752318\\
126	0.378533560154328\\
127	0.368258831151544\\
128	0.284498097135192\\
129	0.346161889094206\\
130	0.228295723318194\\
131	0.521836616745441\\
132	0.233367170797761\\
133	0.257029037807813\\
134	0.45005074089605\\
135	0.68251994146922\\
136	0.335795803211817\\
137	0.522031950542434\\
138	1.04713095261343\\
139	2.12977944240943\\
140	0.88152637481886\\
141	0.679790864076372\\
142	0.629378596924395\\
143	1.79454114115391\\
144	1.35175264397417\\
145	0.66810268784896\\
146	0.875308784250288\\
147	0.994112893454721\\
148	0.726798380650674\\
149	0.675797937425354\\
150	1.09916500467936\\
151	3.87781450463247\\
152	1.14623360464541\\
153	0.834458834676481\\
154	0.712612917763248\\
155	1.40854100244133\\
156	0.699117243717328\\
157	0.667155259850269\\
158	1.25368605362812\\
159	0.827652746885463\\
160	0.905633936972729\\
161	0.771144985528316\\
162	3.20152163968038\\
163	0.869330187617811\\
164	0.804328206753923\\
165	0.915474293362962\\
166	0.555380431100929\\
167	0.761633634262389\\
168	0.560474710124587\\
169	0.809345013742106\\
170	1.03122345004699\\
171	0.686663084369161\\
172	0.455874754829534\\
173	2.00681102625661\\
174	0.664153547838971\\
175	1.03265407940102\\
176	1.05672841029077\\
177	1.26670441380915\\
178	1.14276978720961\\
179	1.56889445893636\\
180	0.536646386203205\\
181	0.713925628161083\\
182	1.44520078609158\\
183	0.858117760097729\\
184	0.699326399504322\\
185	1.16366213806511\\
186	1.103021005809\\
187	1.82382385636917\\
188	2.81634799529464\\
189	1.1196021217575\\
190	0.517624168290533\\
191	0.589509796506088\\
192	0.778161521362024\\
193	0.82099352170167\\
194	0.752438590897534\\
195	0.521530521515\\
196	0.606187062060052\\
197	0.45480476836962\\
198	0.478429814816231\\
199	0.510994719898059\\
};
\addlegendentry{AA(4)-$X_i$}

\addplot [color=blue, mark=o, mark options={solid, blue}]
  table[row sep=crcr]{%
0	0.939332707222954\\
1	0.775878018730737\\
2	0.52218978005057\\
3	0.34576181166399\\
4	0.239386810460419\\
5	0.179760329904193\\
6	0.135399905167853\\
7	0.100078293211898\\
8	0.0741415318293827\\
9	0.0567516444081243\\
10	0.0450400639797682\\
11	0.0365625764352043\\
12	0.0300653734711954\\
13	0.0249790487920783\\
14	0.0209894519894769\\
15	0.0178549525466922\\
16	0.0153538041075848\\
17	0.0132928735556167\\
18	0.0115332051279932\\
19	0.00999914733071849\\
20	0.00866150582882674\\
21	0.00750952705181597\\
22	0.00653086085275068\\
23	0.00570589009210827\\
24	0.00501101509860637\\
25	0.0386080110350201\\
26	0.0108368953710853\\
27	0.00790445566942571\\
28	0.00105300667154167\\
29	0.0026140395885353\\
30	0.000478396412753886\\
31	0.000260755478237408\\
32	5.13872993320747e-05\\
33	9.54676512790972e-05\\
34	3.19899404443756e-05\\
35	7.62833058687295e-06\\
36	1.26598383087129e-05\\
37	6.30436074361326e-06\\
38	2.54208820252076e-06\\
39	8.74177283028626e-07\\
40	1.05178496941788e-06\\
41	5.08376750451218e-07\\
42	1.84937247678515e-07\\
43	1.41580043804873e-07\\
44	8.98869766174184e-08\\
45	1.16773833987936e-07\\
46	2.95398534125587e-08\\
47	2.74727734073017e-08\\
48	1.15384567882773e-08\\
49	5.65676998700185e-09\\
50	3.00075572687428e-09\\
51	1.51300779230017e-09\\
52	2.00582426924874e-09\\
53	1.43570505323951e-09\\
54	7.83259745837035e-10\\
55	4.17293303060352e-10\\
56	4.36828425141268e-10\\
57	9.07993487335089e-11\\
58	5.34433747735894e-11\\
59	3.27938116864602e-11\\
60	2.36692921226e-11\\
61	1.65918854458801e-11\\
62	6.95041655598228e-12\\
};
\addlegendentry{AA(4)-$\rho_i$}

\addplot [color=myalign, mark=o, mark options={solid, myalign}]
  table[row sep=crcr]{%
0	0.939332707222953\\
1	0.775878018730734\\
2	0.522189780050574\\
3	0.34576181166402\\
4	0.239386810460415\\
5	0.179760329904197\\
6	0.135399905167883\\
7	0.100078293211852\\
8	0.0741415318293632\\
9	0.0567516444081557\\
10	0.0450400639797869\\
11	0.0365625764352153\\
12	0.0300653734712139\\
13	0.0249790487920716\\
14	0.0209894519894563\\
15	0.0178549525466948\\
16	0.0153538041075966\\
17	0.0132928735556334\\
18	0.0115332051279951\\
19	0.00999914733072333\\
20	0.00866150582882401\\
21	0.00750952705180929\\
22	0.00653086085275207\\
23	0.00570589009209725\\
24	0.00501101509859894\\
25	0.0760529439028124\\
26	0.109777831777269\\
27	0.00818820411899477\\
28	0.040983311013599\\
29	0.00245990472813928\\
30	0.0087273110757865\\
31	0.0118113399067921\\
32	0.00262570373326139\\
33	0.00392071891328479\\
34	0.000939776398478871\\
35	0.000776895961210679\\
36	7.85704190105531e-05\\
37	4.89483473441605e-05\\
38	2.60752523352547e-05\\
39	1.52541167651507e-05\\
40	1.17110823945485e-05\\
41	4.84490374469956e-06\\
42	2.72359102739857e-06\\
43	1.57477318671165e-06\\
44	9.31768955774695e-07\\
45	1.83439361993033e-07\\
46	8.21669510209386e-08\\
47	4.33437374982447e-08\\
48	2.91254676734709e-08\\
49	2.46703444026753e-08\\
50	2.02127364354861e-08\\
51	1.4629957061567e-08\\
52	8.09225666594473e-09\\
53	4.7501473481478e-09\\
54	1.54934706932386e-09\\
55	9.08396625850869e-10\\
56	4.11417684169783e-10\\
57	3.27437352448869e-10\\
58	2.03445035103763e-10\\
59	1.27703164146907e-10\\
60	1.01692924537382e-10\\
61	6.48776003550244e-11\\
62	2.71267689973564e-11\\
63	1.56290800871884e-11\\
64	8.66022926307678e-12\\
};
\addlegendentry{AA(4)-$\ALIGN{X_i}$}

\addplot [color=mygrassaa, mark=o, mark options={solid, mygrassaa}]
  table[row sep=crcr]{%
0	0.939332707222954\\
1	0.775878018730737\\
2	0.52218978005057\\
3	0.34576181166399\\
4	0.239386810460419\\
5	0.179760329904193\\
6	0.135399905167853\\
7	0.100078293211898\\
8	0.0741415318293827\\
9	0.0567516444081243\\
10	0.0450400639797682\\
11	0.0365625764352043\\
12	0.0300653734711954\\
13	0.0249790487920783\\
14	0.0209894519894769\\
15	0.0178549525466922\\
16	0.0153538041075848\\
17	0.0132928735556167\\
18	0.0115332051279932\\
19	0.00999914733071849\\
20	0.00866150582882674\\
21	0.00750952705181597\\
22	0.00653086085275068\\
23	0.00570589009210827\\
24	0.00501101509860637\\
25	0.0472861872162543\\
26	0.0369133557665361\\
27	0.0243077286884981\\
28	0.00339358864894402\\
29	0.004006149544998\\
30	0.00272795886017933\\
31	0.00145182355618793\\
32	0.000963141972401886\\
33	0.000626931339894855\\
34	0.000382406260279323\\
35	0.000249578050432455\\
36	0.000125771834792358\\
37	7.83468311350858e-05\\
38	4.87209680127145e-05\\
39	2.70726264339664e-05\\
40	1.8073684454531e-05\\
41	1.24014330382696e-05\\
42	8.31515705871774e-06\\
43	5.9074435223904e-06\\
44	3.84269603600576e-06\\
45	2.3515051513209e-06\\
46	1.67736262435263e-06\\
47	1.00885600812651e-06\\
48	6.0552135241441e-07\\
49	3.04396419513559e-07\\
50	1.91540025573079e-07\\
51	1.15108882823498e-07\\
52	8.04252784222734e-08\\
53	5.4322499747552e-08\\
54	2.97040946676448e-08\\
55	1.82105636200619e-08\\
56	1.31887201075807e-08\\
57	9.2507451489442e-09\\
58	5.66700207391725e-09\\
59	4.29701994871042e-09\\
60	2.9747236797146e-09\\
61	1.7505099055633e-09\\
62	1.09097814454633e-09\\
63	8.14732644263881e-10\\
64	4.88040020018802e-10\\
65	2.8477001086397e-10\\
66	2.06525244816436e-10\\
67	1.47914738968308e-10\\
68	9.95786390846205e-11\\
69	5.65482590181822e-11\\
70	2.58371463156479e-11\\
71	1.11894083079744e-11\\
72	5.42892836169947e-12\\
};
\addlegendentry{AA(4)-Gr}

\addplot [color=mygnt, mark=diamond, mark options={solid, mygnt}]
  table[row sep=crcr]{%
0	0.939332707222954\\
1	0.775878018730737\\
2	0.52218978005057\\
3	0.34576181166399\\
4	0.239386810460419\\
5	0.179760329904193\\
6	0.135399905167854\\
7	0.100078293211898\\
8	0.0741415318293827\\
9	0.0567516444081244\\
10	0.0450400639797682\\
11	0.0365625764352043\\
12	0.0300653734711954\\
13	0.0249790487920784\\
14	0.0209894519894769\\
15	0.0178549525466922\\
16	0.0153538041075849\\
17	0.0132928735556167\\
18	0.0115332051279933\\
19	0.00999914733071855\\
20	0.00866150582882681\\
21	0.007509527051816\\
22	0.00653086085275066\\
23	0.00570589009210822\\
24	0.00501101509860637\\
25	0.00442343001493456\\
26	0.00120223476881573\\
27	0.000169537361299914\\
28	5.58552017622082e-06\\
29	1.5643272755996e-08\\
30	1.98820589741501e-12\\
};
\addlegendentry{AlignNT}


\end{axis}
\end{tikzpicture}%

%% file: Figs/tikz/citeseer_iter.tex
%
%
\definecolor{mycolor1}{rgb}{1.00000,0.00000,1.00000}%
\definecolor{mycolor2}{rgb}{0.49400,0.18400,0.55600}%
\definecolor{mycolor3}{rgb}{0.00000,1.00000,1.00000}%
\definecolor{mycolor4}{rgb}{0.63500,0.07800,0.18400}%
\begin{tikzpicture}

\begin{axis}[%
width=3.65cm,
height=3.7cm,
scale only axis,
tick label style={font=\footnotesize},
xmin=0,
xmax=110,
xlabel style={font=\footnotesize\color{white!15!black},yshift=.2cm},
xlabel={Iteration},
ymode=log,
ymin=1.1e-16,
ymax=2e-1,
yminorticks=true,
ytickten = {-14, -11, -8, -5, -2},
ylabel style={font=\footnotesize\color{white!15!black},yshift=-.2cm},
ylabel={Residual norm},
axis background/.style={fill=white},
legend style={at={(axis cs:93,8e-2)}, anchor=north east, legend cell align=left, align=left, font = \tiny, draw=white!15!black}
]
\addplot [color=myscf, mark=x, mark options={solid, myscf}]
  table[row sep=crcr]{%
0	0.0298583785928137\\
1	0.00691394934823653\\
2	0.00325114377704337\\
3	0.00189120278022078\\
4	0.00128724787973461\\
5	0.000981728735250281\\
6	0.000812162527359635\\
7	0.000713299887622297\\
8	0.000653123014901443\\
9	0.00061442777511767\\
10	0.000587729804475816\\
11	0.000567750217464731\\
12	0.000551526759668849\\
13	0.000537385471441801\\
14	0.000524377772912909\\
15	0.000511968830711459\\
16	0.000499862354406867\\
17	0.000487900503377274\\
18	0.000476005883196061\\
19	0.000464147549637544\\
20	0.00045232092372788\\
21	0.000440535872815238\\
22	0.000428809634316773\\
23	0.000417162633307449\\
24	0.000405616038647598\\
25	0.000394190366914706\\
26	0.000382904718431768\\
27	0.0003717763939327\\
28	0.000360820739170175\\
29	0.00035005112448142\\
30	0.00033947900258307\\
31	0.000329114009959393\\
32	0.000318964090713013\\
33	0.000309035630025647\\
34	0.000299333589464853\\
35	0.000289861639509494\\
36	0.000280622286605281\\
37	0.00027161699326396\\
38	0.000262846290466587\\
39	0.000254309882096505\\
40	0.00024600674141681\\
41	0.000237935199786076\\
42	0.000230093027915117\\
43	0.000222477510031757\\
44	0.000215085511356015\\
45	0.000207913539304095\\
46	0.000200957798843047\\
47	0.00019421424241294\\
48	0.000187678614822259\\
49	0.000181346493507481\\
50	0.000175213324530258\\
51	0.000169274454666489\\
52	0.000163525159921483\\
53	0.000157960670784935\\
54	0.000152576194518815\\
55	0.000147366934751141\\
56	0.000142328108628687\\
57	0.000137454961762783\\
58	0.000132742781183899\\
59	0.000128186906503509\\
60	0.000123782739465066\\
61	0.000119525752050675\\
62	0.000115411493295351\\
63	0.000111435594947355\\
64	0.00010759377610051\\
65	0.000103881846912748\\
66	0.000100295711514354\\
67	9.68313701996466e-05\\
68	9.34849209864553e-05\\
69	9.02525606198375e-05\\
70	8.71305850884041e-05\\
71	8.41153897150878e-05\\
72	8.1203468877578e-05\\
73	7.83914154080796e-05\\
74	7.56759197167131e-05\\
75	7.30537686781655e-05\\
76	7.05218443169998e-05\\
77	6.8077122323043e-05\\
78	6.57166704248409e-05\\
79	6.34376466460263e-05\\
80	6.12372974665987e-05\\
81	5.91129559084537e-05\\
82	5.70620395624883e-05\\
83	5.50820485721882e-05\\
84	5.31705635870316e-05\\
85	5.13252436972886e-05\\
86	4.95438243603135e-05\\
87	4.78241153271288e-05\\
88	4.61639985768805e-05\\
89	4.45614262657834e-05\\
90	4.30144186961515e-05\\
91	4.15210623102577e-05\\
92	4.00795077131932e-05\\
93	3.86879677280194e-05\\
94	3.73447154860786e-05\\
95	3.60480825547682e-05\\
96	3.47964571046239e-05\\
97	3.3588282117201e-05\\
98	3.24220536348054e-05\\
99	3.1296319052904e-05\\
100	3.02096754557782e-05\\
101	2.91607679956406e-05\\
102	2.81482883153996e-05\\
103	2.71709730148964e-05\\
104	2.62276021604737e-05\\
105	2.53169978374089e-05\\
106	2.44380227448435e-05\\
107	2.35895788325709e-05\\
108	2.2770605978997e-05\\
109	2.19800807097056e-05\\
110	2.12170149556611e-05\\
111	2.04804548504088e-05\\
112	1.97694795652892e-05\\
113	1.90832001819022e-05\\
114	1.84207586008417e-05\\
115	1.7781326485844e-05\\
116	1.71641042423742e-05\\
117	1.65683200298284e-05\\
118	1.59932288063068e-05\\
119	1.54381114051383e-05\\
120	1.49022736421994e-05\\
121	1.43850454531474e-05\\
122	1.38857800596337e-05\\
123	1.34038531636824e-05\\
124	1.29386621693546e-05\\
125	1.24896254308069e-05\\
126	1.20561815259543e-05\\
127	1.16377885549455e-05\\
128	1.12339234626162e-05\\
129	1.08440813841509e-05\\
130	1.04677750132313e-05\\
131	1.01045339919144e-05\\
132	9.75390432150803e-06\\
133	9.4154477937946e-06\\
134	9.08874144185788e-06\\
135	8.77337700993329e-06\\
136	8.46896044158792e-06\\
137	8.1751113856413e-06\\
138	7.89146271919644e-06\\
139	7.61766008722523e-06\\
140	7.35336145817461e-06\\
141	7.0982366949537e-06\\
142	6.85196714086292e-06\\
143	6.61424521991032e-06\\
144	6.38477405106241e-06\\
145	6.16326707587157e-06\\
146	5.94944769910421e-06\\
147	5.74304894186367e-06\\
148	5.54381310678873e-06\\
149	5.351491454961e-06\\
150	5.16584389397947e-06\\
151	4.9866386770254e-06\\
152	4.81365211230182e-06\\
153	4.64666828267237e-06\\
154	4.48547877503957e-06\\
155	4.32988241917859e-06\\
156	4.17968503565186e-06\\
157	4.03469919252609e-06\\
158	3.89474397058717e-06\\
159	3.75964473668997e-06\\
160	3.62923292506765e-06\\
161	3.50334582621866e-06\\
162	3.38182638318811e-06\\
163	3.26452299490816e-06\\
164	3.15128932645563e-06\\
165	3.04198412584579e-06\\
166	2.93647104727074e-06\\
167	2.83461848047895e-06\\
168	2.73629938608549e-06\\
169	2.64139113663926e-06\\
170	2.54977536321709e-06\\
171	2.46133780737236e-06\\
172	2.37596817822851e-06\\
173	2.29356001455187e-06\\
174	2.21401055163126e-06\\
175	2.1372205927996e-06\\
176	2.0630943854268e-06\\
177	1.99153950119629e-06\\
178	1.92246672061872e-06\\
179	1.85578992149833e-06\\
180	1.79142597131824e-06\\
181	1.72929462336855e-06\\
182	1.66931841647313e-06\\
183	1.61142257821573e-06\\
184	1.55553493153179e-06\\
185	1.50158580456651e-06\\
186	1.44950794362724e-06\\
187	1.39923642921505e-06\\
188	1.35070859493557e-06\\
189	1.30386394924548e-06\\
190	1.25864409993626e-06\\
191	1.21499268121778e-06\\
192	1.17285528334047e-06\\
193	1.13217938468681e-06\\
194	1.09291428620317e-06\\
195	1.05501104811488e-06\\
196	1.01842242884415e-06\\
197	9.83102826051864e-07\\
198	9.49008219729814e-07\\
199	9.16096117251057e-07\\
200	8.84325500359111e-07\\
201	8.53656773991329e-07\\
202	8.24051716870582e-07\\
203	7.95473433809767e-07\\
204	7.6788630968931e-07\\
205	7.41255965023245e-07\\
206	7.15549213085177e-07\\
207	6.90734018489685e-07\\
208	6.66779457240672e-07\\
209	6.43655678176848e-07\\
210	6.21333865699619e-07\\
211	5.99786203876879e-07\\
212	5.78985841718664e-07\\
213	5.58906859729944e-07\\
214	5.39524237531419e-07\\
215	5.20813822698822e-07\\
216	5.02752300655781e-07\\
217	4.8531716555848e-07\\
218	4.68486692261569e-07\\
219	4.5223990922382e-07\\
220	4.36556572362635e-07\\
221	4.21417139812933e-07\\
222	4.06802747556497e-07\\
223	3.92695185910451e-07\\
224	3.79076876821606e-07\\
225	3.65930851947646e-07\\
226	3.53240731527265e-07\\
227	3.4099070392685e-07\\
228	3.29165505950902e-07\\
229	3.1775040381383e-07\\
230	3.06731174772432e-07\\
231	2.9609408939637e-07\\
232	2.85825894438374e-07\\
233	2.75913796357373e-07\\
234	2.66345445321482e-07\\
235	2.57108919857886e-07\\
236	2.48192711953392e-07\\
237	2.39585712758126e-07\\
238	2.31277198693247e-07\\
239	2.23256818120551e-07\\
240	2.15514578414651e-07\\
241	2.08040833544502e-07\\
242	2.00826272005216e-07\\
243	1.93861905249917e-07\\
244	1.8713905647764e-07\\
245	1.80649349812651e-07\\
246	1.74384699873698e-07\\
247	1.68337301698004e-07\\
248	1.62499621010548e-07\\
249	1.56864384835476e-07\\
250	1.51424572449102e-07\\
251	1.46173406595638e-07\\
252	1.41104345079341e-07\\
253	1.36211072591617e-07\\
254	1.31487492836763e-07\\
255	1.2692772095874e-07\\
256	1.22526076187925e-07\\
257	1.18277074764884e-07\\
258	1.14175423124911e-07\\
259	1.10216011259404e-07\\
260	1.06393906407591e-07\\
261	1.02704346862283e-07\\
262	9.91427360610314e-08\\
263	9.57046368498723e-08\\
264	9.23857659477681e-08\\
265	8.91819886358847e-08\\
266	8.60893135710279e-08\\
267	8.31038878403607e-08\\
268	8.0221992149153e-08\\
269	7.74400361783227e-08\\
270	7.47545541351918e-08\\
271	7.21622003987859e-08\\
272	6.96597454008593e-08\\
273	6.72440715509407e-08\\
274	6.49121693827007e-08\\
275	6.26611338010429e-08\\
276	6.04881604494022e-08\\
277	5.8390542233661e-08\\
278	5.63656659311601e-08\\
279	5.44110089560626e-08\\
280	5.25241361896831e-08\\
281	5.07026969726786e-08\\
282	4.89444221613162e-08\\
283	4.72471212977893e-08\\
284	4.5608679901793e-08\\
285	4.40270568022471e-08\\
286	4.25002816335435e-08\\
287	4.10264523460995e-08\\
288	3.96037328555854e-08\\
289	3.82303507585661e-08\\
290	3.6904595100749e-08\\
291	3.56248142803893e-08\\
292	3.43894139519267e-08\\
293	3.31968550809666e-08\\
294	3.20456519932215e-08\\
295	3.09343705298878e-08\\
296	2.98616262719792e-08\\
297	2.88260828239553e-08\\
298	2.78264501149581e-08\\
299	2.68614828285468e-08\\
300	2.59299788140698e-08\\
301	2.50307776296094e-08\\
302	2.41627590663383e-08\\
303	2.33248417628418e-08\\
304	2.25159818630084e-08\\
305	2.17351717021462e-08\\
306	2.09814385643489e-08\\
307	2.02538434696667e-08\\
308	1.95514799947476e-08\\
309	1.8873473146353e-08\\
310	1.821897829203e-08\\
311	1.75871800643099e-08\\
312	1.6977291394989e-08\\
313	1.63885524900185e-08\\
314	1.58202299264415e-08\\
315	1.52716156900108e-08\\
316	1.47420263421304e-08\\
317	1.42308021281742e-08\\
318	1.37373061828407e-08\\
319	1.32609237225969e-08\\
320	1.28010612775908e-08\\
321	1.2357145968821e-08\\
322	1.19286247774305e-08\\
323	1.1514963864249e-08\\
324	1.1115647905254e-08\\
325	1.07301794431006e-08\\
326	1.03580782727789e-08\\
327	9.99888084355112e-09\\
328	9.65213967105719e-09\\
329	9.31742280043002e-09\\
330	8.99431325686249e-09\\
331	8.68240851236779e-09\\
332	8.38132000645037e-09\\
333	8.09067265789734e-09\\
334	7.81010438560529e-09\\
335	7.5392656649343e-09\\
336	7.27781909351684e-09\\
337	7.02543897194112e-09\\
338	6.78181089308137e-09\\
339	6.54663135342565e-09\\
340	6.31960737047824e-09\\
341	6.10045612440735e-09\\
342	5.88890461028836e-09\\
343	5.68468928050531e-09\\
344	5.48755573067541e-09\\
345	5.29725837743342e-09\\
346	5.11356016016282e-09\\
347	4.93623222355797e-09\\
348	4.76505367148649e-09\\
349	4.59981124816567e-09\\
350	4.44029909560899e-09\\
351	4.28631850714773e-09\\
352	4.13767765598022e-09\\
353	3.99419137722606e-09\\
354	3.85568090930523e-09\\
355	3.72197370613067e-09\\
356	3.59290320000827e-09\\
357	3.46830860180924e-09\\
358	3.34803469126312e-09\\
359	3.23193163760676e-09\\
360	3.11985480430658e-09\\
361	3.01166457135067e-09\\
362	2.90722615533806e-09\\
363	2.80640945719963e-09\\
364	2.70908887644985e-09\\
365	2.61514317618612e-09\\
366	2.52445532318966e-09\\
367	2.43691234171648e-09\\
368	2.35240517412562e-09\\
369	2.27082854174284e-09\\
370	2.19208082362485e-09\\
371	2.11606391887924e-09\\
372	2.0426831241216e-09\\
373	1.97184702702009e-09\\
374	1.90346738318942e-09\\
375	1.83745900847215e-09\\
376	1.77373967085411e-09\\
377	1.71222999028163e-09\\
378	1.65285334063942e-09\\
379	1.59553575168299e-09\\
380	1.54020581927411e-09\\
381	1.48679461892243e-09\\
382	1.43523560737955e-09\\
383	1.38546456256643e-09\\
384	1.3374194726176e-09\\
385	1.29104048844775e-09\\
386	1.24626983370783e-09\\
387	1.20305173766595e-09\\
388	1.16133235474645e-09\\
389	1.12105971520939e-09\\
390	1.08218365233528e-09\\
391	1.04465572835992e-09\\
392	1.00842919828781e-09\\
393	9.73458928922672e-10\\
394	9.39701358690641e-10\\
395	9.07114433866724e-10\\
396	8.75657554199426e-10\\
397	8.45291536946031e-10\\
398	8.15978549551886e-10\\
399	7.87682079082903e-10\\
400	7.60366872040464e-10\\
401	7.33998899611903e-10\\
402	7.08545318199509e-10\\
403	6.83974412214824e-10\\
404	6.60255574568571e-10\\
405	6.3735925717174e-10\\
406	6.15256940037507e-10\\
407	5.93921084442652e-10\\
408	5.73325116658565e-10\\
409	5.53443368375431e-10\\
410	5.34251086010967e-10\\
411	5.15724349111696e-10\\
412	4.97840080587467e-10\\
413	4.80576005240186e-10\\
414	4.63910610299886e-10\\
415	4.47823138732497e-10\\
416	4.32293547108832e-10\\
417	4.1730248728948e-10\\
418	4.02831288505638e-10\\
419	3.88861923418978e-10\\
420	3.75376987961101e-10\\
421	3.62359679741118e-10\\
422	3.49793785758125e-10\\
423	3.37663648483368e-10\\
424	3.25954168422547e-10\\
425	3.14650740340568e-10\\
426	3.03739300911869e-10\\
427	2.93206240264384e-10\\
428	2.83038445279603e-10\\
429	2.73223250274856e-10\\
430	2.63748429321647e-10\\
431	2.54602166902816e-10\\
432	2.45773086836677e-10\\
433	2.37250181302491e-10\\
434	2.29022830385898e-10\\
435	2.210807877846e-10\\
436	2.13414154752728e-10\\
437	2.06013388270004e-10\\
438	1.98869264380278e-10\\
439	1.91972886191105e-10\\
440	1.85315654942995e-10\\
441	1.78889288253725e-10\\
442	1.72685774940839e-10\\
443	1.666973874155e-10\\
444	1.6091666251001e-10\\
445	1.55336402601546e-10\\
446	1.49949652795872e-10\\
447	1.44749704668111e-10\\
448	1.39730079570975e-10\\
449	1.34884525250187e-10\\
450	1.30207008671908e-10\\
451	1.25691694594341e-10\\
452	1.21332962147058e-10\\
453	1.17125385885734e-10\\
454	1.13063717048083e-10\\
455	1.09142895031859e-10\\
456	1.0535804501581e-10\\
457	1.0170444123339e-10\\
458	9.81775420579411e-11\\
459	9.47729459530573e-11\\
460	9.14864126242176e-11\\
461	8.83138506510541e-11\\
462	8.5251306039297e-11\\
463	8.22949654280429e-11\\
464	7.94411416434083e-11\\
465	7.66862883294745e-11\\
466	7.40269625607288e-11\\
467	7.14598621753731e-11\\
468	6.89817776951322e-11\\
469	6.65896315190824e-11\\
470	6.42804386284499e-11\\
471	6.20513285579282e-11\\
472	5.98995166854256e-11\\
473	5.78223200739815e-11\\
474	5.58171595950427e-11\\
475	5.38815350632141e-11\\
476	5.20130344427109e-11\\
477	5.02093289034799e-11\\
478	4.84681719038907e-11\\
479	4.67873970189409e-11\\
480	4.51649029921817e-11\\
481	4.35986792841458e-11\\
482	4.20867662465711e-11\\
483	4.06272810475745e-11\\
484	3.92184126188027e-11\\
485	3.78583962270692e-11\\
486	3.65455467120933e-11\\
487	3.5278223721175e-11\\
488	3.4054843284303e-11\\
489	3.28738929473386e-11\\
490	3.17338927775554e-11\\
491	3.06334248911666e-11\\
492	2.9571120342799e-11\\
493	2.85456548328974e-11\\
494	2.75557475194192e-11\\
495	2.66001730746616e-11\\
496	2.5677734133166e-11\\
497	2.47872815240876e-11\\
498	2.39277105814736e-11\\
499	2.3097948524093e-11\\
500	2.22969576664303e-11\\
501	2.15237425550062e-11\\
502	2.07773443827913e-11\\
503	2.00568289022224e-11\\
504	1.93612981418638e-11\\
505	1.86898857970831e-11\\
506	1.80417634995639e-11\\
507	1.74161100976833e-11\\
508	1.68121562256892e-11\\
509	1.62291449911425e-11\\
510	1.56663521576134e-11\\
511	1.51230721162027e-11\\
512	1.45986365729981e-11\\
513	1.40923879529427e-11\\
514	1.3603690172853e-11\\
515	1.31319427397569e-11\\
516	1.26765584817662e-11\\
517	1.22369610155832e-11\\
518	1.18126019759646e-11\\
519	1.14029695828238e-11\\
520	1.1007535970645e-11\\
521	1.06258196717826e-11\\
522	1.02573346930038e-11\\
523	9.9016303663157e-12\\
524	9.55826670362239e-12\\
525	9.22680203813383e-12\\
526	8.90683539337161e-12\\
527	8.59796512471882e-12\\
528	8.2998100612416e-12\\
529	8.01198888925192e-12\\
530	7.73414988303254e-12\\
531	7.46594272193631e-12\\
532	7.2070399868026e-12\\
533	6.95711395844444e-12\\
534	6.71585746459074e-12\\
535	6.48296195539972e-12\\
536	6.25814526190224e-12\\
537	6.04112739987895e-12\\
538	5.83163483757521e-12\\
539	5.62940539541112e-12\\
540	5.4341913236414e-12\\
541	5.2457473082448e-12\\
542	5.06383273286429e-12\\
543	4.88823201895248e-12\\
544	4.71871488130246e-12\\
545	4.55508177621308e-12\\
546	4.39712033867991e-12\\
547	4.24463645982344e-12\\
548	4.09744002436018e-12\\
549	3.95534688778822e-12\\
550	3.81818665328087e-12\\
551	3.68577801757243e-12\\
552	3.55796544347492e-12\\
553	3.43458120277244e-12\\
554	3.3154788456555e-12\\
555	3.20050010216157e-12\\
556	3.08951525377966e-12\\
557	2.98237893532162e-12\\
558	2.87895313068468e-12\\
559	2.77911955069575e-12\\
560	2.68274272235934e-12\\
561	2.58971427824495e-12\\
562	2.4999102159607e-12\\
563	2.41321427891072e-12\\
564	2.32952749328355e-12\\
565	2.24874649747044e-12\\
566	2.17076367279862e-12\\
567	2.09548684414021e-12\\
568	2.02281946944619e-12\\
569	1.9526719422603e-12\\
570	1.88495691002168e-12\\
571	1.81958949425874e-12\\
572	1.75648965857403e-12\\
573	1.69557820708161e-12\\
574	1.636782040587e-12\\
575	1.5800188014819e-12\\
576	1.52522892716145e-12\\
577	1.4723364248218e-12\\
578	1.42128212762336e-12\\
579	1.37199132077707e-12\\
580	1.32441445700907e-12\\
581	1.27848724203681e-12\\
582	1.23414897897355e-12\\
583	1.19135108674432e-12\\
584	1.15004031453138e-12\\
585	1.11015964026862e-12\\
586	1.07165983748241e-12\\
587	1.0344975361986e-12\\
588	9.98621242045271e-13\\
589	9.63991305862343e-13\\
590	9.30559717604142e-13\\
591	8.98290446484782e-13\\
592	8.67140722483704e-13\\
593	8.37069586752697e-13\\
594	8.08042330820228e-13\\
595	7.80020126495914e-13\\
596	7.52976195940196e-13\\
597	7.26861712917187e-13\\
598	7.01651220920819e-13\\
599	6.77320232544248e-13\\
600	6.53835293063816e-13\\
601	6.3116210677371e-13\\
602	6.09275325177005e-13\\
603	5.88142759312606e-13\\
604	5.67750380758708e-13\\
605	5.48062161730917e-13\\
606	5.29055693239367e-13\\
607	5.10709984178044e-13\\
608	4.92997622474936e-13\\
609	4.75899599176969e-13\\
610	4.59400202092366e-13\\
611	4.43462932869771e-13\\
612	4.28088780360231e-13\\
613	4.13246764058278e-13\\
614	3.98917971253623e-13\\
615	3.85078127721607e-13\\
616	3.71724660022793e-13\\
617	3.58833652086802e-13\\
618	3.46390367249925e-13\\
619	3.34377472495476e-13\\
620	3.22783889356172e-13\\
621	3.11590627798397e-13\\
622	3.00781744776852e-13\\
623	2.90355925360537e-13\\
624	2.80286368490685e-13\\
625	2.70568596506757e-13\\
626	2.61183396560309e-13\\
627	2.52125362053208e-13\\
628	2.43380927856493e-13\\
629	2.34941732907645e-13\\
630	2.26797506563481e-13\\
631	2.18933203540453e-13\\
632	2.11336284162316e-13\\
633	2.04009995474242e-13\\
634	1.96935697711018e-13\\
635	1.9010268305769e-13\\
636	1.83514314252303e-13\\
637	1.77149858275843e-13\\
638	1.71006692581163e-13\\
639	1.65074528043799e-13\\
640	1.59346907059941e-13\\
641	1.53825329386061e-13\\
642	1.48490871958866e-13\\
643	1.43340393039171e-13\\
644	1.38369604866556e-13\\
645	1.33572365206278e-13\\
646	1.28938355431934e-13\\
647	1.24467863060753e-13\\
648	1.20152506731049e-13\\
649	1.15984156112678e-13\\
650	1.11960489294839e-13\\
651	1.08079680477546e-13\\
652	1.04331853929634e-13\\
653	1.00714826963863e-13\\
654	9.72183405176039e-14\\
655	9.38500486715886e-14\\
656	9.05962349604012e-14\\
657	8.74512251092445e-14\\
658	8.44217025135948e-14\\
659	8.14931988755344e-14\\
660	7.86698426421087e-14\\
661	7.59390603448485e-14\\
662	7.3305044218289e-14\\
663	7.0764574251438e-14\\
664	6.83111012536915e-14\\
665	6.59438184637408e-14\\
666	6.36526391975499e-14\\
667	6.1446499989055e-14\\
668	5.93183843522551e-14\\
669	5.7259284683578e-14\\
670	5.52731960866568e-14\\
671	5.33598796057108e-14\\
672	5.15090020867205e-14\\
673	4.97220457623338e-14\\
674	4.79957873565083e-14\\
675	4.63303451538531e-14\\
676	4.47297461224699e-14\\
677	4.31755451832332e-14\\
678	4.16777540807238e-14\\
679	4.0232558499927e-14\\
680	3.88369038309303e-14\\
681	3.74898545258158e-14\\
682	3.61925194905172e-14\\
683	3.49364242302668e-14\\
684	3.37232075603821e-14\\
685	3.25554350019726e-14\\
686	3.14251525198936e-14\\
687	3.03356299436735e-14\\
688	2.92840684449275e-14\\
689	2.82683649067571e-14\\
690	2.72890411939399e-14\\
691	2.6339666144688e-14\\
692	2.54268375074378e-14\\
693	2.45467301758767e-14\\
694	2.36971493362516e-14\\
695	2.28735188710113e-14\\
696	2.20801103488935e-14\\
697	2.13127469243817e-14\\
698	2.05761020026818e-14\\
699	1.98641365641784e-14\\
700	1.91714032353096e-14\\
701	1.85095409098906e-14\\
702	1.78635281053174e-14\\
703	1.72493284243495e-14\\
704	1.66485312640041e-14\\
705	1.6070513465413e-14\\
706	1.55135829521469e-14\\
707	1.49788464003254e-14\\
708	1.4457787991396e-14\\
709	1.39563303402157e-14\\
710	1.34720785872969e-14\\
711	1.30072724175446e-14\\
712	1.25544458991688e-14\\
713	1.21217756677241e-14\\
714	1.16983051274507e-14\\
715	1.1296823305989e-14\\
716	1.08999734631222e-14\\
717	1.05229606837349e-14\\
718	1.01567561880946e-14\\
719	9.80610317895993e-15\\
};
\addlegendentry{SCF}

\addplot [color=black, mark=+, mark options={solid, black}]
  table[row sep=crcr]{%
0	0.0298583785928137\\
1	0.00691394934823653\\
2	0.00703332722234696\\
3	0.00822724164695911\\
4	0.00694601138222456\\
5	0.00676539958892092\\
6	0.00377601390995812\\
7	0.0043854000841559\\
8	0.0132644557756803\\
9	0.0142956038755613\\
10	0.0194621541403991\\
11	0.0196699631839548\\
12	0.036826090049191\\
13	0.0122858520483585\\
14	0.0202682753970267\\
15	0.0209325137976128\\
16	0.0200220126728228\\
17	0.00979938695949391\\
18	0.00918392529936148\\
19	0.0217887800971711\\
20	0.0215376462815349\\
21	0.0154221012018656\\
22	0.0101076080176334\\
23	0.0100802794620977\\
24	0.0148107182550049\\
25	0.00831102673761895\\
26	0.026073630035408\\
27	0.00844915438654848\\
28	0.00885346222242682\\
29	0.0101180478598908\\
30	0.0227142023333413\\
31	0.021675977162582\\
32	0.00974798928601747\\
33	0.0203386637223933\\
34	0.0189762686039977\\
35	0.0117144501219268\\
36	0.027789992057338\\
37	0.0542436669125952\\
38	0.0522300121820968\\
39	0.0258020840674312\\
40	0.0476563438593203\\
41	0.0317667229589853\\
42	0.0296155863602633\\
43	0.0394151091396785\\
44	0.0263463816440784\\
45	0.0359518052378593\\
46	0.0380702053351531\\
47	0.0284328409013097\\
48	0.0219430194048259\\
49	0.023618221312842\\
50	0.0401580679774315\\
51	0.0223227153377658\\
52	0.0263503011305991\\
53	0.0266646610606255\\
54	0.0231101890259028\\
55	0.0348358852417438\\
56	0.0219006126632093\\
57	0.0275714489883383\\
58	0.0304575982822305\\
59	0.0218950460923166\\
60	0.0217412520376837\\
61	0.0238185480599373\\
62	0.0248975540626074\\
63	0.0275284622411405\\
64	0.0346137477896085\\
65	0.0291659543520605\\
66	0.0478865665936783\\
67	0.0292294313872083\\
68	0.0482388188791972\\
69	0.0182420770548971\\
70	0.0282248923008449\\
71	0.0248933058238504\\
72	0.0301974108747007\\
73	0.025976787866579\\
74	0.0392014451607189\\
75	0.0146088162666373\\
76	0.0136460561191993\\
77	0.0178195742227578\\
78	0.0216751474684405\\
79	0.0181058277578621\\
80	0.017473441719061\\
81	0.0139312143285603\\
82	0.015623396895954\\
83	0.0128030336479359\\
84	0.0180159109677132\\
85	0.019538350615525\\
86	0.0191715362707348\\
87	0.0168413532505892\\
88	0.0150055357784006\\
89	0.0180264889341458\\
90	0.016866589017788\\
91	0.0127633311665993\\
92	0.0174191856834753\\
93	0.0155172690419887\\
94	0.0177657040693173\\
95	0.0192674409846984\\
96	0.0188411371576982\\
97	0.0234454590561026\\
98	0.0227941300988236\\
99	0.0262820101482419\\
100	0.0174034723621309\\
101	0.0247675243307078\\
102	0.0183473207576597\\
103	0.0214678655005345\\
104	0.0110648342538409\\
105	0.0127533651924497\\
106	0.0145833020803897\\
107	0.0206584697624672\\
108	0.0412840657652045\\
109	0.0471907828185381\\
110	0.0165894355771647\\
111	0.037850226016381\\
112	0.0285082633965449\\
113	0.0390960189863723\\
114	0.0246886920457098\\
115	0.0412014853310326\\
116	0.0328860992344496\\
117	0.0299153218009755\\
118	0.0277546637340904\\
119	0.0387272573951502\\
120	0.0264492762316285\\
121	0.0194816134408299\\
122	0.023698556816804\\
123	0.0179959662872622\\
124	0.0207718412766423\\
125	0.01621736496747\\
126	0.022170894059866\\
127	0.0194591998965976\\
128	0.0148752057001519\\
129	0.0138463296830081\\
130	0.0141577271354296\\
131	0.0125121445560535\\
132	0.0249658660683804\\
133	0.0174054777156712\\
134	0.0201733162763426\\
135	0.0157779141306312\\
136	0.0176242733979757\\
137	0.0181159683489291\\
138	0.0216258057669863\\
139	0.0230607843666951\\
140	0.0179559916799423\\
141	0.0172885609608883\\
142	0.0147025236338371\\
143	0.0171395653820653\\
144	0.0337774234307431\\
145	0.0587983192191488\\
146	0.0549085957610115\\
147	0.0275530493144046\\
148	0.0198125725681187\\
149	0.0225295957437037\\
150	0.0239193431721176\\
151	0.019776861752053\\
152	0.0313488898172163\\
153	0.0210749505976595\\
154	0.0326208829647557\\
155	0.0383333597433237\\
156	0.0389225607442333\\
157	0.0219620512851878\\
158	0.0223621739999232\\
159	0.0178684859625777\\
160	0.014649278929151\\
161	0.0137799478392238\\
162	0.0217187635549377\\
163	0.0165863838197711\\
164	0.0164275076894199\\
165	0.0169964673651523\\
166	0.0167106838811446\\
167	0.0145560982604014\\
168	0.0138476156854949\\
169	0.0110177885724375\\
170	0.0113136492198486\\
171	0.0126891290698966\\
172	0.0139653393780689\\
173	0.0151929292807985\\
174	0.0199346722303972\\
175	0.018401029397513\\
176	0.0173705099097654\\
177	0.0147110152218485\\
178	0.0189179766110276\\
179	0.0187883322923107\\
180	0.0100533738017922\\
181	0.0106545612841273\\
182	0.00876020684400837\\
183	0.0200483784656806\\
184	0.0243194503416446\\
185	0.00988128705317116\\
186	0.00843790295183966\\
187	0.016499702747773\\
188	0.00843665271831675\\
189	0.00710549446673169\\
190	0.00618404756430548\\
191	0.0428492488934502\\
192	0.0196816782431067\\
193	0.0137422389999176\\
194	0.0130146248036251\\
195	0.0120423920170503\\
196	0.015855443801536\\
197	0.0184923833455888\\
198	0.011340569173741\\
199	0.00867027375665\\
200	0.00883863196564314\\
201	0.0121723562843233\\
202	0.011998296403133\\
203	0.0214489860080925\\
204	0.0187461054379498\\
205	0.0214901540470397\\
206	0.0221468877700777\\
207	0.0250568493642028\\
208	0.0105352424208645\\
209	0.0121737352569254\\
210	0.0233018805813124\\
211	0.0105788493854492\\
212	0.0117342397232104\\
213	0.00787943315881869\\
214	0.0159999055218711\\
215	0.0117363007684046\\
216	0.0145948462457087\\
217	0.00986244238950464\\
218	0.0197109788591016\\
219	0.0111937328123276\\
220	0.00820673484291018\\
221	0.0272603278864807\\
222	0.00875862125198519\\
223	0.00798440977002834\\
224	0.00860180283819727\\
225	0.0096814464572469\\
226	0.0204107279088381\\
227	0.0148484364022731\\
228	0.0139663160255474\\
229	0.0143667672387193\\
230	0.0456131560517192\\
231	0.054090213554432\\
232	0.0245974430140508\\
233	0.0198360102871166\\
234	0.0193488344174791\\
235	0.0175218232407111\\
236	0.0210609187209064\\
237	0.0127128730017668\\
238	0.0111165782929566\\
239	0.0124091978108891\\
240	0.0136995138594535\\
241	0.0119271898312228\\
242	0.0138152411632363\\
243	0.0106607024947361\\
244	0.0211558649148066\\
245	0.03321082643966\\
246	0.0225814384054356\\
247	0.0110616367213795\\
248	0.0243501295161544\\
249	0.0137551737799873\\
250	0.0236996989196894\\
251	0.0143745406910308\\
252	0.0132302324546752\\
253	0.0158712826451447\\
254	0.022057332421755\\
255	0.00957606950725802\\
256	0.0158435120891655\\
257	0.0157697385290604\\
258	0.0213148222079201\\
259	0.0163803201109225\\
260	0.016507176342521\\
261	0.0158769397719163\\
262	0.0175428710616045\\
263	0.0143246780906405\\
264	0.0134217777569554\\
265	0.0125150644679718\\
266	0.0158130175402442\\
267	0.0130803348779344\\
268	0.012913817026249\\
269	0.020334806361636\\
270	0.0115625622085297\\
271	0.00931979198025513\\
272	0.018236623568248\\
273	0.0106828710049341\\
274	0.00951689234485671\\
275	0.0113866040450122\\
276	0.0172310663828915\\
277	0.0167200882757166\\
278	0.0200049215455789\\
279	0.0145296305493181\\
280	0.0208120557675343\\
281	0.017263049253107\\
282	0.0137635022904781\\
283	0.0150116192758517\\
284	0.0105570929479265\\
285	0.0104984772201068\\
286	0.0103523141307451\\
287	0.0207462107744711\\
288	0.0110813021337763\\
289	0.011293203060948\\
290	0.0145748356638071\\
291	0.0106944792894996\\
292	0.0106154100272672\\
293	0.0184441325859749\\
294	0.00980317163261027\\
295	0.00681922178227243\\
296	0.00386929519421929\\
297	0.00686262126974612\\
298	0.00730623365378069\\
299	0.00430148583558404\\
300	0.00420442857315153\\
301	0.00290828021745826\\
302	0.0028101806655397\\
303	0.00819414901118907\\
304	0.0177220469725286\\
305	0.00773583628656384\\
306	0.0584417900482416\\
307	0.0313325426197849\\
308	0.0293381871512379\\
309	0.0227047161028273\\
310	0.0502944870298194\\
311	0.049231075375051\\
312	0.0428093701127451\\
313	0.0330455687400749\\
314	0.0340083586781939\\
315	0.023212026208914\\
316	0.0303987650759534\\
317	0.0303046461133009\\
318	0.0413402507269486\\
319	0.0257552091747087\\
320	0.0235967602204703\\
321	0.0176816236480984\\
322	0.0220727936374033\\
323	0.0155400078314516\\
324	0.0179728137944999\\
325	0.0199547296676161\\
326	0.0141736966835028\\
327	0.0136167791166298\\
328	0.00957582323157429\\
329	0.00793479521104081\\
330	0.0104360923733356\\
331	0.026274428780031\\
332	0.0161048300861191\\
333	0.0136404414073936\\
334	0.0108472372877087\\
335	0.0210664077579865\\
336	0.0151975087568182\\
337	0.018252656226182\\
338	0.0217859693471019\\
339	0.0203808671057919\\
340	0.0153323749364299\\
341	0.0168017949967617\\
342	0.0160818485595909\\
343	0.017699066174123\\
344	0.018026831076415\\
345	0.0168406794515215\\
346	0.0215646926586995\\
347	0.0152013522702367\\
348	0.0133415374961992\\
349	0.00973757283499403\\
350	0.00948296280982757\\
351	0.0136185388666012\\
352	0.0497394713436386\\
353	0.0388125711964659\\
354	0.0290709431691382\\
355	0.0173760291056547\\
356	0.0397141692679894\\
357	0.0215634974917575\\
358	0.0196385276740807\\
359	0.0149092939379857\\
360	0.018684080463918\\
361	0.0208770523373258\\
362	0.016562657591299\\
363	0.0184471648872689\\
364	0.0189795223634358\\
365	0.019485722921052\\
366	0.0153986127121485\\
367	0.00990813443884068\\
368	0.0142228307021234\\
369	0.0197168015616841\\
370	0.0136456065492285\\
371	0.0169914793177903\\
372	0.0165407853499292\\
373	0.0171873762613158\\
374	0.0224126719852762\\
375	0.0105477539447029\\
376	0.00661384019569006\\
377	0.0115875496596909\\
378	0.0146671120659714\\
379	0.00891811252587501\\
380	0.00692267942343706\\
381	0.0062996363167425\\
382	0.0114700087377136\\
383	0.00717369334135311\\
384	0.0125401494083248\\
385	0.0225860814661634\\
386	0.0142135106169096\\
387	0.00799863780128012\\
388	0.0086620279062169\\
389	0.0100833353987458\\
390	0.00832016519863919\\
391	0.00750785720463444\\
392	0.0138308184178632\\
393	0.0149349998324229\\
394	0.0093787793414216\\
395	0.0174243043794457\\
396	0.0087687626520259\\
397	0.0130962348115396\\
398	0.00991353456072231\\
399	0.0206960466018654\\
400	0.0111076482485929\\
401	0.0102408601621633\\
402	0.0118881358548398\\
403	0.0235610804653682\\
404	0.0181781646244374\\
405	0.0155110383028292\\
406	0.0211336467797868\\
407	0.0225692921977141\\
408	0.0384358431462906\\
409	0.0339347071458661\\
410	0.0206224031967684\\
411	0.0179658885485052\\
412	0.0165856433155373\\
413	0.0158163276800707\\
414	0.0105317918451096\\
415	0.0155419355540064\\
416	0.0141634429047032\\
417	0.0147316657066555\\
418	0.0144300838780983\\
419	0.00752059439525524\\
420	0.00886440570024512\\
421	0.00971194625919749\\
422	0.0138432619779509\\
423	0.00871862665473319\\
424	0.0124965143452478\\
425	0.0176082445617217\\
426	0.0103718346002204\\
427	0.0103666201976113\\
428	0.0118855815286625\\
429	0.0107589180928488\\
430	0.0171449642374723\\
431	0.0246587079041752\\
432	0.0187748578910339\\
433	0.0154908229011917\\
434	0.0208665436632889\\
435	0.0179705816611942\\
436	0.0254130177938343\\
437	0.0187443024048469\\
438	0.0160264731083726\\
439	0.0209374864109834\\
440	0.0256983001012258\\
441	0.0216801175731523\\
442	0.0226720391142151\\
443	0.0200844090072798\\
444	0.0150834516473579\\
445	0.0174243366258314\\
446	0.0174897318629465\\
447	0.0199171666185184\\
448	0.0297815059872436\\
449	0.0262854788005292\\
450	0.0208379979165468\\
451	0.0226484498945158\\
452	0.0191898899927424\\
453	0.0170633054261726\\
454	0.013049578617096\\
455	0.0273428801409526\\
456	0.0130607716601436\\
457	0.0206092793990618\\
458	0.014196172521008\\
459	0.0201672039519989\\
460	0.0181839962819994\\
461	0.0130569841710321\\
462	0.0144997301226252\\
463	0.0133880243746023\\
464	0.0168703420019242\\
465	0.00781473868627876\\
466	0.0182203187987877\\
467	0.00888578980009203\\
468	0.0229525945094468\\
469	0.0186587384312939\\
470	0.0147796845471613\\
471	0.012352080851007\\
472	0.0209454632088239\\
473	0.0156477910250643\\
474	0.0177281170930076\\
475	0.021817992528409\\
476	0.0171887153727653\\
477	0.0178127200939058\\
478	0.0187440799312918\\
479	0.0239802395642033\\
480	0.0146904310944559\\
481	0.0213849833053148\\
482	0.0165247995775991\\
483	0.0187954079756004\\
484	0.0254231995530209\\
485	0.0195588970863876\\
486	0.0227035026854994\\
487	0.0196291194659189\\
488	0.0226882211370889\\
489	0.0570202679506139\\
490	0.0227785738144452\\
491	0.0160587445218142\\
492	0.0136751437668504\\
493	0.0162113436916596\\
494	0.0123724012322006\\
495	0.0142422351898747\\
496	0.0211913933889936\\
497	0.0178252112421996\\
498	0.0186277728729792\\
499	0.0199039240837956\\
500	0.0146029630737996\\
501	0.0218579745296727\\
502	0.0187427817780701\\
503	0.0208623915473819\\
504	0.0162641564829022\\
505	0.00714714279528311\\
506	0.01029485274094\\
507	0.00792678076578627\\
508	0.00966487547686353\\
509	0.00948528209553281\\
510	0.0062366070398443\\
511	0.0146006217451499\\
512	0.0218039203702377\\
513	0.0195569442825199\\
514	0.0226136503317939\\
515	0.0228385098551034\\
516	0.0334202897288014\\
517	0.0197862454128681\\
518	0.0263127304095064\\
519	0.0274258923537552\\
520	0.0179165365921224\\
521	0.0176908632894398\\
522	0.0184269499945659\\
523	0.0242035370487897\\
524	0.0288354157816975\\
525	0.0295140821133458\\
526	0.0159295927620068\\
527	0.0251484383580024\\
528	0.0228871894573101\\
529	0.0244594117134392\\
530	0.0213972886741235\\
531	0.0278972408096138\\
532	0.0318494147391678\\
533	0.0264947253572341\\
534	0.0388506882911206\\
535	0.0400682663352772\\
536	0.0216310967979779\\
537	0.021474684550911\\
538	0.0229107192959155\\
539	0.0215913161380471\\
540	0.0176698538648167\\
541	0.0166680112094166\\
542	0.0208592309117594\\
543	0.021089515072082\\
544	0.0195005104801098\\
545	0.0203994600130202\\
546	0.0224745013372269\\
547	0.0187592452133193\\
548	0.0258412198286945\\
549	0.0229644170078927\\
550	0.0207501530444935\\
551	0.0186029435364364\\
552	0.0147097035915828\\
553	0.0167139926331528\\
554	0.0123081946791221\\
555	0.0238115815717971\\
556	0.0240625671018637\\
557	0.0126893330185459\\
558	0.0199136613052968\\
559	0.0169157271307747\\
560	0.0126323400105513\\
561	0.0210289421745827\\
562	0.0101369581498892\\
563	0.0188306444952866\\
564	0.011651657898982\\
565	0.0201760102081055\\
566	0.00806054361928902\\
567	0.0178713203162624\\
568	0.0142654901938164\\
569	0.023689269144751\\
570	0.0146686117406766\\
571	0.0131838316883404\\
572	0.018621235042401\\
573	0.018183873324182\\
574	0.0185636765639708\\
575	0.00921443723322366\\
576	0.0129189125933888\\
577	0.0116935339045509\\
578	0.0149423808583736\\
579	0.0185860444297192\\
580	0.0215135217777737\\
581	0.0159312699341228\\
582	0.0153739384828228\\
583	0.0140861164849723\\
584	0.0215268154394398\\
585	0.0131074303072042\\
586	0.0210651581521964\\
587	0.0202516339128311\\
588	0.0110539979635871\\
589	0.0237617817365394\\
590	0.0177223369749728\\
591	0.0182803237927638\\
592	0.0200875596733539\\
593	0.0197809346659745\\
594	0.013875734128543\\
595	0.0195153632015019\\
596	0.0169356400010865\\
597	0.00932480808467753\\
598	0.0103342173615313\\
599	0.011062760886983\\
600	0.0157201914716533\\
601	0.0163534135413644\\
602	0.0169374166867936\\
603	0.0160702335614922\\
604	0.0174359027252144\\
605	0.0115733732180039\\
606	0.0115372540076984\\
607	0.00960833480641308\\
608	0.0195456136066727\\
609	0.0201293338175764\\
610	0.0165909566843625\\
611	0.0169658096769439\\
612	0.0194740673898102\\
613	0.0233523476318579\\
614	0.0216310553254285\\
615	0.0218023544010086\\
616	0.020514151993013\\
617	0.0190516939159599\\
618	0.0218835026157539\\
619	0.0156835069170738\\
620	0.0174008378770355\\
621	0.0159115650760287\\
622	0.0143165172553608\\
623	0.0182962863885978\\
624	0.0129472860833688\\
625	0.0101251916676568\\
626	0.0207306132163794\\
627	0.0192153504696795\\
628	0.0172693531455379\\
629	0.0157001726655381\\
630	0.0167580153083539\\
631	0.02581323550975\\
632	0.0178600199736853\\
633	0.0163402757992929\\
634	0.0373871694553487\\
635	0.00709828928131109\\
636	0.0144745748272472\\
637	0.00968538245961986\\
638	0.0341079032697181\\
639	0.0293908664987973\\
640	0.012358779071415\\
641	0.0315725216950817\\
642	0.0305988427846227\\
643	0.0476935716453668\\
644	0.0391128424197245\\
645	0.033928828577668\\
646	0.0267804748212915\\
647	0.0276148625115574\\
648	0.0331764562787935\\
649	0.0331024121297855\\
650	0.0224010664482308\\
651	0.0235649318201026\\
652	0.0405323239905292\\
653	0.0307278547722302\\
654	0.0461944215246737\\
655	0.0294027113803754\\
656	0.0322185864080428\\
657	0.0199516853374524\\
658	0.0230221629343936\\
659	0.0216096523925272\\
660	0.0202857082231338\\
661	0.0186995370497817\\
662	0.0238414063321467\\
663	0.0180372535291686\\
664	0.014815993911099\\
665	0.0169593250917193\\
666	0.0182963011189976\\
667	0.0188142532285598\\
668	0.0224671985240397\\
669	0.0146184731849492\\
670	0.0241109236281467\\
671	0.0245982632830846\\
672	0.0204688446982117\\
673	0.0203568592141964\\
674	0.0190599055924618\\
675	0.0147458506866349\\
676	0.0175591106658226\\
677	0.0191386432484946\\
678	0.0183863690263796\\
679	0.00964550287198259\\
680	0.0141025545944855\\
681	0.0131758326449471\\
682	0.009964931024742\\
683	0.0106231063144138\\
684	0.0100295379089189\\
685	0.00936837112893026\\
686	0.0137352058317919\\
687	0.00709686148523814\\
688	0.00651778894837729\\
689	0.00580211082899315\\
690	0.00827463267463308\\
691	0.0144411988489005\\
692	0.0169200136058626\\
693	0.0273237175602023\\
694	0.0368718110589912\\
695	0.0318387963999334\\
696	0.0552084663313054\\
697	0.0177563661120538\\
698	0.0216377350686979\\
699	0.0188589647141767\\
700	0.0190961005849836\\
701	0.016270179600943\\
702	0.0106655648358283\\
703	0.0195871623843602\\
704	0.0172384046470511\\
705	0.0144086017838298\\
706	0.0255970636739927\\
707	0.0170745709534841\\
708	0.0159823139733246\\
709	0.0164140738523109\\
710	0.0244160981503383\\
711	0.0259196654284021\\
712	0.0227191019145469\\
713	0.0119321640085956\\
714	0.0133747360292954\\
715	0.0157888307927473\\
716	0.0201058602250807\\
717	0.0184197799229802\\
718	0.0205747420436903\\
719	0.0139957339919126\\
720	0.0179264218730295\\
721	0.0129302981279324\\
722	0.0155070966964172\\
723	0.0185753226938148\\
724	0.0170269234784657\\
725	0.0111594354218862\\
726	0.00922215606840823\\
727	0.00870566558306153\\
728	0.00994400857668187\\
729	0.00872110777550862\\
730	0.0161496695201742\\
731	0.0187751562387921\\
732	0.0129020530025897\\
733	0.00875699195851434\\
734	0.0210339020906599\\
735	0.020154263484199\\
736	0.024493630411324\\
737	0.0199918507910258\\
738	0.0214507050568795\\
739	0.0224544318606723\\
740	0.0183997704350701\\
741	0.0156993397526923\\
742	0.0235611499578458\\
743	0.043460849159282\\
744	0.0530628219177137\\
745	0.0187643372501662\\
746	0.0115207608088372\\
747	0.0179547889740276\\
748	0.039810556619708\\
749	0.0221325768857735\\
};
\addlegendentry{AA(4)-$X_i$}

\addplot [color=myalign, mark=o, mark options={solid, myalign}]
  table[row sep=crcr]{%
0	0.0298583785928137\\
1	0.00691394934823653\\
2	0.00273575596213325\\
3	0.00130557749936355\\
4	0.000858511908241395\\
5	0.000693999367187228\\
6	0.00062914413805795\\
7	0.000532904027594333\\
8	0.000809575346770122\\
9	0.000914086021442422\\
10	0.000981380108714267\\
11	0.00155918135872962\\
12	0.000286623887047403\\
13	0.00123588572861402\\
14	0.00100806457148617\\
15	0.000101646949464084\\
16	4.07386279817153e-05\\
17	4.10140444313258e-05\\
18	2.25186633832672e-05\\
19	3.69343041565974e-05\\
20	1.22373823983838e-05\\
21	8.07283415016821e-06\\
22	4.15705202690127e-06\\
23	2.44585101252136e-06\\
24	1.93304655916737e-06\\
25	1.62151308304847e-06\\
26	1.22873711248277e-06\\
27	7.58436046937501e-07\\
28	3.78372541862724e-07\\
29	1.79542566923068e-07\\
30	7.61705892851578e-08\\
31	4.53802043535772e-08\\
32	1.99132238660885e-08\\
33	1.16072555053835e-08\\
34	6.62671668286089e-09\\
35	4.96739941898522e-09\\
36	4.09220377282137e-09\\
37	3.24036414401467e-09\\
38	2.27076683923174e-09\\
39	1.55890579006137e-09\\
40	1.14534913751581e-09\\
41	9.43295301635863e-10\\
42	5.65359741761929e-10\\
43	2.65193208739515e-10\\
44	1.49893341269161e-10\\
45	1.12141635875591e-10\\
46	9.30148343229162e-11\\
47	7.69811820803385e-11\\
48	5.3438932043945e-11\\
49	3.76024686528912e-11\\
50	3.13530136534147e-11\\
51	2.38069264937755e-11\\
52	1.87704006564663e-11\\
53	1.54691453030084e-11\\
54	1.20119527386776e-11\\
55	8.95628069645944e-12\\
56	5.93485988770909e-12\\
57	4.10925704304756e-12\\
58	3.29570292846206e-12\\
59	2.87247937678511e-12\\
60	2.19572535185854e-12\\
61	1.8218884151856e-12\\
62	1.43157577774388e-12\\
63	9.94821348152285e-13\\
64	6.65725657035183e-13\\
65	4.36779118355253e-13\\
66	3.26869124659099e-13\\
67	2.74986404624836e-13\\
68	2.0443326905853e-13\\
69	1.361486851057e-13\\
70	1.05720141187389e-13\\
71	7.8728660311299e-14\\
72	5.79813607621555e-14\\
73	4.02472333384106e-14\\
74	2.60738464285406e-14\\
75	2.08711193159579e-14\\
76	1.59562400032277e-14\\
77	7.97174782611626e-15\\
};
\addlegendentry{AA(4)-$\ALIGN{X_i}$}

\addplot [color=mygrassaa, mark=o, mark options={solid, mygrassaa}]
  table[row sep=crcr]{%
0	0.0298583785928137\\
1	0.00691394934823653\\
2	0.00417113618527657\\
3	0.001400420272979\\
4	0.00112561972604951\\
5	0.000654710863870351\\
6	0.000577814092889828\\
7	0.000538830309500599\\
8	0.000475310128649527\\
9	0.000310102187854579\\
10	0.000144848864711353\\
11	5.75232843267e-05\\
12	4.37566174277152e-05\\
13	3.11693029715212e-05\\
14	2.29772314214601e-05\\
15	1.78644619665345e-05\\
16	9.42016661932596e-06\\
17	4.02349817467909e-06\\
18	2.05714966296973e-06\\
19	9.85482936317182e-07\\
20	5.63460413496433e-07\\
21	3.233764714249e-07\\
22	2.47500712349495e-07\\
23	2.19463590669659e-07\\
24	1.81505184230062e-07\\
25	1.42868896346026e-07\\
26	1.00215397672743e-07\\
27	7.30034557673001e-08\\
28	3.56578931825789e-08\\
29	1.93104094775477e-08\\
30	1.43314705899211e-08\\
31	1.16988098111579e-08\\
32	1.00363747337296e-08\\
33	8.12656302977217e-09\\
34	6.5564329244981e-09\\
35	4.93513719278168e-09\\
36	3.52868508595399e-09\\
37	2.42247661396327e-09\\
38	1.78331715859215e-09\\
39	1.34261181164791e-09\\
40	9.57496160748986e-10\\
41	6.36866681675104e-10\\
42	4.58589560561412e-10\\
43	4.12601330602325e-10\\
44	3.63670983519071e-10\\
45	2.97799904283496e-10\\
46	2.05800077161226e-10\\
47	1.4248480632639e-10\\
48	1.01159622213237e-10\\
49	7.7369448649794e-11\\
50	5.72718847893016e-11\\
51	4.62202448229797e-11\\
52	3.89602162943659e-11\\
53	2.96322571517031e-11\\
54	2.25110051224456e-11\\
55	1.84591225078559e-11\\
56	1.63781821300344e-11\\
57	1.37489502003989e-11\\
58	9.93101579281484e-12\\
59	6.33392848619039e-12\\
60	4.63880479784421e-12\\
61	3.84284537758575e-12\\
62	3.22477939186411e-12\\
63	2.58775649750958e-12\\
64	1.9959875132466e-12\\
65	1.63647889955153e-12\\
66	1.29694010729437e-12\\
67	1.0599051965964e-12\\
68	8.96863104419904e-13\\
69	7.45216373096504e-13\\
70	5.7739078388343e-13\\
71	4.42117391985176e-13\\
72	3.45889308782183e-13\\
73	2.94022588033348e-13\\
74	2.51803867395426e-13\\
75	1.87013952892102e-13\\
76	1.37199336093788e-13\\
77	1.17152065859856e-13\\
78	1.06179962877014e-13\\
79	9.38145863871243e-14\\
80	7.35777843719159e-14\\
81	5.32165205033488e-14\\
82	4.38302129820394e-14\\
83	3.89204131424074e-14\\
84	3.21203196384954e-14\\
85	2.37695199322388e-14\\
86	1.72468490742471e-14\\
87	1.43013104851185e-14\\
88	1.27252888196456e-14\\
89	1.12004864169689e-14\\
90	9.32403946621811e-15\\
};
\addlegendentry{AA(4)-Gr}

\addplot [color=mygnt, mark=diamond, mark options={solid, mygnt}]
  table[row sep=crcr]{%
0	0.0298583785928137\\
1	0.00691394934823653\\
2	0.00325114377704337\\
3	0.00064059006326909\\
4	0.000304232808501008\\
5	4.60494380803669e-05\\
6	4.20789532395008e-06\\
7	5.208017380764e-08\\
8	1.85444019862073e-11\\
9	1.96753225833899e-16\\
};
\addlegendentry{AlignNT}
\legend{}
\end{axis}
\end{tikzpicture}%

%% file: Figs/tikz/cora_iter.tex
%
%
\definecolor{mycolor1}{rgb}{1.00000,0.00000,1.00000}%
\definecolor{mycolor2}{rgb}{0.49400,0.18400,0.55600}%
\definecolor{mycolor3}{rgb}{0.00000,1.00000,1.00000}%
\definecolor{mycolor4}{rgb}{0.63500,0.07800,0.18400}%
\begin{tikzpicture}

\begin{axis}[%
width=3.65cm,
height=3.7cm,
scale only axis,
tick label style={font=\footnotesize},
xmin=0,
xmax=110,
xtick = {0, 25, 50, 75, 100},
xlabel style={font=\footnotesize\color{white!15!black},yshift=.2cm},
xlabel={Iteration},
ymode=log,
ymin=1.1e-16,
ymax=2e-1,
yminorticks=true,
ytickten = {-14, -11, -8, -5, -2},
ylabel style={font=\footnotesize\color{white!15!black},yshift=-.2cm},
ylabel={Residual norm},
axis background/.style={fill=white},
legend style={at={(axis cs:1,5e-16)}, anchor=south west, legend cell align=left, align=left, font = \tiny, draw=white!15!black}
]
\addplot [color=myscf, mark=x, mark options={solid, myscf}]
  table[row sep=crcr]{%
0	0.0307613903451636\\
1	0.0217693601554288\\
2	0.0113544570677127\\
3	0.0056010740058595\\
4	0.00402895938457537\\
5	0.00316367157048876\\
6	0.00253159170867862\\
7	0.00203660694775321\\
8	0.00164253661481789\\
9	0.00132725542324835\\
10	0.00107435664180063\\
11	0.000871066450191869\\
12	0.000707336921719126\\
13	0.000575241377455123\\
14	0.000468511595019705\\
15	0.000382176418468468\\
16	0.000312280528451849\\
17	0.000255666383954162\\
18	0.000209805178226085\\
19	0.000172665427500159\\
20	0.000142610325559272\\
21	0.00011831709136202\\
22	9.87132037739538e-05\\
23	8.29257241726856e-05\\
24	7.02409335969324e-05\\
25	6.00723097804127e-05\\
26	5.19354491345717e-05\\
27	4.54288712781446e-05\\
28	4.02197154189757e-05\\
29	3.60332228694816e-05\\
30	3.26447802552285e-05\\
31	2.98733687885908e-05\\
32	2.75755864432314e-05\\
33	2.5639859515949e-05\\
34	2.39808428852745e-05\\
35	2.25342047916996e-05\\
36	2.12520082905646e-05\\
37	2.00988170733648e-05\\
38	1.90485484222672e-05\\
39	1.80820161686808e-05\\
40	1.71850635705422e-05\\
41	1.63471735728871e-05\\
42	1.55604503921377e-05\\
43	1.48188817598867e-05\\
44	1.41178089208948e-05\\
45	1.34535481027582e-05\\
46	1.28231212119805e-05\\
47	1.22240646682765e-05\\
48	1.1654293819229e-05\\
49	1.11120067312779e-05\\
50	1.05956158012453e-05\\
51	1.01036989903797e-05\\
52	9.63496488645417e-06\\
53	9.18822750914249e-06\\
54	8.76238798426657e-06\\
55	8.35642106658985e-06\\
56	7.96936509211324e-06\\
57	7.60031436333485e-06\\
58	7.24841326759546e-06\\
59	6.91285163671653e-06\\
60	6.5928610019859e-06\\
61	6.28771150087472e-06\\
62	5.99670926354041e-06\\
63	5.71919415746647e-06\\
64	5.45453780387786e-06\\
65	5.20214180431565e-06\\
66	4.96143613320548e-06\\
67	4.73187766445695e-06\\
68	4.51294880883224e-06\\
69	4.30415624483403e-06\\
70	4.105029730309e-06\\
71	3.91512098494292e-06\\
72	3.7340026360421e-06\\
73	3.56126722166252e-06\\
74	3.39652624615155e-06\\
75	3.23940928415758e-06\\
76	3.08956312966465e-06\\
77	2.9466509871642e-06\\
78	2.81035170237359e-06\\
79	2.68035903022727e-06\\
80	2.5563809380741e-06\\
81	2.43813894217796e-06\\
82	2.32536747580313e-06\\
83	2.21781328724777e-06\\
84	2.11523486634522e-06\\
85	2.01740189797564e-06\\
86	1.92409474132144e-06\\
87	1.83510393355634e-06\\
88	1.75022971682003e-06\\
89	1.66928158734354e-06\\
90	1.5920778656764e-06\\
91	1.51844528698577e-06\\
92	1.44821861050265e-06\\
93	1.38124024718496e-06\\
94	1.31735990473638e-06\\
95	1.25643424918592e-06\\
96	1.19832658221181e-06\\
97	1.14290653348803e-06\\
98	1.09004976735734e-06\\
99	1.03963770312568e-06\\
100	9.91557248390502e-07\\
101	9.45700544728584e-07\\
102	9.01964725217252e-07\\
103	8.60251683215142e-07\\
104	8.20467851884444e-07\\
105	7.82523993918755e-07\\
106	7.46335001066582e-07\\
107	7.11819702932954e-07\\
108	6.78900684648897e-07\\
109	6.47504113012817e-07\\
110	6.17559570669043e-07\\
111	5.88999897998698e-07\\
112	5.61761042324792e-07\\
113	5.35781914103013e-07\\
114	5.11004249782374e-07\\
115	4.87372481029938e-07\\
116	4.64833609989902e-07\\
117	4.4333709033748e-07\\
118	4.22834713840953e-07\\
119	4.03280502173196e-07\\
120	3.84630603729769e-07\\
121	3.66843195236556e-07\\
122	3.49878387903634e-07\\
123	3.336981379163e-07\\
124	3.18266161089501e-07\\
125	3.03547851461523e-07\\
126	2.89510203651869e-07\\
127	2.76121738819906e-07\\
128	2.63352434045626e-07\\
129	2.51173654986929e-07\\
130	2.39558091634499e-07\\
131	2.28479697066511e-07\\
132	2.17913629014892e-07\\
133	2.07836194143535e-07\\
134	1.98224794906227e-07\\
135	1.89057878861571e-07\\
136	1.80314890336315e-07\\
137	1.7197622431541e-07\\
138	1.64023182475689e-07\\
139	1.56437931249761e-07\\
140	1.49203461827977e-07\\
141	1.42303552015618e-07\\
142	1.35722729844707e-07\\
143	1.2944623888311e-07\\
144	1.23460005137547e-07\\
145	1.17750605489178e-07\\
146	1.12305237598485e-07\\
147	1.07111691186052e-07\\
148	1.02158320652473e-07\\
149	9.74340189691916e-08\\
150	9.29281927578842e-08\\
151	8.86307385533895e-08\\
152	8.45320201265394e-08\\
153	8.06228468896308e-08\\
154	7.68944532803811e-08\\
155	7.33384791064929e-08\\
156	6.99469508025201e-08\\
157	6.67122635411211e-08\\
158	6.36271641909579e-08\\
159	6.06847350380973e-08\\
160	5.78783782869606e-08\\
161	5.52018012620514e-08\\
162	5.26490022892884e-08\\
163	5.02142572508674e-08\\
164	4.78921067368082e-08\\
165	4.56773438141698e-08\\
166	4.35650023438999e-08\\
167	4.15503458487948e-08\\
168	3.9628856888536e-08\\
169	3.77962269321793e-08\\
170	3.60483467019622e-08\\
171	3.43812969518866e-08\\
172	3.27913396807831e-08\\
173	3.12749097522638e-08\\
174	2.98286069043101e-08\\
175	2.84491881143874e-08\\
176	2.71335603331062e-08\\
177	2.5878773555253e-08\\
178	2.46820141936204e-08\\
179	2.35405987764991e-08\\
180	2.24519679314642e-08\\
181	2.14136806404307e-08\\
182	2.04234087755942e-08\\
183	1.94789318666576e-08\\
184	1.8578132135933e-08\\
185	1.7718989735921e-08\\
186	1.68995782288154e-08\\
187	1.61180602660495e-08\\
188	1.53726834671524e-08\\
189	1.46617764893654e-08\\
190	1.39837452809518e-08\\
191	1.33370695044128e-08\\
192	1.2720299136614e-08\\
193	1.21320512088348e-08\\
194	1.15710066980887e-08\\
195	1.10359075924596e-08\\
196	1.05255540499823e-08\\
197	1.00388017146685e-08\\
198	9.57455915365205e-09\\
199	9.13178540741777e-09\\
200	8.70948765484302e-09\\
201	8.3067189818696e-09\\
202	7.92257627499118e-09\\
203	7.55619817843806e-09\\
204	7.20676317042683e-09\\
205	6.87348772068597e-09\\
206	6.55562453377863e-09\\
207	6.25246087052342e-09\\
208	5.96331695788491e-09\\
209	5.68754445231486e-09\\
210	5.42452499564919e-09\\
211	5.17366882563732e-09\\
212	4.93441345424955e-09\\
213	4.70622240228381e-09\\
214	4.48858400491729e-09\\
215	4.28101025459961e-09\\
216	4.08303571687301e-09\\
217	3.89421647304327e-09\\
218	3.71412914268032e-09\\
219	3.54236991534333e-09\\
220	3.37855366376943e-09\\
221	3.22231306484278e-09\\
222	3.0732977817332e-09\\
223	2.93117368661546e-09\\
224	2.79562209422071e-09\\
225	2.66633906089205e-09\\
226	2.54303469555008e-09\\
227	2.42543251838999e-09\\
228	2.31326883372736e-09\\
229	2.20629213733351e-09\\
230	2.10426255685305e-09\\
231	2.00695131762189e-09\\
232	1.91414021875899e-09\\
233	1.82562115166932e-09\\
234	1.74119563291527e-09\\
235	1.66067435659316e-09\\
236	1.58387677256524e-09\\
237	1.51063067878377e-09\\
238	1.44077183774395e-09\\
239	1.37414360414413e-09\\
240	1.3105965812509e-09\\
241	1.24998828035319e-09\\
242	1.19218279965996e-09\\
243	1.13705052237828e-09\\
244	1.0844678279601e-09\\
245	1.03431680875723e-09\\
246	9.86485016810482e-10\\
247	9.40865195137117e-10\\
248	8.97355056356882e-10\\
249	8.55857033586634e-10\\
250	8.16278078907308e-10\\
251	7.78529445276467e-10\\
252	7.42526489948717e-10\\
253	7.08188487445059e-10\\
254	6.75438438843885e-10\\
255	6.44202910761024e-10\\
256	6.14411862405416e-10\\
257	5.85998496518693e-10\\
258	5.58899102390933e-10\\
259	5.33052913831229e-10\\
260	5.0840197636195e-10\\
261	4.84891018607694e-10\\
262	4.62467321727957e-10\\
263	4.41080604301406e-10\\
264	4.20682911435156e-10\\
265	4.01228505775932e-10\\
266	3.82673766850467e-10\\
267	3.64977087291728e-10\\
268	3.48098788464098e-10\\
269	3.32001023164722e-10\\
270	3.16647695081431e-10\\
271	3.02004378606728e-10\\
272	2.88038241401272e-10\\
273	2.74717965693838e-10\\
274	2.62013682099065e-10\\
275	2.49896904881215e-10\\
276	2.38340467365716e-10\\
277	2.27318455798419e-10\\
278	2.16806155363268e-10\\
279	2.0677999423866e-10\\
280	1.97217491170612e-10\\
281	1.88097203540702e-10\\
282	1.79398682777073e-10\\
283	1.71102423615414e-10\\
284	1.63189822213121e-10\\
285	1.55643138984004e-10\\
286	1.48445450237701e-10\\
287	1.4158061699197e-10\\
288	1.35033245879698e-10\\
289	1.28788659022175e-10\\
290	1.22832850145411e-10\\
291	1.17152468834302e-10\\
292	1.11734775074638e-10\\
293	1.06567619773087e-10\\
294	1.01639418585183e-10\\
295	9.6939121428631e-11\\
296	9.24561875768016e-11\\
297	8.81805687648032e-11\\
298	8.41026758225284e-11\\
299	8.02133638122189e-11\\
300	7.65039114783491e-11\\
301	7.29660007879022e-11\\
302	6.95917023620464e-11\\
303	6.63734468286361e-11\\
304	6.33040177440846e-11\\
305	6.03765343304863e-11\\
306	5.75844313595231e-11\\
307	5.49214503018406e-11\\
308	5.23816170105979e-11\\
309	4.99592399196694e-11\\
310	4.76488817601339e-11\\
311	4.54453692660159e-11\\
312	4.33437553074881e-11\\
313	4.13393308550552e-11\\
314	3.94276009414853e-11\\
315	3.76042788707888e-11\\
316	3.58652750746455e-11\\
317	3.42066916424239e-11\\
318	3.26248074797545e-11\\
319	3.11160791799545e-11\\
320	2.96771232900335e-11\\
321	2.83047107199933e-11\\
322	2.69957632808774e-11\\
323	2.57473482447683e-11\\
324	2.45566690310659e-11\\
325	2.34210503446527e-11\\
326	2.23379482111352e-11\\
327	2.13049327486063e-11\\
328	2.03196898837197e-11\\
329	1.93800105405841e-11\\
330	1.8483785531156e-11\\
331	1.76290059016738e-11\\
332	1.68137568327977e-11\\
333	1.60362078984896e-11\\
334	1.52946158892666e-11\\
335	1.45873188127717e-11\\
336	1.39127314196177e-11\\
337	1.3269340050781e-11\\
338	1.26557017367619e-11\\
339	1.2070442644822e-11\\
340	1.15122450404803e-11\\
341	1.09798645793011e-11\\
342	1.047210240979e-11\\
343	9.98782238298719e-12\\
344	9.52593851937634e-12\\
345	9.08541226562202e-12\\
346	8.6652593948062e-12\\
347	8.26453713963481e-12\\
348	7.88234572463668e-12\\
349	7.5178276668468e-12\\
350	7.17016714572088e-12\\
351	6.83858346918262e-12\\
352	6.52233422281132e-12\\
353	6.2207112867684e-12\\
354	5.93303384025639e-12\\
355	5.65866309546495e-12\\
356	5.39697985940063e-12\\
357	5.14739636510885e-12\\
358	4.90935683186063e-12\\
359	4.68232519876434e-12\\
360	4.46579195451474e-12\\
361	4.25927029756186e-12\\
362	4.06230217736914e-12\\
363	3.87444262775232e-12\\
364	3.69526787135709e-12\\
365	3.52438029248578e-12\\
366	3.36139503166297e-12\\
367	3.20594832365857e-12\\
368	3.05769120167784e-12\\
369	2.91628941308741e-12\\
370	2.78142673874691e-12\\
371	2.65279976549283e-12\\
372	2.53012134561171e-12\\
373	2.41311626468885e-12\\
374	2.30152376743913e-12\\
375	2.19509004387377e-12\\
376	2.09357763555963e-12\\
377	1.9967608134604e-12\\
378	1.90442140218109e-12\\
379	1.8163504281946e-12\\
380	1.73235515242594e-12\\
381	1.65224092762618e-12\\
382	1.57583400439486e-12\\
383	1.50296016984227e-12\\
384	1.43345729609451e-12\\
385	1.36716651399325e-12\\
386	1.30394212795781e-12\\
387	1.24364116274713e-12\\
388	1.18612889690902e-12\\
389	1.13127697967895e-12\\
390	1.0789614985473e-12\\
391	1.02906338666112e-12\\
392	9.8147553658954e-13\\
393	9.36085809213074e-13\\
394	8.92797304086706e-13\\
395	8.51509741807921e-13\\
396	8.12129102838962e-13\\
397	7.74571894727875e-13\\
398	7.38752692974657e-13\\
399	7.04587706482624e-13\\
400	6.72005766351948e-13\\
401	6.40926874926939e-13\\
402	6.11285843815554e-13\\
403	5.83017394421744e-13\\
404	5.56054808365116e-13\\
405	5.30340667304558e-13\\
406	5.05815014312214e-13\\
407	4.82425130383318e-13\\
408	4.60116137724079e-13\\
409	4.38838845943947e-13\\
410	4.18546421255003e-13\\
411	3.99190878807802e-13\\
412	3.80731157026591e-13\\
413	3.63124518861239e-13\\
414	3.46333965160308e-13\\
415	3.30317074779767e-13\\
416	3.15040755777273e-13\\
417	3.00471798129827e-13\\
418	2.86577940903378e-13\\
419	2.73323665344483e-13\\
420	2.60684732018117e-13\\
421	2.48627758216663e-13\\
422	2.37129593571011e-13\\
423	2.26163984760053e-13\\
424	2.15704674892912e-13\\
425	2.05729817722178e-13\\
426	1.96217056578218e-13\\
427	1.87143377435981e-13\\
428	1.78490410646667e-13\\
429	1.70234549593755e-13\\
430	1.62362684073404e-13\\
431	1.5485293492576e-13\\
432	1.47692075108333e-13\\
433	1.40863992356402e-13\\
434	1.34351433973012e-13\\
435	1.28137368367431e-13\\
436	1.22211265723964e-13\\
437	1.16559736827932e-13\\
438	1.11170476782714e-13\\
439	1.06029729150524e-13\\
440	1.01127200275492e-13\\
441	9.64519344648725e-14\\
442	9.19916251961512e-14\\
443	8.77379599448886e-14\\
444	8.36789659211277e-14\\
445	7.98094783091025e-14\\
446	7.61178576523069e-14\\
447	7.25975300828159e-14\\
448	6.9239497595489e-14\\
449	6.60378042721679e-14\\
450	6.29874441856306e-14\\
451	6.00746915792026e-14\\
452	5.72965002526636e-14\\
453	5.46473576059009e-14\\
454	5.21201611011825e-14\\
455	4.97116592170364e-14\\
456	4.74113752201115e-14\\
457	4.52203864488333e-14\\
458	4.31313195720303e-14\\
459	4.11352414620821e-14\\
460	3.9232963689329e-14\\
461	3.74249909123057e-14\\
462	3.56910248791302e-14\\
463	3.4038021777906e-14\\
464	3.24694493089622e-14\\
465	3.09647138689877e-14\\
466	2.95346518945649e-14\\
467	2.81679486758084e-14\\
468	2.68678198179704e-14\\
469	2.56249337578232e-14\\
470	2.44416343656129e-14\\
471	2.33115343205992e-14\\
472	2.22320980355127e-14\\
473	2.12073811993422e-14\\
474	2.0223717621199e-14\\
475	1.92883210487839e-14\\
476	1.83964174973791e-14\\
477	1.75516512921318e-14\\
478	1.67389734128347e-14\\
479	1.59644355344577e-14\\
480	1.52370531796253e-14\\
481	1.45267580450735e-14\\
482	1.38479993170824e-14\\
483	1.32082802273329e-14\\
484	1.26002826921208e-14\\
485	1.20167044069671e-14\\
486	1.1463419873915e-14\\
487	1.09378443539778e-14\\
488	1.04272498198644e-14\\
489	9.94602828132793e-15\\
};
\addlegendentry{SCF}

\addplot [color=black, mark=+, mark options={solid, black}]
  table[row sep=crcr]{%
0	0.0307613903451636\\
1	0.0217693601554288\\
2	0.0115739226569024\\
3	0.00954579040840216\\
4	0.0287301643042926\\
5	0.0154325634919183\\
6	0.0108727851420687\\
7	0.0280561694554486\\
8	0.0330553300034759\\
9	0.046305695459487\\
10	0.0242638055666314\\
11	0.0317565872960166\\
12	0.0336325318422033\\
13	0.0311335153657411\\
14	0.020895716588482\\
15	0.0288783701820181\\
16	0.0376786778953678\\
17	0.0344133230736374\\
18	0.0396852125675934\\
19	0.0374136778132266\\
20	0.0266213437867668\\
21	0.0349330384752714\\
22	0.0345710873645746\\
23	0.0620698645002473\\
24	0.0399064854582861\\
25	0.0479075955288444\\
26	0.033452289062937\\
27	0.0283851918693263\\
28	0.0372054715266925\\
29	0.046017822757537\\
30	0.0290624472710783\\
31	0.0369289642054486\\
32	0.0431174014560015\\
33	0.028893562843502\\
34	0.0294734377618375\\
35	0.0414321600703184\\
36	0.0414019027999781\\
37	0.0408660030122686\\
38	0.0492840111529391\\
39	0.0382050069583942\\
40	0.0326182562176386\\
41	0.037976058642193\\
42	0.0395810316904311\\
43	0.0216432709066934\\
44	0.0211846537557473\\
45	0.0189285886491355\\
46	0.0191156818009467\\
47	0.0385519870112018\\
48	0.0301251024338112\\
49	0.0522260526317732\\
50	0.040721095707823\\
51	0.045379995250514\\
52	0.046708479702707\\
53	0.0550076861521107\\
54	0.0558426870523648\\
55	0.0434378192754349\\
56	0.051654258783843\\
57	0.0360027850421006\\
58	0.0466948921232027\\
59	0.0569000852837289\\
60	0.0494761821705033\\
61	0.0559201377402861\\
62	0.0327480272382401\\
63	0.051230357631006\\
64	0.0450507454891427\\
65	0.0397130042664524\\
66	0.050726392632891\\
67	0.128691361781586\\
68	0.0481932856237691\\
69	0.0534929072153011\\
70	0.158375836521013\\
71	0.107647314036394\\
72	0.0607345487173408\\
73	0.15319231069966\\
74	0.0520921421325816\\
75	0.164790893123622\\
76	0.135867431118865\\
77	0.155851130174156\\
78	0.113668607172709\\
79	0.177128951947085\\
80	0.142983935795593\\
81	0.0716830807988406\\
82	0.050534123384629\\
83	0.0805300527107119\\
84	0.0851727790784572\\
85	0.0776879800378318\\
86	0.112435016699995\\
87	0.130740830350874\\
88	0.0619131692510733\\
89	0.0545316990614636\\
90	0.169532905974913\\
91	0.101110372441459\\
92	0.0822770619524108\\
93	0.0487300723854129\\
94	0.0450887779280501\\
95	0.0348288254128747\\
96	0.0576502565958447\\
97	0.0248516706947647\\
98	0.0273047717690063\\
99	0.0282312815587899\\
100	0.0324724833515821\\
101	0.0508778671071781\\
102	0.0424557736473682\\
103	0.0429186240809766\\
104	0.0389008901211063\\
105	0.0486120858183629\\
106	0.0422773786672798\\
107	0.046309332977609\\
108	0.0233926785067738\\
109	0.039771560129635\\
110	0.0306054169955443\\
111	0.0382430340041535\\
112	0.0310898302630279\\
113	0.0504707363125847\\
114	0.0393272067545962\\
115	0.0553107463110642\\
116	0.0428121424108488\\
117	0.0578012998363137\\
118	0.0407890913736805\\
119	0.0435439028916443\\
120	0.0398160005095201\\
121	0.0356940373854687\\
122	0.0439859135985562\\
123	0.0350170995541003\\
124	0.048324428912229\\
125	0.0495370598245474\\
126	0.0525843214418626\\
127	0.0493527836787819\\
128	0.0373169473465236\\
129	0.0499529991081209\\
130	0.0302201148793205\\
131	0.0403028125147888\\
132	0.0346085977523881\\
133	0.0460458522100978\\
134	0.0433671663515316\\
135	0.0854286067034665\\
136	0.0577886782130655\\
137	0.0463276045222158\\
138	0.0572623119153969\\
139	0.0508400904177963\\
140	0.0437107860442626\\
141	0.0434963495084718\\
142	0.0504982420377803\\
143	0.0252453516720887\\
144	0.0411491963217376\\
145	0.0563628294681802\\
146	0.0454823698525725\\
147	0.0293623163625726\\
148	0.0383303574379482\\
149	0.0373995642135298\\
150	0.0425438715875674\\
151	0.0360998985865018\\
152	0.0402275998135038\\
153	0.0243411631659304\\
154	0.026888060397054\\
155	0.0536564915014443\\
156	0.0258278436230344\\
157	0.0766390055917152\\
158	0.0428272196807932\\
159	0.0228700783037746\\
160	0.0326421357084441\\
161	0.0387302048108786\\
162	0.0234346928050217\\
163	0.026532042904136\\
164	0.030961723408834\\
165	0.0323887703085275\\
166	0.0216967548921799\\
167	0.029133728961054\\
168	0.0395034249507914\\
169	0.0421317005230353\\
170	0.0336887910171961\\
171	0.0314563837761163\\
172	0.0275260563518241\\
173	0.0377537074700532\\
174	0.0294765351860304\\
175	0.0313853825499725\\
176	0.0258295187650797\\
177	0.0352301831940147\\
178	0.0251391468255713\\
179	0.034568753158698\\
180	0.0307690665601123\\
181	0.0572061068147471\\
182	0.0265740949887549\\
183	0.053737902140602\\
184	0.0286208576767036\\
185	0.0257474475029844\\
186	0.0296573632675327\\
187	0.0264860528660831\\
188	0.0275059321109001\\
189	0.0349516799673662\\
190	0.0301129662824467\\
191	0.0261918285077061\\
192	0.0314522632463082\\
193	0.0365123587726655\\
194	0.0307312168194526\\
195	0.0284608975515433\\
196	0.0297751366019071\\
197	0.0244645696572825\\
198	0.0251446770203358\\
199	0.0270111929224929\\
200	0.0297636155058431\\
201	0.0342715251800786\\
202	0.025472958213043\\
203	0.0275356694124874\\
204	0.0302201705289376\\
205	0.025776393058234\\
206	0.0284588038193998\\
207	0.0291915758188064\\
208	0.0321291809105082\\
209	0.0296739964277617\\
210	0.0296498021440705\\
211	0.0417047724843744\\
212	0.0350407627912563\\
213	0.0293881955187774\\
214	0.0214617711388713\\
215	0.0313818830437998\\
216	0.0276088271546912\\
217	0.0360669601740455\\
218	0.0339269052188519\\
219	0.0286045867067032\\
220	0.0256192544227252\\
221	0.0297981445993043\\
222	0.0327573586870581\\
223	0.0355705143768056\\
224	0.0401077066236\\
225	0.0292777177356792\\
226	0.0255148194402216\\
227	0.0313929970725005\\
228	0.0530170096256651\\
229	0.0445914478508232\\
230	0.0455824206425125\\
231	0.0323287950983367\\
232	0.0405186807826619\\
233	0.0310804069495422\\
234	0.0405177825676647\\
235	0.0405159455009876\\
236	0.0323223944355953\\
237	0.0319200714152132\\
238	0.0337869063652632\\
239	0.0358559396595833\\
240	0.0316061447852755\\
241	0.0372480038944587\\
242	0.0243100108619428\\
243	0.0343337996834239\\
244	0.0377317638598897\\
245	0.032385356595792\\
246	0.029733589665719\\
247	0.0304665204156477\\
248	0.0253147488491289\\
249	0.0198917694089168\\
250	0.018560993564867\\
251	0.0254544141246152\\
252	0.0329835258080661\\
253	0.0221295977576432\\
254	0.0234710383524116\\
255	0.0272168216972337\\
256	0.0231529125514496\\
257	0.0235797041007535\\
258	0.0184933656976886\\
259	0.0279348815731502\\
260	0.0390961631722176\\
261	0.0365329558761994\\
262	0.0360460673335834\\
263	0.0343247245707127\\
264	0.0376822362594417\\
265	0.0664728177531966\\
266	0.0310879204477869\\
267	0.0356304337016001\\
268	0.0534078942328059\\
269	0.0587890515077658\\
270	0.0403914452065965\\
271	0.0330237554074364\\
272	0.03155904118282\\
273	0.046817873008187\\
274	0.0353809666829758\\
275	0.0376582657409829\\
276	0.061637787454881\\
277	0.0451596488014613\\
278	0.030873504811788\\
279	0.0356865619073607\\
280	0.0343368642990679\\
281	0.0484421313464626\\
282	0.0469434594768215\\
283	0.0432549553549801\\
284	0.0284892425563617\\
285	0.0360588438837685\\
286	0.0431361614733175\\
287	0.0350447279055789\\
288	0.0323946412510392\\
289	0.0251949671513899\\
290	0.0375782976099884\\
291	0.0372828385307728\\
292	0.0292537141423347\\
293	0.0355166919223304\\
294	0.0307884006718931\\
295	0.0294596912373172\\
296	0.0313276843347889\\
297	0.0315191217817398\\
298	0.0402689088167582\\
299	0.0378600294809243\\
300	0.0248847463312076\\
301	0.027867353850668\\
302	0.0291207103505055\\
303	0.021779569119256\\
304	0.042202085938569\\
305	0.0338176838864977\\
306	0.0312263140864054\\
307	0.0333201361277315\\
308	0.0305325855830358\\
309	0.0385455074164309\\
310	0.0312903210543897\\
311	0.0306582622130462\\
312	0.0199276625089829\\
313	0.0265638005281855\\
314	0.0220229363848836\\
315	0.0760225201410439\\
316	0.0251819943269084\\
317	0.0229136214438977\\
318	0.0259274275790453\\
319	0.0360790984814403\\
320	0.0284947907336692\\
321	0.0259018042706107\\
322	0.0442810262614371\\
323	0.0349040731475921\\
324	0.0393865816959\\
325	0.0395581087416281\\
326	0.116462087722111\\
327	0.049423753625197\\
328	0.0269287241095992\\
329	0.0265602910886289\\
330	0.0232885683624594\\
331	0.028885041006873\\
332	0.0217252142776215\\
333	0.0199656617325792\\
334	0.0247319527266258\\
335	0.0246437111745871\\
336	0.0293294495161277\\
337	0.0279340116209659\\
338	0.0214862273739881\\
339	0.0377620042227611\\
340	0.0300610610799539\\
341	0.0270154207912971\\
342	0.0291387383313226\\
343	0.0238194861917156\\
344	0.018847959423312\\
345	0.0221608812040795\\
346	0.033739413017498\\
347	0.0397878175796254\\
348	0.031346999519849\\
349	0.0252850606919036\\
350	0.0273883397433659\\
351	0.0398878789668896\\
352	0.0534698200148033\\
353	0.0324628825452622\\
354	0.0358613209813202\\
355	0.0253276593095081\\
356	0.0596492454449071\\
357	0.0347585599229741\\
358	0.0370600590743539\\
359	0.0361223779111621\\
360	0.0304412733923596\\
361	0.122573260045348\\
362	0.0485697610277911\\
363	0.0563205855389615\\
364	0.0454591557410934\\
365	0.0448520945658063\\
366	0.047682445384126\\
367	0.0200289607447088\\
368	0.0190190018317722\\
369	0.026167194279278\\
370	0.0196655971442813\\
371	0.0282676631657151\\
372	0.0258023572454547\\
373	0.0347028833084577\\
374	0.0313040033285424\\
375	0.0331451349146969\\
376	0.0311577154576703\\
377	0.0231332069817509\\
378	0.0277908026983918\\
379	0.0258293793178356\\
380	0.0255332378449179\\
381	0.0285849757709604\\
382	0.0260186454753199\\
383	0.0288039271819204\\
384	0.0282214259776648\\
385	0.0284386004940118\\
386	0.0315972361507571\\
387	0.0279269483665437\\
388	0.0183409159304008\\
389	0.0191377171548889\\
390	0.0408259094982885\\
391	0.0339080586232584\\
392	0.0495466879271443\\
393	0.0377522635847741\\
394	0.023756021268146\\
395	0.0408606793753813\\
396	0.0319499949096455\\
397	0.0660812417282577\\
398	0.0591952573491652\\
399	0.0447163027969842\\
400	0.0307785740525799\\
401	0.050586514003775\\
402	0.0342296775527575\\
403	0.0422638333171973\\
404	0.0520731415093656\\
405	0.0536179074472483\\
406	0.0465705220453521\\
407	0.0424444996232124\\
408	0.0373571366552669\\
409	0.0539556435350327\\
410	0.0570267776336183\\
411	0.0580965860508704\\
412	0.0602946989912297\\
413	0.0464491546559218\\
414	0.0500265922053421\\
415	0.035203852712712\\
416	0.0398141928203502\\
417	0.0368749283990573\\
418	0.046345187820576\\
419	0.0400719222382847\\
420	0.0491951239277907\\
421	0.0305939045757583\\
422	0.02593691609755\\
423	0.0332702529674723\\
424	0.0308548142841704\\
425	0.0255723083656409\\
426	0.0338761997028166\\
427	0.0383282086856096\\
428	0.0335025406808992\\
429	0.0473967003788561\\
430	0.0401823877972348\\
431	0.038636985509099\\
432	0.0426504351469912\\
433	0.0544600303515256\\
434	0.036692706309976\\
435	0.0460673719561213\\
436	0.040793029419293\\
437	0.0361715366203636\\
438	0.0563786068298739\\
439	0.0581455218863648\\
440	0.061285930872589\\
441	0.0440250066914319\\
442	0.0474049951209837\\
443	0.0344551149031255\\
444	0.0557469839939481\\
445	0.0548890543452562\\
446	0.0527871903450485\\
447	0.05709538836635\\
448	0.0379787144104077\\
449	0.0311061410491187\\
450	0.0371810611618674\\
451	0.0275866957376927\\
452	0.0284108820646547\\
453	0.0400983645221438\\
454	0.0273016978613313\\
455	0.038565552796013\\
456	0.0359062045672873\\
457	0.0394046418347454\\
458	0.0762009454574382\\
459	0.0435425123803701\\
460	0.0298493273409305\\
461	0.0302038738136212\\
462	0.0312040882284126\\
463	0.0230912799424877\\
464	0.0336882203211509\\
465	0.0340670113075505\\
466	0.0287432879216262\\
467	0.0258775443839424\\
468	0.0193649542520148\\
469	0.013698910641725\\
470	0.0232394907555982\\
471	0.016854225773563\\
472	0.020167469027529\\
473	0.0191704586845105\\
474	0.0172139123519347\\
475	0.020803444846654\\
476	0.0157974996621458\\
477	0.0224055072592819\\
478	0.0256967887034825\\
479	0.0226470010283654\\
480	0.0245675473649097\\
481	0.0274561494827388\\
482	0.0237421629805664\\
483	0.011900905563007\\
484	0.0201741472361424\\
485	0.0184931522480912\\
486	0.0255062289620697\\
487	0.0168297391502586\\
488	0.024902606988345\\
489	0.0278369585249344\\
490	0.0279044510322014\\
491	0.0319421240459697\\
492	0.0322968217743301\\
493	0.0292134856934263\\
494	0.0267946564407123\\
495	0.0245629925802122\\
496	0.0302100267528873\\
497	0.0408877140515299\\
498	0.0341810180513695\\
499	0.0361187507058937\\
};
\addlegendentry{AA(4)-$X_i$}

\addplot [color=myalign, mark=o, mark options={solid, myalign}]
  table[row sep=crcr]{%
0	0.0307613903451636\\
1	0.0217693601554288\\
2	0.0136288280168109\\
3	0.00810898334208754\\
4	0.00424620600489084\\
5	0.00236685147749986\\
6	0.000940213836041929\\
7	0.000481938954281058\\
8	0.0003265518495882\\
9	0.000140731017937997\\
10	6.72788386895646e-05\\
11	5.5788438261228e-05\\
12	4.84932889280531e-05\\
13	3.49347109279623e-05\\
14	2.06988902983366e-05\\
15	1.30457041683436e-05\\
16	8.51536692723278e-06\\
17	6.8143792865022e-06\\
18	4.41893536744606e-06\\
19	2.63392083389389e-06\\
20	1.85440607635234e-06\\
21	1.59307302690343e-06\\
22	1.08389982171374e-06\\
23	6.11863195253191e-07\\
24	3.19413710970556e-07\\
25	2.69970155704644e-07\\
26	1.85157628624552e-07\\
27	1.54474920387145e-07\\
28	1.31445339787773e-07\\
29	7.17921353256442e-08\\
30	4.30752969809992e-08\\
31	3.55522180678194e-08\\
32	2.7857236124087e-08\\
33	2.12875943951112e-08\\
34	1.42668418701366e-08\\
35	1.03835796866154e-08\\
36	8.99573281457381e-09\\
37	7.30947300673691e-09\\
38	5.99858696404491e-09\\
39	4.57408136127339e-09\\
40	3.3190557806629e-09\\
41	2.88842684016789e-09\\
42	2.16998798795818e-09\\
43	1.54721914256153e-09\\
44	1.13878100833522e-09\\
45	8.35267812150951e-10\\
46	6.54146710566359e-10\\
47	5.6260276154694e-10\\
48	3.59250772695696e-10\\
49	2.66433461724377e-10\\
50	2.18971247863222e-10\\
51	1.66311028272688e-10\\
52	1.35935258278264e-10\\
53	1.22404112335382e-10\\
54	6.63021989093317e-11\\
55	4.5539876483546e-11\\
56	3.20618894198422e-11\\
57	2.78714351866621e-11\\
58	2.57153625774919e-11\\
59	1.80799388824693e-11\\
60	9.23266206059961e-12\\
61	7.9315976644633e-12\\
62	6.58737039194677e-12\\
63	4.91205548457821e-12\\
64	3.63355203879477e-12\\
65	2.2157043531886e-12\\
66	1.69632762853154e-12\\
67	1.39433822192827e-12\\
68	8.73453076010032e-13\\
69	6.59282682371879e-13\\
70	5.23343567193401e-13\\
71	4.35386730070989e-13\\
72	3.89993694897826e-13\\
73	3.08103425092701e-13\\
74	1.66823875003537e-13\\
75	1.17283083400033e-13\\
76	8.65012595234037e-14\\
77	7.83667028440799e-14\\
78	7.06448115581105e-14\\
79	5.18455887741107e-14\\
80	3.49608904994726e-14\\
81	3.06284836107491e-14\\
82	2.56229213880663e-14\\
83	2.28187310557734e-14\\
84	1.60163633805424e-14\\
85	1.06786517028703e-14\\
86	7.14094730860978e-15\\
};
\addlegendentry{AA(4)-$\ALIGN{X_i}$}

\addplot [color=mygrassaa, mark=o, mark options={solid, mygrassaa}]
  table[row sep=crcr]{%
0	0.0307613903451636\\
1	0.0217693601554288\\
2	0.0175045928037686\\
3	0.00542356650632339\\
4	0.005792203951672\\
5	0.00243163690102525\\
6	0.000658560558110899\\
7	0.000394982415898035\\
8	0.00014093915760806\\
9	6.51678442582559e-05\\
10	5.17404035652784e-05\\
11	3.99806703095825e-05\\
12	3.55945494361999e-05\\
13	2.12305152439285e-05\\
14	1.03704441785617e-05\\
15	7.73483349850454e-06\\
16	5.6514398263922e-06\\
17	3.53425610793622e-06\\
18	2.17397054124721e-06\\
19	1.32542784982366e-06\\
20	1.01286114790005e-06\\
21	8.07264929980917e-07\\
22	6.92471419248256e-07\\
23	5.35395062105956e-07\\
24	4.32954763649244e-07\\
25	3.55827501693089e-07\\
26	2.81483296034486e-07\\
27	2.13642821884275e-07\\
28	1.65550289231926e-07\\
29	1.29357353200349e-07\\
30	8.57881172811049e-08\\
31	6.32872971211653e-08\\
32	4.6358224776729e-08\\
33	3.99309759576377e-08\\
34	3.46417014167916e-08\\
35	2.96403600082778e-08\\
36	2.52054597384781e-08\\
37	2.21617729105497e-08\\
38	1.90323874783308e-08\\
39	1.49598552703348e-08\\
40	1.06072368502253e-08\\
41	6.56116299140249e-09\\
42	4.70579011036366e-09\\
43	3.65708551752282e-09\\
44	3.20199551044343e-09\\
45	2.78170376566037e-09\\
46	2.4623689731233e-09\\
47	2.30343655229687e-09\\
48	2.06264870312514e-09\\
49	1.62422072814083e-09\\
50	1.1729196377963e-09\\
51	7.76714795475768e-10\\
52	4.82960175956066e-10\\
53	4.06182991150684e-10\\
54	3.49070612383466e-10\\
55	3.03593242967222e-10\\
56	2.60075096536513e-10\\
57	2.34998112240221e-10\\
58	2.11045764860555e-10\\
59	1.85047376923164e-10\\
60	1.3863825518045e-10\\
61	9.04359176766247e-11\\
62	6.78225585192348e-11\\
63	5.75648617179711e-11\\
64	4.93939196087119e-11\\
65	4.00054318690595e-11\\
66	3.22011339265165e-11\\
67	2.7018997647646e-11\\
68	2.5284746408578e-11\\
69	2.21279661802303e-11\\
70	1.81923357233696e-11\\
71	1.28275572980474e-11\\
72	1.07520426503615e-11\\
73	9.49371381708119e-12\\
74	8.13425685969685e-12\\
75	6.30648194892926e-12\\
76	4.69807005891575e-12\\
77	3.58218560642938e-12\\
78	3.16123164987048e-12\\
79	2.84476902079857e-12\\
80	2.34144487352903e-12\\
81	1.88879474726174e-12\\
82	1.64594994775838e-12\\
83	1.51327378101345e-12\\
84	1.34995378680753e-12\\
85	1.07024796360905e-12\\
86	7.11940137239962e-13\\
87	4.98491167042383e-13\\
88	4.15719032893466e-13\\
89	3.66599224909598e-13\\
90	3.26571688891506e-13\\
91	2.67289000414606e-13\\
92	2.3655581969942e-13\\
93	2.24956718105809e-13\\
94	2.04246085248198e-13\\
95	1.55411722401587e-13\\
96	1.09983657122759e-13\\
97	7.17267403775916e-14\\
98	5.86616901355723e-14\\
99	5.25494548507157e-14\\
100	4.55423755385402e-14\\
101	3.72954493601991e-14\\
102	3.30239362975597e-14\\
103	3.12771271310716e-14\\
104	2.84645605076258e-14\\
105	2.29215854894367e-14\\
106	1.63281380797925e-14\\
107	1.22687731023284e-14\\
108	1.0791976427216e-14\\
109	9.897125688969e-15\\
};
\addlegendentry{GrassAA(4)}

\addplot [color=mygnt, mark=diamond, mark options={solid, mygnt}]
  table[row sep=crcr]{%
0	0.0307613903451636\\
1	0.0217693601554288\\
2	0.0113544570677127\\
3	0.00136659310732321\\
4	6.68063226504466e-05\\
5	3.53331114857243e-07\\
6	3.89599433434472e-11\\
7	1.44414265229457e-15\\
};
\addlegendentry{AlignNT}
\legend{}
\end{axis}
\end{tikzpicture}%

%% file: Figs/tikz/reuters_iter.tex
%
%
\definecolor{mycolor1}{rgb}{1.00000,0.00000,1.00000}%
\definecolor{mycolor2}{rgb}{0.49400,0.18400,0.55600}%
\definecolor{mycolor3}{rgb}{0.00000,1.00000,1.00000}%
\definecolor{mycolor4}{rgb}{0.63500,0.07800,0.18400}%
\begin{tikzpicture}

\begin{axis}[%
width=3.65cm,
height=3.7cm,
scale only axis,
tick label style={font=\footnotesize},
xmin=0,
xmax=110,
xlabel style={font=\footnotesize\color{white!15!black},yshift=.2cm},
xlabel={Iteration},
ymode=log,
ymin=1.1e-16,
ymax=2e-1,
yminorticks=true,
ytickten = {-14, -11, -8, -5, -2},
ylabel style={font=\footnotesize\color{white!15!black},yshift=-.2cm},
ylabel={Residual norm},
axis background/.style={fill=white},
legend style={at={(axis cs:106,8e-4)}, anchor=north east, legend cell align=left, align=left, font = \tiny, draw=white!15!black}
]
\addplot [color=myscf, mark=x, mark options={solid, myscf}]
  table[row sep=crcr]{%
0	0.0197076762812499\\
1	0.0048303764003738\\
2	0.00176920140683371\\
3	0.00093947884106448\\
4	0.000643490626987808\\
5	0.000498436177018471\\
6	0.000409577099987213\\
7	0.0003479136212361\\
8	0.000301947990948283\\
9	0.000266126118724531\\
10	0.000237374266122214\\
11	0.000213811722739761\\
12	0.000194195057939178\\
13	0.000177652253988232\\
14	0.00016354464647823\\
15	0.0001513897801991\\
16	0.000140815286828557\\
17	0.000131529544716868\\
18	0.000123301944379884\\
19	0.00011594895712223\\
20	0.000109323900407978\\
21	0.000103309183996036\\
22	9.78103044874043e-05\\
23	9.2751125956952e-05\\
24	8.80701391666261e-05\\
25	8.37174835666113e-05\\
26	7.96525730660453e-05\\
27	7.58422037463334e-05\\
28	7.22590476346738e-05\\
29	6.8880455850839e-05\\
30	6.56875093111311e-05\\
31	6.26642670466481e-05\\
32	5.97971718286683e-05\\
33	5.70745806632229e-05\\
34	5.44863941454634e-05\\
35	5.20237638991324e-05\\
36	4.9678861569461e-05\\
37	4.7444696258637e-05\\
38	4.5314970037644e-05\\
39	4.32839633593002e-05\\
40	4.13464439395951e-05\\
41	3.94975940545476e-05\\
42	3.77329522896618e-05\\
43	3.60483666375197e-05\\
44	3.44399565135405e-05\\
45	3.29040817892805e-05\\
46	3.14373173569486e-05\\
47	3.00364320632118e-05\\
48	2.86983711037265e-05\\
49	2.74202411678514e-05\\
50	2.61992977773446e-05\\
51	2.5032934383648e-05\\
52	2.39186728821809e-05\\
53	2.28541552756005e-05\\
54	2.18371362751688e-05\\
55	2.08654766740322e-05\\
56	1.9937137361123e-05\\
57	1.90501738716187e-05\\
58	1.82027313911589e-05\\
59	1.73930401477155e-05\\
60	1.66194111380162e-05\\
61	1.58802321457122e-05\\
62	1.51739640164076e-05\\
63	1.44991371611645e-05\\
64	1.38543482648709e-05\\
65	1.32382571800024e-05\\
66	1.26495839893434e-05\\
67	1.20871062237292e-05\\
68	1.15496562229072e-05\\
69	1.10361186291889e-05\\
70	1.05454280048978e-05\\
71	1.00765665656396e-05\\
72	9.62856202233053e-06\\
73	9.20048552560877e-06\\
74	8.79144970683819e-06\\
75	8.40060681050995e-06\\
76	8.02714691314163e-06\\
77	7.67029622423383e-06\\
78	7.32931546518429e-06\\
79	7.00349832221564e-06\\
80	6.69216996979141e-06\\
81	6.39468566104784e-06\\
82	6.11042938209782e-06\\
83	5.83881256718534e-06\\
84	5.57927287181615e-06\\
85	5.33127300119636e-06\\
86	5.09429959140137e-06\\
87	4.86786214084022e-06\\
88	4.65149198975543e-06\\
89	4.4447413454902e-06\\
90	4.24718235149037e-06\\
91	4.05840619803527e-06\\
92	3.878022272821e-06\\
93	3.70565734953121e-06\\
94	3.54095481277235e-06\\
95	3.38357391766038e-06\\
96	3.23318908249895e-06\\
97	3.08948921311103e-06\\
98	2.95217705733539e-06\\
99	2.8209685883735e-06\\
100	2.69559241567874e-06\\
101	2.57578922213879e-06\\
102	2.46131122640801e-06\\
103	2.35192166918626e-06\\
104	2.24739432247977e-06\\
105	2.14751302071277e-06\\
106	2.05207121276907e-06\\
107	1.96087153400762e-06\\
108	1.87372539736098e-06\\
109	1.79045260267267e-06\\
110	1.71088096342517e-06\\
111	1.63484595013128e-06\\
112	1.56219034959794e-06\\
113	1.49276393938073e-06\\
114	1.42642317674107e-06\\
115	1.36303090145149e-06\\
116	1.30245605185698e-06\\
117	1.24457339355825e-06\\
118	1.18926326019857e-06\\
119	1.13641130577645e-06\\
120	1.08590826799742e-06\\
121	1.03764974215575e-06\\
122	9.91535965097085e-07\\
123	9.47471608757767e-07\\
124	9.05365582950352e-07\\
125	8.65130846878218e-07\\
126	8.26684229061669e-07\\
127	7.89946255280772e-07\\
128	7.54840984158344e-07\\
129	7.2129585006662e-07\\
130	6.89241513027224e-07\\
131	6.58611715295119e-07\\
132	6.2934314429464e-07\\
133	6.01375301684734e-07\\
134	5.74650378238661e-07\\
135	5.491131342777e-07\\
136	5.24710785434321e-07\\
137	5.01392893507443e-07\\
138	4.79111262129317e-07\\
139	4.57819837142172e-07\\
140	4.37474611326785e-07\\
141	4.18033533431526e-07\\
142	3.99456421216316e-07\\
143	3.81704878359895e-07\\
144	3.64742215084495e-07\\
145	3.48533372287995e-07\\
146	3.33044849062752e-07\\
147	3.18244633423884e-07\\
148	3.04102136144539e-07\\
149	2.90588127490509e-07\\
150	2.7767467681484e-07\\
151	2.65335094795201e-07\\
152	2.53543878278273e-07\\
153	2.42276657541279e-07\\
154	2.3151014592313e-07\\
155	2.21222091679644e-07\\
156	2.11391231993233e-07\\
157	2.01997249009898e-07\\
158	1.93020727841088e-07\\
159	1.84443116436276e-07\\
160	1.76246687230083e-07\\
161	1.6841450049972e-07\\
162	1.60930369347067e-07\\
163	1.53778826240139e-07\\
164	1.4694509104145e-07\\
165	1.40415040466841e-07\\
166	1.34175178873339e-07\\
167	1.28212610385071e-07\\
168	1.22515012220973e-07\\
169	1.17070609235406e-07\\
170	1.11868149574702e-07\\
171	1.06896881427028e-07\\
172	1.0214653079154e-07\\
173	9.76072802476282e-08\\
174	9.32697486553667e-08\\
175	8.91249717895019e-08\\
176	8.51643837798631e-08\\
177	8.1379799428743e-08\\
178	7.77633972920341e-08\\
179	7.43077034972519e-08\\
180	7.10055763238185e-08\\
181	6.78501914165732e-08\\
182	6.48350276975986e-08\\
183	6.19538538774666e-08\\
184	5.92007155895626e-08\\
185	5.65699230677338e-08\\
186	5.4056039406732e-08\\
187	5.16538693018878e-08\\
188	4.93584483333324e-08\\
189	4.71650326903662e-08\\
190	4.50690893786938e-08\\
191	4.30662868471067e-08\\
192	4.11524860290432e-08\\
193	3.93237317973007e-08\\
194	3.75762447933548e-08\\
195	3.59064136001711e-08\\
196	3.43107872930886e-08\\
197	3.27860683034773e-08\\
198	3.13291056040235e-08\\
199	2.99368881959181e-08\\
200	2.86065388809886e-08\\
201	2.73353083251355e-08\\
202	2.6120569374552e-08\\
203	2.49598116182781e-08\\
204	2.38506362023104e-08\\
205	2.27907508834504e-08\\
206	2.17779652772487e-08\\
207	2.08101863387104e-08\\
208	1.98854140360435e-08\\
209	1.90017372140939e-08\\
210	1.81573296499656e-08\\
211	1.73504462716122e-08\\
212	1.65794195617313e-08\\
213	1.58426561022028e-08\\
214	1.51386332794074e-08\\
215	1.44658961535159e-08\\
216	1.38230544283048e-08\\
217	1.32087795990613e-08\\
218	1.26218022016068e-08\\
219	1.20609091656741e-08\\
220	1.15249413538926e-08\\
221	1.1012791118774e-08\\
222	1.05234000430186e-08\\
223	1.00557567468055e-08\\
224	9.60889479253264e-09\\
225	9.18189068664044e-09\\
226	8.77386197711311e-09\\
227	8.38396542713488e-09\\
228	8.01139526945974e-09\\
229	7.65538154466477e-09\\
230	7.31518851428351e-09\\
231	6.99011312749059e-09\\
232	6.67948357966088e-09\\
233	6.38265792131657e-09\\
234	6.0990227301935e-09\\
235	5.82799184141896e-09\\
236	5.56900513506472e-09\\
237	5.32152739352952e-09\\
238	5.08504716987783e-09\\
239	4.85907575412921e-09\\
240	4.64314615258185e-09\\
241	4.4368121203454e-09\\
242	4.23964724809617e-09\\
243	4.05124406826874e-09\\
244	3.87121323079459e-09\\
245	3.69918267642129e-09\\
246	3.53479688859737e-09\\
247	3.37771614270184e-09\\
248	3.22761581859388e-09\\
249	3.08418571243859e-09\\
250	2.94712941383063e-09\\
251	2.81616368164404e-09\\
252	2.69101785696549e-09\\
253	2.57143331286001e-09\\
254	2.45716291791797e-09\\
255	2.34797051668778e-09\\
256	2.24363045068392e-09\\
257	2.14392709065453e-09\\
258	2.04865439345448e-09\\
259	1.95761546334531e-09\\
260	1.87062215455035e-09\\
261	1.78749469097242e-09\\
262	1.70806127896168e-09\\
263	1.63215776371345e-09\\
264	1.55962727704986e-09\\
265	1.49031993272287e-09\\
266	1.42409249403155e-09\\
267	1.36080809701306e-09\\
268	1.30033595776902e-09\\
269	1.24255110666822e-09\\
270	1.18733412032523e-09\\
271	1.13457088743431e-09\\
272	1.08415237032536e-09\\
273	1.03597437278012e-09\\
274	9.8993732713967e-10\\
275	9.45946093280471e-10\\
276	9.03909761728911e-10\\
277	8.63741456539782e-10\\
278	8.25358166807225e-10\\
279	7.88680567902732e-10\\
280	7.53632862352966e-10\\
281	7.20142626126983e-10\\
282	6.88140637404624e-10\\
283	6.57560767972891e-10\\
284	6.28339817889271e-10\\
285	6.00417400602586e-10\\
286	5.73735811411713e-10\\
287	5.48239907430417e-10\\
288	5.23877002660769e-10\\
289	5.00596747451787e-10\\
290	4.78351026833863e-10\\
291	4.57093872063546e-10\\
292	4.36781352763493e-10\\
293	4.17371485648268e-10\\
294	3.98824166317938e-10\\
295	3.81101058792846e-10\\
296	3.64165539300714e-10\\
297	3.47982601743585e-10\\
298	3.32518810403501e-10\\
299	3.17742208398466e-10\\
300	3.03622250431046e-10\\
301	2.90129764935613e-10\\
302	2.77236861870231e-10\\
303	2.64916898296602e-10\\
304	2.53144414212881e-10\\
305	2.41895084073482e-10\\
306	2.31145651590297e-10\\
307	2.20873908342373e-10\\
308	2.11058625930863e-10\\
309	2.01679518333945e-10\\
310	1.9271720366478e-10\\
311	1.84153158147367e-10\\
312	1.75969686983949e-10\\
313	1.68149875926328e-10\\
314	1.6067756359841e-10\\
315	1.53537309025943e-10\\
316	1.46714356030843e-10\\
317	1.40194606613625e-10\\
318	1.33964581406905e-10\\
319	1.28011410622131e-10\\
320	1.22322788763336e-10\\
321	1.16886959880872e-10\\
322	1.11692689753342e-10\\
323	1.06729243911779e-10\\
324	1.01986367466669e-10\\
325	9.74542567310439e-11\\
326	9.312354510878e-11\\
327	8.89852846233719e-11\\
328	8.50309210280778e-11\\
329	8.12522813422573e-11\\
330	7.76415601471889e-11\\
331	7.41912944072101e-11\\
332	7.08943512964098e-11\\
333	6.77439210389176e-11\\
334	6.47334878430468e-11\\
335	6.18568366621542e-11\\
336	5.9108019285335e-11\\
337	5.64813503132623e-11\\
338	5.39714102094385e-11\\
339	5.15730089950835e-11\\
340	4.92811866684067e-11\\
341	4.70912104089948e-11\\
342	4.4998552290899e-11\\
343	4.29988892752325e-11\\
344	4.10880887936952e-11\\
345	3.92621979654563e-11\\
346	3.75174482073705e-11\\
347	3.58502319672251e-11\\
348	3.42571071028953e-11\\
349	3.2734774084957e-11\\
350	3.12800943173596e-11\\
351	2.98900553057228e-11\\
352	2.85617896061324e-11\\
353	2.72925495726395e-11\\
354	2.60797098536665e-11\\
355	2.49207723861582e-11\\
356	2.38133317632157e-11\\
357	2.27551062450539e-11\\
358	2.17439060333849e-11\\
359	2.07776412208393e-11\\
360	1.98543182757412e-11\\
361	1.89720213516951e-11\\
362	1.81289355516791e-11\\
363	1.73233135433954e-11\\
364	1.65534944188093e-11\\
365	1.58178830150552e-11\\
366	1.5114962192652e-11\\
367	1.44432778508635e-11\\
368	1.38014406886876e-11\\
369	1.31881265438894e-11\\
370	1.26020672915608e-11\\
371	1.20420531078125e-11\\
372	1.15069219421026e-11\\
373	1.09955746435434e-11\\
374	1.05069480977743e-11\\
375	1.00400345019475e-11\\
376	9.59386996914244e-12\\
377	9.16753438190873e-12\\
378	8.76014519234364e-12\\
379	8.37085814077761e-12\\
380	7.99886911050138e-12\\
381	7.64341422563868e-12\\
382	7.30375524570794e-12\\
383	6.97918591483956e-12\\
384	6.6690426574492e-12\\
385	6.37268015559298e-12\\
386	6.08948903246092e-12\\
387	5.81888285137761e-12\\
388	5.56030115554827e-12\\
389	5.3132103086119e-12\\
390	5.07709836639359e-12\\
391	4.85148171187578e-12\\
392	4.63588974389322e-12\\
393	4.42987730957092e-12\\
394	4.2330209285449e-12\\
395	4.04491174306619e-12\\
396	3.86516260604225e-12\\
397	3.69340150437427e-12\\
398	3.52927145574734e-12\\
399	3.37243896591814e-12\\
400	3.22257252656935e-12\\
401	3.07936789668654e-12\\
402	2.94252368222893e-12\\
403	2.81176248749944e-12\\
404	2.68681334373219e-12\\
405	2.56741570679089e-12\\
406	2.45332278698032e-12\\
407	2.34430161272528e-12\\
408	2.24012554191582e-12\\
409	2.1405773272189e-12\\
410	2.04545354626287e-12\\
411	1.95455704470381e-12\\
412	1.86769671628046e-12\\
413	1.7847030962122e-12\\
414	1.70539354798151e-12\\
415	1.62960926142844e-12\\
416	1.55718939182369e-12\\
417	1.48799123764041e-12\\
418	1.42186828838871e-12\\
419	1.35868091951961e-12\\
420	1.29830429214709e-12\\
421	1.24060956965511e-12\\
422	1.18547881908837e-12\\
423	1.13279679710338e-12\\
424	1.08245644477461e-12\\
425	1.03435645102693e-12\\
426	9.88391172969851e-13\\
427	9.44465366007886e-13\\
428	9.02496616511109e-13\\
429	8.62391574408993e-13\\
430	8.24067548985244e-13\\
431	7.87447785127666e-13\\
432	7.52456126125507e-13\\
433	7.19019978366648e-13\\
434	6.87065398371893e-13\\
435	6.56533338770749e-13\\
436	6.27358843606601e-13\\
437	5.99481111138793e-13\\
438	5.72840782688842e-13\\
439	5.47383964236768e-13\\
440	5.23061153243594e-13\\
441	4.99817017079894e-13\\
442	4.77606615234314e-13\\
443	4.56382482084786e-13\\
444	4.36099710247425e-13\\
445	4.16719430919818e-13\\
446	3.98201837336849e-13\\
447	3.8050705575617e-13\\
448	3.63598403098389e-13\\
449	3.4744122703129e-13\\
450	3.3199995834328e-13\\
451	3.17246552223417e-13\\
452	3.03148761606662e-13\\
453	2.89677007026468e-13\\
454	2.76803038042674e-13\\
455	2.645045291783e-13\\
456	2.52748374738586e-13\\
457	2.41517453294905e-13\\
458	2.30784793801653e-13\\
459	2.20528159724429e-13\\
460	2.10728946725195e-13\\
461	2.01365124566161e-13\\
462	1.92415335003129e-13\\
463	1.8386546205527e-13\\
464	1.75694359257218e-13\\
465	1.6788635412444e-13\\
466	1.60426253236744e-13\\
467	1.53297299047188e-13\\
468	1.46486084011586e-13\\
469	1.39977616535914e-13\\
470	1.33758038848824e-13\\
471	1.27811721462993e-13\\
472	1.22132755955645e-13\\
473	1.16705542993333e-13\\
474	1.11518536957055e-13\\
475	1.06562158090055e-13\\
476	1.01828674073438e-13\\
477	9.73019632209057e-14\\
478	9.29766408086748e-14\\
479	8.88450068890435e-14\\
480	8.48978050461535e-14\\
481	8.11279377855114e-14\\
482	7.75212575472727e-14\\
483	7.40752717385543e-14\\
484	7.07867802728824e-14\\
485	6.76397052910613e-14\\
486	6.46335729675441e-14\\
487	6.17607175123941e-14\\
488	5.9016370275405e-14\\
489	5.6398906821922e-14\\
490	5.38906257403863e-14\\
491	5.149363380982e-14\\
492	4.92043365232164e-14\\
493	4.70190767078654e-14\\
494	4.49284290593926e-14\\
495	4.29326436364517e-14\\
496	4.10270690010378e-14\\
497	3.92029030769871e-14\\
498	3.74609622309594e-14\\
499	3.57944827529705e-14\\
500	3.42072371623481e-14\\
501	3.26857395698013e-14\\
502	3.12325106489769e-14\\
503	2.98456761404365e-14\\
504	2.8522592003036e-14\\
505	2.72542256884796e-14\\
506	2.6040848183247e-14\\
507	2.48841062274367e-14\\
508	2.37805096228002e-14\\
509	2.27206611941626e-14\\
510	2.17088088259186e-14\\
511	2.07542069346482e-14\\
512	1.98313954979589e-14\\
513	1.8942741159676e-14\\
514	1.81034963337985e-14\\
515	1.73008992125242e-14\\
516	1.65303020561541e-14\\
517	1.57931601542449e-14\\
518	1.50949004444903e-14\\
519	1.44223974197031e-14\\
520	1.3788273535632e-14\\
521	1.31679092955158e-14\\
522	1.25930617183749e-14\\
523	1.20255087908175e-14\\
524	1.149148808105e-14\\
525	1.09868523257251e-14\\
526	1.05064829034162e-14\\
527	1.00332761724926e-14\\
528	9.57828152478631e-15\\
};
\addlegendentry{SCF}

\addplot [color=black, mark=+, mark options={solid, black}]
  table[row sep=crcr]{%
0	0.0197076762812499\\
1	0.0048303764003738\\
2	0.00742687790917993\\
3	0.00753177094665772\\
4	0.0189599935582662\\
5	0.0440901442577332\\
6	0.0153798529571821\\
7	0.033545944359719\\
8	0.0315220951534906\\
9	0.0213741266134102\\
10	0.0173898059417846\\
11	0.0170358565344111\\
12	0.0466648066516468\\
13	0.0109773785530258\\
14	0.010310204134706\\
15	0.0134626510353376\\
16	0.00873447732652388\\
17	0.0148108106746437\\
18	0.0112985295349889\\
19	0.0119912914204147\\
20	0.0155557699836672\\
21	0.0321968379033493\\
22	0.0455097680986564\\
23	0.0809133704849428\\
24	0.0248883230505314\\
25	0.0260833860448093\\
26	0.0433550035592337\\
27	0.0282641790520371\\
28	0.0195231079563872\\
29	0.014957415601429\\
30	0.0180741568608515\\
31	0.0217172701441693\\
32	0.0150412740917873\\
33	0.0351390728173136\\
34	0.031740836395962\\
35	0.0329198925159917\\
36	0.0332730378280264\\
37	0.0250104539945182\\
38	0.0284655652029749\\
39	0.0189264204715074\\
40	0.0119676292271196\\
41	0.01581727739653\\
42	0.0117081636584477\\
43	0.0083945505383744\\
44	0.00792141348060808\\
45	0.00738955388710154\\
46	0.00940339000186432\\
47	0.00903390033426581\\
48	0.00895998078558805\\
49	0.00755326059452474\\
50	0.0208403469751815\\
51	0.0191600835613943\\
52	0.01607485196075\\
53	0.0235484567517658\\
54	0.0203999076857642\\
55	0.0261794983873165\\
56	0.0266426691725491\\
57	0.0380866726084683\\
58	0.0380926892903886\\
59	0.0258857887664812\\
60	0.0259764907625612\\
61	0.0244757297837716\\
62	0.0263039683211942\\
63	0.0200737176742773\\
64	0.0307798095003294\\
65	0.0279852981080062\\
66	0.0262608353075304\\
67	0.0269692855252649\\
68	0.0389927473558601\\
69	0.0258539706754934\\
70	0.0124179082460969\\
71	0.0165226803846425\\
72	0.00824454546494848\\
73	0.00934830908081985\\
74	0.0135905378375178\\
75	0.0217519390953603\\
76	0.0234031825353955\\
77	0.0190328752484171\\
78	0.0184689110140021\\
79	0.0240013992563033\\
80	0.022441032622915\\
81	0.0210030644273953\\
82	0.0270600543770493\\
83	0.0275423964647565\\
84	0.0248871164155765\\
85	0.0191492973432748\\
86	0.0209093385336142\\
87	0.021187083513357\\
88	0.0304561435813225\\
89	0.0265412840755715\\
90	0.0252786274612514\\
91	0.0210607777342726\\
92	0.0171107379708316\\
93	0.0254484656005008\\
94	0.0251416242644356\\
95	0.0408461802315024\\
96	0.0857990171467961\\
97	0.0251029358342181\\
98	0.0321263991731963\\
99	0.0943910532650046\\
100	0.0640633858865776\\
101	0.0881363521424925\\
102	0.0480965499805968\\
103	0.0430854102890448\\
104	0.118823981955632\\
105	0.0487071428987826\\
106	0.0976190458844148\\
107	0.035945444863129\\
108	0.0423791021174199\\
109	0.0421207435159736\\
110	0.0697218028401365\\
111	0.081085967005454\\
112	0.0255767089047775\\
113	0.0283152184397063\\
114	0.050717575327372\\
115	0.0449165620281973\\
116	0.120762256939007\\
117	0.121030702659187\\
118	0.0710677538176578\\
119	0.0658901926217391\\
120	0.0805955639486998\\
121	0.0521494238755014\\
122	0.0656755890475546\\
123	0.0951828913861836\\
124	0.0783717095924099\\
125	0.0801147409992999\\
126	0.114497030692002\\
127	0.0994985135546004\\
128	0.0914109513524508\\
129	0.0710831328634021\\
130	0.0626488766473327\\
131	0.0421525478569628\\
132	0.0446414495855246\\
133	0.0764986485086704\\
134	0.0612429284581953\\
135	0.0352131470013827\\
136	0.0423750184851866\\
137	0.0309040765157537\\
138	0.0391337346252847\\
139	0.0327687272452992\\
140	0.0180226610167567\\
141	0.0326255382614969\\
142	0.0248275578341183\\
143	0.0430470614131107\\
144	0.0238855452675378\\
145	0.0397290454960567\\
146	0.0397683364419222\\
147	0.0268794200020337\\
148	0.0324989307212599\\
149	0.0126177972478031\\
150	0.0175245634603773\\
151	0.0540183940263291\\
152	0.0231374998028716\\
153	0.0231918024873098\\
154	0.020979038421582\\
155	0.033610575573542\\
156	0.0503691809734114\\
157	0.0186487256875703\\
158	0.0323091230683558\\
159	0.0228024578950297\\
160	0.0288643678997568\\
161	0.0121066295387771\\
162	0.0284649622730596\\
163	0.0318122348254158\\
164	0.0230673293329763\\
165	0.0194614871959855\\
166	0.0211119643364505\\
167	0.0210918032974736\\
168	0.0151900212733168\\
169	0.0159831302158456\\
170	0.0112324962739895\\
171	0.0232996613492075\\
172	0.0275527599799157\\
173	0.0133509368139807\\
174	0.0180742740481007\\
175	0.0197179334300162\\
176	0.0245531836167899\\
177	0.0169898259921928\\
178	0.0126560714894208\\
179	0.0568382509729688\\
180	0.0439129173935681\\
181	0.0302303489124713\\
182	0.0361260529652793\\
183	0.042169004828316\\
184	0.0669981421570004\\
185	0.0538244973343697\\
186	0.0528796420399201\\
187	0.0343971242750707\\
188	0.0250719375939607\\
189	0.0341373692457568\\
190	0.0181317160145844\\
191	0.0743470482144548\\
192	0.0555704919527203\\
193	0.0505537461138352\\
194	0.0801658591835485\\
195	0.0572654336922251\\
196	0.0385674319728377\\
197	0.029609637123951\\
198	0.0465731050386063\\
199	0.0682013888910961\\
200	0.0437872438239269\\
201	0.0512794006174414\\
202	0.0226509981749291\\
203	0.0273623333492046\\
204	0.0267969069236998\\
205	0.0265347784224315\\
206	0.0257697268530509\\
207	0.0245081236944329\\
208	0.0190803178255161\\
209	0.0251734280837542\\
210	0.0187798450124429\\
211	0.0215510951287689\\
212	0.0228897584038034\\
213	0.0261427835387904\\
214	0.0293207640748177\\
215	0.0635672865005374\\
216	0.0465802639769023\\
217	0.0326475803436713\\
218	0.0476412924615845\\
219	0.0462758068029096\\
220	0.0551698640605899\\
221	0.0216025999011976\\
222	0.0256409794436935\\
223	0.0284411086702549\\
224	0.0280484197850654\\
225	0.0258238634662452\\
226	0.0184481625540438\\
227	0.0238945572674453\\
228	0.0203970613250811\\
229	0.0348729571622596\\
230	0.0175120965410838\\
231	0.0186489486671837\\
232	0.0762006285808658\\
233	0.0421260934544176\\
234	0.0448757845452372\\
235	0.0566066599314405\\
236	0.0558104704055609\\
237	0.0581876979738886\\
238	0.0661296129470341\\
239	0.0700571734546662\\
240	0.0623857056336269\\
241	0.0327427894314458\\
242	0.0285954802033736\\
243	0.0450322919574197\\
244	0.0730448025760479\\
245	0.0770160292673368\\
246	0.0478102319133137\\
247	0.0656794407772073\\
248	0.0778599900366019\\
249	0.0685571479955646\\
250	0.0528291438220519\\
251	0.0546365712986816\\
252	0.0315057740913007\\
253	0.032224926180734\\
254	0.0391294309521745\\
255	0.0550832324970464\\
256	0.02343213624831\\
257	0.030832961199795\\
258	0.0183784757835218\\
259	0.0485952490926959\\
260	0.0273280270245648\\
261	0.0257914048870978\\
262	0.0237091162412216\\
263	0.0239656095107005\\
264	0.0413096213739964\\
265	0.0554740151582899\\
266	0.0732798953975564\\
267	0.0585571213623493\\
268	0.0354492166015841\\
269	0.0346854327509942\\
270	0.0427131702124821\\
271	0.0520598150080223\\
272	0.0605594124787677\\
273	0.0328353999889582\\
274	0.0325777513141336\\
275	0.036417208429647\\
276	0.0497113867290909\\
277	0.0487336345175197\\
278	0.0272961736617039\\
279	0.0344751995536942\\
280	0.0522886632682844\\
281	0.027605760405242\\
282	0.0329522562743846\\
283	0.0351704117299958\\
284	0.0356223435205621\\
285	0.0333326900078035\\
286	0.0272202728263373\\
287	0.0278077938835974\\
288	0.0299631131531948\\
289	0.028687751334242\\
290	0.0277566313283242\\
291	0.0226185182456124\\
292	0.028081949221965\\
293	0.0271228249080646\\
294	0.0241056511700268\\
295	0.0257674628063936\\
296	0.0713021235883022\\
297	0.0479586954300402\\
298	0.0215147617675809\\
299	0.0234782728981574\\
300	0.025065057519142\\
301	0.0257108877053414\\
302	0.0267787266969121\\
303	0.0242205643073242\\
304	0.0320204318522606\\
305	0.0278742600200613\\
306	0.0270195056551145\\
307	0.0461802164042266\\
308	0.0301792936835158\\
309	0.0328927401602599\\
310	0.0332549403048036\\
311	0.0295702118179225\\
312	0.0205245393740788\\
313	0.028698976313544\\
314	0.0304615920618247\\
315	0.0284892177092847\\
316	0.0184282362633848\\
317	0.0209766925068472\\
318	0.0226116508946144\\
319	0.0210803750152082\\
320	0.0231859724227807\\
321	0.0245054632713216\\
322	0.0302611558991458\\
323	0.0396415698513691\\
324	0.0297614361773722\\
325	0.0271091673956856\\
326	0.0285607626672712\\
327	0.0252575402267785\\
328	0.0257636750531139\\
329	0.0246878474162199\\
330	0.0113620412755009\\
331	0.0240380811392764\\
332	0.0289912742332738\\
333	0.0251377427279313\\
334	0.0454439862525057\\
335	0.0495412341702317\\
336	0.0574031187083157\\
337	0.0563642741542061\\
338	0.0315599540417729\\
339	0.0253409152219031\\
340	0.0750409132079398\\
341	0.0800503874642337\\
342	0.0616447206362265\\
343	0.0312826623523421\\
344	0.0437552192158672\\
345	0.0227529894693295\\
346	0.0248618754911513\\
347	0.021446233544719\\
348	0.029285148207813\\
349	0.0296748511489658\\
350	0.0209437724981655\\
351	0.0235201920621837\\
352	0.020986272376366\\
353	0.0168344716400783\\
354	0.0199233221165359\\
355	0.011568173722966\\
356	0.0105775269607898\\
357	0.0106141039396432\\
358	0.00910402235637057\\
359	0.00881000678013518\\
360	0.018177923263757\\
361	0.0224209757528021\\
362	0.0224959553817016\\
363	0.0160530573542341\\
364	0.0114391447811568\\
365	0.0141026449785885\\
366	0.015030966224717\\
367	0.011202455137487\\
368	0.0148563608698321\\
369	0.0116561631803191\\
370	0.0715533608356821\\
371	0.0320939259339948\\
372	0.0225327132735645\\
373	0.0770612245632748\\
374	0.0755530207289185\\
375	0.0499941128451628\\
376	0.0556961066195664\\
377	0.0797337305804431\\
378	0.0700720359382514\\
379	0.0462859404258531\\
380	0.0239129790408113\\
381	0.0606053549149042\\
382	0.0741068192423697\\
383	0.0187000790414176\\
384	0.0194330634864514\\
385	0.0193615357757721\\
386	0.0344449659327289\\
387	0.0477842372038321\\
388	0.0688823476348999\\
389	0.0174055645026509\\
390	0.0169660652848807\\
391	0.0671016148558754\\
392	0.0828936921729691\\
393	0.0919886489162617\\
394	0.107516385393223\\
395	0.0792086799518532\\
396	0.0792296067555059\\
397	0.100663554176674\\
398	0.0885046363205504\\
399	0.0440651901207479\\
400	0.112873104040998\\
401	0.0939728475812312\\
402	0.0676731637853271\\
403	0.0389110909812631\\
404	0.124534287415443\\
405	0.118165320319045\\
406	0.107123256163972\\
407	0.0621108496624363\\
408	0.0994304866829359\\
409	0.0480294864306307\\
410	0.118806643508493\\
411	0.0768685988106078\\
412	0.0662176337900905\\
413	0.0775490142702165\\
414	0.102276738305438\\
415	0.0658513018415356\\
416	0.0887830103991738\\
417	0.0984572845858534\\
418	0.0450415655104273\\
419	0.0391391395508234\\
420	0.025580595615264\\
421	0.0226665057079896\\
422	0.0394473350853916\\
423	0.0233462543624364\\
424	0.0348712709007805\\
425	0.0202533809734953\\
426	0.0187095018989557\\
427	0.0264727870383009\\
428	0.0235877749879292\\
429	0.0280682424505156\\
430	0.0251729209769555\\
431	0.0186800943730548\\
432	0.018390668103189\\
433	0.0167998485368785\\
434	0.0190085095985253\\
435	0.0195820482254566\\
436	0.0101914504168059\\
437	0.0259013488795082\\
438	0.0250329974758115\\
439	0.023675532664811\\
440	0.0208900554822278\\
441	0.0276077770796958\\
442	0.0289207463766068\\
443	0.027887796345986\\
444	0.0280153365619471\\
445	0.0302534778591844\\
446	0.0287291844215045\\
447	0.0258407898020461\\
448	0.0286855212734794\\
449	0.025550332460273\\
450	0.0290948367526231\\
451	0.0285368461711996\\
452	0.0184619960496609\\
453	0.0244466441191987\\
454	0.0270026492444269\\
455	0.0233888397538303\\
456	0.0301073466746957\\
457	0.013745479843475\\
458	0.0277614206601048\\
459	0.0186231837255717\\
460	0.0219649288413672\\
461	0.0269634486396151\\
462	0.0278620758759737\\
463	0.0246531463576717\\
464	0.0313750171667807\\
465	0.0301454326132831\\
466	0.0238807834002502\\
467	0.0238433240799445\\
468	0.0218223782223473\\
469	0.0149511213932303\\
470	0.0186773430804\\
471	0.0177938492215755\\
472	0.0193899257293223\\
473	0.0277994360886709\\
474	0.0214474577815295\\
475	0.0302499953103118\\
476	0.00812591214112742\\
477	0.0285566990940151\\
478	0.0685810799222259\\
479	0.0635822256884503\\
480	0.0589131011231567\\
481	0.0680341565322879\\
482	0.0732877675583536\\
483	0.0738358665925784\\
484	0.0621403887201985\\
485	0.0557150936515415\\
486	0.036636445566337\\
487	0.0420640513322866\\
488	0.0449434575491108\\
489	0.0667547766475215\\
490	0.0185098040542241\\
491	0.015867684770825\\
492	0.0350753066968033\\
493	0.0375960496410062\\
494	0.0539069926484277\\
495	0.0667887922245587\\
496	0.0634147275130771\\
497	0.0648891455701001\\
498	0.0692970034695354\\
499	0.0420590613438279\\
500	0.0620253997675912\\
501	0.0602754998701226\\
502	0.0419721699821985\\
503	0.0660171357151709\\
504	0.065068644983486\\
505	0.0564837996244835\\
506	0.068068172548536\\
507	0.0404919768961539\\
508	0.0570903790874672\\
509	0.0406617062181506\\
510	0.0562278508006554\\
511	0.0432156624628136\\
512	0.0238970209423581\\
513	0.0538096307400832\\
514	0.0319707625710713\\
515	0.0284623423365495\\
516	0.0291121480389033\\
517	0.0307258334153798\\
518	0.0231345849755867\\
519	0.0201336719370002\\
520	0.0153526434335101\\
521	0.0176564409231039\\
522	0.0168813880103548\\
523	0.0144238231783759\\
524	0.024464461548052\\
525	0.0215523553793798\\
526	0.0138670833207344\\
527	0.0155723857778054\\
528	0.0333598501463297\\
529	0.0227313451706278\\
530	0.0233712433469042\\
531	0.0205707806460728\\
532	0.0154868354576012\\
533	0.0226630692406947\\
534	0.0308224300401835\\
535	0.0393160187783615\\
536	0.053518169131623\\
537	0.0315758297786693\\
538	0.0304531060137094\\
539	0.0415564045819543\\
540	0.0472359865528026\\
541	0.0374521244908481\\
542	0.0248071783789445\\
543	0.0627639737173423\\
544	0.0205261666499682\\
545	0.0542268498387671\\
546	0.0257799877681428\\
547	0.0241452151002164\\
548	0.028226729022361\\
549	0.0294546024414065\\
};
\addlegendentry{AA(4)-$X_i$}

\addplot [color=myalign, mark=o, mark options={solid, myalign}]
  table[row sep=crcr]{%
0	0.0197076762812499\\
1	0.00483037640037381\\
2	0.00115635673435068\\
3	0.000615450737686731\\
4	0.000352389832734282\\
5	0.00023023731528612\\
6	0.000165829179136708\\
7	0.00012305019345685\\
8	9.94167258063601e-05\\
9	8.17729317563097e-05\\
10	5.12838278301741e-05\\
11	2.68669335076819e-05\\
12	1.46592642549867e-05\\
13	7.25796593297269e-06\\
14	5.09386302558343e-06\\
15	2.80809133461928e-06\\
16	2.31306612514742e-06\\
17	1.26294931658087e-06\\
18	6.89653530450414e-07\\
19	3.86481985160285e-07\\
20	2.68607437212865e-07\\
21	1.63607588391817e-07\\
22	8.8792228488331e-08\\
23	6.75268448066867e-08\\
24	3.87077784160192e-08\\
25	2.8073570034194e-08\\
26	1.25142409468627e-08\\
27	9.76678729284032e-09\\
28	8.93123276478508e-09\\
29	6.13469193941282e-09\\
30	3.08126541268779e-09\\
31	1.72372323793371e-09\\
32	1.18423537893445e-09\\
33	7.4356863032404e-10\\
34	4.72784643153261e-10\\
35	2.75435594209215e-10\\
36	2.28530985560005e-10\\
37	1.89263391556401e-10\\
38	1.42807042402164e-10\\
39	7.13029692748014e-11\\
40	4.87803791171068e-11\\
41	3.36954307643962e-11\\
42	2.05994385885872e-11\\
43	9.87257339006441e-12\\
44	7.61136578152844e-12\\
45	5.97854322922746e-12\\
46	2.67634484130533e-12\\
47	1.54120431030072e-12\\
48	8.38440787580112e-13\\
49	6.0936579972507e-13\\
50	4.49565731488308e-13\\
51	1.96328984349744e-13\\
52	1.45466431480316e-13\\
53	1.00501919968147e-13\\
54	5.08984717008365e-14\\
55	2.96534110051984e-14\\
56	1.40478421572945e-14\\
57	8.6849593349565e-15\\
};
\addlegendentry{AA(4)-$\ALIGN{X_i}$}

\addplot [color=mygrassaa, mark=o, mark options={solid, mygrassaa}]
  table[row sep=crcr]{%
0	0.0197076762812499\\
1	0.0048303764003738\\
2	0.0014950773087173\\
3	0.00068550222858486\\
4	0.000381612731204518\\
5	0.000237745379403511\\
6	0.000155506293984075\\
7	0.000126112943470184\\
8	9.61900716376118e-05\\
9	7.38053635718767e-05\\
10	5.22723940170956e-05\\
11	3.24710299721186e-05\\
12	1.62607434474677e-05\\
13	1.11198430073369e-05\\
14	7.34785488426341e-06\\
15	5.10549812724597e-06\\
16	3.77784371666667e-06\\
17	3.07679877200718e-06\\
18	2.68680191215479e-06\\
19	2.33496588997364e-06\\
20	1.82445514334902e-06\\
21	1.27364621352257e-06\\
22	8.71681956072617e-07\\
23	5.91163272393687e-07\\
24	4.04914758967433e-07\\
25	2.63785451290532e-07\\
26	1.88244783112776e-07\\
27	1.4399773266389e-07\\
28	1.18148793057276e-07\\
29	1.02521706995811e-07\\
30	8.93810720808053e-08\\
31	7.17036257364619e-08\\
32	5.64548505535506e-08\\
33	4.1754656241282e-08\\
34	2.71846806070033e-08\\
35	1.80733351217442e-08\\
36	1.23070287139807e-08\\
37	8.56693719624602e-09\\
38	6.77554645615919e-09\\
39	5.6534825663081e-09\\
40	4.85297595075026e-09\\
41	4.29585216018142e-09\\
42	3.59591619281361e-09\\
43	2.79001726075246e-09\\
44	2.07292174117764e-09\\
45	1.49708770965708e-09\\
46	9.76804674157298e-10\\
47	6.92272212180529e-10\\
48	5.0791963087582e-10\\
49	4.00919977142523e-10\\
50	3.22416404172871e-10\\
51	2.78071125689007e-10\\
52	2.46742367731282e-10\\
53	2.15485267569402e-10\\
54	1.69613560107386e-10\\
55	1.29460311723754e-10\\
56	9.09449183917209e-11\\
57	6.57088668003819e-11\\
58	4.81632492121899e-11\\
59	3.45182137659824e-11\\
60	2.64253977598381e-11\\
61	2.18471346411546e-11\\
62	1.85127467736869e-11\\
63	1.63101576297241e-11\\
64	1.43834350944833e-11\\
65	1.14058606545734e-11\\
66	8.91416375953044e-12\\
67	6.85236124439715e-12\\
68	4.88746766185948e-12\\
69	3.50983789173186e-12\\
70	2.52850401276609e-12\\
71	1.82141642210922e-12\\
72	1.51558795332128e-12\\
73	1.32524700082353e-12\\
74	1.16345807535378e-12\\
75	1.02390640288805e-12\\
76	8.62036476086104e-13\\
77	6.88377456221854e-13\\
78	5.20008494181808e-13\\
79	3.49378311016184e-13\\
80	2.40571365524166e-13\\
81	1.74508404568826e-13\\
82	1.42941020361557e-13\\
83	1.1817368242989e-13\\
84	9.998680510161e-14\\
85	8.56524536608672e-14\\
86	7.77517010243809e-14\\
87	6.50515698504562e-14\\
88	5.00289437495528e-14\\
89	3.63126268162703e-14\\
90	2.70786568213422e-14\\
91	2.05332999074086e-14\\
92	1.58689739733947e-14\\
93	1.17697871373807e-14\\
94	9.9700535590412e-15\\
};
\addlegendentry{AA(4)-Gr}

\addplot [color=mygnt, mark=diamond, mark options={solid, mygnt}]
  table[row sep=crcr]{%
0	0.0197076762812499\\
1	0.0048303764003738\\
2	0.00176920140683371\\
3	0.000208048867573209\\
4	3.6679609515277e-06\\
5	2.82523715315625e-09\\
6	9.67361772299757e-15\\
};
\addlegendentry{AlignNT}
\end{axis}
\end{tikzpicture}%